\documentclass[12pt, reqno]{article}

\usepackage{braket}
\usepackage{hyperref}
\usepackage{theoremref}
\usepackage{amsmath}
\usepackage{amsthm}
\usepackage{amssymb}
\usepackage{amsfonts}
\usepackage{color}
\usepackage[all]{xy} \xyoption{all}
 \usepackage{wrapfig}
\usepackage{enumerate}
\usepackage{enumitem}
\usepackage{slashed}
\usepackage{ytableau}

\usepackage{tikz, tikz-3dplot, pgfplots}
\usetikzlibrary{decorations.markings, patterns, decorations.pathmorphing}

\tikzset{->-/.style={decoration={
markings,
mark=at position #1 with {\arrow{>}}},postaction={decorate}}}

\def\AAA{\mathbb{A}}

\def\NN{\mathbb{N}}

\def\CC{\mathbb{C}}
\def\DD{\mathbb{D}}

\def\GG{\mathbb{G}}

 \def\LL{{\mathbb{L}}}
  \def\MM{{\mathbb{M}}}
\def\PP{\mathbb{P}}
\def\RR{\mathbb{R}}
\def\ZZ{\mathbb{Z}}
\def\QQ{\mathbb{Q}}

\def\gen{\mathfrak{g}}
\def\hen{\mathfrak{h}}

\def\ken{\mathfrak{k}}
\def\len{\mathfrak{l}}
\def\men{\mathfrak{m}}
\def\nen{\mathfrak{n}}
\def\oen{\mathfrak{o}}
\def\pen{\mathfrak{p}}
\def\qen{\mathfrak{q}}
\def\sen{\mathfrak{s}}
\def\ten{\mathfrak{t}}
\def\uen{\mathfrak{u}}

\def\Den{\mathfrak{D}}

\def\Hen{\mathfrak{H}}

\def\Len{{\mathfrak L}}

\def\Nen{\mathfrak{N}}

\def\Ren{\mathfrak{R}}

\def\Ten{\mathfrak{T}}

\def\Ac{\mathcal{A}}

\def\Cc{\mathcal{C}}
\def\Kc{\mathcal{K}}
\def\Dc{\mathcal{D}}
\def\Ec{\mathcal{E}}
\def\Fc{\mathcal{F}}
\def\Gc{\mathcal{G}}

\def\Jc{{\mathcal{J}}}
\def\Kc{{\mathcal{K}}}
\def\Lc{\mathcal{L}}
\def\Mc{\mathcal{M}}
\def\Nc{\mathcal{N}}
\def\Hc{\mathcal{H}}
\def\Oc{\mathcal{O}}

\def\Qc{\mathcal{Q}}
\def\Rc{\mathcal{R}}
\def\Sc{\mathcal{S}}
\def\Tc{\mathcal{T}}

\def\Vc{\mathcal{V}}
\def\Wc{\mathcal{W}}

\def\Zc{\mathcal{Z}}

\def\Gb{\mathbf{G}}

\def\bee{{\mathbf{e}}}
\def\bg{{\mathbf{g}}}
\def\bh{{\mathbf{h}}}

\def\bm{{\mathbf{m}}}
\def\bn{{\mathbf{n}}}
\def\bp{{\mathbf{p}}}

 \def\Aut{{\on{Aut}}}
 
 \def\be{\begin{equation}}
 \def\bef{\begin{figure}}
 \def\bem{\begin{matrix}}
 
 \def\bpm{\begin{pmatrix}}
  
 \def\btp{\begin{tikzpicture}}
  \def\Bun{{\on{Bun}}}

  \def\Ch{{\on{Ch}}}
  \def\CO{{\on{CO}}}

  \def\Coker{\operatorname{Coker}\nolimits}
  \def\CSpin{{\on{CSpin}}}
 \def\CYD{{\CC\hskip -0.07cm\on{YD}}}

 \def\del{{\partial}}
 \def\Det{{\on{Det}}}
 \def\diag{{\on{diag}}}
 \def\dirac{{\slashed{D}}}

  \def\ee{\end{equation}}
 \def\End{\on{{End}}}
 \def\enf{\end{figure}}
 \def\enm{\end{matrix}}
 \def\eord{{\on{eord}}}
 \def\epm{\end{pmatrix}}
 \def\eps{{\varepsilon}}
 \def\etp{\end{tikzpicture}}

\def\Ext{{\on{Ext}}}

 \def\F{{\on{F}}}
 \def\FT{{\on{FT}}}

 \def\G{{\on{G}}}
 \def\gl{{\gen\len}}
 \def\GL{{\on{GL}}}

 \def\Heis{{\on{Heis}}}

 \def\hn{{({1\over 2})^n}}
  \def\Hom{\operatorname{Hom}\nolimits}

   \def\hra{\hookrightarrow}

 \def\i{{\bf {i}}}
 \def\Id{\operatorname{Id}\nolimits}
 \def\Ind{{\on{Ind}}}

 \def\Im{{\on{Im}}}

 \def\K{{\on{K}}}
 \def\kappa{\varkappa}
   \def\Ker{\operatorname{Ker}\nolimits}

 \def\LG{{\on{LG}}}
  \def\Lie{{\on{Lie}}}
   
 \def\lla{\longleftarrow}
   \def\lra{\longrightarrow}

 \def\Mat{{\on{Mat}}}
 
  \def\Mod{{\text{-}\on{Mod}}}
  
 \def\Mp{{\on{Mp}}}
 \def\MF{{\on{MF}}}
 \def\mer{{\on{mer}}}
  
\def\O{{\on{O}}}
\def\Ob{{\on{Ob}}}
 
\def\on{\operatorname}
\def\ol{\overline}
\def\O{{\on{O}}}
\def\OMp{{\on{OMp}}}
\def\oo{{\infty}}
\def\op{{\on{op}}}
\def\ord{{\on{ord}}}
\def\OSp{{\on{OSp}}}
\def\osp{{\mathfrak{osp}}}

 \def\Perv{{\on{Perv}}}
 \def\PGL{{\on{PGL}}}
 \def\phi{{\varphi}}
 \def\PK{{\on{PK}}}
 \def\pt{{\on{pt}}}
 \def\PT{{\on{PT}}}

\def\qis{\on{{qis}}}

\def\red{{\on{red}}}
\def\reg{{\on{reg}}}
\def\Rep{{\on{Rep}}}
\def\rk{{\on{rk}}}

 \def\SF{{\on{SF}}}
 \def\SG{{\on{SG}}}
  
  \def\Sh{{\on{Sh}}}
  \def\SK{{\on{SK}}}
  \def\SL{{\on{SL}}}
  \def\ssl{\mathfrak{{sl}}}
  \def\SLG{{\on{SLG}}}
  
  \def\SMF{{\on{SMF}}}
  
  \def\so{{\sen\oen}}
  \def\SO{{\on{SO}}}
  \def\Sol{{\on{Sol}}}
  
  \def\Sp{{\on{Sp}}}
  \def\sp{{\mathfrak{sp}}} 
  \def\Spec{{\on{Spec}}}
  \def\Spin{{\on{Spin}}}
  \def\ssl{{\sen\len}}
  \def\st{{\on{st}}}
 \def\Stab{{\on{Stab}}}
 \def\Supp{{\on{Supp}}}
 \def\SV{{\on{SV}}}

   \def\Sym{{\on{Sym}}}
   \def\sym{{\on{sym}}}
   
 \def\tr{{\on{tr}}}

\def\uKer{{\ul\Ker}}
\def\ul{\underline}

\def\V{{\on{V}}}

\def\wh{ \widehat}
\def\wt{\widetilde}

\def\ZYD{{\ZZ\hskip -0.07cm\on{YD}}}
\def\Z2YD{{{1\over 2}\ZZ\hskip -0.07cm\on{YD}}}

\def\0{{\ol{0}}}
\def\1{{\ol{1}}}
\def\(({(\hskip -1mm (}
\def\)){)\hskip -1mm )}
\def\-{{\setminus}}
\def \= {{\,\, \simeq \,\,}}
 \def\be{\begin{equation}}
\def\ee{\end{equation}}
\def\ed{\end{document}}

\newtheorem{thm}[equation]{Theorem}
\newtheorem{cor}[equation]{Corollary}
\newtheorem{lem}[equation]{Lemma}
\newtheorem{prop}[equation]{Proposition}

\theoremstyle{definition}

\newtheorem{defi}[equation]{Definition}
\theoremstyle{remark}

\newtheorem{rem}[equation]{Remark}
\newtheorem{rems}[equation]{Remarks}

\newtheorem{ex}[equation]{Example}
\newtheorem{exas}[equation]{Examples}

\numberwithin{itemcounter}{subsection}
\numberwithin{equation}{subsection}
\usepackage[toc,page]{appendix}

\usepackage{etoolbox}
\appto\appendix{\addtocontents{toc}{\protect\setcounter{tocdepth}{1}}}

\begin{document}

\title{ Supersymmetry, differential operators of infinite order and theta functions}

\author{Mikhail Kapranov} 

\maketitle
\abstract{
In 1972, M. Sato proposed an approach to proving modularity  of forms like Thetanullwerte
by characterizing them via certain differential operators
of infinite order (DOI) in the modular variable(s) alone. A DOI is an infinite series in
derivatives decreasing so fast that  it acts on holomorphic functions by a sheaf morphism.
This approach was developed by several authors including Kashiwara, Kawai, Takei
and Yoshida.

 We give an interpretation of this approach using supersymmetry which provides
 a natural source of DOIs: the naive exponential of any odd
 supersymmetry generator is a DOI. The case of  the Riemann theta function of genus $n$ 
 is governed by the supergroup $\OSp(1|2n)$ (and its metaplectic cover)
 acting on a natural super-thickening of the Siegel plane. For $n=2$ this is the
 $3$-dimensional $\Nc=1$ superconformal group and the structure at hand is precisely
 the free massless scalar supermultiplet (combining the Laplace and Dirac equations). For $n>2$
 we get a superextension of the generalized conformal structure existing
 on the Lagrangian Grassmannian as on any Hermitian symmetric space. 
 
 An additional interesting feature here is that the odd supersymmetry generators
 acting ``on-shell'' (i.e., in the space of solutions of the equations of motion) satisfy
 even-style Heisenberg commutation relations. These equations of motion
 upgrade  to a complex of differential operators corresponding to a natural BGG-type
 resolution of the super-Weil representation of $\osp(1|2n)$. 
 
}

% \today	
 \tableofcontents

\vfill\eject
 
 \addtocounter{section}{-1}
 
\numberwithin{equation}{section}
 
  \section{Introduction} 
  
  \paragraph{Differential operators of infinite order.}
  
  One of the remarkable but still somewhat mysterious achievements of M. Sato's
  school of infinite analysis is the  new ``local''  approach 
   \cite {sato-theta, kawai-theta-jap,  kawai-theta, sato-KK-theta} to the
  proof of the modular behavior
  of Thetanullwerte and other special values of theta functions,
   most importantly
  the Riemann Thetanullwert
  \[
  \theta(T)  = \sum_{\bn\in \ZZ^n} e^{\pi \i \, \bn^t T \bn}, \quad T\in \Hen_n \,=\,\bigl\{ T\in \Sym_{n\times n}(\CC)\bigl|
  \Im(T) >0\bigr\}
  \]
  as a function on the Siegel upperhalf plane $\Hen_n$. 
  It proceeds by  characterizing such functions, in their dependence
  on the modular variable(s) alone, by  modular invariant systems
  of  linear differential equations  {\em  of infinite order}, involving
  infinite series like $\sum_{n=0}^\oo u_n(\tau) (d/ d\tau)^n$
  with decay conditions ensuring that the result defines a {\em local
  operator} on holomorphic functions, i.e., an endomorphism
  of the sheaf $\Oc$ of such functions. Such series can be defined
  on any complex manifold $X$ and  form a sheaf of rings
  $\Dc^\oo=\Dc^\oo_X$ acting on $\Oc_X$.

  \vskip .2cm
  
  To give an example, the series $e^{d/d\tau}$ (the shift operator)
  is not a section of $\Dc^\oo_\CC$, as it changes the domain of definition,
  but the expression 
  \[
  \cos(\sqrt{d/d\tau})=\sum_{n=0}^\oo {(-1)^n\over (2n)!}\cdot  {d^n\over d\tau^n}
\] 
is. More generally, if $P$ is any (matrix) differential operator of 
(effective) order $<1$
in the natural sense (see  Definition \ref{def:eord}
 below), 
for example,  $P$  can be an $r$th root of $d/d\tau$ in some matrix algebra 
$\Mat_N(\Dc_\CC)$, then 
$e^P$ lies in $\Mat_N (\Dc^\oo_\CC)$, i.e.,  acts locally. It is from such operators
that Sato, Kashiwara and Kawai \cite {kawai-theta, sato-KK-theta} construct   systems characterizing
 Thetanullwerte. 
 
 \vskip .2cm
 
 \paragraph{Supersymmetry: source of operators of order $1/2$.} 
 The goal of this paper is to relate Sato's approach with the
 concept of {\em supersymmetry} so important in physics. Indeed, supersymmetry
 provides square roots of spacetime translations (and other vector fields)
   and such roots  are 
 (matrix) differential operators
 of order $1/2$.  So:
 \[
 \text{  \bf {The exponential  $e^D$ of any odd supersymmetry
 generator $D$ is a section of} $\Dc^\oo$.  }
 \]
 The simplest example is the operator (``spinor derivative'')
 \[
 D_\xi = {\del\over\del\xi} + \xi {\del\over\del \tau } \,=\,
 \begin{pmatrix}
 0 & 1 \\ {\del /  \del \tau } & 0
 \end{pmatrix}, \quad D_\xi^2 = {\del\over\del \tau}
 \]
 acting in $\Oc_{\CC^{1|1}}$, the sheaf of 
 holomorphic functions on the $(1|1)$-dimensional
 superspace $\CC^{1|1}$ with the  even coordinate $\tau$ and odd coordinate $\xi$.
 Its matrix form comes from  analyzing 
 the decomposition $\Oc_{\CC^{1|1}}= \Oc_\CC \oplus  \xi \Oc_\CC$,
 the procedure known in supergeometry as {\em component analysis}. 
 This matrix form  of $D_\xi$ is in fact one of the operators used in  
 \cite{kawai-theta, sato-KK-theta}.

  \vskip .2 cm

  \paragraph{Even Heisenberg relations for odd supersymmetry generators.} 
  The main point of Sato's approach to modularity $\tau \mapsto -1/\tau $ (in one dimension) 
   can be expressed by the following 
  nonstandard realization of the Heisenberg  commutation relation:
   \be\label{eq:heis-1}
  \biggl[ \sqrt{\del/  \del \tau }, \sqrt{\del/ \del(-1/\tau)} \biggr] = 1/2
  \ee
(in the action on forms of weight $1/2$). 
  It implies that  whenever we can make sense of the operators
  involved, the exponentials of their  appropriate multiples are {\em commuting}
  operators of infinite order. 

\vskip .2cm
  
 Supersymmetry shines a natural light  on this phenomenon 
 and on  its higher-dimensional versions.  They   turn out to come from  the
  action of the Lie superalgebra $\osp(1|2n)$ on $S\Hen_n$, the natural super thickening
  of the Siegel plane $\Hen_n$. The odd part  $\osp(1|2n)_\1$ is the
  standard symplectic space $\CC^{2n}$ and Sato's operators can be interpreted
  as the supersymmetry generators $D_v, v\in\CC^{2n}$. These
  {\em odd operators}  satisfy the {\em even Heisenberg
  relations} (i.e., relations involving usual, not anti-commutators)
   when acting on a natural subspace
   % (Weil representation)
  of the space of holomorphic (super) forms of weight $1/2$ on $S\Hen_n$. 
  
  \paragraph{The main points of the paper.}\label{par:mainpoints}
  More precisely, our ``super-interpretation'' of the theory of Sato and his collaborators
  can be summarized as follows.
  
  \begin{itemize}
  
  \item[(1)]  One  should  ``superize'' the Siegel plane $\Hen_n$ to $S\Hen_n$, a super-domain of
  dimension $({n (n+1)\over 2} | n)$ as discussed above
  and view $\theta(T)$ as a super-form of weight $1/2$ on $S\Hen_n$. The Lie superalgebra $\osp(1|2n)$
  acts on $S\Hen_n$. 
  
  \item[(2)]  $S\Hen_n$ is an open part of the Lagrangian super-Grassmannian 
  %$\SLG(\CC^{2n})$
  parametrizing $n$-dimensional isotropic subspaces in $\CC^{2n|1}$ and so
  carries a certain holomorphic  differential-geometric structure. 
  
   For  $n=2$  this is precisely the $3$-dimensional $\Nc=1$
  superconformal structure of physicists 
  \cite{park, kuzenko-park-etc, kuzenko, nizami}.  Here we recall that the usual Lagrangian Grassmannian $\LG(\CC^4)$
  is the $3$-dimensional projective quadric and so has a natural conformal structure. 
  
  \item[(3)] The structure in (2) gives rise to  a  $\osp(1|2n)$-invariant super-analog $\DD$ of the conformal Laplacian 
  acting from super-forms (not necessarily modular) of weight $1/2$ to tensor-valued superforms of type
  \ytableausetup{centertableaux, boxsize=0.7em}
     $\begin{ytableau}
  {} \\{}
   \end{ytableau}$  (again twisted with $1/2$). 
   
   For $n=2$ this gives the  $3$d  scalar  super-Laplacian
   with values with forms of weight $3/2$ well known in supersymmetry: 
   $\DD \Phi =0$ is the equation of motion for the  {\em  free massless scalar supermultiplet}, see
   \cite[\S2.3]{gates}
    \cite[\S4.1]{deligne-freed}.  It combines the usual Laplace equation for the bosonic scalar
    component and the Dirac equation for its fermionic spinor superpartner. 
   
   \item[(4)]  The supersymmetry generators $D_{p_i},  D_{q_i}$ corresponding to the standard symplectic
   basis $p_i, q_i \in \CC^{2n}=\osp(1|2n)_\1$    satisfy, when acting on $\Ker(\DD)$,
 the even Heisenberg relations  
   \[
   [D_{p_i}, D_{q_j}]_- = {1\over 2} \delta_{ij}, \quad [D_{p_i}, D_{p_j}]_- =  [D_{q_i}, D_{q_j}]_- = 0,
   \]
 and so the operators of infinite order
 \[
 A_i = e^{\sqrt{4\pi \i} D_{p_i}}, \quad  B_i = e^{\sqrt{4\pi \i} D_{q_i}}
 \]
  preserving  $\Ker(\DD)$, all commute with each other on this kernel. 
  
  \item[(5)] The Riemann  Thetanullwert $\theta(T)$ considered as a super-function on $S\Hen_n$  independent on
  the odd variables,  is characterized uniquely, up to a constant factor, by the
  (super) differential equations 
  \be\label{eq:theta-intro}
  \DD\theta = 0, \,\,\, A_i\theta =\theta, \,\,\, B_i\theta = \theta. 
  \ee
  This system  of equations 
  is manifestly modular  invariant  (in the same sense in which
  $\theta(T)$ is a modular form).  Further, this system (or, rather the Koszul complex
  to which it gives rise) is $\RR$-{\em holonomic} in the sense of  \cite {sato-KK-micro}. 
 
  \end{itemize} 
  
  \paragraph { Role of the Weil representation.} 
  As with any discussion of theta functions, Sato's theory is greatly clarified
  by invoking the concept of  the Weil representation.
  In its classical form \cite {lion-vergne, weil}  this is a representation of the double cover of
  $\Sp(2n)(\RR)$ in the space $\Wc=L_2(\RR^n)$. As such,
  it splits into direct sum $\Wc=\Wc_\0 \oplus \Wc_\1$ of irreducible
  representations consisting of functions that are even or odd in the naive
  sense: $f(-q) = \pm f(q)$. By the  ``Gaussian transform''
  (integration against a variable Gaussian,
  see \S \ref {subsec:super-weil-1}\ref {par:sgauss-1}  below and 
   \cite{Kashiwara-vergne, guillemin, gonch-hyper}),   $\Wc_\0$ is embedded 
  into $\F_{1/2}$,
  the space of holomorphic forms of weight $1/2$  on the Siegel
  upperhalf plane $\Hen_n$. It is defined inside $\F_{1/2}$ by
  a system of $2^{\text{nd}}$ order holomorphic ``Laplace type" equations
  which reflects the generalized conformal structure carried by $\Hen_n$
  as an open part of the Lagrangian Grassmannian \cite{gonch-symm, gonch-selecta}.

    \vskip .2cm

   The space $\Wc_\1$ can be similarly embedded into the
  space of vector valued holomorphic forms of weight $1/2$ on $\Hen_n$
  and characterized by Dirac-type equations there.

     \vskip .2cm
  
  The super point of view unites $\Wc_\0$ and $\Wc_\1$ into a single
  irredicuble representation of the double cover $\OMp(1|2n)$ of the supergroup
  $\OSp(1|2n)(\RR)$ and embeds it (by a super-analog of the Gaussian transform,
  see  \S \ref {par:super-lap-n} \ref{par:sgauss-n}) into the space $\SF_{1/2}$
  of holomorphic forms of weight $1/2$ on the
    super-thickening $S\Hen_n$. The image is annihilated by
 the super-Laplacian  $\DD$
  mentioned  in \S \ref{par:mainpoints}. see also 
   Proposition \ref {prop:SG=BW-n} below. 
   
  For $n=2$ this $\DD$  is precisely the ``physical'' super-Laplacian  
  associated to the $3$-dimensional $\Nc=1$ superconformal symmetry. 
  The action of
  $\osp(1|2n)_\1 = \CC^{2n}$ on $\SF_{1/2}$ satisfies
  Heisenberg relations only modulo the ideal of the components of $\DD$
  (which is tensor-valued  for $n>2$)
   and this provides an example of a ``Jacobi structure
  with respect to a $\Dc$-module"
  from \cite{sato-KK-theta}. Using the physical terminology (justified for $n=2$),
  one can say that the Heisenberg relations are satisfied {\em  on-shell}, i.e.,
  on the space of solutions of the equations of motion (in our case, of the
  free massless scalar supermultiplet). 
  
    \vskip .2cm
    Using the  super-Gaussian transform, we can relate the system  \eqref {eq:theta-intro} 
    of local equations  with the more
    familiar system of two groups of equations (one algebraic, one difference)
    on a vector $u=u(x), x=(x_1,\cdots, x_n)\in \RR^n$  in  the distributional completion of
    $\Wc$:
    \be\label{eq:familiar-intro}
    e^{2\pi \, \i \, x_\nu} u(x) = u(x), \quad e^{\del/\del x_\nu} 
    u(x) \bigl( = u(x+\bee_\nu)\bigr) = u(x).
    \ee
   These equations  exhibit $u$, up to a constant factor,  as the  sum of shifted delta-functions
    $u(x)= \sum_{\bm\in \ZZ^n} \delta(x-\bm)$, 
   a classical object in the theory whose Gaussian transform (super or not)  is $\theta(T)$. 
   But \eqref {eq:familiar-intro} cannot be interpreted in terms of the classical
 Weil representation as the operators in
   the exponents  do not belong to the Lie algebra $\sp(2n)$. Instead, they  belong to 
   $\osp(1|2n)$ (in its super-Weil representation)   and this can be seen as one of the
   reasons for introducing 
   supersymmetry into the theory.

      \vskip .2cm
   A right Bernstein-Gelfand-Gelfand
  (BGG) type resolution of the image of $\Wc_\0$ in $\F_{1/2}$  by an equivariant  complex of vector bundles
  and differential operators can be obtained from the ``standard'' parabolic BGG
  resolution of an irreducible highest weight module in a regular block of the category $\Oc$
  by applying the Enright-Shelton equivalence
  \cite {enright-shelton, enright-hunzuker-pruett}.  It is closely related to
  the work of Jozeffiak - Pragacz 
  -Weyman \cite {weyman-resol} on the minimal left resolution of the ideal of symmetric matrices of rank
  $1$, see  \cite {enright-hunziker}. 
   A right   resolution   of the image of $\Wc_\1$ can also  be extracted from
   \cite {enright-shelton, enright-hunzuker-pruett}.  It  corresponds 
   to the minimal  left 
     resolution of a natural
   module over the coordinate ring of the variety of symmetric matrices of rank $1$
   given in   \cite{reiner-roberts}.

  Building on this, we construct a full BGG-type resolution of the image of  $\Wc$  in
  $\SF_{1/2}$  as an $\osp(1|2n)$-module
  by applying the results of Gorelik \cite {gorelik-typical} 
  relating the categories of representations of $\osp(1|2n)$ and its even part
  $\sp(2n)$.

 \paragraph{Possible  further directions. } 
 First, the appearance of supersymmetry in the study of theta functions needs to be understood conceptually.
 In a sense, it exhibits  the elliptic variables (of Jacobi forms) as  ``spinors'' for the modular variables. 
 One wonders about the natural generality of this phenomenon. For example, we can look at
 ``nonabelian theta functions'' and their heat equations  appearing in the theory of moduli of 
 vector bundles on Riemann surfaces. 
 
 \vskip .2cm
  
  In this paper we consider only the  most classical Riemann Thetanullwert but it seems clear that one can
  similarly analyze more general theta functions (with characteristics etc.). This would lead to 
  $\RR$-holonomic
  systems of differential equations of infinite order which are invariant under congruence subgroups
  in $\Sp(2n, \ZZ)$. 
 
 \vskip .2cm
 
 Intuitively, the reason that the odd generators of $\osp(1|2n)$ satisfy unusual ``even'' Heisenberg relations
 in the super-Weil representation is that this representation is ``small'', i.e., realized
 in the space of functions of relatively small number of variables.
  In the classical (non-super) theory
such  small representations correspond to  minimal coadjoint orbits with  the classical Weil
 representation being  a prime example. They are known under the names of Wallach representations
  \cite{enright-hunziker},  or generalized Weil representations \cite {gonch-weil},  or representations annihilated
  by  the Joseph ideal in the enveloping algebra \cite{joseph}.
    For groups  appearing in physics,
  spaces of solutions of conformally invariant equations of motion provide examples of small representations \cite{kobayashi}. 
  It would be interesting to investigate super-analogs of small representations 
  from the point of view of new unusual relations among odd generators. 
  For example, the
 flag superspace $F(2|0, 2|1; 4|1)$  thickening the Grassmannian $G(2,4)$,
 is known as the complex  $4$d  $\Nc=1$ compact  conformal
  superspace \cite{fioresi-lledo-V} and one can look at the induced super-thickening of the matrix ball in $G(2,4)$
  as an analog of $S\Hen_2$ with its superconformal structure. 
  
  \vskip .2cm
  
  From the  physical point of view, we considered only the simplest $3$-dimensional
  $\Nc=1$ superconformal  theory: that of the  free  massless scalar superfield.  One can
  look at more general
  such theories. At the classical level,  the simplest (self-) interacting $3$d $\Nc=1$ superconformal
  theory is that of the $\Phi^4$ scalar superfield  on $\RR^{3|2}$
  with the equation of motion
   $\DD\Phi = \Phi^3$ (in the nonsuper setting, the $\phi^4$ theory is  classically
   conformal
    in $4$
   dimensions, not $3$). In $3$d $\Nc=1$ superconformal field theories, 
   the exponentials of  odd supersymmetry generators
  seem important operators which may also satisfy interesting relations,
  in the classical or quantum setting.

  \vskip .2cm
  
  Finally, one can look into extended supersymmetry ($\Nc=r>1$) analogs of the 
  approach we present,
  involving the  supergroup $\OSp(r|2n)$ and the corresponding $\Nc=r$  thickening of the Siegel upperhalf plane.

    \paragraph{Organization of the paper.} Physics thrives on exceptions. Exceptional groups, 
    exceptional (almost-) isomorphisms between groups such as $SU(2)\to SO(3)$ and other 
    ``coincidences''  all 
    have important physical significance.  I owe this philosophical point to Y. I. Manin. 
    
      \vskip .2cm
    
    In our case, the exceptional isomorphism of complex Lie algebras $\sp(4) = \so(5)$
    connects Siegel modular forms of genus $2$ to $3$-dimensional conformal geometry
    and this extends to the natural super-versions of both theories. 
    For this reason the paper is organized in a ``genetic'' fashion, emphasizing such
    connections. 
    
      \vskip .2cm
    
    Chapter 
    \ref {sec:susy-source} 
     is dedicated to general properties of differential operators of infinite order. 
    Here we emphasize the role of supersymmetry in producing interesting operators
    of order $1/2$ whose exponentials give examples of differential operators of infinite order. 
    
      \vskip .2cm
    
    In Chapter \ref {sec:gen1} 
     we treat  $1$-dimensional (Jacobi) theta functions and their natural place in
    the  supersymmetric picture governed by the Lie supergroup $\OSp(1|2, \RR)$ and its metaplectic
    cover $\OMp(1|2)$ acting on the super-projective line $\PP^{1|1}$
    and the super-thickening of the Lobachevsky plane there. 
    
      \vskip .2cm
    
    Chapter \ref {sec:gen2} is specially devoted to the case of Siegel forms of genus $2$
    and their supersymmetric generalizations which are directly related to $3$-dimensional $\Nc=1$
    superconformal symmetry.  Here, forms of various weights are directly interpreted as
    (super) conformal densities   and  properties of the super-Laplacian
    have natural physical meaning. We treat this particular case in considerable detail
    so as to enable a more brief treatment of the general case later. 
    
      \vskip .2cm
     
    Finally, in Chapter \ref {sec:higher-gen} we develop the theory in full generality. Here, the conformal
    geometry is no longer the governing structure. Instead, we have a certain ``generalized conformal
    structure'' present on any Hermitian symmetric space and studied in 
    \cite{gonch-symm, gonch-selecta}. It seems that such
    structures do not have a direct physical interpretation. 
    
    Our approach
    uses a supersymmetric extension of this structure for the case of the  Lagrangian Grassmanian.
    In this generality, many simple aspects of the genus $2$ case become nontrivial
    representation-theoretic issues. The most important is the BGG resolution of the super-Weil
    representation of $\osp(1|2n)$ which we establish using the 
    Enright-Shelton reduction    \cite {enright-shelton, enright-hunzuker-pruett}
    of  the corresponding singular block of the category $\Oc$ for $\sp(2n)$ to a regular block for
    a different Lie algebra. This reduction is  followed by the Gorelik equivalence \cite {gorelik-typical}
    relating  some representations of $\osp(1|2n)$ with those of its even part $\sp(2n)$.

  \paragraph{Notations and conventions.} 
  We work in the complex analytic category, so all manifolds, vector bundles,
  sections, differential operators etc. will be assumed holomorphic
  unless specified otherwise. 
  We use the boldface notation $\i=\sqrt{-1}$ for the imaginary unit,
  to be distinct from   the ordinary letter $i$
 that  is free for other use.  By $1_n$ we denote the unit $n\times n$ matrix. 
 We often write $s\in\Fc$  to signify that $s$ is a global section of a sheaf
 $\Fc$ on a space $X$, i.e., that $s\in H^0(X,\Fc)$.  By $(\eps_{i_1\cdots i_p})$ we denote
 the ``Levi-Civita symbol'', the totally antisymmetric expression
 equal to $+1$ when $i_1 <\cdots < i_p$ (and to $0$ when
 $i_1, \cdots, i_p$ are not distinct).

   \paragraph{Acknowledgements.} I would like to thank K. Coulembier, M. Gorelik, A. Grekov, 
   A. Polishchuk, V. Serganova, 
    C. Stroppel, K. Vilonen and T. Xue for useful discussions and
    correspondence. 
     This work was supported by the  World Premier International Research Center Initiative (WPI Initiative), 
 MEXT, Japan and  by the JSPS  KAKENHI grant 20H01794. 
 \vfill\eject

 \numberwithin{equation}{subsection}
 
 \section{Supersymmetry as a source of differential operators of infinite order}\label{sec:susy-source}

 \subsection{Differential operators and differential bundles}\label{subsec:DODB}
 
 \paragraph{Differential operators.}\label{par:DO}
  Let $X$ be a complex manifold. 
 Let $E,F$  be 
 holomorphic vector bundles  
 on $X$  which we identify with their sheaves of holomorphic sections. 
 We denote by $\Dc(E,F) = \Dc_X(E,F)$ the sheaf
 of holomorphic differential operators from   $E$ to   $F$. 
 Recall that we have the increasing filtration $\Dc(E,F) = \bigcup_{m=0}^\oo \Dc^{\leq m}(E,F)$,
 where $ \Dc^{\leq m}(E,F)$ consists of operators of order $\leq m$, with $\Dc^{\leq 0}(E,F) = \ul\Hom_{\Oc_X}(E,F)$
 consisting of morphisms of vector bundles. 
 
 For any three holomorphic bundles $E,F,G$ we have 
 the composition 
 \be\label{eq:DO-comp}
 \Dc(F,G)\otimes_\CC\Dc(E,F) \lra \Dc(E,G).
 \ee
 Such composition makes
 each $\Dc(E,E)$ into a sheaf of associative rings. We denote 
 $\Dc_X = \Dc(\Oc_X, \Oc_X)$. For any $E$ as above the sheaf
 $\Dc(\Oc_X,E)$ is a right $\Dc_X$-module,  $\Dc(E, \Oc_X)$
 is a left $\Dc_X$-module and we have
 \[
 \Dc(E,F)  =  \ul\Hom_{\Mod_{\Dc_X}}\bigl(\Dc(\Oc_X, E), \Dc(\Oc_X, F)\bigr) =
 \ul\Hom_{_{\Dc_X}\Mod} \bigl( \Dc(F, \Oc_X), \Dc(E, \Oc_X)\bigr). 
 \]
 In the sequel we will omit the adjective ``holomorphic'' before bundles, differential
 operators etc. 
 
 \paragraph{Differential bundles.}  Cf. \cite[\S 2.3F]{hennion-K}. 
 
 \begin{defi}   Let $E$ be a sheaf of $\CC$-vector spaces  on $X$.
 
 \vskip .2cm
 
(a)  By a {\em differential chart} on $E$ we will mean a datum of:
 \begin{itemize}
\item  An open covering $\{U_a\}_{a\in A}$ of $X$; we write $U_{ab}=U_a\cap U_b$. 

 \item A vector bundle $E_a$ on $U_a$ for each $a\in A$. 
 
 \item
 Isomorphisms of sheaves $\phi_a: E|_{U_a} \to E_a$ such that the resulting
 transition functions
 \[
 \phi_{ab} = \phi_b \circ\phi_a^{-1}: E_a|_{U_{ab}} \lra E_b|_{U_{ab}}
 \]
 are invertible differential operators. 
 \end{itemize} 
 (b) By a {\em differential bundle} (a {\em $\Dc$-bundle} for short) on $X$ we will
 understand a sheaf $E$ of $\CC$-vector spaces together with a differential chart.
 
 \vskip .2cm
 
 \noindent (c) By a {\em differential operator} between two  $\Dc$-bundles $E,F$
 we understand a morphism of sheaves which on each $U_a$ is given,
 in the differential charts, 
 by an (invertible)  differential operator $E_a\to F_a$. The sheaf of such operators is
 denoted $\Dc(E,F) = \Dc_X(E,F)$. 
 \end{defi} 
 
 In particular, for each $\Dc$-bundle $E$ we have a right $\Dc_X$-module
 $\Dc(\Oc_X, E)$ and a left $\Dc_X$-module $\Dc(E, \Oc_X)$. These
 modules are locally free (over $\Dc_X$). 
 
 \begin {prop}
 The correspondence $E\mapsto \Dc(\Oc_X, E)$ resp. $E\mapsto \Dc(E, \Oc_X)$
 induces an equivalence resp. anti-equivalence of the category of
 $\Dc$-bundles and differential operators between them with the
 category of locally free right resp. left $\Dc_X$-modules. \qed
 \end{prop}
 
 \begin{ex}\label{ex:inf-neigh}
 Let $i: X\hookrightarrow Y$ be a closed embedding and 
  and $I_X\subset \Oc_Y$ the ideal of $X$. The sheaf $\Oc_{X^{(d)}} = 
 \Oc_Y/I_X^{d+1}$ is known as  the structure sheaf of the 
  $d$th infinitesimal neighborhood $X^{(d)}\subset Y$. In general, it does not
  have a structure of a sheaf of $\Oc_X$-modules. One way
  to give such a structure  (in fact, a structure of an $\Oc_X$-algebra) 
  is to give  a projection $p:Y\to X$, $p\circ i=\Id_X$.
 Such a projection always exists locally on $X$ and $Y$ but may not exist globally.
   Nevertheless, $\Oc_{X^{(d)}}$ always has a natural
  structure of a $\Dc$-bundle on $X$. Indeed,   although
  two (locally defined) projections
give   two possibly different structures
  of  an $\Oc_X$-module on $\Oc_{X^{(d)}}$, 
   the identity map of $\Oc_{X^{(d)}}$ is given, with respect
  to these two structures,  by an invertible differential operator
  as one can easily verify in coordinates.  
    \end{ex}
    
    \vfill\eject
    
    \subsection{Differential operators of infinite order}\label{subsec:DO-inf}
    
    \paragraph{ Definitions.} Let $\dim_\CC X= n$ and $E,F$ be (holomorphic)
    vector bundles on $X$. In 1959 M. Sato explained (see \cite{SKK, KK-Kimura}) how to
    complete $\Dc(E,F)$ to a sheaf $\Dc^\oo(E,F)$ whose sections, called
    {\em differential operators of infinite order}, are given by possibly  infinite series
    in  higher derivatives but act on holomorphic sections by sheaf morphisms. 
    The intrinsic definition of Sato is as the sheaf of cohomology with support
    in the diagonal $\Delta\subset X\times X$:
    \be\label{eq:D-inf-coh}
    \Dc^\oo(E, F) \, = \, \ul H^n_\Delta\bigl(X\times X, (E^*\otimes\Omega^n_X) \boxtimes F
    \bigr). 
   \ee
 We refer to  \cite{SKK, KK-Kimura} for cohomological definition of the 
 action of sections of $\Dc^\oo(E,F)$
 on sections of $E$ (to give sections of $F$) as well as of the composition
  $\Dc^\oo(F,G)\otimes_\CC\Dc^\oo(E,F) \to \Dc^\oo(E,G)$. We denote $\Dc^\oo_X = \Dc^\oo(\Oc_X, \Oc_X)$. 
The sheaf $\Dc(E,F)$ is embedded into $\Dc^\oo(E,F)$ as the subsheaf of
{\em meromorphic cohomology with support} formed by classes of cocycles
in appropriate coverings which have only polar singularities. 

If we choose a local coordinate system $z_1,\cdots, z_n$ and local trivializations
$E\to\Oc_X^p$, $F\to\Oc_X^q$, then sections of $\Dc^\oo(E,F)$ are given by
series
\be\label{eq:D-oo-expl}
A \,\, =
\sum_{I=(i_1,\cdots, i_n)\in \ZZ_+^n} u_I(z)\,  \del^I, \quad u_I(z)\in \Mat_{p\times q}(\Oc_X), \quad \del^I = {\del^{i_1+\cdots +i_n}
\over \del z_1^{i_1} \cdots \del z_n^{i_n}}, 
\ee
 satisfying the  following {\em overconvergence condition}: the series 
 \be\label{eq:overconv}
 \sum_{I\in\ZZ_+^n} \,  i_1^! \cdots i_n^!\,  u_I(z)\,  \lambda^I, \quad \lambda^I=
  \lambda_1^{i_1} \cdots \lambda_n^{i_n} 
 \ee
 represents an entire (matrix-valued)  function of $\lambda_1,\cdots, \lambda_n$. This condition is invariant under  changes
 of coordinates and trivializations and can be used as an alternative explicit definition of
 $\Dc^\oo(E,F)$. 
 
 \paragraph{Examples of operators of infinite order.} 
 
 \begin{ex}
 Let $X=\CC$. 
  The space of global sections of $\Dc^\oo_\CC$ can be written as
 \[
 H^0(\CC, \Dc^\oo_\CC) \,  =\,  H^1_\Delta(\CC\times \CC, \Omega^1\boxtimes \Oc)\, =\, 
 {(\Omega^1\boxtimes \Oc) (\CC^2-\Delta) \over (\Omega^1\boxtimes \Oc) (\CC^2)}. 
 \]
 An element of $(\Omega^1\boxtimes \Oc) (\CC^2-\Delta)$ is a 1-form $K(w,z) dw$ holomorphic on $\CC^2-\Delta$. 
 The operator $A\in H^0(\CC, \Dc^\oo_\CC)$ corresponding to such a form acts on locally defined
   holomorphic functions $f\in H^0(U, \Oc_\CC)$ by
 \[
f\mapsto Af, \quad  (Af)(z) \,=\,\int_{|w-z|=\eps (z)} K(w,z) f(w) dw, 
 \]
 where $\eps(z)$ is any sufficiently  small  positive real number, the integral being independent of such choice. 
 More precisely, we choose $\eps(z)$
 so that the disk  $\{ w: |w-z|\leq \eps(z\}$ lies in $U$ and  we see that $Af$ is again defined on all of $U$, so $A: \Oc_\CC\to\Oc_\CC$
 is a sheaf homomorphism.

 A series $\sum_{i=0}^\infty u_i(z) \del^i$ of the form \eqref{eq:D-oo-expl} satisfying \eqref{eq:overconv} gives the $1$-form
 \[
 K(w,z) dw \,=\, \sum_{i=0}^\oo { i! \,  u_i(z) \over (w-z)^{i+1} } dw\,\in \, (\Omega^1\boxtimes \Oc)(\CC^2-\Delta)
 \]
 and the corresponding operator $A\in H^0(\CC, \Dc^\oo_\CC)$ acts on functions by
 \[
 (Af)(z) \,=\, \sum_{i=0}^\oo u_i(z) {d^n f \over dz^n},
 \]
 as follows from the Cauchy formula and the Taylor expansion of $f$ near $z$. 
 The meromorphic part of the local cohomology consists of classes of $K(w,z)dw$
 which have only a pole at $\Delta$, in which case $A$ is of finite order. 
 \end{ex}
 
 \begin{ex}\label{ex:d/dx}
 (a) Continuing with  the case $X=\CC$ and denoting
 $\del=d/dz$,  the series $e^\del = \sum_{i=0}^\oo \del^i/i^!$ does not define a section of $\Dc^\oo_\CC$,
 as it does not satisfy  \eqref{eq:overconv}. Conceptually, $e^\del$ is the shift operator
 $(e^\del f)(z) = f(z+1)$
which shifts the domain
 of definition instead of preserving it. 
 \vskip .2cm
 
 (b) On the other hand, the series $\cos (\sqrt{\del})  = \sum_{i=0}^\oo (-1)^i \del^i /(2i)!$  does satisfy \eqref{eq:overconv} and so defines a section of $\Dc^\oo_\CC$. 
 \end{ex}
 
 In the same  vein, we have the following \cite{sato-KK-theta}. 
 
 \begin{defi}\label{def:eord}
  Let $E$ be a vector bundle on $X$, 
 let $D\in  \Dc(E,E)$ be a differential operator and $a=p/q\geq 0$
  be a rational number. We say that $D$ has {\em effective order} $\leq a$
  and write $\eord(P)\leq a$, if $D^q \in \Dc^{\leq p}(E,E)$. 
 \end{defi}
 
 \begin{prop}\label{prop:eord<1}
 In the situation of Definition \ref{def:eord}
  suppose that $a<1$ and $\eord(D)\leq a$.  Then the series $e^D = \sum_{i=0}^\oo D^i/i!$ defines a section of $\Dc^\oo(E,E)$. \qed
 \end{prop}

 \begin{ex}
 Let
 \[
 D = \begin{pmatrix} 0&1 \\ \del & 0 \end{pmatrix} \,\in \, \Dc (\Oc_\CC^2, \Oc_\CC^2), \quad \del = d/dz. 
 \]
 Then $D^2 = 1_2 \cdot \del$ is of order $1$, so $\eord(D)\leq 1/2$
 and   $e^D\in\Dc^\oo (\Oc_\CC^2, \Oc_\CC^2)$. Matrix elements of $e^D$ are operators similar to that in
 Example \ref{ex:d/dx}(b). 
 \end{ex}
 
 \paragraph{Characteristic varieties for operators of infinite order: abstract definition via
 the sheaf  $\Ec_X^\RR$.} 
 \label{par:char-inf}
 Classically,  the {\em characteristic variety} $\Ch(A)\subset T^*X$ of a (scalar) differential operator
 $A\in \Dc_X$ of finite
 order 
  is defined as the zero locus of its highest symbol 
  which has an invariant meaning as a homogeneous function on $T^*X$.

 \vskip .2cm
 
 For an  operator $A\in\Dc_X^\oo$ of infinite order  there is no highest symbol
 to take the zero locus of. A natural abstract definition of $\Ch(A)$ uses the sheaf $\Ec_X^\RR$
 on $T^*X$ known as the sheaf of {\em holomorphic microlocal operators}. We don't give here
 the precise definition which, similarly to \eqref{eq:D-inf-coh}, proceeds via cohomology with support,
 but now  support is taken in various cone-type regions. This definition  
 can be found in  \cite[\S8.9.5]{bjork} \cite [Def.11.4.2(ii)]{KaSha}
 \cite[\S1.4]{kash-micro}, see also a detailed discussion of the case $n=\dim(X)=1$ in
 \cite[\S3.2]{kash-micro}. 
 
 \vskip .2cm
 
 As explained in these references, $\Ec_X^\RR$ is an $\RR$-conic 
 (constant on orbits of $\RR_{>0}$ in $T^*X$) sheaf of rings on $T^*X$. 
 Denoting $\pi: T^*X\to X$ the projection, we have a monomorphism $\Dc^\oo_X \hra \pi_* \Ec_X^\RR$
 which is an isomorphism for $n>1$.
 Among sections of $\Ec_X^\RR$ on appropriate domains in $T^*X$, one finds familiar pseudo-differential
 operators such as $(\del/\del x_i)^{-1}$ or $\bigl(\sum \del^2/\del x_i^2\bigr)^{-1}$ but also fractional
 powers such as $(\del/\del x_i)^s$, $s\in\CC$. 
 
 The above embedding induces an embedding   $\pi^{-1}\Dc_X^\oo \hra \Ec_X^\RR$ of sheaves on $T^*X$.
 In particular, for any (left) $\Dc_X^\oo$-module $\Mc$ we have an $\Ec_X^\RR$-module
 $\Ec_X^\RR \otimes_{\pi^{-1}\Dc_X^\oo} \Mc$. Extending  the standard terminology, call
 a complex   of left $\Dc_X^\oo$-modules {\em strictly perfect}, if it is locally quasi-isomorphic to
 a finite complex of  free modules of finite rank.

 \begin{defi}
 Let $\Mc^\bullet$ be a strictly perfect complex of   $\Dc_X^\oo$-modules.
  Its {\em characteristic variety}  $\Ch(\Mc^\bullet)\subset T^*X$ is defined as the union of the supports
 of the $\Ec_X^\RR$-modules $\ul H^i(\Ec_X^\RR \otimes^L_{\pi^{-1}\Dc_X^\oo} \Mc^\bullet)$. 
 \end{defi}
 
 \begin{ex}\label{ex:Ch-A}
 Let $A\in \Dc_X^\oo$ be a single differential operator of infinite order. Considering it
 as a $2$-term complex $\Dc_X^\oo \buildrel A \over \to \Dc_X^\oo$, we arrive at the following definition:
 $\Ch(A)$ is the set of codirections $(x,\eta)\in T^*X$ such that $A$ is not invertible in the
 ring $\Ec_{X, (x,\eta)}^\RR$,  the stalk of $\Ec_X^\RR$ at $(x,\eta)$. Similarly for
 matrix operators $A\in\Mat_N(\Dc_X^\oo)$. 
 \end{ex}
 
 As the support of any $\Ec_X^\RR$-module, $\Ch(\Mc^\bullet)$ is a closed $\RR$-conic subset in $T^*X$. 
% involutive in the sense of \cite[Def.6.5.1]{KaSha}. 
Following \cite[Def.1.1]{sato-KK-micro} we say that a strictly perfect complex $\Mc^\bullet$ is {\em $\RR$-holonomic},
if $\Ch(\Mc^\bullet)$ is contained in a (real) subanalytic Lagrangian  subset $T^*X$.
Theorem 1.3 of 
  \cite{sato-KK-micro} then implies:
  
  \begin{thm}\label{thm:Char-SS}
  If $\Mc$ is an $\RR$-holonomic strictly perfect complex of $\Dc_X^\oo$-modules, then the solution complex
  $\Sol(\Mc) :=R\ul\Hom_{\Dc_X^\oo} (\Mc, \Oc_X)$ is {\em $\RR$-constructible}, i.e., it is constructible with
  respect to a certain (real) subanalytic stratification $\{X_\alpha\}$ of $X$, with all stalks of each 
  $\ul H^j(\Sol(\Mc^\bullet))$
  being finite-dimensional $\CC$-vector spaces. Further,  $\SS(\Sol(\Mc^\bullet))$, the micro-support of the complex $\Sol(\Mc^\bullet)$
  in the sense of \cite[\S5.1]{KaSha}, 
   is contained in $\Ch(\Mc^\bullet)$. 
  \qed
  
    \end{thm}
    
    \paragraph{Characteristic variety via zeroes of the total symbol.} For operators
    of finite order the  above definition of the characteristic
    varieties is of course compatible with the standard one using highest symbols and filtrations. 
    For instance,  in the situation of Example \ref{ex:Ch-A}, an operator $A\in\Dc_X$ of order $r$
    is invertible in $\Ec_{X,(x,\eta)}^\RR$ if and only if  its highest symbol $\sigma^r_A$,   a function on $T^*X$ homogeneous of degree $r$, does not vanish at $(x,\eta)$.
    For operators of infinite order one would like to have some more concrete
    procedure of understanding or at least estimating $\Ch(A)$. 
    
    \vskip .2cm
    
    A remarkable idea going back to the 1970 paper \cite{kawai-1970}  of Kawai,
    is to use the {\em set of asymptotic directions for the zero locus of the (coordinate dependent) total symbol} as a
  useful substitute for the zero locus of the highest symbol.
  That is, suppose $X$ is an open subset of $\CC^n$, with coordinates $z_1,\cdots, z_n$. 
 For a (square) matrix operator $A\in \Mat_N(\Dc_X^\oo)$ as in \eqref{eq:D-oo-expl}, its total symbol
 is the matrix function
 \[
\sigma_A(z, \lambda) = \sum_I u_I(z) \lambda^I, \quad z\in X, \lambda\in (\CC^n)^*
 \]
 entire in $\lambda\in  (\CC^n)^*$; unlike in \eqref{eq:overconv}, we do not multiply by the factorials. 
 
 \begin{defi}
 Let $A\in\Mat_N(\Dc_X^\oo)$ as before. 
 A codirection $(z,\eta) \in T^*X = X\times(\CC^n)^*$ with $\eta\neq 0$ is called {\em sym-non-characteristic}
 (that is, {\em symbol-non-characteristic}) for $A$, if there is an open $\RR$-conic neighborhood $V\supset\eta$
 in $(\CC^n)^*$ such that $\sigma_A(z, \lambda)$ is invertible, as a matrix,
   for all $\lambda\in V$ with $\|\lambda\|$ sufficiently
 large. The {\em sym-characteristic variety}  $\Ch^\sym(A)\subset T^*X$ is defined as the set of
 all codirections that are not sym-non-characteristic. 
 \end{defi}
 
 Informally, $\Ch^\sym(A)$ is the set of codirections which are asymptotic
 directions  of  the non-invertibility locus of $\sigma_A(z,\lambda)$. 
 For  $N=1$  (scalar operators), non-invertibility  simply  means  vanishing, so we are talking about asymptotic directions of the zero locus.   
  
 \vskip .2cm
 
 As the total symbol is coordinate dependent, it is not clear whether  the definition of 
 $\Ch^\sym(A)$
 is, in general,  independent on the choice of coordinates. Nevertheless, it is
 quite   appealing and useful
 as the following examples and results show. 
 
 \begin{exas}
 (a) Let $A$ be a scalar operator of  finite order $r$, so that $\sigma_A(z,\lambda)$ is a polynomial in $\lambda$ of degree $r$
 and the highest symbol $\sigma_A^r(z,\lambda)$ is the part homogeneous of degree $r$. 
 The hypersurface $Z_{\sigma_A}\subset T^* X$,  the zero locus of $\sigma_A$,  has asymptotic directions, where
 it approaches infinity in $\lambda$; these directions are precisely given by vanishing of $\sigma_A^r$.
 So in this case $\ \Ch^\sym(A) = \Ch(A)$ gives the classical concept. 
 Similarly for matrix operators of finite order. 
 
 \vskip .2cm
 
 (b) Let $X=\CC$ and   $A=\cos\sqrt{d/dz}$ from Example \ref{ex:d/dx}(b). Then 
 $\sigma_A(z,\lambda) = 
 \cos\sqrt{\lambda}$
 which vanishes at $\lambda=\pi^2(n+(1/2))^2$, $n\in \ZZ$. Thus the real positive direction $\RR_+$ is the
 only asymptotic direction of zeroes of $\sigma_A$. This means that $\Ch^\sym(A) = \CC\times\RR_+\subset
 \CC\times\CC = T^*\CC$. In this example $\Ch^\sym$ for an operator in
 one variable has real dimension $3$, while for a nonzero operator of finite order 
 the dimension  would be always $2$ (as we would have a
  Lagrangian subvariety in $T^*\CC$). 
 This is related to the fact that the space of solutions of $Af=0$  (in any 
 open $U\subset \CC$)
 is infinite-dimensional
 and the sheaf of solutions is not constructible. 
  \end{exas}
  
  The following result is proved by Aoki, Kashiwara and Kawai
  \cite[Th.1]{aoki-kash-kawai} by applying the work of Aoki \cite{aoki-exp-2}.
  
  \begin{thm}\label{thm:AKK-inv}
  Let $X\subset \CC^n$ and $A\in\Dc_X^\oo$ be a scalar differential operator
  of intinite order. Then $\Ch(A)\subset \Ch^\sym(A)$,
  i.e., $A$ is invertible in $\Mat_N(\Ec^\RR_{X, (z,\eta)})$ for any sym-non-characteristic direction 
  $(z,\eta)$ for $A$. \qed
  \end{thm}
  
  Renewed attention to these questions at the present time would likely bring new results
  relating $\Ch^\sym(A)$ and $\Ch(A)$ for matrix operators. We will need one
  particular result in this direction.  
  
  \vskip .2cm

  It deals with the situation when $A = e^C-1_N$ is related to
   the exponential of a matrix
  differential operator of finite order $C: E\to E$  in a vector bundle $E$
    on an open $X\subset \CC^n$. Assume that 
  $E=\bigoplus_{p=1}^N E_p$
 is decomposed into a direct sum of subbundles and let 
  $C_{pq}: E_q\to E_p$  be the $(p,q)$th matrix element of $C$. 
  To ensure that $e^C \in\Mat_N(\Dc_{X}^\oo)$
  we assume that $C$ has effective order $<1$ and, more precisely,
  $\ord(C_{pq}) \leq r_p-r_q+\rho$ where the $r_p\in\QQ$, $p=1,\cdots, N$,
   are rational and $\rho <1$.

  \begin{prop}   \label{thm:aoki-inv} 
  In the situation described, suppose that  $(z,\eta)\in T^*X$ is such
  that  no eigenvalue of the block  matrix
  $\|\sigma_{C_{pq}}^{r_p-r_q+\rho}(z,\eta)\|_{p,q=1}^N$ lies in $\i \RR$. 
  Then $e^C-1_N$ is invertible in $\Mat_N(\Ec^\RR_{X, (z,\eta)})$. 
   \end{prop}
   
  \noindent{\sl Proof:} This is  \cite[Lem.2.4]{yoshida} as well as
   \cite[Lem.2.2] {KK-Takei}. 
   More precisely, the statement cited concerns the case when each $E_p=\Oc$,
   so the $C_{pq}$ are scalar differential operators. The general case reduces to
   that by choosing a basis in each $E_p$ and repeating each $r_p$
   the number of times equal to $\rk(E_p)$. \qed
   
   \vskip .2cm
   
   If $C$ has constant coefficients  and when the highest symbols of the
   matrix elements coincide with the total symbols
   (which is the case we will eventually need), then
   Proposition \ref  {thm:aoki-inv}  can be seen as an instance of a matrix analog of
   Theorem \ref {thm:AKK-inv}, as $\i\RR$ is the set of asymptotic directions
   of the zero locus of the function $e^z-1$. 
 
 \vfill\eject

    \subsection{Orders and exponentials for endomorphisms of $\Dc$-modules}
    \label{subsec:orders}

    A differential operator from a vector bundle $E$ to itself and be seen as an endomorphism of the corresponding induced $\Dc$-module. 
 We would like to extend the exponential formalism to endomorphisms of more general $\Dc$-modules. Our approach is  based on  \cite[\S 2]{sato-KK-theta} 
and \cite[Eq.(1.14)] {KK-Takei}.

    \paragraph{The meaning of the exponential.} 
        Let $X$ be a complex manifold as before.  For any coherent left $\Dc_X$-module   $\Mc$   we denote $\Mc^\oo = \Dc_X^\oo\otimes_{\Dc_X}\Mc$. 
    
    Suppose given such $\Mc$ and  $\phi\in\End_{\Dc_X}\Mc$. It would be interesting  to  find conditions guaranteeing the existence of 
     ``the exponential''   $e^\phi\in \End_{\Dc^\oo_X}\Mc^\oo$. Defining $e^\phi$ as the sum of the series $\sum \phi^m/m^!$
     (and studying its convergence) is awkward as it requires working with the topology on $\Mc^\oo$ which is not very explicit. 
     It is more convenient to work with the $1$-parameter family $e^{z\phi}$, $z\in \CC$ and characterize it by a differential equation. 
     
     More precisely, cf. \cite[\S1]{sato-KK-theta}, we consider the product $\CC\times X$ with $z$ being the coordinate on $\CC$. 
      and put 
     \[
     \wt\Mc = \Oc_\CC\wh\otimes_\CC \Nc, \quad \wt\Mc^\oo = \Oc_\CC \wt\otimes_\CC \Nc^\oo. 
     \]
     Then $\End_{\Dc_X^\oo} \wt\Mc^\oo$ consists of $1$-parameter families $U(z), z\in\CC$ of endomorphisms of $\Mc^\oo$,
     holomorphically depending on $z$ and carries the structure of $\Dc_\CC$-module via the structure on $\Oc_CC$.
     
     \begin{defi}
     The exponential $e^{z\phi}$ is an endomorphism $u(z)\in \End_{\Dc_X^\oo} \wt\Mc^\oo$ satisfying the differential equation 
     $\del u(z)/\del z  = \phi\cdot u(z)$ and the initial condition $u(0)=\Id$. 
          \end{defi}
          
          It is clear that the solution in question is unique, if it exists so in the latter case the  notation $e^{z\phi}$  makes clear sense.

    \paragraph{The order of an endomorphism with respect to a filtration.}  Let $F= (F_i)_{i\in\ZZ}$ be a good filtration on $\Mc$,
    in particular $F_i=0$ for $i\ll 0$ and $F_{i+1} = \Dc_X^{\leq 1} F_i$ for $i\gg 0$. 
    
    \vskip .2cm

    We say that the  $\ord_F(\phi) \leq a$,  if $\phi(F_i)\subset F_{i+a}$ for $i\gg 0$. 
    For $a,b\in \ZZ_{>0}$ we  say  that    $\eord_F(\phi)\leq a/b$, if  $\ord_F(\phi^b) \leq a$. 
    We refer to $\ord_F$ and $\eord_F$ as the {\em order} and {\em effective order} with respect to $F$. 
    These concepts depend on the choice of $F$, although a shift of $F$ leaves them unchanged. 
    The following is an adaptation of  \cite[Lemma 1.2]{sato-KK-theta}.
    
    \begin{prop}\label{prop:exp-endo}
    If $\eord_F (\phi)<1$, then $e^{z\phi}$ exists. 
    \end{prop}

    \noindent{\sl Proof:} The statement can be proved locally, so in this proof we assume that $X$ is a Stein manifold
    and think of $\Dc_X, \Mc$ etc. as the spaces of global sections of the corresponding sheaves.

     Let $a<b$ be positive integers. Denote  $\psi=\phi^b$ and suppose
    $\ord_F(\psi)\leq a$.  Let $i\gg 0$, so that $\Dc_X F_i = \Mc$. A basis $m_1,\cdots, m_r$ of $F_i$
    is then a set of generators of $\Mc$, so we have a surjection $\Dc_X^r = \bigoplus_{\mu=1}^r \Dc_X e_\mu \buildrel\pi\over \to \Mc$
    with the basis vector $e_\mu$ going to $m_\mu$. Further, as $i\gg 0$, we have $\psi(F_i)\subset F_{i+a} = \Dc_X^{\leq a} F_i$
    and so there is an $r\times r$  matrix $P=\|P_{\mu\nu}\|$ of differential operators $P_{\mu\nu}\in\Dc_X^{\leq a}$ such that the 
    diagram commutes
    \[
    \xymatrix{
    \Jc\ar[d] \ar[r] & \Dc_X^r \ar[r]^\pi \ar[d]^P & \Mc\ar[d]^\phi \ar[r]& 0
    \\
      \Jc \ar[r] & \Dc_X^r \ar[r]^\pi & \Mc \ar[r]& 0
    }
    \]
 where $\Jc=\Ker(\pi)$, a left ideal in $\Dc_X^r$; in particular, $\Jc  P \subset \Jc$. In order to make sense of
 \[
 e^{z\phi} \,= \sum_{m=0}^\oo {z^m \phi^m\over m!} \,=\, \sum_{k=0}^{b-1}\phi^k  u_k(z), \quad u_k(z) := \sum_{l=0}^\oo { \psi^l z^{k+lb}\over
 (k+lb)!}
   \]
 via the differential equation, we make sense of each $u_k(z)$, which should satisfy the equations
 \be\label{eq:eqs-u_k}
 \begin{cases}
 \del u_k/\del z = u_{k-1}(z), \,\,\, u_k(0)=0, \quad k=1,\cdots, b-1; 
 \\
 \del u_0/\del z = \psi \cdot u_{k-1}(z), \,\, u_0(0) = \Id. 
 \end{cases}
 \ee
 Having the $u_k(z)$ satisfying these, it is immediate that $u(z) = \sum_{k=0}^{b-1} \phi^ k u_k(z)$ satisfies the conditions for $e^{z\phi}$.

 To construct $u_k(z)$, we consider the series
 \[
 U_k(z) =  \sum_{l=0}^\oo { P^l z^{k+lb}\over
 (k+lb)!} \,\in \, \Mat_r(\wt\Dc_X^\oo), 
 \]
 whose belonging to $\Mat_r(\wt\Dc_X^\oo)$ follows from the fact that $\ord(P) \leq a < b$  which gives the hyperfactorial decay
 of the denominators. These series satisfy the analogs of \eqref{eq:eqs-u_k} with $P$ instead of $\psi$. 
 We claim that the $U_k(z)$ descent to   $u_k(z) \in \End_{\Dc_X^\oo}(\wt\Mc^\oo)$, i.e., 
  that
 $\wt\Jc^\oo U_k(z) \subset \wt\Jc^\oo$.  
 
 For this, put  $\Jc^{\leq i} = \Jc \cap (\Dc_X^{\leq i})^r$. This gives  a   good filtration on $\Jc$, in particular
 $\Jc^{\leq m_0+a} = \Dc_X^{\leq a} \Jc^{\leq m_0}$ for $m_0 \gg 0$.  Choose  such $m_0$ and choose a
 set $R_1,\cdots, R_d\in\Jc^{\leq m_0}$ generating $\Jc$ as a left ideal. 
 As $\Jc P\subset \Jc$, we have
 \[
  \Jc^{\leq m_0} P \,\subset\,  \Jc^{m_0+a} \,= \, \Dc_X^{\leq a} \Jc^{\leq m_0}.
  \]
 So there exists a $d\times d$ matrix $A = \|A_{\mu\nu }\|$ over $\Dc_X$ with $\ord(A_{ll'})\leq a$ such that
 \[
 R_\mu  P = \sum_{\nu }  A_{\mu\nu } R_{\nu} \text { and therefore } R_\mu P^l = \sum_\nu \| A^l\|_{\mu\nu} R_\nu\text{ for any } n\geq 1.
 \]
 This implies that for any power series with complex coefficients $f(x,z) = \sum_{l,m=0}^\oo a_{lm} x^l z^n$ 
 such that the matrix series   $f(P, z)$ and $f(A,z)$   converge over $\Dc_X^\oo$ for any $z$, 
 we have
 \be\label{eq:comm-R}
 R_\mu f(P,z) = \sum \|f(A,z)\|_{\mu\nu} R_\nu,  
 \ee
 an equality of matrices over $\wt\Dc_X^\oo$. Let us take
 \[ 
 f(x,z) = \sum_{l=0}^\oo { x^l z^{k+lb}\over
 (k+lb)!}. 
 \]
 Then both $f(P,z)$ and $f(A,z)$ make sense as matrices over $\wt\Dc_X^\oo$, since the orders of
 $P$ and $A$ are $\leq a$ and in the denominators we have $(k+lb)!$ with $b>a$. 
 The equality \eqref{eq:comm-R} then implies that $R_\mu U_k(z) \subset \wt\Jc^\oo$.
 Since the $R_\mu$ also generate $\wt\Jc^\oo$ over $\wt \Dc_X^\oo$, we conclude that
 $\wt\Jc^\oo U_k(z) \subset \wt\Jc^\oo$ as claimed. The proposition is proved. \qed
  
   \vfill\eject
  \subsection{Differential operators on supermanifolds} \label{subsec:DO-super}
  
  \paragraph{Generalities on  super-algebra. Center and anticenter. }
  \label{par:super-general}
   We use \cite{manin} as a general reference. 
    Accordingly, 
  all linear algebra objects (vector spaces, modules, sheaves etc.) will be 
  assumed  $\ZZ/2$-graded: $V=V_\0\oplus V_\1$.  By $|v|\in \{\0, \1\}$
  we denote the parity of  a homogeneous $v\in V$. 
  We denote by $\Pi$
  the parity change functor. 
  
  \vskip .2cm
  
   For any super-vector space $V$  the {\em parity operator}
  $(-1)^\F = (-1)^\F_V: V\to V$
 is the automorphism acting by $+1$ on $V_\0$ and by $-1$ on $V_\1$.
 The  notation comes from the fact that in many cases there is a natural
 ``fermionic number operator''  $\F$ with integer eigenvalues whose
 parity gives the super-grading.  The correspondence $V\mapsto (-1)^\F_V$
 is a tensor automorphism of the identity functor on the tensor category
 of super-vector spaces, i..e,
 \[
 (-1)^\F_{V\oplus W} = (-1)^\F_V \oplus (-1)^\F_W, \quad 
  (-1)^\F_{V\otimes W} = (-1)^\F_V \otimes (-1)^\F_W. 
 \]
 It follows that  if $R$ is any associative (or Lie, or any other kind of)
 Lie superalgebra, then $(-1)^\F_R: R\to R$ is an algebra
 automorphism..
 
 \vskip .2cm

 Let  $R$ be an associative superalgebra and $M,N$ two left $R$-modules.
  By a  {\em skew morphism} $f: M\to N$ 
  we will mean  an even $\CC$-linear map such that
  $f(am) = (-1)^{|a|} af(m)$. For example, $(-1)^\F : M\to M$ is
  a skew morphism. Alternatively, a skew morphism is just a
  $(-1)^\F_R$-semilinear morphism: $f(am) = (-1)^\F (a) f(m)$. 
   As the composition of two skew morphisms is
  a  morphism of $R$-modules in the usual sense, $f$ is a skew
   morphism if and only if 
   $(-1)^\F \circ  f = f \circ (-1)^\F$ is a usual morphism.   
  
  \vskip .2cm
  
    By $\Zc(R)$ we denote the (super) center of $R$, formed by $a\in R$ such that
  $ab = (-1)^{|a|\cdot |b|} ba$ for all $b\in R$. This is a commutative superalgebra.
Following  \cite{gorelik-ghost},  by $\Ac(R)$ we denote the {\em anticenter} 
of $R$ formed by $a$ such that $ab = (-1)^{(|a|+1)|b|} ba$ for all $b$.
Thus an even element of $\Ac(R)$ anticommutes with odd elements
and commutes with even ones. 
  Inside $R$, we have 
  \[
  \Zc(R)\cdot  \Ac(R) = \Ac(R)\cdot \Zc(R) \,\,\subset\,\, \Ac(R), \quad \Ac(R) \cdot \Ac(R)\,\, \subset\,\,  \Zc(R).
  \]
  The sum $\wt\Zc(R) = \Zc(R) + \Ac(R)\subset R$ is called the {\em ghost center} of $R$.
  The sum is direct if no nonzero $a\in\Zc(R)$ is a zero-divisor in $R$,
  see  \cite{gorelik-ghost}. 
  
  \begin{ex}
  Let $R=\End(V)$ for a super-vector space $V$. Then $(-1)^\F\in\Ac(R)$
  and $\wt \Zc(R) = \CC\cdot 1 \oplus \CC\cdot (-1)^\F$. 
  \end{ex}

  \paragraph{Operators of finite order.}    
   Let $X$ be a (complex analytic, as always)  supermanifold of
  dimension $n|m$. Thus   $X$ is a ringed space
  $X=(X_\red, \Oc_X)$ where $X$ is a manifold of dimension $n$ and
  $\Oc_X$ is a $\ZZ/2$-graded sheaf of commutative superalgebras
  locally isomorphic to $\Oc_{X_{red}}\otimes\Lambda[\xi_1, \cdots, \xi_m]$.

  A vector bundle of rank $r|s$ on $X$ is a  sheaf $E$
   of ($\ZZ/2$-graded) $\Oc_X$-modules which is locally isomorphic to
   $\Oc_X^r \oplus (\Pi \Oc_X)^s$. In particular, we have the
   {\em Berezinian bundle}  or the {\em bundle of volume forms}
   $\omega_X$ which is a vector bundle of rank $1|0$ on $X$.

   Extending the framework of \S   \ref{subsec:DODB} \ref{par:DO}, 
  for any two vector bundles $E,F$ we define  the sheaf $\Dc(E,F) = \Dc_X(E,F)$
   of differential operators $E\to F$,  the composition maps
   as in \eqref{eq:DO-comp} and denote $\Dc_X=\Dc(\Oc_X, \Oc_X)$,
   cf. \cite{penkov}. 
   
 Similarly to Example \ref{ex:inf-neigh}, a vector bundle $E$ on $X$ can
   be seen as a $\ZZ/2$-graded $\Dc$-bundle on $X_\red$, and
   a section of $\Dc_X(E,F)$ can also be seen as a differential operator
   between the corresponding $\Dc$-bundles on $X_\red$. 
   That is,  we have a morphism of sheaves on $X_\red$
   \be\label{eq:D-Dred}
  \phi:  \Dc_X(E,F) \lra \Dc_{X_\red}(E,F),
   \ee
   where on the right $E,F$ are considered as $\Dc$-bundles on $X_\red$.
   
   \begin{prop}
  $\phi$  is an isomorphism of $\ZZ/2$-graded sheaves. 
   \end{prop}
   
   \noindent {\sl Proof:} It is enough to consider $X=\CC^{n|m}$ and work in coordinates
   $z_1,\cdots, z_n, \xi_1,\cdots, \xi_m$. Further, it is enough to consider the case when
    $E = \Oc_X^{r|s}, F=\Oc_X^{u|v}$ are trivial. For simplicity 
    of notation consider the case   $E=F=\Oc_{\CC^{n|m}}$
    the general case is similar. 
    
    The choice of coordinates identifies $\Oc_{\CC^{n|m}}$ with  
    $\Oc_{\CC^m} \otimes \Lambda$,
    $\Lambda =\Lambda[\xi_1,\cdots, \xi_m]$
    with a trivial vector bundle with fiber $\Lambda$. A differential operator in 
    $\Oc_{\CC^{n|m}}$
considered as a super-line bundle on $\CC^{n|m}$ 
    has, in the obvious multi-index notation,   the form
    \[
    \sum_{I,J,K} a_{IJK}(z)\xi^I \del_\xi^J \del_z^K, \quad a_{IJK}(z)\in\Oc_{\CC^n}.
    \] 
 A differential operator in $\Oc_{\CC^{n|m}}$ considered as a trivial bundle 
 $\Oc_{\CC^n}\otimes\Lambda$ on $\CC^n$, has the form
 \[
 \sum_K u_K(z) \del_z^K, \quad u_K \in \End_\CC(\Lambda) \otimes \Oc_{\CC^n}. 
 \]
 To  compare  these two type of objects,  we notice that $\End_\CC(\Lambda)$ is 
 identified with the space of
 differential operators on $\CC^{0|m}$, i.e., of expressions of the form
 $\sum_{I,J} c_{IJ} \xi^I \del_\xi^J$, $c_{IJ}\in\CC$. \qed
 
 \paragraph{Example: the parity operator $(-1)^\F$.}

  Given a vector bundle $E$ on a supermanifold $X$, the operator
 $(-1)^F$ on local sections of $E$ is a canonical section of $\Dc_X(E,E)$.
 To write it in coordinates, we first consider the case
 $X=\CC^{0|m}$ and $E=\Oc$, i.e., $H^0(X,E)=\Lambda[\xi_1,\cdots, \xi_m]$. 
 In this case we find directly that
 \be\label{eq:-1^F}
 \begin{gathered}
 (-1)^\F = \sum_{p=0}^n (-2)^p \sum_{1\leq i_1 < \cdots < i_p\leq m}
 \xi_{i_1}\cdots \xi_{i_p} {\del\over\del \xi_{i_p}} \cdots {\del\over\del\xi_{i_1}} = 
 \\
 = 1-2\sum_i \xi_i {\del\over\del\xi_i} + 4 \sum_{i<j} \xi_i\xi_j{\del\over\del\xi_j}
 {\del\over\del\xi_i} - \cdots
 \end{gathered}
 \ee
 If $X$ has arbitrary dimension $n|m$, then $(-1)^\F: \Oc_X \to \Oc_X$
 is represented by the
 expression \eqref{eq:-1^F} 
   in any local
 coordinate system $(z_1,\cdots, z_n, \xi_1,\cdots, \xi_m)$.   
 In the case of a vector bundle $E$, we  represent $E$ locally as
 $\Oc_X\otimes V$ for a super-vector space $V$ and   tensor
   \eqref{eq:-1^F}   with the operator  $(-1)^\F_V$  acting on $V$.

   \paragraph{Penkov's equivalence.}  Note that $X_\red$ an be seen as a supermanifold of dimension $n|0$
  equipped with a canonical closed embedding  of supermanifolds $\eps: X_\red\hookrightarrow X$.  In \cite{penkov} I. B. Penkov used the direct and inverse
  images under $\eps$ to identify the categories of $\Dc$-modules on $X_\red$ and
  $X$. More precisely, we denote
  \[
    \Dc_\to = \Oc_{X_\red}\otimes_{\Oc_X}\Dc_X, 
    \quad \Dc_\leftarrow = \omega_{X_\red}\otimes _{\Oc_X} \Dc_\to \otimes _{\Oc_X} \omega_X^*.
       \]
  Then $\Dc_\to$ is a $(\Dc_{X_\red}, \Dc_X)$-bimodule and
  $\Dc_\leftarrow$ is a $(\Dc_X, \Dc_{X_\red})$-bimodule. 
  
  \begin{thm}[\cite{penkov}] \label{thm:penkov}
  The functors
  \[
  \eps_*: \Mc \mapsto \Mc \otimes_{\Dc_{X_\red}} \Dc_\leftarrow, 
  \quad \eps^*: \Nc \mapsto \Dc_\to \otimes_{\Dc_X} \Nc
  \]
  define mutually quasi-inverse equivalences between the categories
  of quasicoherent left $\Dc_{X_\red}$-modules $\Mc$  and  of quasicoherent 
  left $\Dc_X$-modules $\Nc$. \qed
  \end{thm}
  
 For future reference recall  the isomorphisms implying Theorem 
 \ref{thm:penkov}:
 \be
 \Dc_\leftarrow \otimes_{\Dc_{X_\red}} \Dc_\to\,\, \= \,\,\Dc_X, \quad
 \Dc_\to \otimes_{\Dc_X} \Dc_\leftarrow \,\,\= \,\,\Dc_{X_\red}.
 \ee
 
 \paragraph{Operators of infinite order.}  Let $X$ be a supermanifold of dimension $n|m$ and $E,F$ be
 two vector bundles on $X$. The sheaf $\Dc^\oo(E,F) = \Dc_X^\oo(E,F)$
 of differential operators of infinite order on $X_\red$  can be defined
 in one of the three equivalent ways:
 
 \vskip .2cm
 
 \noindent{\sl (i) Cohomological defininion.}  Assuming that $\dim X = n|m$, 
 \[
 \Dc_X^\oo(E,F) \,=\, 
 \ul H^n_\Delta \bigl(X_\red\times X_\red, (\omega_X\otimes E^*)\boxtimes F\bigr). 
 \]
 
 \noindent{\sl (ii) Coordinate definition.}  A section $P\in\Dc_X^\oo(E,F)$ is 
 a datum , for any local  coordinate system $(z_1,\cdots, z_n, \xi_1,
 \cdots, \xi_m)$ on $X$ and any local trivializations $E\to \Oc_X^{p|q}$,
 $F\to \Oc_X^{r|s}$, of a series
  \[
    \sum_{I,J,K} a_{IJK}(z)\xi^I \del_\xi^J \del_z^K, \quad a_{IJK}(z)\in
    \Hom_\CC(\CC^{p|q}, \CC^{r|s}) \otimes \Oc_{\CC^n}.
    \] 
    Such series are required to satisfy the overconvergence condition:
    \[
    \forall I,J \text{ the series } \sum_K k_1 ! \cdots k_n ! \,\,  a_{IJK}(z) p_1^{k_1} \cdots p_n^{k_n}
    \text{ converges } \forall (p_1,\cdots, p_n)\in\CC^n
    \]
and be connected  with each other via standard formulas for transformations
of total symbols, see, e.g., \cite{SKK}. 

\vskip .2cm

\noindent{\sl (iii) Reductionist definition.}  
\[
\Dc^\oo_X (E,F) = \Dc^\oo_{X_\red}(E,F),
\]
where on the right we consider $E,F$ as differential bundles on $X_\red$. 

\vskip .2cm

We omit the straightforward proof of equivalence of these three definitions.   
     \vfill\eject
     
    \subsection{(Even-style) exponentials of odd vector fields} \label{subsec:exp-odd-vect}
    
    \paragraph{Odd vector fields.}
    Let $X$ be a supermanifold and $v\in T_{X, \1}$ be an odd vector field on $X$.
    Denote by $D_v\in\Dc_X$ the corresponding differential operator.
    
    \begin{prop}\label{prop:e^D_v}
    The series $e^{D_v} = \sum_{i=0}^\oo D_v^i /i!$ defines a section of $\Dc_X^\oo$. 
    \end{prop}
    
    Note that $e^{D_v}$ is {\em not}  an instance of  the familiar procedure exponentiating a Lie superalgebra
    to a Lie supergroup. In that procedure, one multiplies an odd vector field
    by an odd (square-zero) coefficient which results in exponential series terminating.
    Here, we do not do that and the series is infinite. 
    
    \vskip .2cm
    
    \noindent{\sl Proof:} It is enough to assume that $X$ is a domain in $\CC^{n|m}$
    with standard coordinates $z_1,\cdots, z_n$, $\xi_1,\cdots, \xi_m$,
    so $\Oc_{\CC^{n|m}} = \Oc_{\CC^n}\otimes\Lambda$, where we denote
    $\Lambda = \Lambda[\xi_1,\cdots, \xi_m]$. Then we can view $D_v$
    as an   $\End_\CC(\Lambda)$-valued (or, equivalently, $2^m\times 2^m$ matrix)
    differential operator on $\CC^n$ or, put in yet another way,
    as a
section of $\Dc_{\CC^n} (\Oc_{\CC^n}\otimes\Lambda, \Oc_{\CC^n}\otimes\Lambda)$
    and it is enough to show that $e^{D_v}$ understood in this sense,
    lies in $\Dc^\oo_{\CC^n} (\Oc_{\CC^n}\otimes\Lambda, \Oc_{\CC^n}\otimes\Lambda)$.     
    
 For this, note that the square $D_v^2 = (1/2) [D_v, D_v]$ us equal to
    $(1/2) D_{[v,v]}$ where $[v,v] $ is the self (super)commutator of $v$.  
    As $[v,v]\in T_{X, \0}$  is an even vector field, $D_{[v,v]}$ considered as an
    $\End(\Lambda)$-valued differential operator on $\CC^n$, is of order 1. 
    Thus $D_v^2$ is of order 1, so $\eord(D_v)\leq 1/2$ and our statement
    follows from Proposition \ref{prop:eord<1}. \qed
    
    \paragraph{Equivariant bundles.} Let $\GG$ be a (complex) Lie supergroup
    with Lie superalgebra  $\gen = \gen_\0\oplus \gen_\1$.  Let $X$ be a supermanifold
    wih a $\GG$-action and $E$ be a $\GG$-equivariant vector bundle on $X$.
    Each $y\in \gen_1$ given then a differential operator $D_y\in\Dc_X(E,E)$. 
    
    \begin{prop}
    The series $e^{D_y}$ defines a section of $\Dc_X^\oo(E,E)$. 
    \end{prop}
    
    \noindent{\sl Proof:} Smilar to that of Proposition \ref{prop:e^D_v}, as 
    $2D_y^2 = [D_y, D_y]  = 
    D_{[y,y]}$ is an operator of order 1 and so $\eord(D_y) \leq 1/2$.  \qed
    
    \paragraph{The infinite order completion of the enveloping algebra.}
    Let $\GG, \gen$ be as before and $U(\gen)$ be the universal
    enveloping (super)algebra of $\gen$. As well known, $U(\gen)$ can be
    identified with $\Dc(\GG)^\GG$, the algebra of left-invariant differential
    operators on $\GG$. So we denote
    \be
    U^\oo(\gen) := \Dc^\oo(\GG)^\GG
    \ee
    and call the {\em infinite order completion} of $U(\gen)$ the algebra
    formed by left-invariant operators of infinite order on $\GG$.
     It is easy to characterize elements of
    $U^\oo(\gen)$ directly, as infinite series of summands in the Poincar\'e-Birkhoff-Witt
    decomposition of $U(\gen)$ satisfying the overconvergence condition
    similar to that used in \eqref{eq:overconv}.  Proposition \ref{prop:e^D_v}
    implies:
    
    \begin{prop}
    Let $y\in\gen_\1$ be any odd element. Then $e^y = \sum_{i=0}^\oo y^i/i!$ is a well
    defined element of $U^\oo(\gen)$. \qed
    \end{prop}
    
    Note that for  an even element $y\in\gen_\0$ the series $e^y$ is never in $U^\oo(\gen)$ since it
    represents the shift operator on $\GG$. 
    
        \vfill\eject

 \subsection {Twisted bundles and differential operators.} \label{subsec:twisted}
 
 Modular forms coming from theta functions are often of half-integer weight, that is,
 they are sections of fractional (half-integer) powers of basic geometric line bundles. 
 One,  now standard way of handling such ``fractional bundles'' 
 (which may not exist globally or, if exist, not have
 expected equivariance properties) is by using the concept of twisted rings
 of differential operators 
 \cite[\S 2]{BB-jantzen}.  Operationally, it is also convenient to speak about such
 bundles directly, as objects of certain gerbes. So we recall this formalism
 referring to \cite{breen, brylinski-loop} for additional background. 
 
 \paragraph{Stacks and gerbes.} 
 
 Let $S$ be a topological space. Recall the concept of a {\em stack}  
  (``categorical analog of a sheaf'') of categories on $S$. Thus a stack
  $\Gc$ associated to any open $U\subset S$ a category $\Gc(U)$
  with restriction functors $\Gc(U)\to \Gc(U')$ for  $U'\subset U$
  satisfying the natural axioms of descent. We will be only interested
  in stacks of groupoids (categories with all morphisms isomoprhisms).
  
  \vskip .2cm
  
   Let $\Gamma$ be a discrete abelian group. 
 By n  $\Gamma$-{\em stack} we will
 mean a stack  $\Gc$ of groupoids on $S$ with $\Gamma$ 
 acting freely on any  nonempty $\Hom_{\Gc(U)}(x,y)$ such that composition of morphisms
 is bimultiplicative. 
 
   \vskip .2cm
 
 A $\Gamma$-{\em gerbe} is a $\Gamma$-stack $\Gc$  such that for sufficiently
 small $U$ the groupoid $\Gc(U)$ is nonempty and the $\Gamma$-action
 on any $\Hom_{\Gc(U)}(x,y)$ is transitive (making it into a $\Gamma$-torsor). 
 A $\Gamma$-gerbe $\Gc$ gives rise to a class $[\Gc]\in H^2(S, \Gamma)$
 (obstruction to $\Gc(S)$ being nonempty). 
 
 \vskip .2cm
 
 We will be particularly interested in $\Gamma=\CC^*$. 
 
  \paragraph{Twisted  sheaves and bundles.} 
  Let $X$ be a complex supermanifold. We recall that the
   concepts
  of topology, sheaves, gerbes etc. for $X$ refer to the
  space $S=X_\red$. By $\Sh_X$ we denote the stack of (super,
  i.e., $\ZZ/2$-graded) sheaves  of $\CC$-vector spaces
  on   $X$ and their isomorphisms. 
  We denote by $\Bun_X$ the stack formed by holomorphic
  (super) vector bundles on $X$ and their isomorphisms. 
  These are $\CC^*$-stacks with the 
   action of $\CC^*$ given by composing with scalar automorphisms on
 source or target. We have a $1$-morphism of $\CC^*$-stacks
 $\Bun_X\to \Sh_X$ given by taking the sheaf of holomorphic
 sections. 
  
    \vskip .2cm

  Let $\Gc$ be a $\CC^*$-gerbe on $X$. A $\Gc$-{\em twisted sheaf}
  resp.  $\Gc$-{\em twisted bundle}
  is a  ($1$-)morphism  of $\CC^*$-stacks
  $\Fc: \Gc\to\Sh_X$, resp.  $E: \Gc\to\Bun_X$.
  That is, a $\Gc$-twisted bundle gives,  for any open $U\subset X$ and any $x\in \Ob\,  \Gc(U)$
 a bundle  $E(x)$ on $U$ and for any morphism
  $\phi: x\to y$ in $\Gc(U)$ we have an isomorphism of vector bundles
  $E(\phi): E(x)\to E(y)$  such that
   $E(\lambda\phi) = \lambda E(\phi)$, $\lambda\in\CC^*$. 
   The bundles $E(x)$ are called {\em determinations} of $E$. Any two
   determinations are locally  identified by an isomorphism unique up to
   a constant factor. Similarly for a twisted sheaf. 
   
     \vskip .2cm
   
   Given two $\CC^*$-gerbes $\Gc, \Gc'$, we have their
   {\em tensor product} gerbe  $\Gc\otimes\Gc'$
   which is obtained by the ``stackification'' of the pre-stack of groupoids
  $\Gc \boxtimes \Gc'$ with 
  \[
   \begin{gathered}
   \Ob(\Gc\boxtimes\Gc' (U)) = 
   \Ob(\Gc(U)) \times\Ob(\Gc'(U)),
   \\
   \Hom_{(\Gc\boxtimes\Gc')(U)}( (x,x') , (y,y')) = \Hom_{\Gc(U)}(x,y)
   \otimes \Hom_{\Gc';U)}(x', y')
   \end{gathered} 
   \] 
   (tensor product of $\CC^*$-torsors). We also denote by 
   $\Gc^\vee = \Gc^\op$ the {\em dual}, or  {\em opposite stack} of 
   $\Gc$,
   i.e., the stack formed by the opposite categories
   $\Hom_\Gc(x,y)^\op$. 
   
   \vskip .2cm
   
   Given a $\Gc$-twisted bundle $E$ and a $\Gc'$-twisted bundle
   $E'$, we have a $\Gc\otimes\Gc'$-twisted bundle $E\otimes E'$ with
   determinations $(E\otimes E')(x,x') = E(x) \otimes E'(x')$. 
   We also have the $\Gc^\vee$-twisted {\em dual bundle} $E^\vee$
   with determinations $E^\vee(x) = E(x)^\vee$. 
   
   \begin{exas}\label{ex:frac-powers}
   (a) A vector bundle $E$ in the ordinary sense is a $\Gc$-twisted bundle
   for the trivial $\CC^*$-gerbe $\Gc = \pt^{\CC^*}$ with  each 
   $\pt^{\CC^*}(U)$
   having one object $\pt$ with automorphism group $\CC^*$.    
   
   \vskip .2cm 
   
   (b) Let $L$ be a line bundle on $X$ and $N\geq 1$.  All $N$th
   roots of $L$ form a twisted line bundle 
  $L^{1/N}$.  Its  underlying $\CC^*$-gerbe $\Gc$ (sometimes
  also denoted $L^{1/N}$)
  associates to an open $U\subset X$ the groupoid whose objects
  are pairs $(M, \phi)$ where $M$ is a line bundle on $U$ and 
  $\phi: M^{\otimes N}
  \to L|_U$ is an isomorphism. A morphism $(M, \phi) \to (M', \phi')$ is
  a morphism $\psi: M\to M$  such that $\phi'\circ \psi^{\otimes N} = \lambda \phi$
  for some $\lambda\in\CC^*$. The functor $\Gc\to\Bun_X$
  is the forgetful functor $(M,s) \mapsto M$. 
  
  \vskip .2cm
  
  (c) In particular, we can speak  about the  tensor product $E\otimes L^{1/N}$
  where $E$ is a vector bundle and $L$ is a line bundle on $X$.
  
   \end{exas}
   
   \vskip .2cm
  
  \begin{exas}\label{exas:ourTDO}
  (a) 
  Let $K$ be a Lie supergroup and $p: Y\to X$
  be a principal homogeneous  right $K$-space, i.e., a right $K$-torsor. Let 
   \[
 1\to H \lra \wt K \lra K\to 1
 \]
 be a central extension of Lie supergroups with $H$ being
 a  discrete (purely even) group. 
  We have then the $H$-gerbe $\delta(Y)$ on $X$
  whose objects over $U$ are lifts $\wt Y_U \to U$
  of $p^{-1}(U)$ to a $\wt K$-torsor
  over $U$. Denoting by $\ul K, \ul H$ the sheaves on $X$ formed by  
  holomorphic maps with values in $K$ or $H$ (since $H$ is discrete,
  the latter sheaf is constant), the class 
  $[\delta(Y)]\in H^2(X,H) = H^2(X, \ul H)$ is the image of
  the class of $Y$ in $H^1(X, \ul K)$ under the coboundary
  map $\delta: H^1(X,\ul K)\to H^2(X, \ul H)$. 
  Given a character $\sigma H\to \CC^*$, we form the
  {\em pushforward $\CC^*$-gerbe} $\Gc_\sigma := \sigma_* \delta(Y)$ on $X$
  with same objects $\wt Y_U$ as $\delta(Y)$ and
  \[
  \Hom_{\chi_* \delta(Y)}(\wt Y_U, \wt Y'_U) \,=\,
 \sigma_*  \Hom_{\delta(Y)}(\wt Y_U, \wt Y'_U) 
  \]
  (the pushforward of  torsors under $\sigma$). 
  
  \vskip .2cm
  
  (b) In particular (cf. \cite[\S2.5]{BB-jantzen}), let 
  $K=T= (\CC^*)^d$ be an algebraic torus,
  with Lie algebra $\ten$ and let  $s:\ten\to\CC$
 be a  character  (linear functional). We have then the exponential 
 extension
 \[
 1\to H=\Ker(\exp) \lra \ten \buildrel \exp\over\lra T\to 1
 \]
 and the character $\sigma = e^s|_H: H\to\CC^*$. 
 Given a principal $T$-bundle
   $p: Y\to X$, we have the $\CC^*$-gerbe $\chi_*\delta(Y) = 
   (e^s)_*\delta (Y)$ on $X$. We further have the
   $\sigma_*\delta(Y)$-twisted line bundle (equivalently,
   a principal homogeneous $\CC^*$-space)  $Y^s$ which associates
   to any $\ten$-torsor $\wt Y_U\to U$ the 
    pushforward under $e^s: \ten\to\CC^*$
   of  $\wt Y_U$. 
   
   \vskip .2cm
 
 (c)  If $T=\CC^*$, then $s$ can be  viewed as a complex number.
  If $Y$ comes from a line bundle $L$ on $X$
 by removing the zero section, 
  we get the twisted bundle $L^s$, $s\in \CC$,  called the $s$th
  {\em complex power} of $L$. This generalizes the case $s=1/N$
  from Example \ref {ex:frac-powers}(b). 
  \end{exas}

 \paragraph{Morphisms and differential operators.} 
 Let $\Gc_i$ be a $\CC^*$-gerbe on $X$ and $E_i$ be
 a $\Gc_i$-twisted bundle,   $i=1,2$.  Then we have the
 $\Gc_1^\vee\otimes \Gc_2$-twisted sheaves on $X$ 
 \be\label{eq:twisted-mDO}
 \ul\Hom(E_1, E_2), \quad \Dc(E_1, E_2), \quad \Dc^\oo(E_1, E_2)
 \ee
 of $\Oc$-morphisms, differential operators and differential
 operators of infinite order from $E_1$ to $E_2$ respectively. 
Their values on an object $(x_1, x_2)\in \Ob \,(\Gc_1\boxtimes\Gc_2)(U))$ are, respectively,
 \[
 \ul\Hom(E_1(x_1), E_2(x_2), \quad \Dc(E_1(x_1), E_2(x_2) ), \quad \Dc^\oo(E_1(x_1), E_2(x_2))
 \]
 If $\Gc_1 = \Gc_2 =\Gc$, then $\Gc_1^\vee\otimes\Gc_2$ is
 equivalent to the trivial gerbe and we can take $x_1=x_2$
 throughout and arrive at 
 
 \begin{defi}
 Let $E_i$, $i=1,2$,  be $\Gc$-twisted bundles on $X$.
  A {\em morphism} 
 (resp. {\em differential operator}, resp. {\em differential operator
 of infinite order})  $A: E_1\to E_2$ is a system of
 morphisms of vector bundles (resp. differential operators,
 resp. differential operators of infinite order) $A_x: E_1(x) \to E_2(x)$
 given for any open $U\subset X$ and any $x\in \Ob\, \Gc(U)$
 compatible with restrictions to $U'\subset U$ and such that for
 any $\phi: x\to y$ we have $A_y E_1(\phi) = E_2(\phi) A_x$. 
  \end{defi}
  
  In this case \eqref{eq:twisted-mDO} defines actual (non-twisted)
  sheaves. 
  
  \vskip .2cm
  
  \begin{ex}\label{ex:TDO} [Twisted rings of differential operators]
  Let $L$ be a line bundle on $X$ and $s\in \CC$. The sheaf 
  $\Dc(L^s) = \Dc(L^s, L^s)$ is an example of a {\em twisted
  ring of differential operators} in the sense of \cite[\S2.5]{BB-jantzen}. 
  For any actual (non-twisted) bundle $E$ on $X$ the tensor
  product $E\otimes L^s$ as well as $L^s$ are twisted with
  respect to the same $\CC^*$-torsor, so we have a locally free sheaf
  $\Dc(E\otimes L^s, L^s)$ of left $\Dc(L^s)$-modules. For any
  two  such bundles $E_1, E_2$ we have
  \[
  \ul\Hom_{\Dc(L^s)} \bigl( \Dc(E_1\otimes L^s, L^s), \Dc(E_2\otimes L^s, L^s)
  \bigr) \, \= \, \Dc (E_2^\vee\otimes L^{-s}, E_1^\vee\otimes L^{-s}).
  \]
  \end{ex}

          \vfill\eject
          
   \section{Genus 1: $\OSp(1|2)$ acting on $\PP^{1|1}$}\label{sec:gen1}

   \subsection {$\PP^{1|1}$ with its SUSY structure}\label{subsec:P11}
   
   \paragraph{SUSY structures on supermanifolds.} Let $X$ be a (complex analytic) supermanifold of
   dimension $n|m$.  Its tangent bundle $TX$ is a (super) vector bundle on $X$ of rank $n|m$. 
   In the sequel we will drop the adjective ``super''  unless essential for understanding. 
   Any distribution of subspaces in $TX$, i.e., a vector  subbundle
   $\Tc \subset TX$ gives rise to the {\em Frobenius pairing} 
     \[
 \phi_\Tc:   \Lambda^2_{\Oc_X} \Tc \lra TX/\Tc
   \]
   induced by the Lie bracket of sections of $TX$. As well known, vanishing of $\phi_\Tc$ is
   equivalent to the fact that $\Tc$ is integrable to a foliation (Frobenius theorem).

   \begin{defi}\label{def:SUSYstr}
   By a {\em SUSY structure} on $X$ we will mean a vector  subbundle
   $\Tc \subset TX$ of rank $0|m$ such that  $\phi_\Tc$ is a surjective morphism of vector bundles
   (i.e., it is locally a projection onto a direct summand).  A supermanifold with a SUSY strucrure
   will be called a {\em SUSY manifold}. 
   \end{defi} 
   
   In this case the restriction $\Pi\Tc|_{X_\red}$
   is an ordinary vector bundle on the complex manifold $X_\red$, equipped with a surjective
   morphism 
   \[
   \phi_\Tc^\red: S^2_{\Oc_{X_\red}} (\Pi\Tc|_{X_\red})  \lra TX_\red.
   \]

   \begin{rems}
  (a)  In examples motivated  by physics it is often  assumed that $X_\red$ is equipped with a Riemannian metric
   or a conformal structure,  $\Pi\Tc|_{X_\red}$ is a direct sum of some copies of  the  spinor bundle(s),
    and $\phi_\Tc^\red$
   is built out of the spinor pairings (gamma-matrices). 
   
   Thus,  the expression ``$\Nc=p$ supersymmetry''
    means that 
   $\Pi\Tc|_{X_\red}= S^{\oplus p}$ is the direct sum of $p$ copies of the  spinor bundle 
 (in dimension when it is unique)  while
   ``$\Nc=(p,q)$ supersymmetry''  means that $\Pi\Tc|_{X_\red}= S_+^{\oplus p} \oplus S_-^{\oplus q}$
 (in dimension  when there are two spinor bundles 
   $S_\pm$).

    For us, it is convenient to adopt the  definition
   above 
   as our eventual examples will be more general than those coming from  usual conformal structures. 
   
   \vskip .2cm
   
   (b)    Non-triviality of $\phi_\Tc$, i.e., non-integrability of $\Tc$ allows us to view a SUSY structure
    $\Tc$ as a kind of
   ``odd contact structure'' on $X$. 
   \end{rems}

   %\vfill\eject
   
   \paragraph{$\PP^{1|1}$ as a SUSY curve and its automorphism group $\OSp(1|2)(\CC)/\{\pm 1\}$.} 
   \label{par:P11-susy} 
    SUSY curves of dimension $1|1$ are important in
   superstring theory, where they are often called simply super Riemann surfaces \cite{witten}
   or algebraic supercurves \cite{felder} which may cause confusion. Our terminology,
   following that of \cite{manin-SUSY, voronov},   avoids such confusion. 
   
   \vskip .2cm
   
   For a SUSY curve $X$ the line bundle $\Pi\Tc^*|_{X_\red}$ on $X_\red$ is a theta-characteristic
   (square root of  the canonical class), as $\phi_\Tc^\red: (\Pi\Tc|_{X_\red})^{\otimes 2} \to TX_\red$
   is an isomorphism. 
   
   \vskip .2cm
   
   The simplest SUSY curve, that of genus $0$, is the projective superline $\PP^{1|1}$. 
   We recall \cite{manin} that the projective superspace $\PP^{n|m}$ parametrizes
   $1|0$-dimensional subspaces in the linear  superspace $\CC^{n+1|m}$. 
   The automorphism group of $\PP^{n|m}$ as of  a supermanifold is the Lie supergroup
   \[
   \Aut_{\on{Supermfld}} (\PP^{n|m}) = \PGL(n+1|m)(\CC) = \GL(n+1|m)(\CC)/\CC^*,
   \]
   where $\GL(n+1|m)(\CC)$ is the  supergroup formed by automorphisms of $\CC^{n+1|m}$
   as a linear superspace. In general, $\PP^{n|m}$ does not carry a natural SUSY structure.
   However, $\PP^{1|1}$ does carry one which is  suggested by the fact that $T\PP^1 = \Oc_{\PP^1}(2)$
   has a natural square root, namely $\Oc_{\PP^1}(1)$. 
   
   To avoid confusion between the super and non-super cases, we denote
   \be
   \Lc_r = \Oc_{\PP^1}(-r), \quad S\Lc_r = \Oc_{\PP^{1|1}}(-r), \quad r\in \ZZ
   \ee
   and write $\Lc= \Lc_1 = \Oc_{\PP^1}(-1)$ and similarly $S\Lc= S\Lc_{1}$.
   Thus $S\Lc$ is the  tautological rank $1|0$ bundle on $\PP^{1|1}$, embedded
   as a subbundle into $\wt Y = \Oc_{\PP^{1|1}}\otimes Y$.  
   
   \vskip .2cm
   
   To construct the SUSY structure, let us make  $Y=\CC^{2|1}= \Pi\CC \oplus \CC^2 $ into a symplectic (super)   vector space 
   with the super-symplectic form $\eta$  combining
   a symplectic form on $\CC^2$ and a symmetric  bilinear form on $\CC$, i.e., by putting
   \be\label{eq:eta-gen1}
   \eta\bigl( (\theta, u, v), (\theta', u', v',\bigr) \,=\, uv'-vu' + \theta\cdot\theta', 
   \ee
   where $(u,v)$ are linear coordinates on $\CC^2$ and $\theta$ is a linear coordinate on $\CC$. 
   After this, $\Tc\subset TY = T\PP^{1|1}$ appears as a direct analog of the classical
   contact structure present on the projectivisation of any symplectic vector space. 
   
   \vskip .2cm
   
   That is,   all points of $\PP^{1|1}$ 
   (we consider points with values in any commutative superalgebra) are isotropic with
   respect to $\eta$ and so we have 
    we have $S\Lc\subset(S\Lc)^\perp$, where $(S\Lc)^\perp\subset \wt Y$ is the 
   $\eta$-orthogonal
   to $S\Lc$,  of  rank $1|1$.    Now, $T\PP^{1|1} \= \Hom(S\Lc, \wt Y/S\Lc)$ (the so-called
   Euler isomorphism, see \cite[Ex.4.3.12]{manin}) and we put 
   $\Tc = \Hom(S\Lc, (S\Lc)^\perp/S\Lc)$ thus defining the SUSY structure. 
   
   \vskip .2cm
   
   Let 
   \be
  \GG:=  \OSp(1|2)(\CC)\subset \GL(2|1)(\CC)
   \ee
   be the complex  Lie supergroup formed by automorphisms of $\CC^{2|1}$
   preserving $\eta$. Here, to avoid writing $\on{SpO}$, it is convenient tor reverse the order and put the
   odd compoment first.  The underlying Lie group of $\GG$ is $\GG_\red = \on{O}(1)(\CC) \times \Sp(2)(\CC)$,
   where  $\on{O}(1)(\CC)=\{\pm 1\}$ is disconnected and $\Sp(2)(\CC) = \SL(2, \CC)$. 
   The intersection of $\GG$ with $\CC^*\subset \GL(2|1)(\CC)$ is equal to $\{\pm 1\}$. 
   Accordingly, for the automorphism group of $\PP^{1|1}$ as a SUSY curve we find
   \be
   \Aut_{\on{SUSY}}(\PP^{1|1}, \Tc) \,=\, \GG/\{\pm 1\}. 
   \ee
   Putting the odd part first, we represent points of $\CC^{2|1}$ as a supermanifold,
  e.g.,  as  the  superscheme $\AAA^{2|1}= \Spec\,  S^\bullet ((\CC^{2|1})^*)$, 
   by  column vectors  $(\theta, u,v)^t$, similarly to \eqref{eq:eta-gen1}
    with $\theta$ being an odd variable and $u,v$ being even variables. 
    We denote by $H$ the 
 Gram matrix of $\eta$ in the  coordinate system just described: 
\be\label{Gram-H=gen1}
 H = \begin{pmatrix} 1&0&0
 \\
 0&0&1
 \\ 0&-1&0
 \end{pmatrix}.  
 \ee
 Thus,  points of  $\GG$ (with values in any commutative $\CC$-superalgebra)  
 are represented by matrices
  \be\label{eq:OSp-matrix:gen1}
 g=\begin{pmatrix}
 s&\alpha&\beta 
 \\
 \gamma& a & b
 \\
 \delta & c & d
 \end{pmatrix}
 \ee
 with Latin variables being even, Greek ones being odd,   subject to the equations
 $g^T H g = H$. Here $g^T$ is the {\em supertranspose} of $g$, defined  for supermatrices  as
 \be\label{eq:super-trans}
 \begin{pmatrix} A_{\0\0}& A_{\0\1}\\A_{\1\0}  & A_{\1\1}
 \end{pmatrix}^T = \begin{pmatrix} A_{\0\0}^t & -A_{\1\0}^t\\ A_{\0\1}^t & A_{\1\1}^t
 \end{pmatrix}. 
 \ee
 
% \vfill\eject
 
 \paragraph{The super analog of the  Lobachevsky plane.}\label{par:super-lobach}
     Points
 of $\PP^{1|1}$ can be represented in
 homogeneous
 coordinates as $(\theta, u,v)^t$, i.e., by
 points of $\CC^{2|1}$,  see \S \ref{par:P11-susy}.  Consider the affine chart $\AAA^{1|1}\subset \PP^{1|1}$
  given by the condition of $v$ being invertible.  It has inhomogeneous coordinates
 $\tau = u/v$, even and  $\xi = \theta/v$, odd. 
 Consider the odd vector field 
  \be\label{eq:D-xi}
 D_\xi = {\del\over\del\xi } - \xi {\del\over\del\tau}
 \ee
 on $\AAA^{1|1}$, known as the {\em spinor derivation}. The following is well known. 
 
  \begin{prop}\label{prop:D-xi-gen1}
 The SUSY structure $\Tc$ restricted to $\AAA^{1|1}$ has the
 form
 $
 \Tc|_{\AAA^{1|1}} \,=\, \Oc_{\AAA^{1|1}}\cdot D_\xi$.  
  \end{prop}
  
  The surjectivity (in our case, isomorphicity) of the Frobenius pairing $\phi_\Tc$  comes from  the fact that for the super-commutator
  we have 
 $[D_\xi, D_\xi] = \del/\del\tau$. 
 
\vskip .2cm
  
  \noindent{\sl Proof of the proposition:} By definition, 
  $\Tc= \Hom(S\Lc, (S\Lc)^\perp/S\Lc)$ on all of $\PP^{1|1}$.
  Consider a  point $l= (\theta,  u,v)^t = (\xi, \tau, 1)^t$ of $\AAA^{1|1}$. The value of the vector
  field $D_\xi$ at this point is represented by the  odd vector $(1, -\xi, 0)^t$ which is orthogonal to $l$
  with respect to $\eta$ and so spans the $(0|1)$-dimensional space  $l^\perp/l$. \qed

  \vskip .2 cm
   
 The action of $g\in \GG$
 as in \eqref{eq:OSp-matrix:gen1}  in  coordinates $(\tau, \xi)$
 is obtained by normalizing the action of $g$  on the column vector
 $(\theta, u,v)^t$  which gives
 \be\label{eq:osp-act-1}
 g\cdot (\tau, \xi) = \biggl( {\gamma\xi + a\tau + b\over \delta\xi + c\tau + d}, \, 
 {s\xi+ \alpha \tau+\beta \over \delta \xi + c\tau + d}\biggr). 
 \ee
 
 \vskip .2cm 
 
 Let $\Hen = \{ \tau\in \CC| \, \Im(\tau)>0\}\subset \CC$ be the Lobachevsky upperhalf plane and 
  $S\Hen\subset \CC^{1|1}$ the super-thickening of $\Hen$ induced from $\CC^{1|1}$. 
  We call $S\Hen$ the  {\em super Lobachevsky plane},   cf. \cite{baranov}. The action 
  \eqref{eq:osp-act-1} restricted to the real
   Lie supergroup
  \[
 \GG_\RR :=  \OSp(1|2)(\RR)\subset  \GG= \OSp(1|2)(\CC)
  \]
  preserves $S\Hen$. 
  
  \paragraph{Forms and superforms.} 
  
  For $r\in \ZZ$ we denote
  \be
  \F_r = H^0(\Hen, \Lc_r), \quad \SF_r = H^0(S\Hen, S\Lc_r)
  \ee
  and call elements of these spaces {\em forms} (resp. {\em super-forms}) of
  weight $r$ 
   (not necessarily modular or automorphic). Thus $\SF_r$ is a super-vector space
   acted upon by $\GG_\RR$ and $\F_r$ is acted upon by $\SL(2, \RR)$ which is the
   unit connected component of $(\GG_\RR)_\red$. The even part $\SF_{r, \0}$ is
   identified with $\F_r$.

   \iffalse
   with odd and even components  identified as follows
   \[
   \SF_r = \SF_{r, \0} \oplus \SF_{r, \1}, \quad \SF_{r, \0} = \F_r, \quad \SF_{r, \1} = \F_{r-1}.
   \] 
   \fi
  
  \vskip .2cm
  
  The restriction $S\Lc_r$ to $\AAA^{1|1}$ and therefore to $S\Hen$ is canonically trivialized:
  local sections of $S\Lc_r$ are represented by locally defined functions in $(\theta, u,v)$ homogeneous of
  degree $-r$. So putting $(\theta,u,v) = (\xi, \tau,1)$, we identify them with locally defined functions
  $f(\tau, \xi) = f_0(\tau) + f_1(\tau)\xi$. 
   In particular, any $f\in \SF_r$ is represented  in this form 
  where $f_0, f_1$ are holomorphic functions on $\Hen$.

  The action of $g\in \GG_\RR$
   in the form \eqref{eq:OSp-matrix:gen1} (considered as a point with values in some commutative $\RR$-superalgebra)  on  $f\in \SF_r$ as above has the form
   \be\label{eq:osp-2-forms}
   (g^*f) (\tau, \xi) = (\delta\xi + c\tau + d)^{-r}\,  f \biggl( {\gamma\xi + a\tau + b\over \delta\xi + c\tau + d}, \, 
 {s\xi+ \alpha \tau+\beta \over \delta \xi + c\tau + d}\biggr). 
   \ee
   This implies that we can identify
   \be\label{eq:SF-r}
   \SF_{r,\0} \= \F_r, \quad \SF_{r, \1} \= \F_{r-1}, \quad r\in \ZZ
   \ee
   as $\SL(2,\RR)$-modules. 
   
  % \vfill\eject
   
   \paragraph{The gerbe $\Gc_{1/2}$.}\label{par:g-1/2-p11}
    The line bundle $\Lc_1=\Oc_{\PP^1}(-1)$ on $\PP^1$
   does not have a globally defined square root and we denote by $\Gc_{1/2}$ the $\CC^*$-gerbe
   of local determinations of such a root.  That is,  for an open $U\subset \PP^1$
   the objects of the groupoid $\Gc_{1/2}(U)$ are pairs $(M,\phi)$ where $M$ is a line bundle
   over $U$ and $\phi: M^{\otimes 2}\to \Lc_1|_U$ is an isomorphism, 
    see Example  \ref{ex:frac-powers}(b).  So $\Gc_{1/2}(U)\neq\emptyset$
     for any $U\neq \PP^1$ but $\Gc_{1/2}$ is globally nontrivial. 
   The square root $\Lc_{1/2} = \Lc_1^{\otimes{1\over 2}}$ is then well defined as a $\Gc_{1/2}$-twisted
   line bundle, see \S \ref {subsec:twisted}. Further, for any $m\in \ZZ$ we have the
   $\Gc_{1/2}$-twisted line bundle $\Lc_{m+{1\over 2}} = \Lc_m\otimes\Lc_{1\over 2}$, see
   Example  \ref{ex:frac-powers}(c). 
   
   \vskip  .2cm
   
   Recall that $\PP^{1|1}$ is a ringed space with the same underlying space as $\PP^1$. 
   With this understanding, note that  the categories of determinations of $\Lc_1^{\otimes{1\over 2}}$
   and of $S\Lc_1^{\otimes{1\over 2}}$ over any $U\subset \PP^1$ are in equivalence,
   the claim reducing to
 the following obvious fact. Square roots of 
  a nonvanishing  even  function $f\in \Oc_{\PP^{1|1}}(U)_\0^*$ are in bijection 
     with square  roots of the restriction $f|_{\PP^1}\in \Oc_{\PP^1}(U)^*$. 
     Therefore we have the $\Gc_{1/2}$-twisted line bundles $S\Lc_r$, $r\in {1\over 2} + \ZZ$,
     on $\PP^{1|1}$,   extending $\Lc_r$. 
     
     \vskip .2cm
     
     The gerbe $\Gc_{1/2}$ is equivariant under the group $\SL(2,\CC)$ acting on $\PP^1$
     and each $\Lc_r$, $r\in {1\over 2}+\ZZ$, is an $\SL(2,\CC)$-equivariant $\Gc_{1/2}$-twisted line
     bundle.  
     
       \vskip .2cm
     
     Further, $S\Lc_r$ is a $\GG$-equivariant $\Gc_{1/2}$-twisted line
     bundle on $\PP^{1|1}$. Precisely, this means the following.
      Let $A$ be a supercommutative algebra
     and $g\in\GG(A)$ be an $A$-valued point of $\GG$. Let $A_\red = A_\0/(A_\1)^2 A_\0$
     (functions on the underlying even part of the superscheme  $\Spec(A)$). 
     Then $g_\red$, the  specialization of $g$ to $A_\red$, is an $A_\red$-valued point of the
     even part $\GG_\red = \{\pm 1\} \times \SL(2,\CC)$, where the factor $\{\pm 1\}$ acts on
     $\PP^1$ trivially.  The $\SL(2,\CC)$-equivariance of $\Gc_{1/2}$ gives, for each open $U\subset \PP^1$ and each 
     object $x=(M,\phi) \in \Gc_{1/2}(U)$, an object $g_\red(x)\in \Gc_{1/2}(g_\red(U))$.
     The $\GG$-equivariance of $S\Lc_r$ means that for each $g$  as above we have  an identification      
     $\eps_g: S\Lc_{r, x} \to g^* S\Lc_{r, g_\red(x)}$
     of the 
     corresponding (super) line bundles on (the induced super-thickening of) $U$, this system of identifications satisfying
     obvious compatibility conditions. 
     
    % \vfill\eject

   \paragraph{ Super-forms of half-integer weight and the super-metaplectic group.}\label{par;super-MF}
After restriction to $\Hen\subset \PP^{1}$ the gerbe $\Gc_{1/2}$ becomes trivial, as the line
bundle $\Lc_1 = \Oc_{\PP^1}(-1)$ is trivial on $\Hen$ and so has (also trivial) square root $\Lc_{1\over 2}$. The gerbe 
$\Gc_{1/2}|_\Hen$ is  also equivariant with respect to the subgroup $\SL(2, \RR)$ preserving 
$\Hen$ but it is not trivial as an equivariant
gerbe, i.e., it does not have global equivariant objects. 
Such a situation leads, in  a standard way, see, e.g.,  \cite{ADK, brylinski-loop} to a central extension
of the group so that the gerbe  has a global object equivariant with respect to the extension. 
In our case we get  the extensrion
 \be\label{eq:MP2}
  1\to \ZZ/2 \lra \Mp(2) \lra \SL(2,\RR) \to 1
  \ee
  with  $\Mp(2)$ known as the {\em metaplectic
  group}. Recall that the  $\SL(2,\RR)$-action on elements of $\F_1 = H^0(\Hen,  \Lc_1)$ considered,
  in the standard trivialization 
  as functions $f(\tau)$, is given by 
  \[
  g^*f(\tau) = (c\tau+d)^{-1} f\bigl( a\tau+b \bigl/  c\tau+d\bigr), \quad g = \begin{pmatrix}
  a&b\\ c&d
  \end{pmatrix} \in \SL(2,\RR). 
  \]
So $\Mp(2)$ can be defined as a set of pairs $(g,\phi)$ where  $g$ as before
and $\phi$ is a branch of $(c\tau+d)^{1/2}$ in $\Hen$, i.e., a function
  $\Hen\to \CC$ such that $\phi(\tau)^2 = c\tau+d$,
  see, e.g.,  \cite{Kashiwara-vergne, nishiyama}. Thus for any $r = m+{1\over 2} \in{1\over 2}+\ZZ$
  we have a canonical $\Mp(2)$-equivariant line bundle $\Lc_r$ on $\Hen$ whose  sections
  are identified with functions $f(\tau)$ and the action of $\wt g = (g,\phi)\in \Mp(2)$
  is given by
  \[
  (\wt g^* f)(\tau) = (c\tau+d)^{-r} f\bigl(a\tau+b\bigl/c\tau+d\bigr), \quad (c\tau+d)^{-r} := 
  \phi(\tau)^{-2r} = \phi(\tau)^{-1} (c\tau+d)^{-m}. 
  \]
  The space  $\F_r = H^0(\Hen, \Lc_r)$  will be called the space  {\em forms of (half-integer) weight $r$}. 
  It is acted upon by $\Mp(2)$.
  
  \vskip .2cm
  
   Applying the same considerations to the super line bundle $S\Len_1$ on $S\Hen$, 
   we obtain the
   Lie supergroup $\OMp(1|2)$ fitting into a central extension
  \[
  1\to\ZZ/2 \lra \OMp(1|2)\lra G_\RR =  \OSp(1|2)(\RR) \to 1
  \]
  restricting to \eqref{eq:MP2} over $\SL(2,\RR)$. That is, points (in our standard sense)
  of $\OMp(1|2)$ are pairs
  $(g, \phi)$ where $g$ is a point of $G_\RR$ in the form \eqref{eq:OSp-matrix:gen1}
  and $\phi\in H^0(S\Hen, \Oc)$ is such that
  $\phi(\tau,\xi)^2 =  \delta\xi+c\tau+d$, with the obvious law of multiplication of such pairs. 
  We call $\OMp(1|2)$ the {\em super-metaplectic group}. So for any $r\in{1\over 2}+\ZZ$
  as before, we have the line bundle
  $S\Lc_r$ on $S\Hen$, equivariant with respect to $\OMp(1|2)$. Its sections can be
  seen as super-functions $f(\tau, \xi)$ and the action of $\wt g = (g,\phi)\in \OMp(1|2)$
  is given by the same formula as in  \eqref{eq:osp-2-forms} but 
  with $ (\delta\xi + c\tau + d)^{-r}$ understood using $\phi$. 
  The space $\SF_r = H^0(S\Hen, S\Lc_r)$ is acted upon by $\OMp(1|2)$. Its
  elements will be called {\em super-forms of (half-integer) weight}  $r$. 
  Eq. \eqref{eq:SF-r}   extends to an isomorphism of $\Mp(2)$-modules
     \be\label{eq:SF-r-half}
   \SF_{r,\0} \= \F_r, \quad \SF_{r, \1} \= \F_{r-1}, \quad r\in {1\over 2} + \ZZ. 
   \ee

  \vskip .2cm
  
  In particular, we consider the
  action on $S\Hen$ of the subgroup $\SL(2,\ZZ) \subset G_\RR$ given by the conditions
  $s=1$, $\alpha=\beta=\gamma=\delta=0$, $a,b,c,d\in\ZZ$, $ad-bc=1$. 
  \iffalse
  This action has the form
  \[
  \begin{pmatrix}
  a&b\\c&d
  \end{pmatrix} (\tau,\xi) = \biggl( {a\tau+b\over c\tau+d},\,  {\xi\over c\tau+d} \biggr).
  \]
  \fi
 Let $\Mp(2,\ZZ)$ be the central extension of $\SL(2, \ZZ)$ induced by $\Mp(2)\to \SL(2,\RR)$. 
 Suppose given a subgroup $\Gamma\subset \Mp(2,\ZZ)$ and a character
 $\chi: \Gamma\to \CC^*$. 
  For
 $r\in {1\over 2}\ZZ = \ZZ\cup ({1\over 2} + \ZZ)$ we denote 
 \[
 \MF_r(\Gamma, \xi) = \Hom_\Gamma(\CC_\chi, \F_r), \quad \SMF_r(\Gamma, \xi) = 
 \Hom_\Gamma(\CC_\chi, \SF_r)
 \]
 the space of $\Gamma$-modular forms of weight $r$ and character $\chi$ and
 the  
 the (super) space of {\em super-modular forms} of this weight and character. 
 Explicitly,  an element of $\SMF_r(\Gamma, \chi)$
 a function $f = f(\tau,\xi) \in \Oc(S\Hen)$ satisfying 
 \[
 f  \biggl( {a\tau+b\over c\tau+d},\,  {\xi\over c\tau+d} \biggr) = \chi(g) \phi(\tau)^{-2r}  f(\tau,\xi),\quad
 \forall (g, \phi)\in \Gamma, \quad 
 g=\begin{pmatrix} a&b\\c&d
 \end{pmatrix}. 
 \]
 Eqs. \eqref{eq:SF-r}  \eqref{eq:SF-r-half} imply: 
 \be
 \SMF_r(\Gamma,\chi)_\0 = \MF_r(\Gamma,\chi), \quad \SMF_r(\Gamma,\chi)_\1 = \MF_{r-1}(\Gamma,\chi), \quad r\in {1\over 2} \ZZ. 
 \ee

 \vfill\eject
 
 \subsection{Even Heisenberg relations for odd generators of $\osp(1|2)$}\label{subsec:evenheis-1}
 
 \paragraph{The Lie superalgebra $\osp(1|2)$.}
 Consider the complex Lie superalgebra 
 \[
 \gen := \osp(1|2) = \Lie(\GG), \quad \GG= \OSp(1|2)(\CC). 
 \]
 It consists of $3\times 3$ matrices $x$ such that $x^T H + Hx=0$, where $H$ is as in 
 \eqref{Gram-H=gen1}. Explicitly \cite[\S 2.3.1]{musson-book}  this means that $x$ has the form
 \[
 x=\begin{pmatrix}  0&\varkappa &\lambda 
 \\
 -\lambda & k & l 
 \\
 \varkappa & m & n
 \end{pmatrix}, \quad k+n=0. 
 \]
 with $\varkappa, \lambda $ parametrizing $\gen_\1$ and $k,l,m, n=-k$ parametrizing $\gen_\0$.
 That is, $\gen_\0=\ssl(2)$  and $\gen_\1 = \CC^2$
 is the symplectic space  which is the even part of $\CC^{2|1}$. 
 We denote the basis of $\gen$ by
 \[
 \begin{gathered}
  e = 0\oplus \begin{pmatrix} 0&1\\0&0
 \end{pmatrix}, \quad 
 f = 0 \oplus \begin{pmatrix} 0&0\\ 1&0 \end{pmatrix}, 
 \quad h=  0\oplus \begin{pmatrix} 1 & 0\\ 0 & -1 \end{pmatrix} \quad \in\ssl_2= \gen_\0, 
 \\
 p  = \begin{pmatrix} 0&0&1 \\
 -1&0&0\\
0&0&0
 \end{pmatrix}, \quad
 q=\begin{pmatrix}
 0&1&0\\
 0&0&0
 \\
 1&0&0
 \end{pmatrix} \quad \in \gen_\1. 
 \end{gathered} 
 \]
 The (super)commutation law then has the form
 \be\label{ops21-rels}
 \begin{gathered}
 [h,e]_- = 2e, \,\, [h,f]_- = -2f, \,\, [e,f]_-=h \quad \text{(standard  $\ssl(2)$)};
 \\
 \biggl[ 0\oplus \begin{pmatrix} a&b\\ c&d\end{pmatrix}, p\biggr]_-  = ap+bq,  \quad 
 \biggl[ 0\oplus \begin{pmatrix} a&b\\ c&d\end{pmatrix}, q\biggr]_-   = cp+dq \quad 
 \text{($\gen_\0=\ssl(2)$ acts on $\gen_\1$ as on $\CC^2$)}; 
 \\
 [p,p]_+=-2e, \quad [q,q]_+=2f, \quad [p,q]_+ = -h. 
 \end{gathered}
 \ee
 
 \paragraph{The LHS of the even Heisenberg relation as the ghost Casimir.}\label{par:ghostcas-1}
  Let $U(\gen)$
 be the universal enveloping (super) algebra of $\gen$. Consider the element
 \[
 Q\, =\,  [p,q]_-   - {1\over 2} \,\, =\,\,  pq-qp -{1\over 2} \,\,  \in \,\, U(\gen)_\0. 
 \]
 Note that the last relation in \eqref {ops21-rels} involves the super (i.e., as $p,q$ are odd, anti)
 commutator, so in $U(\gen)$ we have the identity $pq+qp=h$ while $Q$ involves the
 ``unnatural'' type of  commutator:  that with the minus sign for two odd arguments. 
 This element is known as the {\em ghost Casimir} (or {\em sCasimir})
  of $\gen$ and is a part of the general theory
 of  ghost centers of universal enveloping superalgebras \cite {pinczon, lesniewski,
 musson-center, arnaudon,  gorelik-ghost}, see \S 
 \ref {subsec:DO-super}
\ref{par:super-general} above for terminology.
 It has the following properties, going back to 
 \cite[Prop.1.2]{pinczon} which we quote in the form of \cite[Prop.II.I(i)-(ii)]{lesniewski}
 but with slightly different normalization of the basis vectors.
 
 \begin{prop}\label{prop:ghost-cas1}
 (a) $Q$ commutes with $U(\gen)_\0$ and anticommutes with $U(\gen)_\1$,
 i.e., $Q$ is an element of the anticenter $\Ac(U(\gen))$. 
 
 \vskip .2cm
 
 (b) We have $Q^2=  C+1/4$, where
 \[
 C= 2(ef+fe) + h^2 + qp-pq
 \]
 is the standard quadratic Casimir generating the center of $U(\gen)$ as an associative superalgebra. 
 \qed
 \end{prop}

 \paragraph{The action of the ghost Casimir on Verma modules.}  For $\lambda\in \CC$ we denote
 by $\SV_\lambda$ the induced Verma module over $\gen$ generated by one vector 
 $\ket{\lambda}$ satisfying 
 \[
 h\ket{\lambda} = \lambda\cdot \ket{\lambda}, \, e\ket{\lambda} = p\ket{\lambda} = 0
 \]
 and no other relations. In other words,
 \be\label{eq:pen-1-2}
 \SV_\lambda = \Ind_\pen^\gen \, \CC_\lambda, \quad\text{where} \quad  \pen = \CC h \oplus \CC e \oplus \CC p
 \ee
 is the natural parabolic subalgebra in $\gen$ and $\CC_\lambda$ is the $1$-dimensional $\pen$-module
 with $e,p$ acting by $0$ and $h$ acting by $\lambda$. It is straightforward that
 \[
 (\SV_\lambda)_\0 = \V_\lambda, \,\, (\SV_\lambda)_\1 = \V_{\lambda-1}
 \]
 as $\gen_\0=\ssl(2)$-modules, cf. \eqref{eq:SF-r}.
  Here $\V_\lambda$ is the standard Verma module
 for $\ssl(2)$.  Similarly to the case of $\ssl(2)$ one sees that 
 $\SV_\lambda$ is irreducible for $\lambda\notin \ZZ_+$ while for $\lambda\in \ZZ_+$ 
 it has the irreducible quotient which is a finite-dimensional representation of $\gen$. 
 The following is proved in the same way as  \cite[Prop.II.1(iii)]{lesniewski}
 (which treats the irreducible quotient for $\lambda\in \ZZ_+$). 
 
 \begin{prop}  
 The action of $Q$  on $(V_\lambda)_\0$ is given by the scalar $\lambda+1/2$ and 
 on $(V_\lambda)_\1$ it is given by its negative $-\lambda-1/2$. \qed
 \end{prop}
 
 \begin{cor}\label{cor:Verma-1/2}
 For $\lambda=-1/2$ (and for no other values of $\lambda$) the odd generators $p,q\in\gen_\1$
 acting on $V_\lambda$, satisfy the even Heisenberg relation $[p,q]_- = 1/2$. 
 \qed
 \end{cor}
 
 \paragraph{Action of $\osp(1|2)$ in forms of weight $1/2$.}  Let $r\in {1\over 2}\ZZ$. 
 The line bundle $S\Lc_r$ on $S\Hen$ being  $\OMp(1|2)$-equivariant, the Lie algebra
  $\gen = \Lie(\OMp(1|
 2))\otimes \CC$ acts in $S\Lc_r$
 by differential operators 
 which we denote 
 $D^{(r)}_y, y\in \gen$, see \S \ref{subsec:exp-odd-vect}. 
 
 Let us find the operators $D^{(r)}_p, D^{(r)}_q$
 corresponding to the odd generators $p,q\in\gen_\1$ explicitly, using the trivialization of $S\Lc_r$ on $S\Hen$
 described in \S \ref {subsec:P11} \ref {par:super-lobach} and \ref{par;super-MF}. That is, a (local)
 section of $S\Lc_r$ is represented by a function $f(\tau,\xi) = f_0(\tau) + f_1(\tau)\xi$ and the $\OMp(2)$-action is given by 
 \eqref{eq:osp-2-forms}, with  $(\delta+c\tau+d)^{-r}$ for half-integer $r$ understood as in 
 \S \ref {subsec:P11}\ref{par;super-MF}. 
 
 \begin{prop}\label{prop:DpDq}
 In the described trivialization the action of $\gen_\1$ in $S\Lc_r$ has the form
 \[
 \begin{gathered}
 D^{(r)}_p =  {\del\over\del \xi} - \xi {\del\over\del\tau} =  \begin{pmatrix}
 0 & 1 \\ {-\del \over  \del \tau } & 0
 \end{pmatrix},
 \\
 D^{(r)}_q =  \tau{\del\over\del\xi} - \tau \xi {\del\over\del\tau} - r\xi 
 =\,\,    \begin{pmatrix}
 0 & \tau \\
 -\tau {\del \over   \del \tau } -r & 0
 \end{pmatrix},
 \end{gathered} 
 \]
 where the matrices on the right describe the component presentation, viewing $f(\tau,\xi)$
 as a pair $(f_0(\tau), f_1(\tau))$. 
 \end{prop}
 
 \noindent{\sl Proof:} To find $D^{(r)}_p$ by differentiation, we take $g = 1+\eps p$ where $\eps$ is an odd parameter
 and take ${d\over d\eps}|_{\eps=0} (g^*f)(\tau, \xi)$. Writing $g$ in the form \eqref{eq:OSp-matrix:gen1}
 and using the Eq. \eqref{eq:osp-2-forms} for the action of $g$ on forms, 
   we find
 \[
 (g^*f)(\tau, \xi) = f(\tau-\eps\xi, \xi+\eps)
 \]
 and the linear term in $\eps$ gives $\del f/\del \xi -\xi\, \del f/ \del\tau $  by the chain rule. To find $D^{(r)}_q$, we
 take $g=1+\eps q$, then
 \[
 \begin{gathered}
 (g^*f)(\tau, \xi) = (\eps\xi +1)^{-r}\cdot f\biggl( {\tau\over\eps \xi +1}, \, \, {\xi+\eps\tau\over \eps\xi + 1}\biggr) = 
(1-r\eps\xi)\cdot   f\bigl( \tau(1-\eps\xi), \xi+\eps\tau\bigr) 
 \end{gathered}
 \]
  and direct differentiation gives the second claim. \qed 
  
  \begin{cor}\label{cor:Heis-1}
  For $r=1/2$ (and no other values of $r$) the operators $D^{(r)}_p, D^{(r)}_q$ satisfy
  the even Heisenberg relation $[D^{(r)}_p, D^{(r)}_q]_- = {1\over 2}$. \qed
 
  \end{cor}
  
  \begin{rem}
  This can be also seen from Corollary \ref {cor:Verma-1/2}, as $\SV_\lambda$ can be identified, 
  as a $\gen$-module, as
  \[
  \SV_\lambda = H^{1, \mer}_{\{0\}} (S\Lc_{\lambda+1}).
  \]
  Here the RHS is the space of meromorphic first cohomology of $S\Lc_{\lambda+1}$ with support in
  the  distinguished point  $0\in \PP^{1|1}$, i.e., the space of polar parts of meromorphic
  sections with a  pole at $0$. 
   Such an identification is a super-analog of the result of Borho-Brylinski
  \cite[Prop. 3.5]{brylinski}. 
  \end{rem}

 %%%BELOW COMMENTED OUT
 \iffalse
  Let 
 \[
 P= \left\{ \begin{pmatrix}
 s&\alpha& 0 \\
 0 & a& 0
 \\
 \delta&c&d
 \end{pmatrix}
 \right\} \subset G=\OSp(1|2)(\CC)
 \]
 be the parabolic subgroup with $\Lie(P)=\pen$ from \eqref{eq:pen-1-2}. It is the stabilizer the $\CC$-point
  $0+0\xi\in \AAA^{1|1}  \subset \PP^{1|1}$ so that we have an identification
 $G/P\= \PP^{1|1}$. 
 \fi
 %%%ABOVE COMMENTED OUT
 
 \vfill\eject
 
 \subsection{The Jacobi theta-function from the point of view of  supersymmetry}\label{subsec:jacobi}
 
 \paragraph{
 $\Theta(z,\tau)$ and its Thetanullwert.} \label{par:thetanull}
 The basic theta-function of Jacobi (denoted by $\theta_3$  in his nomenclature)
 is 
 \be
 \Theta(z,\tau) = \sum_{n\in \ZZ} \exp(\pi \i n^2\tau + 2\pi \i nz), \quad \tau\in\Hen, \,\, z\in \CC.
 \ee
 We will refer to $z$ as the {\em elliptic variable} and to $\tau$ as the
 {\em modular variable}. This function satisfies the heat equation
 \be\label{eq:heat-gen1}
 {\del\Theta\over\del\tau} = { -\i\over 4\pi}\,\,  {\del^2\Theta \over \del z^2}
 \ee
 and the quasi-periodicity conditions (``elliptic behavior")
 \be\label{eq:theta-1-period}
 \Theta(z+1, \tau) = \Theta (z,\tau), \quad
 \Theta (z+\tau, \tau) = \exp(-\pi \i \tau - 2 \pi \i z)\,  \Theta(z,\tau). 
 \ee
 The {\em (Jacobi) Thetanullwert} is the function 
 \[
 \theta(\tau) = \theta_J(\tau) := \Theta(0, \tau) = \sum_{n\in \ZZ} \exp(\pi \i n^2\tau), \quad \tau\in\Hen
 \]
 of the modular variable alone. It satisfies the periodicity and the modular relation
 \be\label{1dtheta-mod}
 \begin{gathered}
 \theta(\tau+2)=\theta(\tau),
 \\
 \theta (-1/\tau) = e^{-\pi \i /4} \sqrt{\tau}\,\,  \theta (\tau),
 \end{gathered}
 \ee
 where $\sqrt{\tau}$ is the branch of the square root in $\Hen$ which is positive on $\RR_+$. 
 In other words,  considered as an (even)  section of $\Lc_{1/2}$, it is invariant under $\wt T_2$
 and quasi-invariant under $\wt S$, where $\wt T_b,  \wt S\in\Mp(2,\ZZ)$
 are elements defined as
 \be\label{eq:tilde-TS}
 \begin{gathered}
 \wt T_b = (T_b, 1), \quad T_b=\begin{pmatrix} 1&b\\0&1\end{pmatrix}, \,\, 1= \text{ determination of } (c\tau+d)^{1/2} = 1^{1/2}
 \text{ equal to } +1, \\
 \wt S = (S, \sqrt{\tau}), \quad S=\begin{pmatrix} 0&-1\\1&0\end{pmatrix}, \,\,
 \sqrt{\tau}= \text{ the above branch of}  (c\tau+d)^{1/2} = \tau^{1/2}.
 \end{gathered} 
 \ee
 More generally, it is known \cite[\S I.1.9, III.8.6]{mumford} that $\theta \in \MF_{1/2}(\Gamma, \chi)$,
 where:
 \begin{itemize}
 \item  $\Gamma\subset \Mp(2,\ZZ)$ is the preimage of the subgroup in
 $\SL(2, \ZZ/2)$ consisting  of % $\bpm 1&0\\0&1\epm$
  $1_2$ and  $\bpm0&1\\1&0\epm$.
  
  \item $\chi: \Gamma\to\CC^*$ is the  character taking values in $8$th roots of $1$,
   uniquely defined by the conditions $\chi(\wt T^2)=1$, $\chi (\wt S) = e^{\pi i/4}$,
  and $\chi(\sigma)=-1$ where $\sigma\in \ZZ/2$ is the generator of the center. 
   \end{itemize}
   
   \vfill\eject
   
   \paragraph{The elliptic variable as spinors for the modular variable.}
    \label{par:ell-spinors}
    The heat equation
   \eqref{eq:heat-gen1} means that $\del/\del z$ acts on $\Theta$ as
   $\sqrt{4\pi \i \del/\del\tau}$, so the shift operators in $z$ can be written formally
   as operators of infinite order involving the modular variable alone: 
   \[
   \begin{gathered}
   \Theta(z+1, \tau)=e^{\del/\del z}\Theta(z,\tau)= e^{\sqrt{4\pi i \del/\del\tau}} \Theta(z,\tau),
   \\
   \Theta(z+\tau, \tau) = e^{\tau \del/\del z} \Theta(z,\tau) = e^{\tau\sqrt{4\pi \i \del/\del\tau}} \Theta(z,\tau).
   \end{gathered} 
   \]
   Because of the square roots these are not honest differential operators
   of infinite order but they are well defined
   sections of
   the sheaf $\Ec^\RR$ of holomorphic microlocal operators \cite{KK-Kimura} on  an appropriate domain
    in $T^*\CC$. The original idea of Sato \cite{sato-theta}, see also 
   \cite[\S2.4]{kawai-theta}, was to use this to rewrite the quasi-periodicity conditions
   \eqref{eq:theta-1-period} characterizing $\Theta(z,\tau)$ as equations
   \be\label{eq:theta-microdiff-1}
   e^{\sqrt{4\pi \i \del/\del \tau}}\,  \theta(\tau) = \theta(\tau), \quad e^{\tau \sqrt{4\pi \i \del/\del\tau} + \pi \i \tau} \,
   \theta(\tau) = \theta(\tau)
   \ee
   of infinite order in $\tau$   that would  characterize $\theta(\tau)$. In order to realize this idea
   it is natural to convert $\sqrt{\del/\del \tau}$ into a $2\times 2$ matrix differential operator as in 
   \cite{sato-KK-theta, kawai-theta}. The conceptual meaning of this procedure
   (apparently unnoticed in {\em loc. cit.}) is that the elliptic variable appears as the
   ``spinors'' for the modular variable which in this case plays
    the role of  ``space-time''  so that
   $\sqrt{\del/\del\tau}$ becomes the spinor derivation $D_\xi$ from 
   \eqref{eq:D-xi} and $\Hen$ becomes extended to $S\Hen \subset \PP^{1|1}$.
   Let us present the analysis of   \cite{sato-KK-theta, kawai-theta} from this angle.
   
   \paragraph{Characterization of $\theta(\tau)$ from SUSY point of view.} 
   Consider $\theta(\tau, \xi) = \theta(\tau)+0\xi$ as an even section of the line bundle
   $S\Lc_{1/2}$ of super-forms of weight $1/2$ on $S\Hen$. The following is
   a   SUSY reformulation of the results of  \cite{kawai-ex}. 
   
    \begin{thm}\label{thm:1d-theta-charact}
     (a)  $\theta(\tau, \xi)$ satisfies the following super-differential equations of infinite order:
   \[
   \begin{gathered}
   e^{\sqrt{-4\pi \i}\,  D^{(1/2)}_p}\theta = \theta, 
   \quad e^{\sqrt{-4\pi \i}\,  D^{(1/2)}_q} \theta = \theta,
   \quad 
  \text{ where}
  \\
   D^{(1/2)}_p = {\del\over\del\xi} -\xi {\del\over\del\tau}, 
   \quad D^{(1/2)}_q=
   \tau{\del\over\del\xi} -\tau \xi {\del\over\del\tau} - {1\over 2} \xi 
    \end{gathered} 
   \]
   are the odd generators of $\osp(1|2)$ acting in $S\Lc_{1/2}$, see
   Proposition \ref{prop:DpDq},  and $\sqrt{-4\pi \i}$ is chosen the same in both
   equations. 
    
    (b) 
   Any local holomorphic section $\Phi=\Phi (\tau, \xi)$ of  $S\Lc_{1/2}$  satisfying the  system  
      \[
     e^{\sqrt{-4\pi \i}\,  D^{(1/2)}_p} \Phi  = \Phi, \quad
      e^{\sqrt{-4\pi \i}\,  D^{(1/2)}_q} \Phi = \Phi, 
\]
is a constant multiple of $\theta(\tau)$. 
   \end{thm}
   
   \noindent 
    We give the proof in the remainder of this section and in the next one. 
   
   \vskip .2cm

 \begin{rems}
 (a)  Since $[D^{(1/2)}_p, D^{(1/2)}_q]_-=1/2$ on $S\Lc_{1/2}$ (Proposition \ref{cor:Heis-1}),
   we have
   \[
   [\sqrt{-4\pi \i} \, D^{(1/2)}_p, \sqrt{-4\pi \i}\,  D^{(1/2)}_q]_- = 2\pi \i, 
   \quad \text{and so} \quad   e^{\sqrt{-4\pi \i} D^{(1/2)}_p} 
    e^{\sqrt{-4\pi \i} D^{(1/2)}_q}
     e^{\sqrt{-4\pi \i} D^{(1/2)}_q}  e^{\sqrt{-4\pi \i} D^{(1/2)}_p}
   \] 
there, i.e. $\theta$ is represented as a common eigenvector of two commuting operators. 

\vskip .2cm

  (b)  Note that the system of Theorem \ref{thm:1d-theta-charact} on 
  $S\Lc_{1/2}$
   is manifestly modular invariant: the elements $p,q\in \osp(1|2)_\1$ are  obtained
   from each other (up to sign)  by adjoint action of the element $\wt S\in \Mp(2, \ZZ)$
   from \S \ref{par:thetanull} effecting the modular transformation. 
   Thus the theorem can be used to deduce the modular behavior
   \eqref{1dtheta-mod}  of $\theta(\tau)$
   from the purely  local conditions given by the above differential equations of infinite order.  Given the origin of these equations discussed in
   \S \ref{par:ell-spinors}, we can say that
   {\em the modular behavior of $\theta(\tau)$  follows from the elliptic behavior 
   of $\Theta(z,\tau)$ in $z$ and the heat
   equation}. 
   
   \vskip .2cm
   
   (c) To compare the above system with \eqref {eq:theta-microdiff-1},
   we recall that $D_p^2 = - \del/ \del\tau$, which explains passing from
   $\sqrt{4\pi \i}$ to $\sqrt{-4\pi \i}$.

 \end{rems}

\paragraph{Super-Gaussians.} \label{par:super-gauss-1}
 For any $x\in\CC$ put
   \be\label{eq:super-gauss-1}
   \theta_x(\tau,\xi) = e^{\pi i x^2\tau} (1+ \sqrt{-\pi i} \, x\, \xi)  = e^{\pi \i x^2 \tau +
    \sqrt{-\pi \i} \, x\,\xi}. 
   %\, e^{\pi \i x^2\tau}\xi.
   \ee
   We call $\theta_x(\tau, \xi)$ as a function of  its three variables, the {\em super-Gaussian}. 
   Its even part $\theta_x(\tau,0)$ is just the variable Gaussian in $x$. 
   
   \begin{prop}\label{prop:super-Gauss-1}
   The super-Gaussian satisfies the equations
   \[
      \begin{gathered}
   \sqrt{-4\pi \i}\,  D^{(1/2)}_p\,  \theta_x(\tau, \xi) = -2\pi \i \,\, x \, \theta_x(\tau, \xi),
   \\
   \sqrt{-4\pi \i} \, D^{(1/2)}_q\,  \theta_x(\tau, \xi) = - {\del\over \del x} \theta_x(\tau, \xi).
\end{gathered}
   \]
   \end{prop}
   
   \noindent{\sl Proof:} Direct computation. NB: as $D^{(1/2)}_p, D^{(1/2)}_q$ are not derivations in the usual sense,
   appliying a formula like  $D(e^f) = D(f) e^f$ 
 to the RHS of  
      \eqref{eq:super-gauss-1} is illegal.  \qed
   \vskip .2cm
   
   \noindent{\sl Proof of part (a) of Theorem \ref{thm:1d-theta-charact} :} 
    Note that
   \be\label{eq:theta-super-sum}
   \theta(\tau) = \sum_{n\in\ZZ} \theta_n(\tau, \xi)
   \ee
   (summation over integer values $x=n$; the terms with $\xi$ cancel). 
    Proposition \ref{prop:super-Gauss-1} implies that 
 for any $n\in \ZZ$ we have
   \be\label{eq:AB-thetan-gen1}
   \begin{gathered}
   e^{\sqrt{-4\pi \i}\,  D^{(1/2)}_p} \, \theta_n = e^{-2\pi \i n} \theta_n = \theta_n, 
   \\
   e^{\sqrt{-4\pi \i} \, D^{(1/2)}_q}\,  \theta_n = \theta_{n-1}
   \end{gathered}
   \ee
   and so the first claimed equation holds term by term  w.r.t. 
    \eqref{eq:theta-super-sum} while the second  one holds
   by shifting the terms in the sum.  \qed
   
   \vskip .2cm
   
   This argument has a natural interpretation in terms of the super-Weil
   representation of $\OMp(1|2)$, see \S
   \ref {subsec:super-weil-1}
 \ref {par:Weil-compar-gen1}. 
   
   \vfill\eject
   
   \subsection {Analysis of the Koszul complex}\label{subsec:koszul-1}

 \paragraph{ Koszul complex:  definitions and statements.} Consider the commuting
 operators
 \[
 A= e^{\sqrt{-4\pi i} D^{(1/2)}_p}-1, \quad B= e^{\sqrt{-4\pi\i} D^{(1/2)}_q}-1
 \]
  in $S\Lc_{1/2}$.  Since they commute, we can form the Koszul complex
  \[
  \Cc^\bullet = 
  \bigl\{ S\Lc_{1/2} \buildrel (A,B)\over\lra  S\Lc_{1/2}^{\oplus 2} 
  \buildrel B-A \over\lra S\Lc_{1/2}\bigr\}
  \]
  of sheaves on $\Hen$ situated in degrees $0,1,2$. . 
  Thus $\ul H^0(\Cc^\bullet)$ is the sheaf of solutions of the system in
  Theorem \ref{thm:1d-theta-charact}(b)  and the  theorem is a part
  of the following more precise result of \cite{kawai-ex}.
  
  \begin{thm}
  \label{thm:Koszul-gen1}
  We have $\ul H^0(\Cc^\bullet) \= \ul\CC_\Hen$, the constant sheaf
  spanned by $\theta(\tau)$, while $\ul H^1(\Cc^\bullet) = \ul H^2(\Cc^\bullet)=0$. 
  \end{thm}
  
 We will give an outline of  the proof. The argument splits into two parts:  
 the study of the characteristic variety and the 
 computation at a single point,
 represented by the following two propositions. which imply the theorem.

 \begin{prop}\label{prop:koszul-lc-1}
 The sheaves $\ul H^i(\Cc^\bullet)$ are locally (and hence globally) constant on $\Hen$. 
 \end{prop}
 
 \begin{prop}\label{prop:koszul-gl-1}
 The stalk $\Cc^\bullet_\i$
 of $\Cc^\bullet$ at $\tau=\i$
  has the zeroth cohomology  space $1$-dimensional and spanned by 
  $\theta(\tau)$
 while other cohomology spaces vanish. 
 \end{prop}

 \paragraph{ Study of the characteristic variety: proof of Proposition \ref{prop:koszul-lc-1}.}
 \label{prop:char-var-gen1}
  As the differentials in $\Cc^\bullet$ are made out of sections of $\Dc^\oo_{S\Hen}$,
 we realize  $\Cc^\bullet$ as the solution complex
  $\Sol(\Nc^\bullet) = R\ul\Hom_{\Dc_{S\Hen}^\oo}(\Nc^\bullet, \Oc_{S\Hen})$
 for the  corresponding (``dual'')  Koszul complex of free $\Dc_{S\Hen}^\oo$-modules:
 \[
 \Nc^\bullet  = \bigl\{ \Dc_{S\Hen}^\oo\buildrel (-A,B)\over  \lra (\Dc_{S\Hen}^\oo)^{\oplus 2}\buildrel A+B \over  \lra \Dc_{S\Hen}^\oo\bigr\}
 \]
 situated in degrees $0, -1, -2$.  Using component analysis, we identify
 \[
 \Oc_{S\Hen} = \Oc_\Hen \oplus \xi \Oc_\Hen = \Oc_\Hen^{\oplus 2}, \quad
 \Dc^\oo_{S\Hen} = \Mat_2 (\Dc^\oo_\Hen). 
 \]
 By the Morita equivalence between modules over $\Dc^\oo_\Hen$ and $\Mat_2(\Dc^\oo_\Hen)$
 we identify 
 \[
 \Cc^\bullet  = R\Hom_{\Dc^\oo_\Hen}(\Nc^\bullet_\red, \Oc_\Hen), \quad\text{where} \quad
 \Nc_\red^\bullet  = \bigl\{ (\Dc_{\Hen}^\oo)^{\oplus 2} \buildrel (-A,B)\over  \lra (\Dc_{S\Hen}^\oo)^{\oplus 4}\buildrel A+B \over  \lra (\Dc_{S\Hen}^\oo\bigr)^{\oplus 2}\}
 \]
 with $A,B$ now  written as elements of $\Mat_2(\Dc^\oo_\Hen)$. 
 
 \vskip .2cm
 
 Proposition \ref{prop:koszul-lc-1} means that $\SS(\Cc^\bullet)$ is contained in
 (and hence equal to)  the zero section
 of $T^*\Hen$. To prove it, it suffices, by Theorem \ref {thm:Char-SS}, to show that the complex
 $\Ec^\RR_\Hen \otimes_ {\Dc^\oo_\Hen}\Nc^\bullet_\red$ is exact outside the said zero section. 
 
 \vskip .2cm
 
 Let $Z_A$ resp. $Z_B\subset T^*\Hen$ be the locus where $A$ resp. $B$ is not invertible in 
 $\Mat_2(\Ec^\RR_\Hen)$. Recall that whenever one of the commuting elements giving rise to
 a Koszul complex is invertible, then the complex is exact, as it becomes realized as the
 cone of an invertible morphism. So
 \[
 \Ch(\Nc^\bullet_\red) \,=\,
 \on{Supp} \, \ul H^\bullet \bigl( \Ec^\RR_\Hen \otimes_{\Dc^\oo_\Hen} \Nc^\bullet_\red\bigr) \,\,\subset \,\, 
 Z_A\cap Z_B. 
 \]
 To estimate $Z_A$, notice that $D^{(1/2)}_p$ is, in components,
  a $2\times 2$ matrix operator with constant coefficients:
   \be\label{eq:sigma-A-gen1}
 D^{(1/2)}_p = \bpm   0&1 \\ -\del_\tau & 0\epm , 
 \quad \sigma_{D^{(1/2)}_p}(\lambda) = \bpm 0&1 \\ -\lambda & 0
 \epm. 
 \ee
  We now  apply Proposition \ref {thm:aoki-inv} to $C= \sqrt{-4\pi \i} \, D_p^{(1/2)}$
  and put $r_1=1/2, r_2=1, \rho=1/2$. 
  Then the matrix $\|\sigma_{C_{ij}}^{r_i-r_j+\rho}(\lambda)\|$ is equal to
  $\sqrt{-4\pi \i} \, \sigma_{D^{(1/2)}_p}(\lambda)$. As the eigenvalues
  of $\sigma_{D^{(1/2)}_p}(\lambda)$ are $\pm\sqrt{-\lambda}$, the condition
  given by Proposition \ref {thm:aoki-inv}  for invertibility of $A=e^C-1$ is
  $\lambda\notin  \i \RR_{\leq 0}$. This means  that 
 \[
 Z_A \,\,\subset \,\, \Hen \times \i \RR_{\leq 0}\,\,\subset \,\, \Hen \times \CC \,\,=\,\, T^* \Hen. 
 \]
 Next, $B$ is obtained from $A$ by the modular transformation $S: \tau \mapsto -1/\tau$, so
 $Z_B = (dS)^*(Z_A)$. But for $\tau\in\Hen$ the direction $(d_\tau S)^* (\i \RR_{\leq 0}) \subset \CC$
 is transversal to $\i \RR_{\leq 0}$, so their intersection is $\{0\}$. This means that $Z_A\cap Z_B$
 is contained in  the zero section of $T^*\Hen$.
 Proposition   \ref{prop:koszul-lc-1} is proved.

   \paragraph{Computation at $\tau=\i$: proof of Proposition  \ref {prop:koszul-gl-1}.}
      We present a 
   version of the argument of \cite {kawai-ex},
   slightly expanded and adapted to our notation.

   Denote $F=S\Lc_{1/2, \i}$ the stalk at $\tau=\i$ of the sheaf
   $S\Lc_{1/2}$, i.e., the space of sections holomorphic near $\i$. 
   We consider $A,B$ as operators on $F$, so
   \[
   \Cc^\bullet_\i = C^\bullet := \bigl\{
   F\buildrel (A.B) \over\lra F^{\oplus 2} 
   \buildrel B-A\over\lra F\bigr\}. 
   \]
   It suffices to prove three statements:
   \begin{itemize}
   \item[(K1)] $A$ is surjective.
   
   \item[(K2)] $B: \Ker(A)\to \Ker(A)$ is sirjective.
   
   \item[(K3)] $\Ker(A)\cap \Ker (B) = \CC\cdot \theta(\tau)$. 
   \end{itemize}
  Here (K1-2) form the dual version (with kernels rather then cokernels)
   of the condition that $A,B$
  form a regular sequence and ensure that $H^{>0}(C^\bullet)=0$,
  while (K3) idenifies $H^0(C^\bullet)$. 
  
  \vskip .2cm
  
  We identify $F$ with $\Oc_{S\Hen, \i} = \Oc_{\Hen, \i}^{\oplus 2}$
  using the trivialization of $S\Lc_{1/2}$. The condition on a (local)
  section of $S\Lc_{1/2}$ to be an eigenfunction of $D^{(1/2)}_p$ is then
   a $2\times 2$ matrix ODE with constant coefficients. Solving it
   in a standard way, we find that such eigensections (local or global)
   are given by the super-Gaussians $\theta_x(\tau, \xi)$ for
   $x\in \CC$. By Proposition \ref {prop:super-Gauss-1} we have
   $D^{(1/2)}_p \theta_x = \sqrt{-\pi \i} \, x \, \theta_x$. 
   
   \vskip .2cm
   
   By the Grothendieck-Martineau
   duality (a version of the Fourier transform in the complex domain,
   see \cite{martineau}),
   the space
    $\Oc_{\Hen, \i}^\vee$,   the topological dual of
   $\Oc_{\Hen, \i} =. \CC\{\{ \tau-\i\}\}$ is  identified with $\Oc(\CC_\lambda)$, the
   space of entire functions in the dual variable $\lambda$.
   Accordingly, $F^\vee$ is identified with $\Oc(\CC_\lambda)^2$.

   As with any kind of Fourier transform, the action of $\del/\del\tau$
   on $\Oc_{\Hen, \i}$ transposes to the multiplication by $\lambda$
   in $\Oc(\CC_\lambda)$. 
   This implies, since $A$ has constant
   coefficients, that the action of $A$ on $F$
   transposes to the multiplication, on $\Oc(\CC_\lambda)^{\oplus 2}$,
    by the matrix
   \[
   \sigma_A(\lambda)^t\, \,= \,\, \exp\biggl( \sqrt{4\pi \i} \bpm 0&-\lambda \\ 1 & 0\epm\biggr) -1,
   \]
   the transpose of the total symbol of $A$, cf. 
   \eqref{eq:sigma-A-gen1}.  This matrix is invertible outside of
    $\lambda\in\ZZ$, so the multiplication by $\sigma_A(\lambda)^t$ on
    $\Oc(\CC_\lambda)^{\oplus 2}$ is injective. This means that
    $A$ is surjective, proving (K1).
    
 \vskip .2cm
    
    Next, $\Ker(A)$ is dual to the cokernel of $\sigma_A(\lambda)^t$ on $\Oc(\CC_\lambda)^{\oplus 2}$ which is  closed and topologically
    spanned by the (identical) $1$-dimensional
    contributions at $\lambda\in \ZZ$,
    where $\sigma_A(\lambda)^t$ acquires $1$-dimensional kernel
    and cokernel. To identify these, note that $A= e^{\sqrt{-4\pi \i} D^{(1/2)}_p} -1$
    commutes with $D^{(1/2)}_p$, so the $1$-dimensional contribution to 
    $\Ker(A)$ from $\lambda=n\in \ZZ$ is spanned by the 
    $D^{(1/2)}_p$-eigensection annihilated by $A$, i.e., by $\theta_n(\tau, \xi)$. 
    That is, $\Ker(A)$ consists of sums
     $\sum_{n\in \ZZ} c_m \theta_n(\tau, \xi)$ convergent near $\tau=i$.
     The condition for convergence is easily seen to be
     $|c_n| \leq Ce^{(1-\eps) \pi n^2}$ for some $\eps>0$. Now,
     $e^{\sqrt{-4\pi i} D^{(1/2)}_q} \theta_n = \theta_{n-1}$, so the action of $B$ on $\Ker(A)$
     is terms of the sequences $(c_n)$ is given by the ``difference
     derivative'': $B((c_n)) = (c_{n+1}-c_n)$. So to show (K2)
     we need to solve the equation $c_{n+1}-c_n = d_n$
     for any $\Phi = \sum d_n \theta_n\in\Ker(A)$. This can be
     done uniquely up to an additive constant (``difference integral'')
     and one sees immediately that such $c_n$ again satisfy the
     growth condition so $\sum c_n\theta_n$ converges and lies
     in $\Ker (A)$. 
     
     \vskip .2cm
     
     Finally, $\Ker(A)\cap\, \Ker(B)$ consists of $\Phi(\tau, \xi) = \sum c_n \theta_n(\tau, \xi) $
     such that $c_{n+1}=c_n$ for all $n$, i.e., 
     $\Phi (\tau, \xi) = c\, \theta(\tau)$. This proves Proposition 
     \ref {prop:koszul-gl-1} and Theorem 
     \ref {thm:Koszul-gen1}. 
   
      \vfill\eject

   \subsection{The (super-) Weil representation of $\OMp(1|2)$}\label{subsec:super-weil-1}
   
   \paragraph{Group level: $\Mp(2)$.} \label{par:weil-group-1}
  
 We start by recalling the classical theory \cite{weil, lion-vergne} for the simplest case of the group $\SL(2,\RR)$.
 The real symplectic space $E_\RR=(\RR^{2}, \omega)$ with the standard
 basis $p, q$,  $\omega(p,q)=1$,  gives rise to the {\em Heisenberg group} 
 $
 \Hc = E_\RR \times \RR
 $
 with the product
 \[
 (v, a) \cdot (v', a') \,=\, \bigl(v+v', a+a' + {1\over 2} \omega(v, v')\bigr),
 \quad v, v'\in E_\RR, \,\, a, a'\in \RR. 
 \]
 It is the integrated  real version of the complex Lie algebra
 $\hen = \CC p + \CC q + \CC 1$ spanned by $p,q,1$ in the Heisenberg algebra
 $\Heis = \CC \langle p,q\rangle$, $pq-qp=1$. 
 The standard representation of
 $\Heis\supset \hen$ in  the space $\CC[q]$  of polynomials
 integrates to a unitary representation  $\rho$ of
 $\Hc$ in the Hilbert space $L_2(\RR)$, where $\RR=\RR_q
 = (\RR\cdot q)^*$ is the $1$-dimensional real vector space with
 coordinate $q$. It is known that $\rho$ is the unique, up to an
 isomorphism, infinite-dimensional irreducible unitary representation
 of $\Hc$. 
 
 \vskip .2cm
 
 The group $\SL(2,\RR)$ acts on $\Hc$ by group automorphisms
 $g\cdot (v,a) = (g(v), a)$. For each $g$  we get a  new representation
 $\rho_g: \Hc \to U(L_2(\RR^))$  sending $(v,a)$ to
 $\rho(g\cdot (v,a))$. By uniqueness of $\rho$ there is a (unique, up to
 a unitary scalar)  unitary operator $\varpi(g)\in U( L_2(\RR))$
 such that $\rho_g(v,a)= \varpi(g)^{-1} \rho(v,a) \varpi(g)$. For any choice of
 the $\varpi(g)$ we then have the equalities
 \[
 \varpi(g_1g_2) \,=\, c(g_1, g_2) \cdot \varpi(g_1)
  \varpi(g_2), \quad \exists c(g_1, g_2)\in \CC^*,
 \,\,\, |c(g_1, g_2)| = 1. 
 \]
 The result of  Segal-Shale-Weil \cite{weil}  (in our particular case) is that one can choose the $\varpi(g)$ so that
 all $c(g_1, g_2) = \pm 1$ and  the $\varpi(g)$ organize into a
 representation $\varpi: \Mp (2) \to U(L_2(\RR))$, now
 known as the {\em Weil representation}.    
 The choice of the $\varpi(g)$ can be fixed, essentially uniquely
 (up to modifying the cocycle $c(g_1, g_2)$ by a coboundary)
 by specifying the values of $\varpi(g)$ for a set of $g$ generating
 the group  $\Mp(2)$.    This is usually done by
  specifying that
 for $u=u(x)\in L_2(\RR)$ 
 and $b\in \RR$  we have
 \be\label{eq:weil-expl-1}
 \varpi \left( \bpm 1&0 \\ c&1 \epm\right)  u \,=\, e^{-\pi  \i c x^2}
 \cdot u,\quad
 \text{while} \quad 
 \varpi  \left(\bpm 0& -1\\ 1&0 \epm\right) u  \,=\, \wh u, 
 \ee
 where $\wh u$ is the Fourier transform of $u$.

 \vskip .2cm
 
  We denote $L_2(\RR)$
 with this action of $\Mp(2)$  by $\Wc$.
 Let $\Wc_\0, \Wc_\1\subset\Wc$ be the subspaces
  of even and odd functions: $u(-x) = \pm u(x)$.
  Formulas  \eqref{eq:weil-expl-1} 
 imply that the decomposition
  $\Wc = \Wc_\0 \oplus \Wc_\1$  
 is $\Mp(2)$-invariant. 
 It is further known that $\Wc_\0$ and $\Wc_\1$ are irreducible
 unitary representations of $\Mp(2n)$. Sometimes
 the term ``Weil representation" is applied to $\Wc_\0$ only. 
 
 \paragraph{Lie algebra level: $\ssl(2)$.}  By differentiating \eqref {eq:weil-expl-1}
 in the standard way,  we get a representation of $\ssl(2)$ in 
 $C^\oo(\RR)$ given by
 \be\label{eq:weil-sl2}
 e\mapsto  {1\over 4\pi \i } {d^2\over dx^2}, \,\,\, f\mapsto - \pi \i  x^2, \,\,\,
 h\mapsto -x {d\over dx} - {1\over 2}. 
 \ee
 (The definitions of $e$ and $f$ are often rescaled to get the coefficients
 $\pm 1/2$ but it is convenient for us not to do this.) 
 We will consider the algebraic version of this representaion in $\CC[x]$, the space
 of polynomials. As before, it splits into  the sum $\CC[x]_\0 \oplus \CC[x]_\1$
 with $\CC[x]_\0 = \CC[x^2]$ and $\CC[x]_\1 = x\, \CC[x^2]$ are the spaces of
 even and odd polynomials. They are  generated by elements $1$ and $x$ which
 are highest  vectors (i.e., annihilated by $e$) with $h$-weights $-1/2$ and $-3/2$
 respectively.  Since they are free over $\CC[f]$, we see that
 \[
 \CC[x]_\0 \= \V_{-1/2}, \quad  \CC[x]_\1 \= \V_{-3/2}
 \]
 are the induced (and irreducible, in this case) Verma modules over $\ssl(2)$ with 
 highest weights
 $-1/2$ and $-3/2$ respectively. 
 
 \paragraph{Lie algebra level: $\osp(1|2)$.}
  We extend \eqref{eq:weil-sl2} to a representation of
 $\gen = \osp(1|2)$ on $\CC[x]$ by adding the action of $p,q\in\gen_\1 = \CC^2$ by
 \be\label{eq:osp(1\2)-weil}
 p\mapsto {1\over\sqrt{-4\pi \i }}\, \,  {d\over dx}, \quad q\mapsto \sqrt{-\pi \i } \,\, x. 
 \ee
 One checks that the relations \eqref{ops21-rels} are satisfied, so $\CC[x]$ (considered
 as a super vector space w.r.t.  the above $\ZZ/2$-grading) becomes a
 $\gen$-module which we call the {\em (super) Weil representation} of $\gen$,
 see \cite {malcom}. Note that $p,q$  (the odd generators of $\gen$)
  satisfy,
 in their action on $\CC[x]$, the even Heisenberg relation
 $
 [p,q]_- = {1/2}. 
 $
 In fact, since $\CC[x]$ is free over  $\CC[q]$, it is identified with the Verma module
 over $\gen$ with highest weight $-1/2$:
 \[
 \CC[x] \= \SV_{-1/2}, 
 \]
so the appearance of the even Heisenberg relation matches Corollary \ref {cor:Verma-1/2}. 

\paragraph{Group level: $\OMp(1|2)$.
The Fourier transform automorphism. }  \label{par:grsweil-1}
 
 The theory of unitary representations of Lie supergroups was initiated in
 \cite{carmeli}. In brief, a (real)  Lie supergroup $G$
 gives a (super) Harish-Chandra pair consisting of the underlying
 ordinary Lie group $G_\red$ and the Lie superalgebra $\gen = \Lie (G)\otimes_\RR\CC$.
  A unitary representation of $G$ is, by definition, a representation of this
 Harish-Chandra pair, i.e., a  $\ZZ/2$-graded unitary representation of $G_\red$ 
 together  with compatible anti-self-adjoint representation of the Lie superalgebra
 $\gen$ on the appropriate space of vectors smooth for the  $G_\red$-action. 
 We refer to \cite{carmeli} for more details. 
 
  \vskip .2cm

We apply this to $G = \OMp(1|2)$ and $\gen = \osp(1|2)$. 
 In this case, Eq. \eqref {eq:osp(1\2)-weil} defines an action
 of $\gen$ on $C^\oo(\RR)$, the space of smooth vectors of $\Wc$. Therefore
 we get a representation  of $\OMp(1|2)$ on $\Wc$ which we call the
 {\em (super) Weil representation} of $\OMp(1|2)$.  It is irreducible, since
 both $\Wc_\0$ and $\Wc_\1$ are irreducible under $\Mp(2)$. 
 
   \vskip .2cm
   
   We will use the inner automorphism
   \be\label{eqFT-auto-1}
 \FT:   \OMp(1|2) \lra \OMp(1|2), \quad g\mapsto \wt S g\wt S^{-1},
   \ee
 where $\wt S \in \Mp(2)$ is  the   preimage of
 $\bpm 0&-1\\1&0\epm\in \SL(2, \ZZ)$ from 
 \eqref {eq:tilde-TS} (the other preimage will give the same conjugation). 
 We call $\FT$
 the {\em Fourier transform automorphism}. We use the same name
 and notation $\FT$ for the induced automorphism of the Lie superalgebra
 $\gen = \osp(1|2)$. The action of $\FT$ on $\gen_\1 = (\CC^2,\omega)$
 takes $p\mapsto q$, $q\mapsto -p$.

 \paragraph{The Gaussian and super-Gaussian transforms.}
\label{par:sgauss-1}
 Recall \eqref {eq:super-gauss-1} the super-Gaussian $\theta_x(\tau, \xi)$, for $x\in \CC$.
 Restricting to $x\in \RR$, we see that $\theta_x(\tau,\xi)$ decays exponentially
 as $x\to \oo$. So for any $u=u(x)\in L_2(\RR)$ the {\em super-Gaussian transform}
  \be\label{eq:SG-1}
  \SG(u) (\tau, \xi) := \int_{x\in\RR} u(x) \theta_x(\tau, \xi)
 dx
 \ee
is a well-defined section of $\Oc_{S\Hen}$. Since $\theta_x(\tau, 0) = e^{\pi \i x^2\tau}$
is an ordinary Gaussian in $x$, the restriction $\G(u)(\tau) = \SG(u)(\tau, 0)$ 
is the familiar Gaussian transform, i.e., integration against a
variable Gaussian   \cite[\S 3]{guillemin}
\[
G(u)(\tau) = \int_{x\in \RR} u(x) e^{\pi \i x^2\tau} dx. 
\]
It vanishes on $\Wc_\1$ since Gaussians are even functions in $x$. 
The $\xi$-component of $\SG$ is the integration against a linear function
times a variable Gaussian and so vanishes on $\Wc_\0$ and picks up $\Wc_\1$. 
Our approach combines the two procedures into a single supersymmetric one. 
 Let us consider $\SG(u)$ as a section of $S\Lc_{1/2}$, i.e., as an element of $\SF_{1/2}$
 using the trivialization of $S\Lc_{1/2}$.

 \begin{prop}
 The  super-Gaussian transform considered as an operator
$ \SG: \Wc \lra \SF_{1/2}$,
  is an $\FT$-twisted  injective morphism of $\OMp(1|2)$-representations,
  i.e.,
  \[
  \SG(\varpi(g) f) \,=\, \FT(g) (\SG(f)). 
  \]
 \end{prop}
 
 This extends the familiar fact that $\G: \Wc \to \F_{1/2}$ is an
 (also $\FT$-twisted, in our conventions)  morphism of 
 $\Mp(2)$-representations  vanishing on $\Wc_\1$ and injective on $\Wc_\0$,
 see \cite {guillemin}. 
 
 \vskip .2cm
 
  \noindent{\sl Proof:}   To prove that $\SG$ is an $\FT$- twisted morphism of $\OMp(1|2)$-representations,
  it is enough  to verify it  at the level of Lie (super) algebras,
 i.e., to show that $\SG$ commutes (after applying $\FT$) with the action of $p,q\in\gen_\1$ in the source
 and target. For this, it suffices to show such commutativity for the $\theta_x$,
 i.e., that
 \be\label{eq:theta_x-comm-gen1}
(p \theta_x) (\tau, \xi) = D_p^{(1/2)} \theta_x(\tau, \xi), \quad (q\theta_x)(\tau, \xi)
= D_q^{(1/2)}\theta_x(\tau, \xi),
 \ee
 where $p,q$ act  by \eqref{eq:osp(1\2)-weil}  w.r.t.  the dependence on $x$ while
 $D_p^{(1/2)}, D_q^{(1/2)}$ act  w.r.t  the dependence on $\tau, \xi$.
 But this folllows from Proposition \ref {prop:super-Gauss-1}. 
 
 Once we know that $\SG$ is an $\FT$-twisted
  morphism of representations, it must be
 injective since $\Wc$ is irreducible. \qed

 \begin{rem}\label{rem:sgauss-BW-1}
 The conceptual meaning of  the relations \eqref{eq:theta_x-comm-gen1}
 is that 
  $\theta_x(\tau,\xi)$ considered as a function of $x$, i.e., as an element
 of $\Wc$, 
 is a highest vector  for  (the $\FT$-image of) the parabolic subalgebra 
 in $\gen$ corresponding to $(\tau, \xi)$. The (super) Gaussian transform is therefore just an infinite-dimensional
 version of the classical Borel-Weil procedure which realizes an
 irreducible   finite-dimensional
 representation $V$ of a complex semisimple group as the space 
 of sections of a line bundle $\Lc_V$ on the flag variety $F$. That is,
 the fibers of $\Lc_V$ are the duals of the  $1$-dimensional spaces of highest vectors
 (functionals) in $V^*$ for the Borel subgroups corresponding to points of $F$
 and a vector $v\in V$ gives a section of $\Lc_V$ by pairing with these
 highest functionals. 
  \end{rem} 
 
 \paragraph{The $\RR$-holonomic system in terms of the super-Weil representation.}
 \label{par:Weil-compar-gen1}
 Like the usual Gaussian transform, $\SG: \Wc\to \SF_{1/2}$ is not
 surjective. Elements of $\SF_{1/2}$ can be loosely interpreted
 as corresponding to some (distribution) completion of $\Wc$
 containing, for example, delta-functions   which give the $\theta_x(\tau,\xi)$
 themselves. Using this informal  interpretation, we can
 reformulate the  system of differential operators of infinite order
 from Theorem \ref {thm:1d-theta-charact}(b) for $\Phi (\tau, \xi) \in \SF_{1/2}$
 as a system of equations on a distribution $u(x)$ on $\RR$:
 \be
 u(x+1) = u(x), \quad e^{2\pi \i x} u(x) = u(x)
 \ee
 which gives $u(x) = c\cdot \sum_{n\in \ZZ} \delta(x-n)$,
 so $\SG(u) = c\cdot \theta(\tau)$. 
 
 \vskip .2cm
 
 Thus the super-Gaussian transform changes the non-local  difference equation
 $u(x+1)=u(x)$ into a local differential equation
 of infinite order  $e^{\sqrt{- 4\pi \i } D_q^{(1/2)}} \Phi = \Phi$. 
 Note that even this informal  comparison  becomes possible
 only after adopting the supersymmetric viewpoint and cannot
 be done within the framework of the standard Weil representation .

 \vfill\eject
 
 \section{Genus 2: Three-dimensional $\Nc=1$ superconformal symmetry}\label{sec:gen2}
 
 \subsection{Conformal spaces and   conformal Laplacians}\label{subsec:conf-lap}
 
 We recall some background material, see, e.g., \cite[\S3.1]{bryant}, in the complex analytic framework.

 \paragraph{ Conformal spaces and conformal densities.}
 Let $V$ be a finite-dimensional complex vector space of dimension $d\geq 2$. A  (linear)  {\em conformal structure} on $V$
 is an equivalence class of non-degenerate scalar products $g\in S^2(V^*)$ modulo the scaling
 equivalence: $g\sim \lambda g$, $\lambda\in \CC^*$.  Such a class is uniquely defined by its null-cone
 $K\subset V$ consisting of $v$ such that $g(v,v)=0$ for any $g$ from the class. Indeed, given $K$, the relevant $g$
 are precisely the  irreducible equations of $K$. 
 
 \begin{defi}\label{def;conf-space}
 Let $X$ be a complex manifold of dimension $d\geq 2$. A (holomorphic)   {\em conformal structure} on $X$ is a datum of
 a  linear conformal structure in each tangent space $T_xX$, depending holomorphically on $x$. 
 A {\em conformal manifold} is a complex manifold with a conformal structure. 
 \end{defi}
 
 \begin{prop} A conformal structure on $X$ can be described in three equivalent ways:
 \begin{itemize}
 \item[(i)] As a datum of  holomorphic Riemannian metrics $g = g(x)$  on $X$ given locally, near each point $x_0\in X$
 and defined up to variable holomorphic rescaling: $g(x) \sim \lambda(x) g(x)$, $\lambda(x) \in \Oc_X^*$.
  Note that there may not be
 a global Riemannian metric on all of $X$ defining a given conformal structure. 
 
 \item[(ii)] As a datum of a line bundle $\LL$ on $X$ and a global holomorphic  symmetric bilinear pairing $g: S^2(TX) \to \LL$,
 which is everywhere non-degenerate in the obvious sense. 
 
 \item[(iii)] As a datum of an analytic  hypersurface $K\subset TX$ such that for each $x\in X$
 the intersection $K_x = K\cap T_xX$ is a non-degerate quadratic cone. 
 \end{itemize}
 \end{prop}
 
  Using the datum (ii) above we will sometime write a conformal manifold $X$ as $(X, \LL, g)$
  and call $g$ the {\em conformal metric}. 
  
  \vskip .2cm
 
 \noindent {\sl Proof:} Indeed, (i) is Definition \ref {def;conf-space} while (iii) comes from the description of linear conformal structures
 in terms of null-cones as above. Finally, given $K$ considered as a family of the $K_x$, we
 recover $\LL$ in (ii) uniquely: the fiber $\LL_x$ of $\LL$ at $x\in X$ is the subspace in $S^2(T^*_xX)$
 formed by $0$ and all the equations of $K_x$. \qed
 
 \vskip .2cm

 By {\em conformal densities of weight $r$} on a conformal manifold $(X,\LL, g)$ we will mean sections (local or global)
 of $\LL^{\otimes r}$. 
 
 \begin{prop}\label{prop:conf-dens-omega}
 For a conformal manifold $(X,\LL, g)$ of dimension $d$ we have an isomorphism of line bundles
 $\LL^{\otimes d} = \omega_X^{\otimes -2}$, where $\omega_X = \Omega^d_X$ is the canonical bundle of $X$.
 In other words, $\LL$ is a determination of $\omega_X^{\otimes -2/d}$.
 \end{prop}
 
 \noindent In other words, 
  conformal densities can be always seen as (possibly fractional) powers of volume forms. 
 
 \vskip .2cm
 
 \noindent{\sl Proof:} This rephrazes the classical construction of the  volume form
$\sqrt{g}\,  dx$ associated to a Riemannian metric $g$. \qed 
 
 \paragraph{Example: the conformal quadric.} \label{ex:conf-quadr}
 Let $q$ be a non-degenerate quadratic form on $\CC^{d+2}$. The equation $q(x)=0$ defines the projective quadric
 $Q^d \subset \PP^{d+1} = \PP(\CC^{d+2})$. It has a natural conformal structure defined as follows. 
 For any $x\in Q$ let $\PT_xQ\subset  \PP^{d+1}$ be the projective subspace (hyperplane) in $\PP^{d+1}$
 tangent to $Q$ at $x$. The intrinsic tangent space $T_x(\PT_xQ)$ is identified with $T_xQ$. 
 On the other hand, the intersection $\PK_x := Q\cap \PT_xQ$ inside $\PP^{d+1}$,  a hypersurface in $\PT_xQ$,
 is a quadratic cone with one singular point, namely $x$. This means that we have a canonical
 quadratic cone $K_x \subset T_xQ$ which corresponds to the tangent cone at $x$ to $\PK_x$
 under our identification $T_x (\PT_xQ) =  T_xQ$. 
 
 We will call $Q$ with the conformal structure given by these $K_x$ the {\em $d$-dimensional conformal
 quadric}. Its conformal structure  is invariant under the orthogonal group $\SO(d+2, \CC)\subset \GL(d+2, \CC)$ preserving $q$.

 Let us choose  one point  $\oo\in Q$ as the  ``infinity''; as a point in $\PP^{d+1}$, it corresponds
 to a $1$-dimensional isotropic subspace $l_\oo\subset \CC^{d+2}$. The  complement
 $Q-K_\oo$ (open Schubert cell) is identified with a $d$-dimensional complex affine
 space by 
  (stereographic) projection
 from $\oo$:
  \be\label{eq:stereo}
 \sigma: Q- K_\oo \buildrel \sim \over\lra \,\, \PP(\CC^{n+2}/l_\oo) - \PP (l_\oo^\perp/l_\oo).  
 \ee
 More precisely, the target of $\sigma$ is an affine space over the vector space
 $\Hom(\CC^{n+2}/l^\perp_\oo, l_\oo^\perp/l_\oo)$. It carries a canonical (constant) conformal
 sructure (since $l_\oo^\perp/l_\oo$ carries a non-degenerate scalar product induced by $q$)
 and $\sigma$ is a conformal isomorphism. 
    In this way we realize $\SO(d+2, \CC)$ as the group of  conformal (birational)  
   transformations of  the flat 
 ``complex Minkowski space'' $\CC^d$.
 
 \begin{prop}\label{prop:conf-dens-quad}
 The line bundle $\LL$ for the conformal structure on $Q^d$ just described, is identified with
 $\Oc_{Q^d}(2)$.
 \end{prop}
 
 \noindent {\sl Proof:} Let $l\subset \CC^{d+2}$ be a $1$-dimensional subspace representing a point
 $[l]\in\PP^{n+d}$. It is standard (``Euler sequence'') that $T_{[l]}\PP^{d+1}$ is canonically identified
 with $l^*\otimes (\CC^{d+2}/l)$. If $[l]\in Q^d$, i.e., $l$ is isotropic, then $T_{[l]}Q^d = l^* \otimes(l^\perp/l)$. 
 Now, $l^\perp/l$ has a nondegenerate quadratic form induced by $q$, so as a representation of the
 stabilizer of $[l]$ in $\SO(d+2,\CC)$, the second symmetric power $S^2(l^\perp/l)$ splits into the sum of  two
 irreducibles as $\CC \oplus S^2_0(l^\perp/l)$,  the trace and the traceless part.
 Thus the only invariant subspace in $S^2(l^* \otimes (l^\perp/l)) = S^2(L^*) \otimes S^2(l^\perp/l)$
 is $S^2(l^*)\otimes \CC\,  \= \, S^2(l^*)$ which is precisely the fiber of $\Oc_{\PP^{d+1}}(2)$ at $[l]$. 
 \qed

 \paragraph{Conformal transformations of $\CC^d$ as Clifford-fractional-linear transformations.}  \label{par:conf=fractional}
 A beautiful explicit formula for general conformal  transformations $\CC^d\to \CC^d$
 was given by  Vahlen, Maass and Ahlfors \cite{vahlen, maass, ahlfors1, ahlfors2, waterman} 
 in terms of   the spinor covering $\Spin(d+2, \CC)$.   
 
   \vskip .2cm

  Let us represent $V=\CC^d$ with the constant scalar product $g$
  as an orthogonal  direct sum  $\CC \oplus \CC^{d-1}$ so that $g = 1 \oplus (-b)$,
  where $b$ is a scalar product on $\CC^{d-1}$. 
    Let
  $\Ac$   be the Clifford algebra of $b$, of dimension $2^{d-1}$ over $\CC$,  generated by
   symbols $E_2,\cdots, E_d$ corresponding to a basis $e_2,\cdots, e_d \in \CC^{d-1}$ and 
   subject to the relations $E_iE_j+E_jE_i = b(e_i, e_j)$. It has:
    \begin{itemize}

   \item  A $\CC$-linear anti-involution
   $x\mapsto \ol x$ defined by $\ol E_i = -E_i$. The space $  \Ac^{\leq 1}$
   of elements of the form $t = t_1 + \sum_{i=2}^d t_i E_i$, $t_i\in \CC$,  is identified with $V=\CC^d$,
   with the scalar product  $g$ recovered as $g(t, t')=   t \cdot \ol {t'} =  t'\cdot \ol t$ which in
   this case lies in $\CC\cdot 1_\Ac$. 
      
     \item  A $\CC$-linear anti-involution
  $x\mapsto x^*$ defined by $E_i^* = E_i$. 
  
     \item A $\CC$-linear projection $\Ren: \Ac\to \CC$ (``Clifford real part'') taking $1_\Ac$ to $1$
    and any totally antisymmetrized product of several $E_i$ to $0$.  The
     {\em (squared) Clifford norm }
    of $a\in\Ac$ is defined as $|a|^2 := \Ren(a\cdot \ol a) = \Ren( \ol a\cdot a)$, which extends
    $g(t,t)$  defined for $t\in V$. 
   
   \end{itemize}

 \noindent 
 The {\em Clifford group} $\Gamma\subset \Ac^*$ consists of $a$ such that the transformation
 $f_a: x\mapsto ax(\ol a^*)^{-1}$ preserves $V$. Such a transformation is orthogonal and
 $\Gamma$ is identified with the group $\Spin^\CC(d, \CC)$ fitting into the central extension 
 \[
1\to \CC^* \lra \Gamma \to \SO(d, \CC) \to 1
\]
 Any non-isotropic ``vector'' $t\in V$ lies in $\Gamma$ and $f_t$ is the reflection w.r.t. $t$.
 It us known that $\Gamma$ is generated by all non-isotropic $v\in V$. The closure
 $\ol\Gamma$ is a semigroup generated by $V$.  The correspondence $t\mapsto |t|^2$
 extends to a homomorphism $\Gamma\to \CC^*$. We put
 \be
 \GL_2'(\Ac) = \biggl\{ \gamma = \bpm a&b\\c&d  \epm\in\Mat_2(\Ac)  \biggl| \, 
 \begin{matrix}
 a,b,c,d\in \ol\Gamma, \,\,\, ad^*-bc^*\,\, \in\,\, \CC^*\cdot 1_\Ac, 
 \\
  ab^*, cd^*, c^*a, d^* b\,\,  \in\,\,  V
 \end{matrix} \biggr\}. 
 \ee
 
 The following is a summary (and a complexified version) of the results of 
  \cite{vahlen, maass, ahlfors1, ahlfors2, waterman}. 
  
  \begin{thm}\label{thm:ahlfors}
  (a) The set $\GL_2'(\Ac)$ is a group and $\Det: \gamma \mapsto ad^*-bc^*$ is a homomorphism
  $\GL_2'(\Ac)\to \CC^*$. 
  
  \vskip .2cm
  
  (b) For $\gamma\in \GL_2'(\Ac)$ the rational  map
  \[
   \Phi_\gamma: V \to \Ac, \,\,\, f_\phi (t)  =  (at+b)(ct+d)^{-1}
  \]
  takes values in $V$ and defines a conformal (birational) map $\Phi_\gamma: V\to V$. This gives a surjective
  homomorphism $\GL'_2(\Ac) \to \SO(d+2, \CC)$ with kernel $\CC^*$ thus identifying
  $\GL'_2(\Ac)$ with $\Spin^\CC(d+2, \CC)$. 
  
  \vskip .2cm
  
  (c) The group $\SL_2'(\Ac)=\Ker(\Det)$ is identified with $\Spin(d+2, \CC)$. 
  
  \vskip .2cm
  
  (d) Conformal densities of weight $r$ transform under $\gamma\in \SL'_2(\Ac)$ as above as
\[
(\gamma^* f)(t) = |ct+d|^{-2r} f\bigl( (at+b)(ct+d)^{-1}\bigr).
\]
 \qed
  \end{thm}

 \paragraph{The conformal Laplacian.} Let $(X,\LL, g)$ be a conformal manifold of dimension $d>1$. The
 {\em conformal Laplacian} of $X$ is a canonical differential operator
 \be\label{eq:conf-lap}
 \Delta: \LL^{\otimes {2-d\over 4} }
   \lra \LL^{\otimes {-2-d \over 4}  } 
 \ee
 defined as follows 
  \cite[p.108]{fefferman}. Choose a local trivialization of $\LL$, 
  so $g$ produces an actual Riemannian
  metric $\wt g$ on (an open set of) $X$. Then, viewing sections of $\LL^{\otimes {2-d\over 4} }$
  and  $\LL^{\otimes {-2-d \over 4}  } $, put 
 \[
  \Delta = \Delta_{\wt g} - {d-2 \over 4(d-1)} R_{\wt g}, 
  \]
  where $\Delta_{\wt g}$ is the usual Laplace-Beltrami operator associated to $\wt g$ and $R_{\wt g}$
  is the scalar curvature of $\wt g$. It is then verified that $\Delta$ thus defined transforms under
   changes of the trivialization as an operator in \eqref{eq:conf-lap}. See also \cite[\S3.1.3]{bryant} for
   a more conceptual way to obtain
 $\Delta$. 

\vskip .2cm
  
 If  $d\not\equiv 2 \text{ mod } 4$,  then $\Delta$ acts between fractional powers of $\LL$
 and should be understood as explained in \S \ref {subsec:twisted}
 (i.e., as acting between any  local determinations
 of such fractional powers).  

 \begin{exas}\label{ex:conf-lap-Q}
 (a) For the flat case ($X=\CC^d$, $\LL=\Oc$, $g = g_\st := \sum (dz_i)^2$ constant), $\Delta$ coincides with the usual
 Laplacian $\Delta_\st = \sum \del^2/\del z_i^2$. In view of Proposition \ref {prop:conf-dens-omega}
 this gives 
  the following classical fact: if $u(z)$ is a (local) holomorphic solution of  $\Delta_\st u=0$
  and $\Phi: \CC^d\to\CC^d$ is a conformal map, then
  $J(\Phi)(z)^{d-2\over 2d} u(\Phi(z))$ is again  annihilated by $\Delta^\st$. Here $J(\Phi)(z)$ is the
  Jacobian of $\Phi$.

 \vskip .2cm
 
 (b)  The conformally flat case when $X=Q=Q^d$
   has been discussed already by Dirac \cite{dirac}, see
   also \cite{eastwood}.  As in this case $\LL = \Oc_{Q}(2)$, the conformal Laplacian acts as 
   \[
   \Delta: \Oc_Q\biggl(1-{d\over 2}\biggr)) \lra \Oc_Q\biggl( -1-{d\over 2}\biggr). 
   \]
   
   (c) The case $d=1$   should formally be excluded from the above discussion
   as the formula for the conformal Laplacian has vanishing denominator.
    Nevertheless, we do have a Laplace-type operator
   on the $1$-dimensional quadric $Q=Q^1$ which acts as 
  $\Delta:\Oc_Q(1/2) \to \Oc_Q(-3/2)$ in accordance with Dirac's  general formula from (b). 
  Indeed, in this case 
    $Q$ is a conic in $\PP^2$, isomorphic to $\PP^1$ as an algebraic
    curve, and $\Oc_Q(1) = \Oc_{\PP^1}(2)$. In this case  we do have an invariant differential operator
    $\Delta: \Oc_{\PP^1}(1) \to \Oc_{\PP^1}(-3)$ whose kernel consists of global sections
    of $\Oc_{\PP^1}(1)$ (given by $d^2u/dz^2=0$ in an affine chart $\CC=\AAA^1\subset \PP^1$). 
    The ``conformal group'' is now $\SO(3,\CC)=\SL(2, \CC)/\{\pm 1\}$ acting by 
   $\Phi(z) = (az+b)/(cz+d)$. 
    The action of Example (a) in this case has a factor $J(\Phi)(z)^{-1/2} = cz+d$, 
 so we get the standard action of $\SL_2(\CC)$ on the
    space of linear functions $c_0+c_1z$.

 \end{exas}
 
 \paragraph{The conformal Dirac operator.} The conformal invariance  properties of the Dirac operator
 were pointed out by Dirac himself \cite{dirac} for the conformal quadric and by Hitchin  
 \cite[Prop.1.3]{hitchin}  and Lott \cite[Prop.2] {lott} in general. Here is a brief adaptation of their
 results to the complex analytic case.
 
 \vskip .2cm
 
 Let $\CO(d, \CC)\subset \GL(d,\CC)$ be the complex Lie group generated by $\SO(d)$ and by dilations
 $\CC^* \cdot 1_d$. Thus a $d$-dimensional conformal manifold has a $\CO(d,\CC)$-structure. 
 Denote by $p: \CSpin(d, \CC) \to \CO(d,\CC)$ be the $2:1$ covering  which induces
  $\Spin(d,\CC)\to \SO(d,\CC) $ over $\SO(d,\CC)$
 and the   $2:1$ covering $\CC^*\to \CC^*$ over $\CC^*\cdot 1_d$.
 As well  known, $\Spin(d,\CC)$ has one irreducible  spinor representation $M$ for $d$ odd and two such representations $M_\pm$ for $d$ even. They  extend to representations of $\CSpin(d,\CC)$
 uniquely by the requirement that $p^{-1}(c\cdot 1_d)$ acts by $\pm\sqrt{c}$. 
 
 \vskip .2cm
 
 A {\em spin structure} on a $d$-dimenisonal conformal manifold $(X,\LL,g)$ is a lift of the structure group from $\CO(d)$
 to $\CSpin(d)$. Assuming a spin structure is given, we form the 
 {\em conformal spinor bundle}  $V$ associated to $M$
 if $d$ is odd and two conformal spinor bundles $M_\pm$ associated to $V_\pm$ if $d$ is even. 
 To simplify notation (and  because we will need only the odd case in  this paper)
 let us also  write $M= M_+\oplus M_-$ and  $V=V_+\oplus V_-$ for $d$ even. 
 The gamma-pairing
 is then $\Gamma: V \otimes V\to TM$.

 \begin{prop}\label{prop:dirac-gen}
 Let $(X,\LL, g)$ be a $d$-dimensional conformal manifold with a spin structure. The Dirac
 operators associated to metrics in the conformal class of $X$  (given by  trivializations of $\LL$), unite to give a canonical
 differential  operator
 \[
 \dirac:  V \otimes \LL^{\otimes{-d\over 4}} \lra V \otimes \LL^{\otimes {-2-d\over 4}}.
 \]
  
 \end{prop}
 
 \noindent {\sl Proof:} By  Let $g(x)$ and $g'(x) = \lambda(x) g(x)$ be two conformally
 equivalent metrics and $\dirac_g, \dirac_{g'}$ be the corresponding Dirac operators. By \cite[Prop.2] {lott}
 (proved in the real $C^\oo$ case but translating to the complex analytic case without changes)
 \[
 \lambda^{d+1\over 4}   \dirac_{g'}(\psi)  = 
 \dirac_g(\lambda^{d-1\over 4} \psi)
 \]
 where $\psi$ is any section of the spinor bundle $V_g$ for $g$  which  is considered  ``the same''
 as the spinor bundle $V_{g'}$. This identification is based of viewing $M$ as a representation not
 of $\CSpin(d,\CC)$ where a central element from $p^{-1}(c)$ acts by a $\sqrt{c}\, $ but
 of $\Spin(d,\CC)\times \CC^*$ where the second factor acts trivially.
 Informally,  this is  the standard ``metric''  version of  the spinor bundle which
 does not have any conformal weight built into it (unlike our $V$ which does). 
 
 \vskip .2cm
  
  In other words, see Rem.1
 in \cite[\S1.4]{hitchin}, the conformally invariant Dirac acts $\wt V \to \wt V \otimes L^*$,
 where:
 \begin{itemize}
 \item  $\wt V$ is associated to the representation of $\Spin(d,\CC)\times \CC^*$ taking
 $(g, \lambda)$ to $\rho(g) \lambda^{d-1\over 2}$, where $\rho$ is the usual spin representation. 
 
 \item $L$ is the line bundle whose sections, denoted $e^\sigma$,  scale the metric as $e^{2\sigma}$,
 so $L = \LL^{\otimes{1\over 2}}$ in our notation. 
  \end{itemize}
  
  \noindent  Comparing these two pictures using the isogeny $\Spin(d,\CC)\times \CC^*\to \CSpin (d,\CC)$,
  we see that Hitchin's  $\wt V$   corresponds to our $V \otimes \LL^{\otimes{-d\over 4}}$ and
 his  $\wt V\otimes L^*$  corresponds to our $V \otimes \LL^{\otimes {-2-d\over 4}}$.
  
   \qed

 \begin{ex}
 Let $d=1$, so $T:=TX$ is a line bundle,  $\LL=T^{\otimes 2}$ and $g=\Id$. The group $\CO(1,\CC)$ is just $\CC^*$, so
 $V=T^{\otimes {1\over 2}}$. In this case the proposition gives the operator
 \[
 \dirac: T^{\otimes{1\over 2}}\otimes T^{\otimes {-1\over 2}} \lra 
 T^{\otimes{1\over 2}}\otimes T^{\otimes {-3\over 2}}
 \]
 in which we recognize  the de Rham operator $d: \Oc_X\to\Omega^1_X$.  
  \end{ex}

 \paragraph{Real structure. The future tube.} \label{par:future}
 Now consider the real Minkowski space
 $\RR^{1, d-1}$ with coordinates $x_1,
 \cdots, x_d$ and  the quadratic form $g(x) = x_1^2 -x_2^2-\cdots-x_d^2$
 which we consider as a constant Riemannian metric of signature $(1, d-1)$. Let 
 \[
 C^+ = \{ x\in \RR^{1, d-1}: g(x) >0, \, x_1 > 0\} \subset \RR^{1, d-1}
 \]
 be the positive future cone. The {\em future tube} of dimension $d$ is the domain
 \[
 \Ten_n = \{ z\in   \CC^n = \RR^{1, n-1}\otimes_\RR \CC \, | \, \Im(z) \in C^+\},
 \]
 see \cite{sergeev}. 
 It can be embedded into the complex quadric $Q^d\subset \PP^{d+1} = \PP(\CC^{d+2})$
 associated to the quadratic form 
 \[
 q(w_0, \cdots, w_{d+1}) = w_0^2  + w_1^2 - w_2^2 -  \cdots  -   w_{d+1}^2 . 
 \]
 of signature $(2,d)$. Let $b(w, w')$ be the symmetric bilinear form associated to $q$
 and $h(w, w') = b(\ol w, w')$ the corresponding Hermitian form. 
 Let
 $\Den_d\subset Q^d$ be the open part consisting of $1$-dimensional isotropic
 complex lines $l$ such that the $h|_l <0$. Then, a version of the  stereographic projection from
 \S \ref{ex:conf-quadr} 
 identifies $\Den_d$ with $\Ten_n$, see \cite[\S A6]{satake} or
 \cite[Prop.6.1.7]{huybrechts} for more details.
 In this way $\Ten_d$ is realized as a bounded symmetric
 domain of type IV. 
 
 \vskip .2cm
 
  The  action of the real group $\SO(2,d)$ on $Q^d$ preserves  $\Den_d$
 and so $\SO(2,d)$  acts on $\Ten_d$.
  In terms of 
 $\Spin(2,d)$ this action can be obtained by taking the real version of
  the discussion of \S \ref{par:conf=fractional}. That is, we consider the
   Clifford algebra $\Ac_\RR$ over $\RR$ generated by anticommuting $E_2, \cdots, E_d$
   subject to $E_\nu^2 = +1$ and form the corresponding groups $\GL'_2(\Ac_\RR)$
   and $\SL'_2(\Ac_\RR)$, the latter being identified with  $\Spin (2,d)$.
   The Clifford-fractional-linear transformations $\Phi_\gamma$ for $\gamma\in \SL_2'(\Ac_\RR)$
   restrict to biholomorphic transformations of 
   $\Ten_d$   and 
  this action identifies
 \be\label{eq:Td-quot}
\Ten_d  \simeq \SO(2,d)/(\SO(2)\times \SO(d)). 
\ee

    \begin{exas}
  (a) For $d=1$  we have $\Ten_1 = \Hen\subset \CC$, the Lobachevsky plane,
  with $\Ac_\RR = \RR$ and $\SL'_2(\Ac_\RR) = \SL(2,\RR)$ acting on 
 $\Hen$ in the standard way. 
 
 \vskip .2cm
 
 (b) For $d=2$ we have $\Ten_2=\Hen\times\Hen$, while $\Ac_\RR = \RR[E_2]/(E^2-1)=\RR\oplus \RR$,
 so $\SL'_2(\Ac_\RR) = \SL(2, \RR) \times\SL(2, \RR)$ acting on  $\Hen\times\Hen$ in the standard
 way.
 
  \vskip .2cm
  
  (c) For $d=3$  we have $\Ac_\RR =\Mat_2(\RR)$; we will discuss this case in the next section. 
 
  \end{exas}

  \noindent For $d>1$ the transformations from $\SO(2,d)$ are, 
   moreover, conformal with respect to the    conformal
  structure induced by the embedding $\Ten_d\subset \CC^d$ from the flat conformal
  structure given by $g$.  By the invariance of the conformal Laplacian we have the
  following classical fact.
  
  \begin{prop} Let $d\geq 2$. 
   The action   of $SO(2,d)$ on $\Ten_d$, $\gamma\mapsto \Phi_\gamma:\Ten_d\to\Ten_d$
    extends to an action on the 
     space  of global holomorphic solutions 
  in $\Ten_n$ of the
  Laplace equation
  \[
  \Delta u=0, \quad \text{where} \quad \Delta = {\del^2\over\del z_1^2} -
  {\del^2\over\del z_2^2} - \cdots - {\del^2\over\del z_d^2}
  \]
  is the standard Laplacian corresponding to the quadratic form $g$. 
  Explicitly, $\gamma\in \SO(2,d)$ acts by the rule
  \[
u\mapsto   \gamma^*u, \quad (\gamma^*u)(z) = J(\Phi_\gamma)(z)^{d-2\over 2d} u(\Phi_\gamma(z)).
  \]
  \qed
  \end{prop}
  
  This action was studied in \cite {jacobsen} for $d=2$ and  in  \cite{kobayashi}  in general. 
 \vfill\eject

 \subsection{The genus $2$ Siegel plane as the $3$d future tube}
 \label{subsec:Siegel2}
 
 The exceptional isomorphism of complex Lie algebras $\so(5) = \sp(4)$ gives a connection between
 conformal geometry in $3$d and symplectic geometry in $4$d. 
 
 \paragraph{The genus 2 Lagrangian Grassmannian as  the  conformal quadric.}\label{par:conf-dens-quad}
  Consider the  $4$-dimensional symplectic vector space $E=\CC^4$ with basis $p_1, p_2, q_1, q_2$
  and the symplectic form $\omega\in \Lambda^2 E^*$ given by 
  \[
   \omega (p_i,q_j) = \delta_{ij}, \quad  \omega (q_i,q_j) = \omega (p_i,p_j) = 0.
   \]
   The linear coordinates dual to this basis will be denoted
  $y_1, y_2, x_1,  x_2\in E^*$. We will arrange these coordinates into
  column vectors $y=(y_1, y_2)^t$, $x= (x_1, x_2)^t$. 
  
  \vskip .2cm
  
   The {\em Lagrangian Grassmannian} $\LG(E)$ is the space of all isotropic $2$-planes
    $L\subset E$.  We denote by $[L]$ the point of $\LG(E)$ represented by $L$.
    By definition, $\LG(E)$ is embedded in $G(2,E)$, the Grassmannian fof all $2$-planes.
   We notice that  $\LG(E)$
     is identified with $Q^3$, the $3$-dimensional
   projective quadric. Indeed, $G(2,E)$ is realized inside $\PP(\Lambda^2 E) = \PP^6$ as
   the Pl\"ucker quadric $Q^4$, and $\LG(E)$ is its hyperplane section given by $\omega=0$. 
   So $\LG(E)$ has a natural conformal structure and 
   the group $\Sp(E)= \Sp(4,\CC)$ of symplectic transformations of $E$
   acts on it by conformal transformations  which, as we have seen,  form the group $\SO(5,\CC)$. 
   This gives  an identification $\Sp(4,\CC) \= \Spin(5,\CC)$. 
   
   \vskip .2cm
   Another way to see the conformal structure is as follows. 
   Let  $\Lambda$  be the tautological rank $2$ bundle on $\LG(E)$ whose fiber over $[L]$
   is $L$. Then $T\LG(E)$, the tangent bundle of $\LG(E)$, is naturally identified with $S^2 \Lambda^*$.
   A fiber of this bundle at $[L]$ is   $S^2L^*$, the space of quadratic forms on $L$ and so contains
   the quadratic cone  $K_{[L]}$ formed by degenerate quadratic forms thus acquiring a conformal structure.
   Note that $\bigwedge^2 \Lambda \, \= \, \Oc_{\LG(E)}(-1)$, the restriction of $\Oc(-1)$ from the ambient
   projective space, and so  by Proposition \ref {prop:conf-dens-quad} the line bundle of conformal densities is 
   $\LL = \Oc_{\LG(E)}(2) = (\bigwedge^2\Lambda^*)^{\otimes 2}$. The canonical $\LL$-valued Riemannian
   metric is the invariant projection
   \[
   g: S^2(T\LG(E)) = S^2(S^2\Lambda^*) \lra \bigl(\bigwedge\nolimits^2 \Lambda^*\bigr)^{\otimes 2} = \LL.  
   \]

   \paragraph{The open Schubert cell as the $3$d complex Minkowski space.}
   
    $\LG(E)$ has  two distinguished points
 \[
 L_q=\bigoplus \CC q_i ,\,  L_p= \bigoplus \CC p_i \,\subset \,  E.
 \]

The {\em open Schubert cell} $\LG^\circ(E)$ consists of $L$ which are transversal to $L_p$. 
 Each $L\in \LG^\circ(E)$ is a graph of a self-adjoint operator 
 $T: L_q \to L_p = L_q^*$ which
  we identify  with its  (symmetric) $2\times 2$
 matrix $T=\|t_{ij}\|$   in the standard bases of $L_q$ and $L_p$. 
 Thus $\LG_n^\circ(E)$ is identified with 
 $\Sym_{2\times 2}(\CC)$, the space of symmetric
 $2\times 2$ matrices.    Using coordinate column vectors  $x,y$ as above, 
 the Lagrangian subspace associated
 to $T\in \Sym_{2\times 2}(\CC)$  given by the equations
 \be\label{eq:space-LT-2}
 L_T \,=\, \bigl\{ (x,y) \bigl| \, y=Tx, \text{i.e., } y_i = \sum_{j=1}^2  t_{ij} x_j \bigr\}. 
 \ee
 The space $\Sym_{2\times 2}(\CC) = \CC^3$ carries the quadratic form $\det(T)$,
 which gives a flat conformal structure. In other words, we represent $\CC^3=S^2(\CC^2)$,
 so that $\CC^2$ is the space of spinors for $(\Sym_{2\times 2}(\CC), \det)$. 
 It is immediate from the discussion in \S \ref {prop:conf-dens-quad}
 that the map $T\mapsto L_T$ is a conformal map from $\Sym_{2\times 2}(\CC)$
 with this flat structure into $\LG(E)$. 
 In fact, this map is an instance for $d=3$ of the inverse of the stereographic 
 projection \eqref{eq:stereo}. 
 
 \vskip .2cm

 As $\Sym_{2\times 2}(\CC)$ is a linear subspace in $\Mat_{2\times 2}(\CC)$,
 we will understand  functions of  the $t_{ij}$ in their restriction to
  $\Sym_{2\times 2}(\CC)$: for example, $t_{12}=t_{21}$ in such restriction.
  As for derivatives, we will use the
 expressions
  \be
  \del_{ij} = {\del\over \del t_{ij}} + {\del\over \del t_{ji}}
  %\quad 1\leq i,j\leq 2 
    \ee
  which are vector fields on  $\Mat_{2\times 2}(\CC)$
tangent to   $\Sym_{2\times 2}(\CC)$, so applying them  to functions on 
$\Mat_{2\times 2}(\CC)$ commutes with restriction.  Note that $\del_{ij}$
is not the same as ``$\del/\del t_{ij}$ where $t_{ij}=t_{ji}$''
   For example,
  $\del_{11}(t_{11})=2$ while $\del_{12}(t_{12})=1$ on $\Mat_{2\times 2}(\CC)$
  as well as on    $\Sym_{2\times 2}(\CC)$.  Obviously, $\del_{ij}=\del_{ji}$.

 \vskip .2cm

 Note that $g\in\Sp(4,\CC)$ whose form in the basis $p_1, p_2, q_1, q_2$
 has the form $g =\begin{pmatrix} A&B\\ C&D\end{pmatrix} $ acts on the coordinate
 vectors $x$ and $y$ by 
 \be
 g : (x,y) \mapsto \bigl(x' = Cy+Dx,\quad  y'=Ay+Bx\bigr). 
 \ee
 So the action on $\LG^\circ(E)$ is given by the standard Siegel fractional-linear formula
 \be\label{eq:act-Siegel}
 g (L_T) = L_{\Phi_g(T)}, \quad \Phi_g (T) = (AT+B)(CT+D)^{-1}. 
 \ee
 This is in fact a $d=3$ instance of the  formula of Theorem \ref{thm:ahlfors}. In this case
 the Clifford algebra $\Ac$, of dimension $4$, is isomorphic to $\Mat_2(\CC)$, the space $V=\Ac^{\leq 1}$,
 of dimension $3$, is identified with $\Sym_{2\times 2}(\CC)$ and the group $\SL'_2(\Ac)$
 is identified with $\Sp(4, \CC)$.

   \paragraph{The genus 2 Siegel plane as the  future tube.} 
   The genus $2$ {\em Siegel upperhalf plane} $\Hen_2$ is defined as the space of symmetric
   $2\times 2$ complex matrices with positive definite imaginary part: 
   \[
   \Hen_2 = \bigl\{ T = \|t_{ij}\|\in \Sym_{2\times 2}(\CC) \bigl| \, \Im(T) > 0 \bigr\} . 
   \]
   The  quadratic form $\det(a)$ on the space $\Sym_{2\times 2}(\RR)$  of real symmetric $2\times 2$ matrices
  has  signature $(1,2)$ and so can be seen as
   the $3$-dimensional Minkowski space $\RR^{1,2}$. In other words, we represent
   $\RR^{2,1}$ as  
    $S^2(\RR^2)$, where  $\RR^2$ is realized as the space of spinors
   for $\RR^{1,2}$.   
   Further, the cone $\Sym_{2\times 2}^+(\RR)$ of positive definite
   symmetric matrices is identified with the future cone $C^+\subset \RR^{1,2}$. 
   This means that $\Hen_2$ is identified with the $3$-dimensional future tube $\Ten_3$. 
   
   \vskip .2cm
   
   The  action \eqref{eq:act-Siegel} of the real group $\Sp(4, \RR)$ on $\Sym_{2\times 2}(\CC)$
   restricts to an action on $\Hen_2$ by biholomorphic conformal transformations thus identifying
   $\Spin(4,\RR)$ with $\Spin(2,3)$, the spin cover of the group of biholomorphic
   automorphisms of $\Ten_3$. This is the $d=3$ instance of the real version of Theorem \ref{thm:ahlfors}
   as discussed in \S \ref{subsec:conf-lap}\ref{par:future}. Indeed, in this case the algebra $\Ac_\RR$
   is isomorphic to $\Mat_2(\RR)$ while $\Ac_\RR^{\leq 1}$ is identified with $\Sym_{2\times 2}(\RR) = \RR^{1,2}$ 
 and $\SL'_2(\Ac_\RR) \= \Spin(4, \RR)$. This gives an identification
 \[
 \Hen_2 = \Sp(4, \RR)/U(2)
 \]
which is equivalent (by passing to  finite isogenies of the groups and subgroups in question)
to  the $d=3$ instance of 
 \eqref {eq:Td-quot}.

 \paragraph{ Siegel forms and conformal densities.} \label{par:sieg-forms-2}
  Let $\ZYD_2$ be the set of pairs 
 $\alpha = (\alpha_1\geq
 \alpha_2)$, $\alpha_i\in\ZZ$. We refer to elements  $\alpha\in\ZYD_2$
 as 
 {\em $\ZZ$-Young diagrams  with $2$ rows}.  They are dominant (integer) weights for
 $\GL(2,\CC)$ and so parametrize irreducible  algebraic representations $\rho_\alpha$
 of $\GL(2, \CC)$ or, what is equivalent. Schur functors $\Sigma^\alpha:  M \mapsto \Sigma^\alpha M$
  on the
 category of $2$-dimensional vector spaces and their isomorphisms. 
 Note the natural identifications
  \be
 \Sigma^\alpha M \, \= \, S^{\alpha_1-\alpha_2}(M) \otimes \bigl(\bigwedge\nolimits^2 M\bigr)^{\otimes \alpha_2}, \quad  \Sigma^\alpha(M^*) \= (\Sigma^\alpha M)^* \= \Sigma^{-\alpha_2, -\alpha_1}M. 
 \ee
 We also denote 
 $\Sigma^\alpha = \Sigma^\alpha \CC^2$ so the irreducible representation associated to $\alpha$ 
 has the form
 $\rho_\alpha: \GL(2, \CC) \to \Aut (\Sigma^\alpha)$. 
 
 \vskip .2cm
 
 Denote  the 
  vector bundle $\Sigma^\alpha(\Lambda)$ on   $\LG(E)$ by   $\Lc_\alpha$. 
  For $\alpha = (r,r)$, $r\in \ZZ$ the functor $\Sigma^\alpha$ takes $M\mapsto \bigl(
 \bigwedge^2M\bigr)^{\otimes r}$  and we write $\Lc_r$ for $\Lc_{(r,r)}$.

   \vskip .2cm
   
   After restriction to $\Hen_2\subset\LG^\circ(E)$, the bundle $\Lc_\alpha$ is 
    trivialized
 with fiber $\Sigma^\alpha$ and is equivariant with respect to 
 $\Sp(4, \RR)$. We denote $\F_\alpha = H^0(\Hen_2, \Lc_\alpha)$
 the space of holomorphic sections which we call 
   {\em Siegel forms (not necessarily automorphic)
  of weight}  $\alpha$. The group $\Sp(4, \RR)$ acts on $\F_\alpha$.
  Explicitly, elements of $\F_\alpha$ are represented by functions
  $f: \Hen_2\to \Sigma^\alpha$ and the action of $g=\begin{pmatrix}
  A&B\\C&D\end{pmatrix}\in\Sp(4, \RR)$ on them is given by
  \be\label{eq;SF-alpha}
  f\mapsto g^*(f), \quad g^* (f)(T) = \rho^{-1}_\alpha(CT+D)( f(\gamma(T))). 
  \ee
 Here passing to  the inverse $\rho_\alpha^{-1}$ is needed as we convert the right action on functions
 (induced by the left action on the space)  into a left action again. 
 
  \begin{exas}
  (a) 
  For  $\alpha = (r,r)$ we have $\rho_\alpha(x) = \det(x)^r$, $x\in \GL(2,\CC)$. 
  Siegel forms of weight $(r,r)$ are traditionally called simply {\em
  (Siegel)  forms
  of weight $r$} and we write   $\F_r = H^0(\Hen_2, \Lc_r)$ for the space of such forms. 
  
    \vskip .2cm
  
  (b) As the bundle $\LL$ of $\LG(E)$ is identified with $\Lc_{-2} = \Oc_{\LG(E)}(2)$,
  conformal densities on $\Hen_2$ are the same as Siegel forms of weight $-2$. 
  
  \vskip .2cm
  
  (c) Since $\Omega^1_{\LG(E)} \= S^2(\Lambda)$,  we have the identification of
  the canonical bundle
  \[
  \omega_{\LG(E)}  =
   \Omega^{3}_{\LG(E)}  
 \,  \= \, 
   \bigl(\bigwedge\nolimits^2 (\Lambda)\bigr)^{\otimes (3)} \, = \, \Lc_3.
  \]
  Therefore volume forms
  on $\Hen_2$ are the same as Siegel forms of weight $3$. 
   \end{exas}
   
   \paragraph{The gerbe $\Gc_{1/2}$, the metaplectic group and forms of half-integer weight. }\label{par:g_1/2-gen2}
   Similarly to  \S \ref{subsec:P11}\ref{par:g-1/2-p11}, we denote by $\Gc_{1/2}$ the $\CC^*$-gerbe
   on $\LG(E)$ formed by square roots  of the line bundle $\Lc_1$. Thus we have a canonical
   such square root $\Lc_{1/2}$ understood as a $\Gc_{1/2}$-twisted line bundle. 
   The gerbe $\Gc_{1/2}$ is $\Sp(4, \CC)$-equivarint and so is $\Lc_{1/2}$.

   \vskip .2cm
   
  Let  $\Z2YD_2$ be the set of pairs $\alpha = (\alpha_1 \geq \alpha_2)$ with  $\alpha_i\in{1\over 2}\ZZ$,
   $\alpha_1-\alpha_2\in\ZZ$, i.e., of integral and half-integral weights for $\GL(2,\CC)$. 
   Let
   \[
   \Z2YD{}^-_2 \,\, = \,\, \Z2YD{}_2 \setminus \ZYD_2 \,\, = \,\, \bigl({1\over 2}, {1\over 2}\bigr) + \ZYD_2
   \]
   be the set of strictly half-integral weights. 
    For  
 $\alpha = ({1\over 2}, {1\over 2}) + \beta \,\in\, \ZYD^-_2$  we define the 
 $\Sp(4,\CC)$-equivariant
 $\Gc_{1/2}$-twisted vector  bundle
 $\Lc_\alpha = \Lc_{1/2}\otimes \Lc_\beta$.      
 
    \vskip .2cm
   
   After restriction to $\Hen_2$, the gerbe $\Gc_{1/2}$ is trivial and $\Sp(4,\RR)$-equivariant
   and so gives rise to the central extension
   \[
   1\to \ZZ/2 \lra \Mp(4) \lra \Sp(4,\RR) \to 1 
   \]
   known as the {\em metaplectic group}. 
     Similarly to \S  \ref{subsec:P11}\ref{par;super-MF}, $\Mp(4)$ is defined as the
   set of pairs $\wt\gamma = (\gamma, \phi)$ 
   where $\gamma\in\Sp(4, \RR)$ is in the block form  as before and
   $\phi = \phi(T)$ is a function in $\Hen_2$ satisfying $\phi(T)^2  = \det(CT+D)$. 
   For $\alpha \in \Z2YD_2^-$,
   the $\Gc_{1/2}$-twisted line bundle $\Lc_\alpha$ can be regarded as an actual
   $\Mp(4)$-equivariant line bundle on $\Hen_2$. 
   
   \vskip .2cm
   
   Extending previous notation, for any $\alpha \in \Z2YD_2$ we denote $\F_\alpha = H^0(\Hen_2, \Lc_\alpha)$ and call
   elements of this space Siegel forms of weight $\alpha$ (or of weight $r\in {1\over 2} \ZZ$, if
   $\alpha = (r,r)$). It is a representation of $\Mp(4)$. If $\alpha$ is not integer, so
   $\alpha = ({1\over 2}, {1\over 2}) +\beta$ as above, then elements of $\F_\alpha$ can be
   seen as functions $f: \Hen_2\to\Sigma^\beta$ with the action of
    $\wt g = (g,\phi)
   \in\Mp(4)$
   as above on them given by 
  \[
  f\mapsto \wt g^*(f), \quad \wt g^* (f)(T) = \phi(CT+D)^{-1} \cdot \rho^{-1}_\beta(CT+D)( f(\gamma(T))). 
  \]
 
  To give a unified treatment for any $\alpha\in\Z2YD_2$, let 
    $\wt \GL(2, \CC)\to \GL(2,\CC)$ be the $2$-fold covering corresponding
 to taking $\det(x)^{1/2}$, i.e., the group fitting into the Cartesian square
 \be\label{eq:wt-GL2}
 \xymatrix{
 \wt \GL(2, \CC) \ar[d] \ar[r] & \CC^*\ar[d]^{z^2}
 \\
 \GL(2, \CC) \ar[r]^{\det} & \CC^*.
 }
 \ee
 An element $\wt g = (g,\phi)\in\Mp(2)$ lifting $\gamma$ as above,
 defines a lift $\wt{CT+D}$ of the $\GL(2,\CC)$-valued function $T\mapsto CT+D$ on $\Hen_2$
 to a $\wt\GL(2,\CC)$-valued function. 
 Now, irreducible representations $\rho_\alpha$ of $\wt\GL(2,\CC)$
 are labelled by $\alpha \in \Z2YD_2$, so we can form the space
 $\Sigma^\alpha$ for such $\alpha$. Then $\F_\alpha$ consists of holomorphic
 $f: \Hen_2\to\Sigma^\alpha$ with $\wt g$ as above taking
 \be\label{eq;SF-alpha-half}
 f\mapsto \wt g^*(f), \quad \wt g^* (f)(T) =  \rho_\alpha^{-1}(\wt{CT+D})( f(\gamma(T))). 
 \ee

 \paragraph{The Siegel Laplacian.} 
   
   Specializing Example \ref{ex:conf-lap-Q} (b) to $d=3$, we obtain the (genus 2) {\em Siegel Laplacian}
   which is an $\Mp(4)$-invariant differential operator
   \be\label{eq:siefel-lap2}
   \Delta\,  = \, 
   \del_{11} \del_{22} - \del_{12} \del_{21}: \,\,
    \Lc_{1/2} \lra \Lc_{5/2}. 
   \ee

   \begin{rem}\label{rem:sing-forms-gen2}
   The operator $\Delta$ is known in the classical theory of Siegel modular forms.
 Sections  $\Lc_{1/2}$  annihilated by $\Delta$  (with further automorphy properties)
    are known as 
 {\em singular forms}    \cite{freitag-zweites, freitag-LNM}. The weight $ r= 1/2$ 
 is known in this case  as the {\em singular weight},  since only    for $r=1/2$ the condition of 
  vanishing of $\Delta$ on $\Lc_r$, if 
 defined explicitly by 
   \eqref{eq:siefel-lap2},  is   $\Mp(4)$-invariant. In fact, it was shown by
Freitag  \cite{freitag-zweites}  (with  a more general result  proved 
by Resnikoff \cite{resnikoff-sing})
that any section of $\Lc_{1/2}$ automorphic with respect
to a congruence subgroup of $\Mp(4,\ZZ)$, is  a singular form. 
The general context for  such invariant operators  is that of morphisms
between Verma modules, to be discussed in Appendix 
\ref{app:A}  below.

   \end{rem}

   \paragraph{The Siegel Dirac operator.} 
   For the $3$d complex Minkowski space $M^3_\CC = \Sym_{2\times 2}(\CC)$
   with coordinates $T=\|t_{ij}\|$ and metric $\det(T)$, spinors are represented
   by column vectors $\psi = (\psi_1, \psi_2)^t$. In these coordinates the
    Dirac operator $\dirac$
   on spinor fields $\psi(T)$ has the well known form, cf. 
   \cite[Eq.(2.30)]{deligne-freed}: 
   \be\label{eq:siegel-dirac}
   \begin{gathered}
   (\dirac \psi)_i = \eps_{jk} \del_{ij} \psi_k,  , \quad \text{i.e.,} 
   \\
   (\dirac\psi)_1 = \del_{11} \psi_2 - \del_{12} \psi_1, \,\,\,  (\dirac\psi)_2 = 
   \del_{21} \psi_2 -  \del_{22} \psi_1.  
   \end{gathered} 
   \ee
   
   \begin{prop}
   Considering $\psi(T)$, $T\in \Hen_2\subset \Sym_{2\times 2}(\CC)$
   as sections of the (trivial rank $2$) vector bundle 
   $\Lc_{({3\over 2}, {1\over 2})}$ on $\Hen_2$, 
   the operator \eqref{eq:siegel-dirac} defines an $\Mp(4)$-invariant
   differential operator
   \[
   \dirac: \Lc_{({3\over 2}, {1\over 2})} \lra \Lc_{({5\over 2}, {3\over 2})}
   \]
    between  equivariant vector bundles on $\Hen_2$. 
   \end{prop}
   We call this $\dirac$ the {\em Siegel Dirac operator}. 
   
   \vskip .2cm
   
   \noindent{\sl Proof:} 
   In dimension $3$ we have $\Spin(3,\CC) = \SL(2, \CC)$ and 
   $\CSpin(3,\CC) = \GL(2, \CC)$. This means that for the conformal quadric $Q^3=\LG(\CC^4)$ 
   the conformal  spin bundle is found as $V= \Lambda^*=\Sigma^{(0, -1)}\Lambda$
   whose restriction to $\Hen_2$ we denoted $\Lc_{(0, -1)}$. 
   So our statement by specializing Proposition \ref {prop:dirac-gen}
   to $d=3$, once we notice that $\LL=\Lc_{(2,2)}$. 
   \qed

 \vfill\eject

 \subsection {The genus $2$   Lagrangian super-Grassmannian as  the $3$d  $\Nc=1$ superconformal space}\label{subsec:SLG-gen2}
 
 Although super-generalizations of the Minkowski space and of the Poincar\'e group
 can be defined quite generally using spinors, natural conformal compactifications
 of such spaces exist only in a limited range of dimensions, essentially
 governed by exceptional isomorphisms of Lie algebras. The $3$-dimensional
 case is one of them.

 \paragraph{The super-thickening of $\LG(\CC^4)$.}  
 Let $E=(\CC^4, \omega)$
 be as before. Similarly to \S\ref{subsec:P11}\ref{par:P11-susy}
  put $Y=\CC^{4|1} =  \CC^{0|1}\oplus E$, 
  and  equip it 
    with the super-symplectic form $\eta$    by putting
   \be\label{eq:eta-gen2}
   \eta\bigl( \zeta, v), ( \zeta', v')\bigr) \,=\, \omega(v,v') + \zeta\cdot\zeta', 
   \quad x, x'\in E, \,\,\, \zeta, \zeta'\in \CC. 
   \ee

 \begin{defi}
 The {\em Lagrangian super-Grassmannian}
 of $E$ is  the supermanifold  $\SLG(E)$ parametrizing isotropic subspaces
 $L\subset Y$ of purely even dimension $2|0$. 
 \end{defi}
 
\noindent Thus,  $\SLG(E)_\red = \LG(E)$. 

\vskip .2cm

  By construction $\SLG(E)$ carries the tautological
 rank $2|0$ bundle 
   \[
   S\Lambda   \subset \wt Y := Y\otimes \Oc_{\SLG(E)}.
   \]
    It extends  the tautological bundle on
 $\LG(E)$  denoted $\Lambda$. As $\SLG(E)$ consists
 of isotropic subspaces, the form $\eta$ gives
  a projection $h: \wt Y/S\Lambda \to (S\Lambda)^*$. 
  The tangent bundle to $\SLG(E)$ is described as follows.

  \begin{prop}
  We have a Cartesian square of vector bundles 
  \[
  \xymatrix{
  T\SLG(E) \ar[r] \ar[d] & \Hom^+(S\Lambda, (S\Lambda)^*)\ar[d]
  \\
  \Hom(S\Lambda, \wt Y/S\Lambda) \ar[r]_{ \ul h} &\Hom(S\Lambda, (S\Lambda)^*),
  }
 \]
 where $\Hom^+(S\Lambda, (S\Lambda)^*) = S^2((S\Lambda)^*)$ is the bundle of  self-adjoint maps
 and $ \ul h $ is the composition with $h$. 
 \end{prop}
 
  \noindent{\sl Proof:} 
 By definition $\SLG(E)$
 is a closed sub-supermanifold in the super-Grassmannian
 $G=G(2|0, Y)$ of all $2|0$-dimensional subspaces which carries its own tautological
 bundle $S\Lambda_G \subset \wt Y_G := Y\otimes \Oc_G$. It is standard 
 \cite[Th.4.3.11]{manin} that
 $TG \= \Hom(S\Lambda_G, \wt Y_G/S\Lambda_G)$, so $TG|_{\SLG(E)}\=
 \Hom(S\Lambda, \wt Y/S\Lambda)$ and $T\SLG(E)$ is cut out in it  by
 the condition of ``preserving  isotropy infinitesimally".  
 This condition simply means (by differentiating the isotropy condition) 
  that the composition of
 a morphism $S\Lambda\to \wt Y/S\Lambda$ with $h$ is self-adjoint. \qed
 
 \vskip .2cm
 
 \noindent The proposition imples that $\SLG(E)$ is a supermanifold of dimension $3|2$. 
 
  \vskip .2cm

 Let $(S\Lambda)^\perp \subset \wt Y$ be the $\eta$-orthogonal to $S\Lambda$,
 a vector bundle
 of rank $2|1$. Note that $S\Lambda \subset(S \Lambda)^\perp$,
 so $(S\Lambda)^\perp/S\Lambda$ is a bundle of rank $0|1$. 
 It follows that
 \[
 \Tc := \Hom (S\Lambda, (S\Lambda)^\perp/S\Lambda) \subset T\SLG(E), \quad \rk(\Tc) = 0|2, 
 \]
 as it is an odd rank subbundle in $\Hom(S\Lambda, \wt Y/S\Lambda)$ and so  projects to $0$
 in the even rank bundle $\Hom(S\Lambda, (S\Lambda^)*)$.
 Note that the restriction of $(S\Lambda)^\perp/S\Lambda$ to $\LG(E)$ is trivial
 with fiber $\CC^{0|1}$, as a $2|0$-dimensional subspace in $Y$ depending
 on even parameters only, must lie in $\CC^4$, the even part. So
 $\Pi \Tc|_{\LG(E)} = \Lambda^*$.
 
 \begin{prop}\label{prop:T-gen2-susy}
 $\Tc\subset T\SLG(E)$ is a SUSY structure (Definition \ref {def:SUSYstr})
 and the corresponding Frobenius pairing  restricts on $\LG(E)$ to
  the canonical isomorphism
 \[
 S^2(\Pi \Tc|_{\LG(E)})  = S^2 \Lambda ^* \lra T\LG(E). 
 \]
 \end{prop}
 The proof will be given in the next paragraph.

   \paragraph{The open cell as  the super-Minkowski space
 $M^{3|2}_\CC$.} \label{par:cell-superMink}
   Let  $\SLG^\circ(E)\subset\SLG(E)$ be  the induced super-thickening of the 
  open part $\LG^\circ(E)$. Recall that $\LG^\circ(E)$ consists of Lagrangian subspaces in
  $E=L_p \oplus L_q$ transversal to $L_p$.
 As in \S  \ref {subsec:Siegel2}
 \ref{par:conf-dens-quad},
 denote by $x_i$ resp. $y_i$, $i=1,2$,  the linear coordinates in $L_q$ resp. $L_p$ dual to the basis $\{q_i\}$
  resp. $\{p_i\}$. A $2|0$-dimensional subspace (not necessarily isotropic) in $Y$ transversal to
  $L_p$ has the form
     \be\label{L-T-xi-gen2}
L_{T, \xi} = \bigl\{ (\zeta, x,y)^t\bigl| \, 
  y_i = \sum_{j=1}^2 t_{ij} x_j, \quad \zeta =   \xi x := \sum_{i=1}^2 \xi_i x_i  
  \bigr\},  
   \ee
 where $T=\|t_{ij}\|$ is any (not necessarily symmetric) $2\times 2$ matrix and
 $\xi = (\xi_1, \xi_2)$ are odd parameters assembled into a row vector. 
 Denote by $\zeta\in\CC^{0|1}$ the standard basis vector with $\theta = 1$. 
 The condition for $L_{T,\xi}$ to be isotropic is
 \be\label{eq:sigma-gen2}
 \sigma(T, \xi) \,  := \, t_{12} - t_{21} + \xi_1 \xi_2 \,=\, 0
 \ee
 as $L_{T,\xi}$ is spanned by
 \[
 \wt q_1 = \xi_1\zeta + t_{11} p_1 + t_{12} p_2 + q_1, \quad
 \wt q_2  =  \xi_2\zeta + t_{21} p_1+t_{22} p_2 + q_2
 \]
 and $\sigma(T,\xi) = \eta(\wt q_1, \wt q_2)$. 
 
 \vskip .2cm
 
 So $\SLG^\circ(E)$ is embedded into the affine superspace $\CC^{4|2}$ with four even coordinates
 $T=\| t_{ij}\|$ and two odd coordinates $\xi = (\xi_1, \xi_2)$ and defined inside it by the
 nonlinear equation $\sigma(T,\xi)=0$.  For many purposes it will be convenient for us
 to work with $\SLG^\circ(E)$ in this implicit form, using:
 \begin{itemize}
 \item Functions $f(T,\xi)$ on $\CC^{4|2}$ understood in their restriction to $\SLG^\circ(E)$,
 i.e., modulo $(\sigma)\subset \Oc_{\CC^{4|2}}$, the ideal generated by $\sigma$. 
 
 \item Derivations (vector fields) $D$ on $\CC^{4|2}$ which are tangent to $\SLG^\circ (E)$, i.e.,
 such that $D(\sigma) \subset (\sigma)$. 
 \end{itemize}
 Among such tangent vector fields we find the {\em spacetime}  and {\em spinor derivatives}
 \be\label{eq:del-D-gen2}
 \del_{ij} = {1\over 2} \biggl({\del\over \del t_{ij}} + {\del\over\del t_{ji}}\biggr), \quad
 D_i = {\del\over\del \xi_i} - \sum_{j=1}^2 \xi_j {\del\over \del t_{ij}}.
 \ee
 Indeed, one sees immediately that $\del_{ij}(\sigma) = D_i(\sigma) = 0$. These derivatives
 satisfy the relations
 \be\label{eq:D_i-gen2}
 \begin{gathered}
 [D_i, D_j]_+ = -2\del_{ij}, \quad [D_i, \del_{jk}]_- = 0,
 \\
 D_i(\xi_j) = \delta_{ij}, \quad D_i(t_{jk}) = -\delta_{ij} \xi_k \end{gathered}
 \ee
 in the first line of which  we recognize the familiar relations of the $3$d $\Nc=1$ SUSY algebra (or the
 super-Galileo algebra), see \cite[\S2]{gates}  \cite[\S2.3]{deligne-freed}. 
 
 \vskip .2cm
 
 For component analysis it is desirable to have a system of independent coordinates in
 $\SLG^\circ(E)$, identifying it with $\CC^{3|2}$. For these, we take
 \be\label{eq:xij-gen2}
 v_{ij} \,=\, {1\over 2} (t_{ij} + t_{ji}) \,=\, v_{ji}
 \ee
 and the same $\xi_i$ as before, So $v= \|v_{ij}\|$ is a symmetric matrix. Note that 
 \be\label{eq:del-ij-gen2} 
  \del_{ij}  =  {1\over 2} \biggl({\del\over \del v_{ij}} + {\del\over\del v_{ji}}\biggr)
 \ee
 in the new coordinates $(v,\xi)$ as well. Further, in these coordinates the derivative $D_i$ has
 the form
 \be\label{eq:D_i-in-xij-gen2}
 D_i = {\del\over \del \xi_i} - \sum_j \xi_j \del_{ij}. 
 \ee
 Note also that such a  formula does not hold in the $(T,\xi)$ coordinates, even after restricting to
 $\SLG^\circ(E)$. This is because the meanings of $\del\over\del \xi_i$ in the two coordinate
 systems are different:
 \[
 t_{12} = v_{12} - {1\over 2} \xi_1 \xi_2, \quad t_{21} = v_{12} + {1\over 2} \xi_1 \xi_2,
 \]
 so  for example,   $\del / \del \xi_1$ in coordinates $(v,\xi)$ does not annihilate $t_{12}$ or $t_{21}$. 
 
 \vskip .2cm
 
 We see  that $\SLG^\circ(E)$ is precisely the  ``toy superspace'' model of \cite[Ch.2]{gates}, denoted $M^{3|2}$
 in \cite[Ch.2]{deligne-freed}. As we are dealing with complex points,
 we use the notation $M^{3|2}_\CC$.
 
 \vskip .2cm
 
Now,  Proposition \ref {prop:T-gen2-susy} is a consequence of the next
 \begin{prop}
 The restriction of the subbundle $\Tc$ to $\SLG^\circ(E)$
 is spanned by the $D_i$. 
 \end{prop}

\noindent{\sl Proof:} This is obtained by  direct calculation using the
parametrization \eqref{L-T-xi-gen2}, similarly to the proof of Proposition
\ref {prop:D-xi-gen1}.\qed

% \vfill\eject

  \paragraph{The supergroup $\OSp(1|4)$  as the $3$d $\Nc=1$ superconformal group.}
  \label{par:OSP14}
  By definition, $\SLG(E)$ is acted upon by the Lie supergroup
  $\OSp(1|4, \CC) = \ul\Aut(Y, \eta)$, $Y= \CC^{0|1}\oplus E = \CC^{4|1}$. 
  Using the basis $\{p_i\}$ of $L_p$ and  $\{q_i\}$ of $L_q$ with  $E=L_p\oplus L_q$, we write points of $\OSp(1|4,\CC)$ as 
   \be\label{eq:OSp4-matrix}
   \begin{gathered}
 g=\begin{pmatrix}
 s&\alpha&\beta 
 \\
 \gamma& A & B
 \\
 \delta & C & D
 \end{pmatrix},  \quad \text{subject to} \quad  g^T H g = H, \quad 
 H = \begin{pmatrix} 1&0&0
 \\
 0&0&-1_2
 \\ 0&1_2&0
 \end{pmatrix}, 
 \\
  s\in\Mat_{1\times 1}, \,\,  A,B,C,D \in\Mat_{2\times 2},\,\, 
 \alpha,\beta\in \Mat_{1\times 2}, \,\, \gamma,\delta\in\Mat_{2\times 1}, 
 \end{gathered} 
 \ee
  with Latin variables being even and Greek ones being odd.      This identifies
 \be\label{eq:SLG-4-quot} 
 \SLG(E) \,=\, \OSp(1|4, \CC)/\PP, \quad 
 \PP := \Stab(L_p)
 =  \left\{\begin{pmatrix}
 s&0&\beta 
 \\
 \gamma & A & B
 \\
 0 & 0 & D
 \end{pmatrix} \right\}. 
\ee
 Using the dependent \eqref {eq:sigma-gen2}
  coordinates $T,\xi$ on the open chart $ \SLG^\circ(E) \subset \CC^{4|2}$,
we write the 
 action of $g\in\OSp(1|4,\CC)$  as above in the form
\be\label{eq:SSiegel-bir-2}
g(T,\xi) = \bigl( (\gamma\xi + AT+B)(\delta\xi +CT+D)^{-1}, \,
(s\xi + \alpha T+\beta)(\delta\xi+CT+D)^{-1}\bigr). 
\ee
Here $\gamma\xi = \|\gamma_i \xi_j\|$ is the product of a $2\times 1$ matrix $\gamma$ and a
 $1\times 2$ matrix $\xi$, and similarly for $\delta\xi$. 
 
 \vskip .2cm
  
  \begin{rems}
  (a)  As the SUSY structure $\Tc$ on $\SLG(E)$ is
  defined naturally in terms  of the form $\eta$, the action of  $\OSp(1|4,\CC)$
  preserves $\Tc$. This realizes $\OSp(1|4,\CC)$ as the (complex)  $3$d $\Nc=1$
  superconformal group, see \cite{park} \cite[\S5.1]{kuzenko}. 
  
  \vskip .2cm
  
  (b) 
  Since $\SLG^\circ(E)$ consists of $L$ transversal to $L_p$,
   $\PP$ 
  acts on $\SLG^\circ(E)=\CC^{3|2}$ by automorphisms of SUSY manifolds,
  as is also obvious from \eqref{eq:SSiegel-bir-2}.
  The commutant $[\PP, \PP]$ 
 is  the super-Poincar\'e group acting on the
  complex Minkowski superspace $M^{3|2}_\CC$ while $\PP$ is obtained
  from it 
  by adding scalar dilations. 
  
  \vskip .2cm
  
  (c) The radical 
$
\NN =  \left\{\begin{pmatrix}
 1&0  &\beta
 \\
 \gamma & 1_2 & B
 \\
 0 & 0 & 1_2
 \end{pmatrix} \right\}$   
 of $\PP$
  is not (super) commutative. This means that
$\SLG(E)$ is not a  super symmetric space in the sense of Serganova
\cite{serganova}. 

  \end{rems}
  
   \paragraph{The Super Siegel plane and super-Siegel forms.}\label{par:ssieg-forms-2}
   We define $S\Hen_2$, the {\em  genus 2 super Siegel
   half plane} to be the super thickening of $\Hen_2$ induced from the thickening
   $\SLG^\circ (E)$ of  $\LG^\circ (E)$ via the open embedding 
   $\Hen_2\subset \LG^\circ(E) = \CC^{3|2}$. 
   Thus the coordinates on $S\Hen_2$ are $(T,\xi)$ as above with the condition
   that  the matrix $v  = \|v_{ij}\|$ from \eqref {eq:xij-gen2} has positive definite
   imaginary part.

   The  real supergroup $\OSp(1|4,  \RR)$ acts on $S\Hen_2$
   by \eqref {eq:SSiegel-bir-2}, identifying 
    \[
 S\Hen_2 \, \= \, \OSp(1|4, \RR)/U(2).
 \] 
 In parallel to \S \ref{subsec:Siegel2}\ref{par:sieg-forms-2}, for each $\alpha\in \ZYD_2$
 we denote the vector bundle $\Sigma^\alpha (S\Lambda)$ on $\SLG(E)$ by $S\Lc_\alpha$
 and abbreviate $S\Lc_{(r,r)}$, $r\in \ZZ$,  to $S\Lc_r$. The space of
 {\em Siegel superforms} of weight $\alpha$ is denoted by
 \[
 \SF_\alpha = H^0(S\Hen_2, S\Lc_\alpha). 
 \]
 It is acted upon by $\OSp(1|4, \RR)$.   
 Note that   $S\Lambda|_{\SLG^\circ(E)}$ 
 and hence $S\Lambda|_{S\Hen_2}$ is trivialized and identified with
 $\Oc_{S\Hen_2}  \otimes \CC^2$. Therefore each $\alpha \in \ZYD_2$
 the  restriction of the bundle $S\Lc_\alpha$ to  $S\Hen_2$ is also trivialized:
 $ S\Lc_\alpha|_{S\Hen_2} \, \= \, \Oc_{S\Hen_2} \otimes V_\alpha$.
 So  elements of $\SF_\alpha$  can be viewed as $\Sigma^\alpha$-valued 
 (super) functions $f = f(T, \xi): S\Hen_2\to \Sigma^\alpha$ and the action of 
 $ g\in \OSp(1|4,\RR)$ in the form \eqref{eq:OSp4-matrix}
 on such functions has the form
 \[
(g^*f )(T,\xi) =  \rho_\alpha(\delta\xi+CT+D)^{-1}   f(g(T, \xi)). 
 \]
 Note that $\F_\alpha$, the space of ordinary Siegel forms of weight
 $\alpha$, is contained in the even part $(\SF_\alpha)_\0$ but does not concide with it because
 the terms with $\xi_1\xi_2$ are even as well. 

 \vskip .2cm
 
 Similarly to \S \ref{subsec:P11}\ref{par:g-1/2-p11}, the gerbe $\Gc_{1/2}$ of determinations
 of square roots of $\Lc_1$ on $\LG(E)$ is identified with the gerbe of
 determinations of  a square root of $S\Lc_1$ on $\SLG(E)$. Thus we have
 a $\OSp(1|4, \CC)$-equivariant $\Gc_{1/2}$-twisted line bundle
 $S\Lc_{1/2}$ on $\SLG(E)$ and further,, for each half-integer 
  $\alpha \in ({1\over 2}, {1\over 2})+\ZYD_2$
 an $\OSp(1|4, \CC)$-equivariant $\Gc_{1/2}$-twisted vector bundle $S\Lc_\alpha$.

 \vskip .2cm 
 
 We define the  {\em super-metaplectic group} $ \OMp(1|4)$ 
 as the 
  $2$-fold
 covering of  the supergroup  $\OSp(1|4, \RR)$,
 whose points are pairs $\wt g = (g,\phi)$ where $g$ is a point of $\OSp(1|4,\RR)$
 and $\phi = \phi(T, \xi)$ is a determination of 
  $\det(\delta\xi+CT+D)^{1/2}$ in $S\Hen_2$. 
  
  \vskip .2cm

  After restriction to $S\Hen_2$, the gerbe $\Gc_{1/2}$ becomes trivial so the
 $S\Lc_\alpha$ become actual trivial line bundles, equivariant with respect to $\OMp(1|4)$. 
 We extend the name
  ``Siegel superforms of weight $\alpha$'' and the notation
   $\SF_\alpha = H^0(S\Hen_2, S\Lc_\alpha)$ to half-integer $\alpha$ as well. 
   Thus $\SF_\alpha$ is acted upon by $\OMp(1|4)$. 
   
   \vskip .2cm
   
   Using the group $\wt\GL(2,\CC)$ from \eqref{eq:wt-GL2}, we can formulate the above
   in a unified fashion for any $\alpha\in \Z2YD_2$ (both integer and half-integer).
   That is, a lift $\wt g = (g,\phi)$ of $g\in \OSp(1|4,\RR)$ to $\OMp(1|4)$ defines a lift
   of the $\GL(2,\CC)$-valued (super) function
   $(T,\xi) \mapsto \delta\xi+CT+D$ on $S\Hen_2$ to a $\wt\GL(2,\CC)$-valued function 
    $(\delta\xi+CT+D)\, \wt{}$. As any $\alpha\in \Z2YD_2$ defines a representation
    $\rho_\alpha: \wt\GL(2,\CC) \to \Aut(\Sigma^\alpha)$, we can say that $f\in \SF_\alpha$
    is a function $f(T,\xi): S\Hen_2\to \Sigma^\alpha$ with $\wt g = (g,\phi)\in  \OMp(1|2)$ acting by 
    \be\label{eq:osp-ac-forms-gen2}
 (\wt g^*f )(T,\xi) =  \rho_\alpha((\delta\xi+CT+D)\, \wt{}\, \,\, )^{-1}   f(g(T, \xi)). 
    \ee
For $\alpha = (r, r)$, $r\in {1\over 2}\ZZ$
 we can identify $\Sigma^\alpha$ with $\CC$ (with nontrivial
 action of $\wt\GL(2, \CC)$, if $r\neq 0$) and  can view  sections of 
 of $S\Lc_\alpha$ as ordinary superfunctions $f(T,\xi)$ on $S\Hen_2$.

 \paragraph{The Lie superalgebra $\osp(1|4)$.} \label{par:osp-1-4}
 Consider the complex Lie superalgebra
 \[
 \gen = \osp(1|4) = \Lie (\OSp(1|4, \CC))
 \]
 Viewing $Y=\CC^{4|1} = \CC^{0|1}\oplus E$ as $\CC^{0|1}\oplus L_p \oplus  L_q$ as in \S 
 \ref {par:OSP14}, we write elements $x\in \gen$ as block matrices characterized by 
 $x^T H = -Hx$ where $H$ is as in \eqref{eq:OSp4-matrix} or explicitly
 \[
 x=\bpm 0&\kappa^t & \lambda^t
 \\ -\lambda  & K&L
 \\
 \kappa &M  &N 
 \epm, \quad \kappa, \lambda \in\Mat_{2\times 1}, \,\, K,L,M,N \in\Mat_{2\times 2}, \,\, L=L^t, M=M^t, K=-N^t. 
 \]
%%%BELOW COMMENTED OUT
 \iffalse
 \[
 x=\bpm 0&\zeta^t & -\eps^t
 \\\eps & Q&R
 \\
 \zeta&S &-Q^t
 \epm, \quad \eps,\zeta\in\Mat_{2\times 1}\,\, Q\in\Mat_{2\times 2}, 
 \quad R,S\in\Sym_{2\times 2}. 
 \]
 \fi
 %%%%ABOVE COMMENTED OUT
 
 Thus $\lambda\in L_p, \kappa \in L_q$ parametrize $\gen_\1$ while $\gen_\0=\sp(4)$. More 
 invariantly, we have the following well known fact. 
 
 \vskip .2cm
 
 Let $\pen=\Lie (\PP)$ and $\nen = \Lie(\NN)$ be the Lie algebra of the parabolic
 subgroup $\PP=\Stab(L_p)$ and its unipotent radical $\NN$. 
 Recall that $E=(\CC^4,\omega)$ is our symplectic vector space. 
 
 \begin{prop}\label{osp4-s2}
 (a)  We have an identification of super vector spaces
  \[
 \bigl\{ \underset{\1} {E} ; \,\, \underset {\0} {S^2E} = \sp(E)\bigr\}
 \buildrel \sim\over  \lra \gen = \osp(1|4).
 \]
 Under this identification,  the superbracket $E\otimes E\to S^2E$
  becomes the symmetrized product.   
 The Lie bracket on $S^2E = \sp(E)$ is the Poisson bracket of quadratic 
 Hamiltonians.  The bracket between $\sp(E)$ and $E$ is the standard action
 of $\sp(E)$ on $E$. 
 
 \vskip .2cm
 
 (b)  Under the identification of (a), we have
 \[
\pen =   \bigl\{ \underset{\1} {L_p}; \,\, \underset{\0} {L_p\cdot E} \bigr\}, 
\quad
\nen =  \bigl\{ \underset{\1} {L_p}; \,\, \underset{\0} {S^2 L_p} \bigr\},
 \]
 with the only nonzero component of the bracket on $\nen$ 
 being the symmetrized product
 $L_p\otimes L_p\to S^2L_p$. 
 \qed
  \end{prop}
  
 As $\gen$ acts on $\SLG(E)$, for any $x\in \gen$ we have a vector field $D_x$ on
 $\SLG(E)$.
 
 \begin{prop}\label{prop:SUSY-LG-2}
 The SUSY structure $\Tc$ of $\SLG(E)$ is spanned by $D_y, y\in\gen_\1$. 
 \end{prop}
 
 Then proposition will follow from a precise formula for the action of
  $D_y$  on $\SLG^\circ(E)$
 in coordinates $(T,\xi)$ for $y$ belonging to the symplectic basis $p_i, q_i$
 of $E=\gen_\1$.  Let us 
  organize various variables
 and derivations into vectors and matrices:
 \[
 \xi=(\xi_1, \xi_2) ,\,\, {\del\over\del\xi}   = 
   \biggl( {\del\over\del\xi_1}, {\del\over\del\xi_2}\biggr), \,\, {\del\over \del T} = 
 \left\|{\del \over  \del t_{ij}}  \right\|
  \]
  and form the vectors of 
  spinor derivations 
  \be \label{eq:Di-D'i-gen2}
 \begin{gathered}
 D = {\del\over\del\xi} - \xi\cdot {\del \over \del T}, \quad \text{i.e.,} \quad
 D_i \,=\,{\del\over\del \xi_i}  -  \sum_j \,\xi_j {  \del\over \del t_{ij} } ,
 \\
 D' =  \biggl({\del\over\del\xi}  -  \xi\cdot {\del \over \del T}\biggr) \cdot T, 
  \quad \text{i.e.,} \quad
  D'_i =\, \sum_j  t_{ji} D_j.  \end{gathered} 
 \ee
 Thus $D_i$ are the spinor derivations from \eqref  {eq:del-D-gen2}. 
 Proposition \ref {prop:SUSY-LG-2} is a consequence of the next
 
 \begin{prop}\label{prop:Di=Dpi-2}
 We have  $D_i = D_{p_i}, D'_i = D_{q_i}$.
 \end{prop} 
 
 \noindent  {\sl Proof:}  Similar to that of Proposition 
 \ref{prop:DpDq} in the $1$-dimensional case, applying 
 \eqref {eq:SSiegel-bir-2}
  to $g=1+\eps p_i$ or $1+\eps q_i$ for an odd parameter $\eps$
 and taking the part linear in $\eps$. 
  \qed 
 
 \paragraph{Action of $\gen_\1$ on  super-Siegel forms.} 
 As on    $S\Hen_2$ each equivariant bundle
 $S\Lc_\alpha$ is trivialized as a bundle 
  the action of $D_y$,  $y\in\gen_\1$  on $S\Lc_\alpha$
  can be written as a matrix (super) differential
 operator  $D_y^{(\alpha)}$. Let us find these operators explicitly
 for $\alpha = (r,r)$, $r\in {1\over 2}\ZZ$, in which case we denote
 them $D_y^{(r)}$. 
 
 \begin{prop}\label{prop:D_x^r-gen2}
 For $r\in {1\over 2}\ZZ$ we have, with respect to the above trivialization, 
 \[
 D^{(r)}_{p_i} = D_i, \quad D^{(r)}_{q_i} = D'_i - r \xi_i. 
 \]
  
 \end{prop}
 
 \noindent{\sl Proof:} Similar to that of Proposition \ref{prop:Di=Dpi-2}, now
 applying  \eqref{eq:osp-ac-forms-gen2} for $\alpha  = (r,r)$
 when  $\rho_\alpha(x) = \det(x)^r$, $x\in\GL(2)$. \qed
 
 \vfill\eject

 \subsection {The super-Laplacian:  on-shell Heisenberg relations} \label{subsec:super-Lap-gen2}

 \paragraph{ The scalar supermultiplet: superpartners and auxiliary fields.} 
 The original understanding of supersymmetry was as ``symmetry between bosons and fermions"
 which arranges particles/fields into pairs of ``superpartners'':  one bosonic, one fermionic. 
 For example, the scalar field/particle, a boson has  as  superpartner  the spinor field/particle,
 a fermion.
 
 \vskip .2cm
 
 The point of view of  {\em superspace} and {\em superfields}, i.e., working with functions or
 other tensor fields on a supermanifold extending the usual spacetime,  is 
 (despite seeming completely self-evident now) a somewhat
 later development \cite{salam}. It does produce superpartners but at the price of introducing additional
 ``parasitic'' fields  (usually called {\em auxiliary fields}) which have no independent
  physical interpretation and need to be eliminated somehow. 
 
 \vskip .2cm
 
 The simplest instance is provided by the so-called {\em scalar supermultiplet}, i.e., the theory of
 a scalar field (function) $\Phi (x, \xi)$ on the super-spacetime. Here $x=(x_1, \cdots, x_d)$
 are even  spacetime coordinates and $\xi=(\xi_1,\cdots, \xi_N)$ are odd ``spinor coordinates''.
 Writing  the component expansion
 \[
 \Phi(x,\xi) = \phi(x) + \sum_{i=1}^N \phi_i(x) \xi_i + \sum_{i<j} \phi_{ij}(x) \xi_i \xi_j + \cdots, 
 \]
  we interpret the $\xi$-constant term $\phi(x)$ as a scalar field in the usual sense and the
  $\xi$-linear term $\psi(x,\xi) = \sum \phi_i(x) \xi_i$ as a spinor field, so $\phi$ and
  $\psi$ are superpartners.  But the terms quadratic and higher in $\xi$
  represent additional fields which are more than we asked for. Therefore,   one typically imposes
  {\em constraints}, i.e., the requirement that some of these additional fields vanish and then look for equations of motion that express the rest through
  other fields, so the physically meaningful  part of  
  $\Phi(x,\xi)$ 
  is no more than linear in $\xi$. Only such $\Phi$ implement the philosophy of superpartners
  in a clean way.
  
  \vskip .2cm
  
  In this section we will see how this phenomenon plays out for Siegel superforms of genus 2. 
  
  \paragraph{The super-Laplacian on $M^{3|2}_\CC$.} 
  \label{par:super-lap-2}
   The complex Minkowski superspace
 $M^{3|2}_\CC = \CC^{3|2}$ carries a (constant)  {\em super-Riemannian metric}, i.e., a super-symmetric
 bilinear form, consisting of:
 \begin{itemize}
 \item  The symmetric  bilinear form (Riemannian metric)  $g$ on $\CC^3 = \Sym_{2\times 2}(\CC)$
 corresponding to the quadratic form $T\mapsto \det(T)$.
 
 \item The symplectic form $\eps$ on the  spinor space $\CC^2$ whose symmetry group $\SL(2,\CC)$ is the
 spinor covering for the  $\SO(3,\CC)$ preserving the metric on $\CC^3$. In our notation, the coordinates
 $\xi_1, \xi_2$ form a symplectic basis of the dual space. We identify $\eps$ with the 
 Levi-Civita symbol $(\eps_{ij})$  with  $\eps_{12}=1$, $\eps_{21}=-1$, $\eps_{11}=\eps_{22}=0$. 
  \end{itemize}
  
  \noindent  We call  the {\em  super-Laplacian} and denote   $\DD$ the quadratic 
     form in the spinor derivatives
     $D_i$ corresponding to the odd part $\eps$ of the super-Riemannian metric, i.e.,
     \be
     \DD  = D_i D^i = \sum \eps_{ij} D_i D_j = D_1D_2 - D_2 D_1 = [D_1, D_2]_- 
     \ee
     acting in functions $\Phi =\Phi (T,\xi)$ on $M^{3|2}_\CC$. 
     Here the second expression $D_i D^i$ is  written using the ``super-Einstein rule'' 
     understood as the third expression.  Thus $\DD$ is the  even-style  commutator of the odd
     derivations $D_1$ and $D_2$. In
     the supersymmetry literature 
     it is typically denoted by $D^2$, see \cite [\S2.2a]{gates} \cite[Eq.(2.41)]{deligne-freed}
     \cite[p.2]{howe}; we prefer a single symbol notation suggestive of two copies of the
     letter $D$. 
     
     \vskip .2cm
     
     We write the ``superfield'' $\Phi$ in components  using the  independent coordinates
     $v=\|v_{ij}\|$ (a symmetric matrix) and $\xi = (\xi_1, \xi_2)$ from \eqref {eq:xij-gen2}
     \be\label{eq:Phi-comp-gen2}
     \Phi (v, \xi) = \phi(v) + \sum_{i=1}^2 \psi_i(v) \xi_i + F(v) \xi_1\xi_2,
     \ee
     so $\phi$ is a scalar on $M^3_\CC$ while $\psi = (\psi_1, \psi_2)^t$ is a spinor
     and $F$ is an auxiliary  field without  independent physical meaning.  Applying 
     \eqref {eq:D_i-in-xij-gen2}, we see that in this form the action of
     $\DD$ is given by
     \be\label{eq:DD-2}
     \DD(\phi, \psi, F) = (-2F, -2\dirac \psi,  2\Delta\phi), 
     \ee
 where $\dirac$ is the Dirac operator \eqref{eq:siegel-dirac} and $\Delta$ is the
 Laplacian \eqref{eq:siefel-lap2}. The equation $\DD \Phi=0$ is the Euler-Lagrange 
 equation for the simplest (free, massless) {\em Wess-Zumino action}
 \[
 S[\Phi] = \int dv d\xi \,(D_i\Phi) (D^i \Phi) \,=\, \int dv d\xi \, \bigl(
 (D_1 \Phi) (D_2 \Phi) - (D_2 \Phi ) (D_1 \Phi) \bigr)
 \]
 (the integration over $dv$ should of course be understood along the real locus). In
 components the equation means that:
 \begin{itemize}
 \item  $\phi$ satisfies the massless wave (Laplace) equation $\Delta \phi=0$;
 
 \item $\psi$ satisfies the massless Dirac equation $\dirac \psi=0$;
 
 \item The auxiliary  field $F$ vanishes identically.
 
 \end{itemize}
 \noindent Thus the constraints $F=0$ are imposed by  the equations of motion,
 i.e., hold ``on shell'' in the physical terminology.  More generally, for the Wess-Zumino model
 with potential, $F$ is expressed through the potential
 \cite[(4.22)]{deligne-freed}. 
 
 Eq. \eqref {eq:DD-2} shows 
 part (a) of the following well known   proposition, see, e.g., \cite[(2.43-44)]{deligne-freed}.   
 \begin{prop}\label{prop:super-lapl-anti-2}
 (a) We have
 \[
 \DD^2 = \Delta  = \del_{11} \del_{22} -\del_{12} \del_{21} 
 \]
 acting along the even variables only, i.e., simultaneously in $(\phi, \psi, F)$.
 
 \vskip .2cm
 
 (b) $\DD$ anticommutes with the spinor derivations: $\DD D_i = - D_i \DD$. \qed
 \end{prop}
 
 \noindent{\sl Proof (b):} As $D_i^2=\del_{ii}$ commutes with all $D_j$, we have
 \[
[D_1, \DD]_+  =  D_1(D_1D_2- D_2D_1)  +  (D_1D_2-D_2 D_1) D_1 = 
D_1^2 D_2 -D_1D_2 D_1 + D_1D_2 D_1 - D_2 D_1^2 = 0
 \]
 and similarly for $[D_2, \DD]_+$. \qed
 
 \paragraph{Super-Laplacian as an element of the anticenter.} \label{par:D-anti}
 The above considerations can be seen as calculations in the enveloping algebra
 $U(\nen^-)$ where $\nen^-= L_q\oplus S^2L_q$ is the opposite of 
 the nilpotent radical from 
 Proposition \ref{prop:super-lapl-anti-2},  acting on $\SLG^\circ(E)$.  
 Thus the basis of $\nen^-$ is formed by
 $p_1, p_2$ (odd) and $p_{ij}=p_{ji}, \, i,j=1,2$ (even)  so that the only superbracket
 is $[p_i, p_j]_+ = 2p_{ij}$. By Proposition \ref {prop:Di=Dpi-2} the operator $D_i$ is
 given by the action of $p_i\in\nen^-_\1$. Therefore 
 the super-Laplacian $\DD$ corresponds to the element
 $Q = [p_1, p_2]_- \in U(\nen^-)$ and Proposition \ref{prop:super-lapl-anti-2}(b)
 with its proof
 lifts to the statement that $Q$
  belongs to the anticenter of $U(\nen^-)$: an even element anticommuting with
  all odd elements and commuting with all  even ones, see
  \S \ref{subsec:DO-super}
 \ref{par:super-general}. For future use, let us denote
 \be\label{eq:LDD}
 ^L\DD = (-1)^F \DD. 
 \ee

 \paragraph{Conformal skew-invariance of the super-Laplacian.} 
 We now restrict from $M^{3|2}_\CC = \SLG^\circ(\CC^4)$ to the super-Siegel plane
 $S\Hen_2$.
 
 \begin{prop}\label{prop:conf-superlapl-gen2}
 Considering functions $f(T,\xi)$, $T\in\Hen_2$,
  as sections of the line bundle $S\Lc_{1/ 2}$, the operator $\DD$
  defines an $\OMp(1|4)$-skew-invariant 
  (see 
  \S \ref{subsec:DO-super}
 \ref{par:super-general}) differential operator
  \[
  \DD: S\Lc_{1/2} \lra S\Lc_{3/2}
  \]
  between equivariant line bundles on $S\Hen_2$.  In other words,  
  $^L\DD$ is a morphism of $\OMp(1|4)$-modules. 

 \end{prop}
 
\noindent{\sl Proof:} It is enough to prove this at the level of the Lie algebra $\gen=\osp(1|4)$.
As $\gen$ is generated by $\gen_\1$, it is enough to show that
\[
D_y^{(3/2)} \, \DD = - \DD \,D_y^{(1/2)}, \,\,\,  \forall\, y\in \gen_1.
\]
For $y=p_i$ this is Proposition 
\ref {prop:super-lapl-anti-2} since $D_{p_i}^{(r)} =. D_i$ independently of $r$.
For $y=q_i$ this means explicitly that
\[
\biggl(\sum_j t_{ji} D_j - {3\over 2} \xi_i\biggr)\,  \DD = - \DD \, 
\biggl( \sum_j t_{ji} D_j - {1\over  2} \xi_i\biggr).
\]
 This latter condition we verify directly using the values of the  $D_i$ given in the second line
 of \eqref {eq:D_i-gen2}. \qed

% \vfill\eject

 \paragraph{Even Heisenberg relations  for $\osp(1|4)_\1$ 
 on harmonic  superforms.}\label{par:even-heis-gen2}
 Let us consider each $D_y^{(r)}$, $y\in\gen$, $r\in {1\over 2}\ZZ$,
 as   operators in $\Oc_{S\Hen_2}$,
 i.e., as sections of $\Dc_{S\Hen_2}$ given by 
 Proposition \ref {prop:D_x^r-gen2},  using the trivialization of
  $S\Lc_r$.  Thus $D_{p_i}^{(r)} = D_i$ is the spinor derivative and 
 so $\DD = [D_1, D_2]_- = [D_{p_1}^{(r)}, D_{p_2}^{(r)}]_-$  independently of $r$. 
 
 \vskip .2cm

 Let $\Dc_{S\Hen_2} \DD \subset \Dc_{S\Hen_2}$ be the left ideal
 generated by $\DD$ and $\Mc_\DD=  \Dc_{S\Hen_2}/  \Dc_{S\Hen_2} \DD$ be the
 quotient  left $ \Dc_{S\Hen_2}$-module.  Thus $\Mc_\DD$
  generated by one section $u$ subject to
 the relation $\DD u=0$.  We call $\Mc_\DD$ the {\em harmonic $\Dc_{S\Hen_2}$-module}. 
 
 \vskip .2cm
 
  Let $\Heis_2 = \Heis(E, {1\over 2} \omega)$ be the Heisenberg algebra
 associated to the symplectic space $E=\CC^4$ and the form ${1\over 2} \omega$.
 So it is generated by $y\in E$ subject to the relations
 $
 yz -zy = {1\over 2} \omega(y,z) \cdot 1
 $.

 \begin{prop}\label{prop:even-Heis-gen2}
(a) The operators $D_y^{(1/2)}$, $y\in\gen_1$, satisfy even Heisenberg relations
 modulo $\DD$, i.e.,
 \[
 [D_y^{(1/2)}, D_z^{(1/2)}]_-   - {1\over 2} \omega(y,z)  \quad\in\quad \Dc_{S\Hen_2} \DD.
 \]
  (b) The right action of each $D_y^{(1/2)}$ on $\Dc_{S\Hen_2}$
  preserves $\Dc_{S\Hen_2} \DD$ and so induces a endomorphism of  $\Mc_\DD$
  as a left $\Dc_{S\Hen_2}$-module. This gives an algebra
  homomorphism $\Heis_2 \to \End_{\Dc_{S\Hen_2}}(\Mc_\DD)$. 
 
 \end{prop}
 
 \noindent{\sl Proof:} (a) The statement means explicitly that 
  \[
 [D_{p_1}^{(1/2)}, D_{p_2}^{(1/2)}]_-, \quad 
 [D_{q_1}^{(1/2)}, D_{q_2}^{(1/2)}]_-,  \quad  [D_{p_i}^{(1/2)}, D_{q_j}^{(1/2)}]_- 
 - {1\over 2} \delta_{ij} \quad \in \quad  \Dc_{S\Hen_2} \DD. 
 \]
 The first commutator is equal to $\DD$, as $D_{p_i}^{(r)} = D_i$ for any $r$. Further, 
 using  \eqref {eq:D_i-gen2}  and  \eqref {eq:sigma-gen2} we verify directly that
 \[
 \begin{gathered}
  [D_{q_1}^{(1/2)}, D_{q_2}^{(1/2)}]_- \, =\,  \det (T) \, \DD, 
  \\
   [D_{p_1}^{(1/2)}, D_{q_2}^{(1/2)}]_- \, =\, t_{22} \, \DD, \quad
    [D_{p_2}^{(1/2)}, D_{q_1}^{(1/2)}]_- \, =\,-t_{11} \DD, 
    \\
     [D_{p_1}^{(1/2)}, D_{q_1}^{(1/2)}]_- \, - {1\over 2}  =\, t_{21} \DD, \quad 
      [D_{p_2}^{(1/2)}, D_{q_2}^{(1/2)}]_- \, - \, {1\over 2} \, =\, - t_{12} \DD.
  \end{gathered} 
 \]

 \noindent (b) To prove the first statement it is enough to show that $\DD D_y^{(1/2)} \in \Dc_{S\Hen_2} \DD$.
 But  this follows from Proposition \ref {prop:super-lapl-anti-2}(b), as 
   $ \DD \,D_{y}^{(1/2)}  =  -   D_{y}^{(3/2)} \, \DD \,  \in \, \Dc_{S\Hen_2} \DD$.
 The second statement follows from part (a).\qed 

\vskip .2cm
  
 Consider the sheaf $\ul\Ker(\DD)\subset S\Lc_{1/2}$
% \be\label{eq:SL^D-gen2}
 %S\Lc_{1/2}^\DD := \ul\Ker \{ \DD: S\Lc_{1/2} \lra S\Lc_{3/2} \}  
 % \ee
  whose sections will be called {\em harmonic superforms} of weight $1/2$ 
  (another name would be  {\em singular superforms}, cf. Remark 
  \ref {rem:sing-forms-gen2}).
  
  \begin{cor}\label{cor:heis-harm-gen2}
  The operators $D_y^{(1/2)}, y\in\gen_\1$,  
  preserve $\ul\Ker(\DD)$ and
  give rise to a morphism of algebras $\Heis_2 \to \End (\ul\Ker(\DD))$. 
  \end{cor}
  
  \noindent {\sl Proof:} Follows from Proposition \ref{prop:even-Heis-gen2},
  as $\ul\Ker(\DD) = \ul \Hom_{\Dc_{S\Hen_2}}(\Mc_\DD, S\Lc_{1/2})$. 
  \qed 
  
 \begin{rem} \label{rem:trunc-1-gen2}
 Written in components, the $D_y^{(1/2)}$ are $4\times 4$ matrix differential operators. 
 From this point of view,  the $3\times 3$ matrix operators  $Q_i, P_i, R_0$ used in
  \cite {KK-Takei, berenstein} are obtained, up to minor changes,  by $3\times 3$ truncation of 
  our $D_{p_i}^{(1/2)}$,
  $D_{q_i}^{(1/2)}$, $\DD$ respectively, i.e.,  by imposing
  the ``constraints'' (vanishing of the component with $\xi_1\xi_2$) by hand in advance instead
  of obtaining them as a consequence of $\DD \Phi=0$ as we do.
   Because of the truncation, the commutation
  relations of $R_0$ with the $Q_i, P_i$ are more complicated than those of $\DD$ with the $D_y^{(1/2)}$
  and additional  iterated commutators need to be considered in order to obtain a counterpart of
  Proposition \ref {prop:even-Heis-gen2}(b). 
  \end{rem}
  
  \paragraph{Exponential solutions of $\DD\Phi=0$: super-Gaussians.} 
  \label{par:super-gauss-2}
  The operator $\DD$,
  written in components as a $4\times 4$ matrix differential operator, has constant coefficients
  in its dependence  on the even coordinates
  $v=\|v_{ij}\|$. We can consider it defined on all of $\CC^3 = \LG^\circ(\CC^4)$. 
   So it is classical to express solutions of $\DD\Phi=0$ in terms of exponentials.
  Let us recall the general formalism, following \cite{palamodov}. 
  
  \vskip .2cm
  
  Let $A = \|a_{ij}\|_{i=1,\cdots,r}^{j=1,\cdots, s}$ be an $r\times s$ matrix   differential operator with
  constant coefficients defined on $\CC^d$ with coordinates $v_1,\cdots, v_d$. Thus each
  each $a_{ij}$ is a polynomial in $\del = (\del/\del v_1, \cdots, \del/\del v_d)$. 
 Replacing $\del$ by a vector of variables $\lambda = (\lambda_1, \cdots, \lambda_d)$
 taking values in $(\CC^d)^*$, we get the {\em formal Fourier transform} (or, equivalently, the total symbol)
\[
\wh A(\lambda) \,=\,\|\wh a_{ij}(\lambda)\| \,=\,\sigma_{\wh A}(\lambda) \, \in \, \Mat_{r\times s}(\CC[\lambda]). 
\]
 We consider the homogeneous system
 \[
 Af=0, \quad f=(f_1, \cdots, f_s)^t, \quad f_i = f_i(v) \in \Oc_{\CC^d}. 
 \]
 Let $\mu \in (\CC^d)^*$. A {\em simple exponential soliution} of $Af=0$ with exponent $\mu$
 is a solution of the form
 \[
 f_i(v) \, =\, e^{\mu v}\cdot c, \quad c= (c_1, \cdots, c_s)^t \in \CC^s, \quad \wh A(\mu) c=0.  
 \]
 More conceptually, the sheaf of  solutions of $Af=0$ is
 \[
 \ul\Hom_{\Dc_{\CC^d}}(\Qc, \Oc_{\CC^d}), \quad
 \text{where}\quad  \Qc = \Coker\bigl\{ \Dc_{\CC^d}^{\oplus r} \buildrel \cdot A \over\lra 
 \Dc_{\CC^d}^{\oplus s}\bigr\}.
 \]
 Applying the formal Fourier transform,we get the $\CC[\lambda]$-module
 \[
 \wh Q \,=\,\Coker \bigl\{\CC[\lambda]^{\oplus r} \buildrel \cdot \wh A(\lambda) \over\lra
 \CC[\lambda]^{\oplus s}\bigr\}
 \]
 and let $\wh\Qc$ be the corresponding coherent sheaf on $\AAA^d = \Spec \, \CC[\lambda]$. 
 For $\mu\in (\CC^d)^*$ a $\CC$-point of $\AAA^d$, the fiber
  $\wh \Qc_\mu :=  \wh Q\bigl/\sum_i (\lambda_i-\mu_i)\wh Q$ of $\wh \Qc$ at $\mu$ is identified with
  $\Coker \{ \CC^r \buildrel \cdot \wh A(\mu) \over\lra \CC^s\bigr\}$. So the space of
  simple exponential solutions with exponent $\mu$ is identified with the dual space
  \[
  {\wh\Qc}_\mu^* \, := \,\Hom_\CC(\wh\Qc_\mu, \CC) \,=\,\Ker \{ \CC^s \buildrel\wh A(\mu) \cdot
  \over\lra 
  \CC^r\bigr\}. 
  \]
  We specialize this to $d=3$, $r=s=4$ and $A=\DD$ written in components but express the answer
  in  the super-language. 
  
  \begin{prop}\label{prop:sGauss-2}
  (a) The simple exponential solutions of $\DD\Phi=0$ are given by the {\em super-Gaussians}
  \[
 \theta_x(T, \xi) =   e^{\pi \i\,\,  x^t T x}  (1+\sqrt{-\pi \i}\, \,  \xi\cdot x)  = 
 e^{ \pi \i\,\,  x^t T x + \sqrt{-\pi \i}\, \,  \xi\cdot x},\quad x = (x_1, x_2)^t \in \CC^2. 
\]
 Here $\xi\cdot x = \sum_{i=1}^2 \xi_i x_i$, as $\xi$ is a row vector.
 The exponent of $\theta_x$ is equal to $\pi \i \, x\cdot x^t = \| \pi \i x_\mu x_\nu\|$ considered
 as a point of $\Sym_{2\times 2}(\CC)$ which we identify with $(\CC^3)^*$. 
 
 \vskip .2cm
 
 (b) In addition to $\DD\theta_x=0$,
 the super-Gaussians satisfy the differential equations
 \[
 \begin{gathered}
 \sqrt{-4\pi \i} \,\,   (D^{(1/2)}_{p_\nu} \, \theta_x)(T, \xi) =- 2 \pi \i\,\,  x_\nu\,  \theta_x(T, \xi),
 \\
\sqrt{-4\pi \i} \, \,  \,  (D^{(1/2)}_{q_\nu}\, \theta_x) (T, \xi) = -{\del\over\del x_\nu} \theta_x(T, \xi).
 \end{gathered} 
 \]
 (c) Let $v = \|v_{ij}\|\in \Sym_{2\times 2}(\CC)$. Let $\Phi$ be a holomorphic solution of $\DD\Phi=0$
 defined near $v$ and let $k\geq 0$. There exists a finite linear combination 
 $\Phi' = \sum_{i=1}^N c_i \theta_{x^{(i)}}(T,\xi)$, $x^{(i)}\in \CC^2$, such that $\Phi-\Phi'$ vanishes
 at $v$ together with all the partial derivatives of order $\leq k$. 

  \end{prop}
  
 Note that part (b) illustrates the Heisenberg relations of Proposition 
 \ref {prop:even-Heis-gen2} in the case of exponential solutions. 
 
 \vskip .2cm
 
 \noindent {\sl Proof:} (a) The equality $\DD\theta_x=0$ is verified directly. The simple exponential
 nature of $\theta_x$ w.r.t. the even coordinates (i.e., the fact that it is an eigenfunction for
 the translation by any $v=\|v_{ij}\|\in\Sym_{2\times 2}(\CC)$) is also clear from the definition.
 Equally clear is the statement about the exponent, as $x^t T x$ is the scalar product of $T$ with
 $x x^t$.
 
  To see that the $\theta_x$ exhaust all simple exponential soliutions of
 $\DD\Phi=0$, we use the component form \eqref {eq:DD-2}  consisting of the Laplace and Dirac
 equations and of the condition $F=0$. It implies that the (even) $\Dc_{\CC^3}$-module $\Qc = \eps^*\Mc_\DD$
 associated to the component form of $\DD\Phi=0$ (i.e., the module corresponding to the
 super-$\Dc$-module $\Mc_\DD$ above by the Penkov equivalence
 $\eps^*$) splits into the direct
 sum $\Qc = \Qc_\Delta \oplus \Qc_\dirac$ of the modules associated to the two equations. 
 
 Let $Z\subset \Sym_{2\times 2}(\CC) = (\CC^3)^*$ be the quadratic cone formed by 
 degenerate matrices. Thus $\CC^3\times Z\subset T^*\CC^3$ is the characteristic variety of
 both $\Qc_\Delta$ and $\Qc_\dirac$ or, equivalently,
 \[
 Z\,=\,\Supp\, \wh \Qc_\Delta \,=\,\Supp\, \wh \Qc_\dirac \,=\,\Supp\,  \wh\Qc. 
 \]
 Consider the squaring map
 \be\label{eq:squaring-2}
 \rho: \CC^2 \mapsto Z, \quad x = \bpm x_1 \\ x_2  \epm \mapsto 
 x\cdot  x^t = \|x_\mu  x_\nu\|,
 \ee
which is $2:1$ outside $0$.
 Then it is straightfoward that $\wh \Qc = \rho_* \Oc_{\CC^2}$ so that $\wh \Qc_\Delta$
 and $\wh\Qc_\dirac$ correspond to the parts in the direct image  which are
 even and odd respectively
 with respect to the involution $x\mapsto -x$. Therefore the $\theta_x$, $x\in \CC^2$, exhaust all the 
 vectors in all the fibers of $\wh\Qc$, i.e., all the simple exponential solutions. 
 
 \vskip .2cm
 
 (b) Direct computation. 
 
 \vskip .2cm
 
 (c) This is a version of the Malgrange approximation theorem
 \cite[\S 2.2]  {malgrange-approx} but formulated for the adic topology on the space
 of formal series at $v$. More precisely, the Malgrance approximation uses, in general,
 polynomial-exponential  rather than simple exponential solutions. Because
 of the identification $\wh\Qc = \rho_*\Oc_{\CC^2}$, such solutions correspond to iterated
 partial derivatives of the $\theta_x$ in $x$, and such derivatives can be approximated
 by finite differences thus reducing  to the finite linear combinations of the $\theta_x$ in the end.
 \qed

  \vfill\eject
 
 \subsection{The Riemann theta function from SUSY point of view}
 
 \paragraph{The genus $2$ Riemann theta.} The Riemann theta function is defined
 \cite{mumford}  by 
 \be
 \Theta(z, T) = \sum_{\bn\in \ZZ^2} \exp\bigl( \pi \i \,  \bn^t T \bn  + 2\pi \i \bn^t \cdot  z\bigr),
 \quad T\in \Hen_2, \,\, z\in \CC^2. 
 \ee
 Here both $\bn$ and $z=(z_1, z_2)^t$ are considered as column vectors. It satisfies the
``heat equations''
\be
\del_{ij} \Theta(z,T) = {-\i \over 2\pi} \,\, {\del^2 \Theta \over \del z_i \del z_j},
\quad 1\leq i,j\leq 2, \quad \del_{ij} := {\del\over \del t_{ij}} + {\del\over \del t_{ji}}, 
\ee
and the quasi-periodicity conditions
\[
\Theta( z+\bm, T) =\Theta(z, T), \quad \Theta (z+T\bm, T) = 
\exp\bigl( - \pi \i \, \bm^t T \bm - 2\pi \i\,  \bm^t z\bigr), \quad \bm\in \ZZ^2. 
\]
The Thetanullwert is the function
\be
\theta(T) = \Theta (0,T) = \sum_{\bn\in \ZZ^2} \exp(\pi \i \, \bn^t T \bn\bigr), \quad T\in\Hen_2. 
\ee
  It  satisfies the   periodicity
 \[
 \theta(T+U) = \theta(T), \quad U = \|u_{ij}\|\in \Sym_{2\times 2}(\ZZ), \,\, u_{11}, u_{22}\in 2\ZZ
 \]
  and modularity 
  \be\label{eq:modul-gen2}
  \theta (-T^{-1}) =  -\i \cdot \det(T)^{1/2}\, \theta(T). 
  \ee
  conditions.  They can be expressed by saying that 
  $\theta(T)\in \F_{1/2}=H^0(\Hen_2, \Lc_{1/2})$ is a Siegel form of weight $1/2$ which is 
  $\chi$-invariant with respect to
  a subgroup $\Gamma\subset \Mp(4,\ZZ)$ and a character $\chi: \Gamma \to \CC^*$
  described in  \cite[\S II.5]{mumford}.

  \paragraph{Characterization of $\theta(T)$ by super-differential equations.}
  Let us consider $\theta(T)$ as  Siegel superform  $\theta(T,\xi)$ independent on $\xi$,
  i.e., a section of $S\Lc_{1/2}$ on $S\Hen_2$. 
  Recall that $\gen =\osp(1|4)$ acts on $S\Hen_2$,  so the basis vectors 
  $p_\nu, q_\nu\in E=\gen_\1$, $\nu=1,2$ give rise to differential operators 
  $D_{p_\nu}^{(1/2)}$,
  $D_{q_\nu}^{(1/2)}$ of order $1/2$ on $S\Lc_{1/2}$.

   \begin{thm}\label{theta:susy-2}
   (a) 
     $\theta(T,\xi)$ is a solution of  the following  system of differential equations:
      \[
 \DD \Phi=0, \quad e^{\sqrt{-4\pi\, \i} \, D^{(1/2)}_{p_\nu}} \Phi = \Phi, \quad e^{\sqrt{-4\pi\, \i} \, D^{(1/2)}_{q_\nu}} \Phi = \Phi, \,\,\,\nu=1,2, 
 \]
 
   (b) 
 Any local holomorphic section  $\Phi=\Phi(T, \xi)$ of  $S\Lc_{1/2}$ on $S\Hen_2$ satisfying the system
 in (a),  is a constant multiple of $\theta(T)$. 
 \end{thm}
 
 \begin{rem} Written in components, the above system is a system of $4\times 4$ matrix PDE
 of infinite order. 
 After dropping  the component with $\xi_1\xi_2$  which vanishes
 {\em a posteriori}  in virtue of $\DD \Phi=0$,
 this  system truncates  to the  system of $3\times 3$ matrix  PDE of infinite order
 considered in 
  \cite{KK-Takei, berenstein}, see Remark  \ref{rem:trunc-1-gen2}. 
 Our approach
 elucidates the supersymmetric meaning of that $3\times 3$ system. 
 \end{rem}
 
\noindent  We  prove part (a) of the theorem right away and 
the proof of part (b)
will be finished in the next section.

 \vskip .2cm
 
    First of all, the heat equations  for $\Theta$ imply that $\theta$ is harmonic, i.e., $\Delta \theta =0$:
 \[
  \del_{11} \del_{22}\, \theta(T)  = \del_{12} \del_{21}\,
 \theta(T) \,\,\, ={-1\over 4 \pi^2} \cdot  {\del^4 \Theta \over \del z_1^2\,  \del z_2^2}\biggl|_{z=0}   (T).  
 \]
 This  implies $\DD\theta=0$
 because of the form \eqref{eq:DD-2} of $\DD$ and  the fact that $\theta$ is independent
 of $\xi$. 
 
 \vskip .2cm

 Further, as in  the genus 1 case,  we can express $\theta(T)$ in terms of the super-Gaussians of Proposition
 \ref {prop:sGauss-2} as
 \be \label {eq:theta-super-sum-2}
 \theta(T)  = \sum_{\bn\in \ZZ^2} \theta_\bn(T, \xi),
 \ee
 and so
  for any $\bn\in \ZZ^2$ we have, by part (b) of the same  proposition,
   \[
   \begin{gathered}
   e^{\sqrt{-4\pi \i}\,  D^{(1/2)}_{p_\nu}} \, \theta_\bn = e^{-2\pi \i n_\nu} \theta_\bn = \theta_\bn, 
   \\
   e^{\sqrt{-4\pi \i} \, D^{(1/2)}_{q_\nu}}\,  \theta_\bn = \theta_{\bn-\bee_\nu},
   \end{gathered}
   \]
   where $\bee_\nu$ is the standard basis vector of $\ZZ^2$. 
   So the  second equation of Theorem  \ref{theta:susy-2}(a) 
    holds term by term  w.r.t. 
    \eqref{eq:theta-super-sum-2} while the third  one holds
   by shifting the terms in the sum.  
   
\vskip .2cm

 We now prepare some background for the proof of   Theorem  \ref{theta:susy-2}(b).

   \paragraph{Diagonal factorization of $\theta(T)$ and of
   super-Gaussians.} 
   
   We consider the obvious diagonal embedding of the two copies of the
   Lobachevsky plane and its  (less obvious) super-extension:
   \be\label{eq:d-sd-2}
   \begin{gathered}
   \delta: \Hen\times\Hen \lra \Hen_2, \quad (\tau_1, \tau_2) \mapsto
   \bpm \tau_1 & 0 \\ 0 &\tau_2\epm, 
   \\
    S\delta: S\Hen \times S\hen \lra S\Hen_2, \quad
 \bigl((\tau_1, \xi_1), (\tau_2, \xi_2)\bigr) \,\mapsto \, \biggl(
 \bpm \tau_1 & 0 \\
 \xi_2 \xi_1 & \tau_2\epm, (\xi_1, \xi_2) \biggr). 
\end{gathered} 
   \ee
   Here we use the dependent coordinates $(T, \xi)$ on $S\Hen_2$ from 
 \S \ref {subsec:SLG-gen2}. The off-diagonal term $\xi_2\xi_1$ 
 ensures that the condition \eqref {eq:sigma-gen2} is satisfied. 
 
 The embedding $\delta$ corresponds to the diagonal embedding of
 groups
 \[
 \SL (2,\RR) \times\SL(2, \RR) \lra \Sp(4, \RR)
 \]
 and is equivariant under this group embedding in the obvious way. 
 The Riemann Thetanullwert  satisfies    the {\em diagonal factorization
  property}
  \be
  (\delta^* \theta) (\tau_1, \tau_2) = \theta_J(\tau_1) \theta_J(\tau_2),
  \ee
  where $\theta_J(\tau)$ is the single variable  Jacobi Thetanullwert from
  \S \ref {subsec:jacobi}. This follows directly from the definition. 
  
  As for the super-thickenings, we have two embeddings of Lie supergroups
  \be\label{eq:emb-G-i-2}
  \Gb^{(1)} = \OMp(1|2) \buildrel \eps_1\over \lra \Gb= \OMp(1|4)
  \buildrel \eps_2\over 
  \lla \Gb^{(2)} = \OMp(1|2), 
  \ee
   whose images {\em do not (super) commute}
   (although their underlying even groups do).  The corresponding
   Lie sub-superalgebras
   \[
   \gen^{(1)} = \osp(1|2) \lra \gen = \osp(1|4) 
   \lla \gen^{(2)} = \osp(1|2)
   \]
   are defined by their odd parts
   \[
   \gen^{(1)}_\1 = \CC p_1 \oplus \CC q_1, \quad \gen^{(2)}_\1 = \CC p_2
   \oplus \CC q_2 \quad\subset \quad E= (\CC^4, \omega) = \gen_\1,
   \]
   and we see that the super-commutator between these parts
   is not $0$. The embedding $S\delta$ is equivariant with respect
   to $\Gb^{(1)}\subset \Gb$ acting on the first factor and to
   $\Gb^{(2)}\subset \Gb$ acting on the second factor separately. 
   Its importance comes from the following factorization property
   of super-Gaussians $\theta_x(T,\xi)$, $x=(x_1, x_2)^t \in \CC^2$
   (see Proposition \ref {prop:sGauss-2} for definition).
   
   \begin{prop}\label{prop:sgauss-fac-2}
   We have
   \[
   \bigl( (S\delta)^* \theta_x\bigr) \bigl( (\tau_1, \xi_1), (\tau_2, \xi_2)
  \bigr) 
  \,=\,
  \theta_{x_1}(\tau_1, \xi_1) \cdot \theta_{x_2}(\tau_2, \xi_2), 
   \]
   where the factors in the RHS are the $1$-dimensional super-Gaussians
   \eqref {eq:super-gauss-1}. 
   \end{prop}
   
   \noindent{\sl Proof:} Direct checking. \qed
   
   \vskip .2cm
   
   Note that the $\theta(\tau_i, \xi_i)$ are inhomogeneous, so their
   product depends on the order, even up to sign. 
   We can obtain factorization
   with respect to the product in the other order, if we re-define
   $S\delta$ using an upper-triangular matrix with $\xi_1\xi_2$
   above the diagonal. 
      
   \vfill\eject
   
   \subsection {Analysis of  the Koszul complex}\label{subsec-koszul-2}
   
   In this section we prove Theorem  \ref{theta:susy-2}(b).

      \paragraph{ Reformulation via the  Koszul complex
      of harmonic superforms.} 
   Recall the sheaf  $\ul\Ker(\DD)\subset S\Lc_{1/2}$
   of harmonic superforms on $\Hen_2$. 
   The Heisenberg relations of Proposition \ref{cor:heis-harm-gen2} 
   imply that the four  differential operators of infinite order
   \be\label{eq:AB-nu-gen2}
 A_\nu = e^{\sqrt{-4\pi \i} D_{p_\nu}^{(1/2)}} -1, \quad
 B_\nu =  e^{\sqrt{-4\pi \i} D_{q_\nu}^{(1/2)}} -1, \quad  \nu=1,2, 
 \ee
acting on $\uKer(\DD)$, commute with each other. 
So we can form  the corresponding Koszul complex  of sheaves
 \[
 \Cc^\bullet = \uKer(\DD) \otimes_\CC \Lambda^\bullet (\CC^4) \,\, = \,\,
 \bigl\{ \uKer(\DD) \buildrel (A_1, A_2, B_1, B_2) \over
 \lra \uKer(\DD)^{\oplus 4} \lra \cdots \lra \uKer(\DD)\bigr\},
 \]
situated in the degrees $[0,4]$.  Its differentials are constructed from
the $A_\nu, B_\nu$ in a standard way. 
Thus $\ul H^0(\Cc^\bullet)$
 is the sheaf of solutions of the system of Theorem \ref {theta:susy-2}. 
 As in  \S \ref {subsec:koszul-1}, 
 we deduce that theorem
 from the next more general
 statement.

 \begin{thm}\label{thm:koszul-LD-2}
 We have $\ul H^0(\Cc^\bullet) \= \ul\CC_{\Hen_2}$, the constant
 sheaf generated by $\theta(T)$ while $\ul H^i(\Cc^\bullet)=0$ for $i>0$. 
 \end{thm}

 Further, similarly  ito  \S \ref {subsec:koszul-1}, we deduce
 Theorem \ref {thm:koszul-LD-2} from the  next two propositions:

  \begin{prop}\label{prop:koszul-lc-2}
 The sheaves $\ul H^i(\Cc^\bullet)$ are locally (and hence globally)
  constant on $\Hen_2$. 
 \end{prop}
 
 \begin{prop}\label{prop:koszul-gl-2}
  Consider the diagonal embedding $\delta$ from \eqref {eq:d-sd-2}
 and the corresponding 
 sheaf-theoretic restriction $\delta^{-1}\Cc^\bullet$, a complex on
 $\Hen\times\Hen$.  Then 
 $\ul H^0(\delta^{-1}\Cc^\bullet) \= \ul\CC_{\Hen\times\Hen}$
 is the constant sheaf spanned by  the restriction $\theta(\delta(\tau_1, \tau_2)) = \theta_J(\tau_1)\theta_J(\tau_2)$,
 while $\ul H^i(\delta^{-1}\Cc^\bullet) = 0$ for $i>0$. 
  \end{prop}
  
  Our arguments for these propositions adapt some steps 
 from \cite{KK-Takei} into the supersymmetric language  which
simplifies them considerably. We  present them in the rest
of this section.

  \begin{rem}    
    In fact,  \cite{KK-Takei}  proves  a more general
  statement. 
  The operators $D_y^{1/2}$, $y\in \gen_\1$
  make sense on the entire
  Minkowski superspace $\CC^{3|2}$, not just on its future tube
 $S\Hen_2$. Therefore the
 complex  of sheaves $\Cc^\bullet$ can be defined on the entire $\CC^3 = \Sym_2(\CC)$, not just on $\Hen_2$. 
 The results of    \cite{KK-Takei} 
 amount to determination of $\ul H^\bullet(\Cc^\bullet)$
 on all of $\CC^3$. We do not address this problem here. 
  \end{rem}

   %\vfill\eject
   
   \paragraph{The Koszul complex at the level of $\Dc^\oo$-modules. }
   \label{par:Koszul-D-gen2}
   Consider the harmonic $\Dc_{S\Hen_2}$-module 
   $\Mc_\DD = \Dc_{S\Hen_2}/\Dc_{S\Hen_2}\DD$ from
   \$ \ref {subsec:super-Lap-gen2} \ref {par:even-heis-gen2}. 
   It follows from explicit matrix form \eqref{eq:DD-2} of $\DD$ that
   the right multiplication by $\DD$ is injective on $\Dc_{S\Hen_2}$
   and so the complex
   \[
   \Rc^\bullet = \bigl\{ \Dc_{S\Hen_2} \buildrel \cdot \DD \over \lra
   \Dc_{S\Hen_2} \bigr\}, 
   \]
 situated in degrees $0, -1$, is a left free resolution of $\Mc_\DD$.

 \vskip .2cm
 
 The module $\Mc_\DD$ inherits a good filtration from $\Dc_{S\Hen_2}$
 and the endomorphisms $D_y^{(1/2)}: \Mc_\DD \to\Mc_\DD$, $y\in \gen_1$,
 have effective order with respect to this filtration, as defined in \S 
 \ref{subsec:orders}, less or equal to $1/2$. 
 Let $\Mc_\DD^\oo = \Dc_{S\Hen_2}^\oo \otimes {\Dc_{S\Hen_2}}\Mc_\DD$. 
 By Proposition \ref{prop:exp-endo} we have well defined
 endomorphisms 
 $e^{z D_y^{(1/2)}} \in \End_{\Dc^\oo_{S\Hen_2}}(\Mc_\DD^\oo)$
 for any $z\in \CC$ and so can define  $A_\nu, B_\nu \in \End_{\Dc^\oo_{S\Hen_2}}(\Mc_\DD^\oo)$
 by   \eqref {eq:AB-nu-gen2}. 
 By Proposition \ref {prop:even-Heis-gen2}, the four endomorphisms $A_\nu, B_\nu$  
 commute with each other  and so we can form the associated Koszul complex
 of $\Dc_{S\Hen_2}^\oo$-modules
 \[
 \Nc^\bullet = \Mc_\DD^\oo \otimes_\CC \Lambda^{-\bullet}(\CC^4),
 \]
 situated in degrees $[-4,0]$.  
 
 \begin{prop}\label{prop:C-as-rhom}
 We have a quasi-isomorphism of sheaves $\Cc^\bullet \= \ul{R\Hom}_{\Dc_{S\Hen_2}^\oo}(\Nc^\bullet, S\Lc_{1/2})$.
 \end{prop}
 
 \noindent{\sl Proof:} This follows from two facts.
 
 (1) The sheaf $\Dc^\oo_{S\Hen_2}$ is
 flat over $\Dc_{S\Hen_2}$. By the Penkov equivalence, this follows from the corresponding
 purely even statement which is discussed in \cite[\S 8.2.15]{bjork}. 
 
 (2) The sheaves
 $\ul \Ext^p_{\Dc_{S\Hen_2}}(\Mc_\DD, S\Lc_{1/2})$  vanish for $p>1$, so
 \[
 \uKer(\DD) \, \= \, \ul{R\Hom}_{Dc_{S\Hen_2}}(\Mc_\DD, S\Lc_{1/2})
 \]
 Indeed, notice first that these  sheaves must be supported on  analytic
 subvarieties of codimension $\geq 1$. See  [51, Th. 2.4.2+Rem.3] for the purely even case
 to which our case is reduced by the Penkov equivalence. But $\Mc_\DD$ is 
 $\OMp(1|4)$-equivariant
 and so are our $\ul\Ext^p$ which must therefore vanish. \qed
 
 \vskip .2cm
 
 Next, we notice that $\Nc^\bullet$ is strictly perfect over $\Dc_{S\Hen_2}^\oo$,
  in fact it is quasi-isomorphic
 to a finite complex of free $\Dc_{S\Hen_2}^\oo$-modules. This is because
 $\Mc_\DD$ is  quasi-isomorphic to $\Rc^\bullet$ above,
 therefore $\Mc_\DD^\oo$ is  quasi-isomrophic to a finite complex of
  $\Dc_{S\Hen_2}^\oo$-modules  and so is 
 each term of $\Nc^\bullet$ is perfect over $\Dc_{S\Hen_2}^\oo$. 
 Our statement then follows from the next standard lemma, cf. \cite[Lem.1.5]{KK-Takei}.
 
 \begin{lem}\label{lem:perf-complex}
 let $R$ be a ring and $C^\bullet$ be a finite complex of $R$-modules.
 If each $C^i$ is quasi-isomorphic to a finite complex if free $R$-modules,
 then so is $C^\bullet$ itself.
 \end{lem}
 
 \noindent{\sl Proof of the lemma}: Induction on the length of $C^\bullet$, using
 the folllowing fact: if $K^\bullet, L^\bullet$ are finite complexes of free $R$-modules,
 then $\Hom(K^\bullet, L^\bullet)$ in the derived category of $R$-modules
 is identified with the set of homotopy equivalence classes of morphisms of
 complexes of $R$-modules $K^\bullet\to L^\bullet$.\qed

 \vskip .2cm
 
   So we can invoke Theorem \ref {thm:Char-SS} and deduce
 Proposition  \ref {prop:koszul-lc-2}   from the next
  
  \begin{prop}\label{prop:Ch(N)-gen2}
  The characteristic variety of $\Nc^\bullet$ is contained in  $\ul 0_{\Hen_2}$, the zero section of $T^*\Hen_2$. 
  \end{prop}
  
 % \vfill\eject
  
  \paragraph{Study of the characteristic variety: Proof of Proposition \ref {prop:Ch(N)-gen2}.  }
  \label{par:char-var-gen2}
 We use the coordinates $(x_{ij}, \xi_i)$ on $S\Hen_2$, the coordinates
 $T=\|t_{ij}\|$ on $\Hen_2$  and denote by $\lambda_{ij}$
 the momenta coordinates on $T^*\Hen_2$ corresponding to the vector fields $\del_{ij}$. 
 We arrange them into a symmetric matrix $\lambda = \|\lambda_{ij}\|\in \Sym_{2\times 2}(\CC)$. 
 As in the proof of Proposition \ref {prop:sGauss-2},
  let $Z\subset \Sym_{2\times 2}(\CC)$ be the hypersurface defined by
   $\det(\lambda)=0$ and
  \[
  \ul Z \,=\, \Hen_2 \times Z \,\,\subset \,\, \Hen_2 \times \Sym_2(\CC) \,=\, T^*\Hen_2.
  \]
 
 \begin{prop}
 
a) $\Ch(\Mc_\DD)$ is equal to $\ul Z$. 
 
 \vskip .2cm
 
 (b) $\Ch(\Nc^\bullet)$ is contained in $\ul Z$. 
 \end{prop}
 
 \noindent{\sl Proof:}  Similarly to what was said in the proof of Proposition 
\ref {prop:sGauss-2},
 (a)  follows from the explicit matrix 
form \eqref {eq:DD-2} of $\DD$ in components,
 where it consists of the Laplace and Dirac operators. Further, 
 (b) holds since $\Nc^\bullet$ is constructed out of copies of $\Mc_\DD^\oo$. 
  \qed
 
 \vskip .2cm
 
 In components, we view $S\Lc_{1/2}$ as $\Oc_{\Hen_2} \otimes_\CC \Lambda[\xi_1, \xi_2]
 =\Oc_{\Hen_2}\otimes \CC^4$ and operators there as $4\times 4$ matrix operators
 in $\Oc_{\Hen_2}$.

 Let $Z_{A_i}, Z_{B_i}\subset T^*\Hen_2$ be the loci of non-invertibility of $A_i, B_i$
 in $\Ec_{S\Hen_2}^\RR = \Mat_4(\Ec_{\Hen_2}^\RR)$. Put $Z_A=Z_{A_1}\cap Z_{A_2}$ and
 $Z_B = Z_{B_1}\cap Z_{B_2}$. As in \S \ref {subsec:koszul-1}\ref {prop:char-var-gen1},
 we notice that  the Koszul complex is  exact whenever one of the commuting operators
giving rise to it is invertible and so
\[
\Ch(\Nc^\bullet) \,\subset \, \ul Z\cap Z_A \cap Z_B. 
\]
To find $Z_{A_i}$, we apply Proposition \ref{thm:aoki-inv}  to $C=\sqrt{4\pi \i}\,  D_i$, so $A_i=e^C-1$. Decomposing $\Lambda = \Lambda[\xi_1, \xi_2]$ as 
$\Lambda^0\oplus \Lambda^1\oplus \Lambda^2$ by degree, we write $D_i$
as a block matrix differential operator
\[
D_i \, = \, {\del\over\del \xi_i} - \sum_j \xi_j \del_{ij}  = \bpm
0& \del/\del\xi_i & 0 
\\
-\sum_j \xi_j\del_{ij} & 0 & \del/\del\xi_i
\\
0& -\sum _j \xi_j\del_{ij} & 0
\epm.
\]
Let $D_{i, {p,q}}:  \Oc_{\Hen_2}\otimes\Lambda^q  \to \Oc_{\Hen_2}\otimes\Lambda^p$
be the block matrix elements of $D_i$. 
These block matrix elements consist of  $0$s on the diagonal, the operators $\del/ \del\xi_i: \Lambda^q\to
\Lambda^{q-1}$, of order $0$, just above the diagonal and the operators 
$-\sum _j \xi_j\del_{ij}: \Lambda^q\to\Lambda^{q+1}$, of order $1$, just below the diagonal. So putting, in the situation of Proposition \ref{thm:aoki-inv},
$r_0=1/2$, $r_1 = 1$, $r_2=3/2$ (corresponding to the summands $\Lambda^0$,
$\Lambda^1$ and $\Lambda^2$) and $\rho=1/2$, we see that the condition
$\ord(D_{i, (p,q)}\leq r_p-r_q + \rho$ is satisfied and, moreover, 
$\sigma_{D_{i, (p,q)}}^{r_p-r_q+\rho}(\lambda)$ coincides with the total symbol
$\sigma_{D_{i, (p,q)}}(\lambda)$. So the block matrix in Proposition 
\ref{thm:aoki-inv} is simply $\sqrt{4\pi \i} \, \sigma_{D_i}(\lambda)$
and the sufficient condition for $A_i$ to be invertible in 
$\Mat_4(\Ec^\RR_{\Hen_2})_{(T,\lambda)}$ given there
 is that no eigenvalue of this matrix lies in $\i\,\RR$. 

\begin{prop}\label{prop:ev-sigma-Di-2}
The only eigenvalues of $\sigma_{D_i}(\lambda) =   {\del/ \del \xi_i} - \sum_j  \lambda_{ij} \xi_j$ (as an endomorphism of $\Lambda = \CC^4$)  are $\pm \sqrt{-\lambda_{ii}}$. 
\end{prop}

\noindent{\sl Proof:} We consider the case $i=1$, the  case $i=2$  is similar. 
 For  a complex number $\alpha$ 
let $\sigma(\alpha) = {\del\over\del \xi_1} - \alpha \xi_1\in \End(\Lambda[\xi_1])$. 
As follows from the explicit  matrix form \eqref {eq:sigma-A-gen1} of
$\sigma(\alpha)$, its eigenvalues are $\pm\sqrt{-\alpha}$.
Let us write $\Lambda[\xi_1, \xi_2] = \Lambda[\xi_1] \otimes \Lambda[\xi_2] = \Lambda[\xi_1] \oplus \xi_2\Lambda[\xi_1]$
and consider $\sigma_{D_1}(\lambda)$ as a block $2\times 2$ matrix
with respect to this decomposition. This matrix has the form
\[
\sigma_{D_1}(\lambda) = {\del\over\del\xi_1} - \lambda_{11} \xi_1 - \lambda_{12} \xi_2 = 
\sigma(\lambda_{11})\otimes \Id  -\lambda_{12}  (\Id\otimes \xi_2) = 
\bpm
 \sigma(\lambda_{11})  & 0 \\ - \lambda_{12}\cdot  \Id &\sigma(\lambda_{11})
\epm,
\]
whence the statement. \qed

\vskip .2cm

As in  \S \ref {subsec:koszul-1} \ref {prop:char-var-gen1}, we now  invoke Proposition  \ref {thm:aoki-inv},
to find
\[
\ul Z \cap Z_A \,\,\subset \,\, \bigl\{ (T, \lambda) \bigl| \,\, \lambda\in Z, \,\, \lambda_{11},
\lambda_{22}\, \in\,  \i \, \RR_{\leq 0}\bigr\}. 
\]
 Recall  the squaring map $\rho: \CC^2\to \ZZ$ from \eqref {eq:squaring-2}.  
  It is $2:1$ everywhere except $0$. If $(x,\lambda) \in \ul Z\cap Z_A$, then 
 $y=\rho^{-1}(\lambda)$,
 defined uniquely up to sign, must have $y_i^2 = \lambda_{ii}\, \in\,  \i \, \RR_{\leq 0}$,
 so $y= \pm \sqrt{-\i} \, z$, $z\in \RR^2$. 
 
 \vskip .2cm
 
 Next, $Z_B$ is obtained from $Z_A$ by the change change of variables given by
 the modular transformation $S: T \mapsto -T^{-1}$. Indeed, $q_i\in\gen_1$
 is obtained from $p_i$ by the adjoint action of the matrix
 $\bpm 0&-1_2\\ 1_2 & 0 \epm \in\Sp(4, \ZZ)$ effecting this transformation which implies
 a similar relation between $B_i$ and $A_i$.  In a more precise notation, we have
 $Z_B = (dS)^t (Z_A)$. Now, a standard matrix calculation shows that
 \[
 dS \,=\, - d(T^{-1}) \,=\, T^{-1} (dT)  T^{-1} \,=\, S (dt) S. 
 \]
 That is,  if we identify the tangent space ot $\Hen_2$ at any point $T$ with
 $\Sym_2(\CC)$, then $d_T(S)$ is the endomorphism
 \[
 X \, \mapsto \, T^{-1} X T^{-1} = S(T) X S(T): \,\,\,\, \Sym_2(\CC) \lra \Sym_2(\CC). 
 \]
 Since $S(T)$ is symmetric, the $3\times 3$ matrix of this endomorphism in the basis
 of matrix units, is symmetric as well, so we can computationally identify $(dS)^t$
 with $dS$. It remains now to establish
 
 \begin{lem}
 We have
 \[
 \ul Z \cap Z_A \cap Z_B \,=\, \ul Z \cap Z_A \cap (dS)(Z_A) \,=\, \ul 0_{\Hen_2}
 \]
 is the zero section of $T^*\Hen_2$. 
 \end{lem}
 
 \noindent {\sl Proof:} Let  $(T, \lambda) \in \ul Z\cap Z_A$. Then
  $\lambda = -\i \, z\cdot z^t$, $z\in \RR^2$. If, at the same time,  $(T, \lambda) \in \ul Z \cap Z_B$, then
  by the above $(d_TS)(\lambda) = -\i\,  w\cdot  w^t$ for some $w\in \RR^2$ as well. But  
  \[
   (d_TS)(\lambda) = -\i\,  S(T) z \cdot z^t S(T) = -\i \, (S(T)z) \cdot (S(T) z)^t
   \]
 which implies, by the $2:1$ property of $\rho$, that $w=\pm S(T) z$, so
 $S(T)z\in \RR^2$. 
This means that the real matrix $\Im (S(T))$  
annihilates $z\in \RR^2$
which implies $z=0$ since $\Im( S(T))$ is positive definite.
This proves the lemma and   Proposition \ref {prop:Ch(N)-gen2}. \qed

 \paragraph{Cauchy problem for the diagonal embedding: Proof of Proposition
 \ref {prop:koszul-gl-2}. } 
 Recall the super-diagonal embedding $S\delta$ from
 \eqref {eq:d-sd-2}. 
   We study the $\Dc$-module pullback of $\Mc_\DD$ under $S\delta$.
 Note that while we  consider  $\Dc$-modules  over supermanifolds,
 their characteristic varieties live in the cotangent bundles of underlying ordinary
 (even) manifolds. With this understanding, we have: 
 
 \begin{prop}\label{prop:cauchy-2}
 (a) The subvariety $\delta(\Hen\times\Hen) \subset \Hen_2$ is non-characteristic
 for $\Mc_\DD$.
 
 \vskip .2cm
 
 (b) The (super) $\Dc$-module pullback $(S\delta)^* \Mc_\DD$  on $S\Hen\times S\Hen$
 is identified with $\Dc_{S\Hen\times S\Hen}$ while the higher (left) derived functors
 $L^m (S\delta)^* (\Mc_\DD)$ vanish for $m \geq 1$. 
 
 \vskip .2cm
 
 (c) With respect to the identification in (b), the endomorphisms $D_{p_1}^{(1/2)}$,
 $D_{q_1}^{(1/2)}$ of $\Mc_\DD$ induce 
 the operators of right multiplications on
 $D_p^{(1/2)}\otimes \Id$, $D_q^{(1/2)}\otimes \Id$ on $\Dc_{S\Hen \times S\Hen}$,
 while 
 $D_{p_2}^{(1/2)}$,
 $D_{q_2}^{(1/2)}$ induce $(-1)^\F \otimes D_p^{(1/2)}$ and 
 $(-1)^\F\otimes D_q^{(1/2)}$
 respectively.
 \end{prop}
 
 Here $D_p^{(1/2)} , D_q^{(1/2)}$  are the $1$-dimensional operators
 considered as sections of $\Dc_{S\Hen}$, see Theorem \ref {thm:1d-theta-charact}. 
 
 \begin{rem}
 The appearance of $(-1)^\F$ in the first tensor factor for the identifications
 of  $D_{p_2}^{(1/2)}$,
 $D_{q_2}^{(1/2)}$ entails that they  {\em commute} 
 with the identifications for $D_{p_1}^{(1/2)}$,
 $D_{q_1}^{(1/2)}$. If $(-1)^\F$ was omitted,   we would have odd operators
 acting on  different tensor factors, which  {\em anticommute}. 
 \end{rem}
 
 \vskip .2cm
 
 \noindent{\sl Proof:} (a) By definition \cite{kash-thesis}, non-characteristicity means that
 the conormal bundle $T^*_{\Hen\times\Hen} \Hen_2 \subset T^*\Hen_2$
 intersects $\Ch(\Mc_\DD) = \ul Z$ 
 only along the zero section. Using the matrix
 coordinates $\lambda$ in $T^*\Hen_2$ as in \S \ref {par:char-var-gen2},
 we find that $T^*_{\Hen\times \Hen}\Hen_2$ is given by $\lambda_{11}= \lambda_{22}=0$, 
 while $\ul Z$ is given by $\det(\lambda)=0$, so the intersection is given by
  $\lambda_{12}^2=0$ which gives $\lambda=0$. 
  
  \vskip .2cm
  
  (b) The statement means that the Cauchy data for the equation $\DD\Phi=0$ on
  the hypersurface $S\Hen \times S\Hen$ assemble naturally into a
  single (super) function on $S\Hen \times S\Hen$ which can be arbitrary
  (satisfies no relations). To see this, let us first
  write $\Phi$ in components
  using the  independent coordinates $x=\|x_{ij}\|$ (a symmetric matrix) and $\xi = (\xi_i)$
  from   \eqref {eq:xij-gen2}
as 
  \[
  \Phi(x,\xi)  = \phi(x) + \sum_{i=1}^2 \psi_i(x) \xi_i + F(x) \xi_1 \xi_2, 
  \]
 see \eqref{eq:Phi-comp-gen2}. Then $\DD\Phi=0$ means 
 \eqref {eq:DD-2} that $\phi$
 satisfies the Laplace equation $\Delta \phi=0$, of second order,
 $\psi  =\sum \psi_i(x) \xi_i$ satisfies the Dirac equation $\dirac \psi=0$,
  of first order, while  $F=0$. 
 So the Cauchy data needed to determine the local solution near the hypersurface
 $x_{12}=0$ are:
 \[
 \phi\bigl(\diag(x_{11}, x_{22})\bigr), \,\,\,
  {\del \phi\over \del x_{12}} \bigl(\diag(x_{11}, x_{22})\bigr), \,\,\, \psi_i \bigl(\diag(x_{11}, x_{22})\bigr), \,\,\, i=1,2, 
 \]
 i.e., $4$ even functions which can be arbitrary.
  To interpret the collection of them together, let us pass
 to the dependent coordinates $(T,\xi)$ used in the definition 
 \eqref {eq:sdelta-11} of $S\delta$.
 Since the image of $S\delta$ is given by the equations
 \[
 t_{12}=0, \quad t_{21} = -\xi_1\xi_2,
 \]  
  the pullback  $\Phi\circ S\delta$
 will have the form
 \[
(\Phi\circ S\delta)
(\tau_1, \tau_2, \xi_1, \xi_2) \,=\,   \phi\bigl(\diag(\tau_1, \tau_2 )\bigr) + \sum_{i=1}^2 \psi_i \bigl(\diag(\tau_1,
 \tau_2)\bigr) \xi_i 
 -  {\del \phi\over \del x_{12}} \bigl(\diag(\tau_1, \tau_2)\bigr) \xi_1 \xi_2,
 \]
 the normal derivative giving the term with $\xi_1\xi_2$ by the Taylor formula. 
 This identifies the Cauchy data with an arbitrary  single superfunction, proving (b). 
 
 \vskip .2cm
 
 (c)  The statement identifies the action of $D_y^{(1/2)}, y\in \gen_\1$,
 on solutions $\Phi$ of $\DD\Phi=0$ in terms of their Cauchy data
 (restriction under $S\delta$). As such, it is a statement about
 (finite) jets. So it can be verified on a subspace $\Vc$ of solutions
 whose elements can realize any finite jet of any solution at any point.
 Let us choose for $\Vc$ the space of finite linear combinations
 $\Phi = \sum_{i=1}^N c_i \, \theta_{x^{(i)}}(T, \xi)$ of super Gaussians,
 see Proposition \ref {prop:sGauss-2} (c). 
 The Cauchy data of such solutions are immediate
 because of the factorization
 \[
   \bigl( (S\delta)^* \theta_x\bigr) \bigl( (\tau_1, \xi_1), (\tau_2, \xi_2)
  \bigr) 
  \,=\,
  \theta_{x_1}(\tau_1, \xi_1) \cdot \theta_{x_2}(\tau_2, \xi_2) 
   \]
(Proposition \ref {prop:sgauss-fac-2}). Therefore the product in $\Oc_{S\Hen_2}$
gives an identification
  $(S\delta)^*\Vc \,\= \, \Vc_1\otimes_\CC \Vc_2$, where 
  $\Vc_\nu, \nu=1,2$, is the space of finite linear combinations of
  $1$-dimensional super-Gaussians $\theta_{x_\nu}(\tau_\nu, \xi_\nu)$,
  $x_\nu \in \CC$. Now the action of the $D_y^{(1/2)}$ 
  on $\theta_x(T, \xi)$ is given by Proposition 
  \ref {prop:sGauss-2} (b). After restricting by $S\delta$, this gives
  the action  of $D_{p_1}^{(1/2)}$ and $D_{q_1}^{(1/2)}$ as
  $D_p^{(1/2)}\otimes 1$ and $D_q^{(1/2)}\otimes 1$ 
  (action on the first factor only), while 
  $D_{p_1}^{(1/2)}$ and $D_{q_1}^{(1/2)}$ act similarly on the
  second factor only, without the additional Koszul sign, which
  gives $(-1)^\F\otimes D_p^{(1/2)}$ and $(-1)^\F\otimes D_q^{(1/2)}$,
  as claimed.  
  
  \vskip .2cm
  
  This finishes the proof of Proposition \ref {prop:cauchy-2} and so
  of Proposition  \ref {prop:koszul-gl-2} and Theorem 
   \ref{theta:susy-2}.

   \vfill\eject
 \subsection{The (super) Weil representation of $\OMp(1|4)$}\label{subsec;super-weil-2}
 
 \paragraph{The group level: $\Mp(4)$.} \label{par:weil-group-2}
  The Weil representation $\varpi$
of $\Mp(4)$ is defined, similarly to 
  the  case of genus 1
 in \S \ref {subsec:super-weil-1},  on the space $\Wc= L_2(\RR^2)$.
 The coordinates in $\RR^2$ are denoted $x=(x_1, x_2)^t$. 
 As a projective representation of $\Sp(4, \RR)$, is defined
 uniquely  by its values on a generating subset:
 \be\label{eq:weil-expl-2}
 \varpi \left( \bpm 1_2&0 \\ C& 1_2 \epm\right)  u \,=\, e^{-\pi  \i \, x^t C x} \cdot u,\quad \text{while} \quad 
 \varpi  \left(\bpm 0& -1_2\\ 1_2&0 \epm\right) u  \,=\, \wh u, 
 \ee
 where  $C\in \Sym_{2\times 2}(\RR)$ and $\wh u$ is the Fourier transform
 of $u$. 
 
 As before, we have a splitting $\Wc=\Wc_\0 \oplus \Wc_\1$ into the
 subspaces of odd and even functions, each of which is an
 irreducible representation of $\Mp(4)$.

 \paragraph{The Lie algebra and supergroup levels: 
 $\sp(4)$, $\osp(1|4)$ and $\OMp(1|4)$.} \label{par:weil-2-lie}
 By differentiating $\varpi$, we get the well known representation of 
 $\gen_\0 = \sp(4)$  in $C^\oo(\RR^2)$  by quadratic differential operators.
 Let us recall the description of this representation and of its
 extension to a representation of $\gen = \osp(1|4)$. Both of them
 are defined on the 
  subspace $\CC[x] = \CC[x_1, x_2] \subset C^\oo(\RR^2)$ of polynomials. 
 
 \vskip .2cm

 For $v\in E = (\CC^4, \omega)$ we denote by $l(v)\in \osp(1|E)_\1$ the corresponding odd generator. 
 Let $S^\bullet E$ be the symmetric algebra of $E$; its commutative product  we denote 
 $\odot$. 
 Recall that we have a natural isomorphism
 \[
 l: S^2E \lra\sp(E) = \gen_\0, \quad l(v\odot w)(x) = {1\over 2} \bigl( \omega(v,x) w + \omega(w,x) v\bigr), \quad v,w,x\in E. 
 \]
 As in \S  \ref {subsec:super-Lap-gen2}
 \ref {par:even-heis-gen2}, let $\Heis_2 = \Heis(E, {1\over 2}\omega)$
 be the Heisenberg algebra corresponding to $E$ and the form 
 ${1\over 2}\omega$. 
 The subspace  $\gen_\0 = S^2E$ is embedded as a Lie algebra into $(\Heis(E), [a,b]=ab-ba)$ as the space of
quadratic Hamiltonians:
%\be\label{eq:S2E-Heis-2}
\[
S^2E \lra \Heis_2 \quad v\odot w \mapsto {1\over 2}  (v\cdot w + w\cdot v). 
\]
Recall also the identification of $\gen$ with $E\oplus S^2E$
(as a vector space) 
in Proposition
\ref  {osp4-s2}. From this, we deduce the following.
 
 \begin{prop}\label{prop:osp-heis-2}
  The  above embeddings of $E$ and $S^2E$ into $\Heis_2$ extend to a surjective homomorphism of associative
 algebras $h: U(\gen) \to \Heis_2$. \qed
 \end{prop}

 The algebra
  $\Heis_2$ is realized as the algebra of polynomial differential
 operators in $\CC[x] $ viia
 % \be\label{eq:osp(1|4)-weil}
 \[
 p_\nu \mapsto {1\over\sqrt{-4\pi \i }}\, \,  {\del \over \del x_\nu}, \quad q_\nu \mapsto \sqrt{-\pi \i } \,\, x_\nu. 
 \]
 The Lie algebra $\gen_\0=\sp(4)$ is then
 spanned by the differential operators
 \be\label{eq:xxdd-2}
 e_{\mu\nu } ={1\over 4\pi \, \i} \cdot  {\del^2\over\del x_\mu \del x_\nu}
 = e_{\nu \mu}, \quad 
 f_{\mu\nu } = - \pi \i \,  x_\mu x_\nu = f_{\nu \mu}, \quad h_{\mu\nu }
 = - x_\mu{\del\over\del x_\nu }
 - {1\over 2} \delta_{\mu\nu}, \quad 1\leq \mu, \nu \leq 2. 
 \ee
 We will call $\CC[x]$ with the action of $\gen = \osp(1|4)$ coming from
 $U(\gen)\to\Heis_2$ the {\em super-Weil representation}
 of $\gen$. The super-grading of $\CC[x]$ is given by the decomposition
 $\CC[x]=\CC[x]_\0 \oplus \CC[x]_\1$ into even and odd polynomials. 
 As a representation of $\gen$, the space $\CC[x]$ is irreducible. 
 
 \vskip .2cm

 Note that the same rules define an action of $\gen$ in $\CC^\oo(\RR^2)$,
 the space of smooth vectors of $\Wc$. 
 The formulas \eqref {eq:xxdd-2} for $\gen_\0$-action
 are compatible (by differentiation) with   the group
   action \eqref {eq:weil-expl-2}. 
 
 \vskip .2cm
 
 Further, as in \S  \ref {subsec:super-weil-1} \ref {par:grsweil-1}, we
 extend the $\ZZ/2$-graded  Weil representation of  $\GG_0 = \Mp(4)$ in 
 $\Wc = L_2(\RR^2) = \Wc_\0 \oplus \Wc_1$
 to a representation of the 
  Lie supergroup
 $\GG = \OMp(1|4)$, using the  above extension of the action of $\sp(4)$
 on the space of smooth vectors to an action of $\gen = \osp(1|4)$. 
 
 We similarly define the {\em Fourier transform automorphism}
 \be\label{eq:FT-2}
 \FT: \OMp(1|4) \lra \OMp(1|4), \quad g\mapsto \wt S g \wt S^{-1},
 \ee
 where $\wt S\in \Mp(4)$ is any preimage of 
 $S = \bpm0 & -1_2\\ 1_2 & 0\epm$. The induced automorphism of $\gen$,
 also denoted $\FT$, takes $p_i\mapsto q_i$ and $q_i \mapsto -p_i$.

 \paragraph{The $2$-term BGG resolution and the super-Laplacian.}
 Recall (\S \ref {subsec:SLG-gen2} \ref {par:osp-1-4}) the paraboliic subalgebra
 $\pen\subset \gen$ and its nilpotent radical $\nen = L_p \oplus S^2 L_p$.
 The Levi quotient/subalgebra of $\pen$ is $\men = \gl_2 = \End(\L_p)$.

 \vskip .2cm
 
 Let $\Oc$ be the parabolic Bernstein-Gelfand-Gelfand (BGG) category consisting of
 finitely generated $\gen$-(super)modules $V$  on which $\men$ acts semisimply
  and such that each $v\in V$ is $\nen$-finite, i.e., $\dim_\CC U(\nen) v <\oo$. 
  For example, the super-Weil representation $\CC[x]$ is an object of $\Oc$. 
  
  Further, for each $\alpha = (\alpha_1 \geq \alpha_2)\in  \Z2YD_2$
  we have the irreducible representation $\Sigma^\alpha$ of $\gl_2$ which we
  consider as a representation of $\pen$ via the projection to $\men$.
  The induced 
  {\em (super) Verma module} corresponding to $\alpha$ is defined as
$
  \SV_\alpha = \Ind_\pen^\gen \Sigma^\alpha. 
 $
 It is also an object of $\Oc$. For $\alpha  = (r,r)$, $r\in {1\over 2}\ZZ$,
 we write $\SV_r$ for $\SV_{(r,r)}$. 
 
 \vskip .2cm
 
 Eq. \eqref {eq:xxdd-2} implies that the vector $1\in \CC[x]$ is annihilated by $\nen$ and the action of any $z\in \men = \gl_2$
 on it is given by multiplication  with ${1\over 2} \tr(z)$. In other words, 
 $\CC\cdot 1\= \Sigma^{(-{1\over 2}, -{1\over 2})}$ as an $\men$-module. 
 Therefore we have a $\gen$-module homomorphism
 $\pi: \SV_{-{1/2}} \to \CC[x]$ which is surjective since $\CC[x]$ is  irreducible. 
 
 \begin{prop}
 The kernel of $\pi$ is isomorphic to $\SV_{-{3\over 2}}$, so that we have a
 $2$-term BGG resolution
 \[
 0\to V_{-{3/2}} \buildrel ^L\DD^\vee \over \lra V_{-{1/2}} \buildrel \pi\over\lra \CC[x] \to 0. 
 \]
 \end{prop}
 
 \noindent{\sl Proof:} By Proposition \ref {prop:conf-superlapl-gen2}
 we have an invariant differential operator given by the twisted super-Laplacian $^L\DD= (-1)^\F \DD: S\Lc_{1/2} \to
 S\Lc_{3/2}$. 
 In virtue of the correspondence between morphisms of
 Verma modules and invariant differential operators (Proposition 
 \ref {prop:harsish-twisted-n} from the 
 Appendix),  $^L\DD$ gives a $2$-term resolution
 \[
 0\to \Dc(S\Lc_{3/2}, S\Lc_{1/2}) 
 \buildrel \cdot ^L\DD\over \lra \Dc_{1/2} = \Dc(S\Lc_{1/2}, S\Lc_{1/2}) \lra 
 \Mc_\DD 
 \to  0
 \]
  of the $\Dc_{1/2}$-module $\Mc_\DD$ whose solution sheaf is 
  $\ul\Ker(\DD) = \ul\Ker(^L\DD)$.
  This latter fact, already mentioned in 
  \S  \ref {subsec-koszul-2} \ref {par:Koszul-D-gen2}, follows from the explicit
  form of $\DD$ in components reducing it to the Laplace and Dirac operators. 
  \qed. 
  
  \vskip .2cm
  
  \begin{rem}\label{rem:left-right-D}
  The necessity of composing $\DD$ with $(-1)^\F$ to get
  an invariant operator can be explained as follows. The super-Minkowski space
  $M^{3|2}_\CC = \SLG^\circ(E)$ can be identified with the supergroup
  $\NN^-$, the unipotent radical opposite to $\NN$. 
  Invariant differential operators on $\SLG(E)$ should give, in particular,
  left $\NN^-$-invariant operators on $\NN^-$. 
  
  Now, 
  the operators $D_i$ given by the infinitesimal action on $ \SLG(E)$ 
  of the odd generators
  $p_i\in \nen^-_\1$ of the Lie algebra of $\NN^-$ are then realized
  as the action  of $\nen$ on $\NN$ by  infinitesimal {\em left} translations   which
   are {\em right invariant} vector fields, but they are not left invariant,
  since $\NN$ is not (super) commutative. The super-Laplacian
   $\DD = [D_1, D_2]_-$ is then the 
  {\em right  $\NN^-$-invariant} operator  on $\NN^-$
  corresponding to $Q=[p_1, p_2]_-\in U(\nen^-)$. As $Q$ lies in the anti-center
  of $U(\nen^-)$, see \S \ref{subsec:super-Lap-gen2}\ref{par:D-anti}, 
  the left-invariant operator corresponding to $Q$ is $^L\DD = (-1)^F\DD$. 

   \end{rem}

 \paragraph{The (super-)Gaussian transform.}\label{par;sgauss-2}
  Similarly to 
 \eqref {eq:SG-1}, we define the {\em genus $2$ super-Gaussian transform}
 to be the map
 \be\label{eq:SG-2}
 \SG=\SG_2: \Wc \lra  \SF_{1/2}, \quad u(x) \mapsto \SG(u)(T, \xi) = \int_{x\in \RR^2} 
 u(x) \theta_x(T, \xi) dx. 
 \ee
 As before, it combines the classical Gaussian transform (integration against a variable Gaussian)
 capturing even functions and the integration against a variable Gaussian times a variable
 linear functions  which captures odd ones. Proposition \ref {prop:sGauss-2}
 implies the following.
 
 \begin{prop}\label{prop:SG=BW-2}
 $\SG$ is an injective $\FT$-twisted morphism of representations of $\OMp(1|4)$ taking values 
   in $\Ker(\DD)\subset \SF_{1/2}$.\qed
 \end{prop}
 
 Similarly to Remark \ref {rem:sgauss-BW-1}, the conceptual reason for this proposition
 is that $\theta_x(T,\xi)$, considered as a function of $x$, i.e., as an element of $\Wc$, 
   is a highest vector with respect to the parabolic subalgebra in $\gen$ corresponding
   (after applying $\FT$) to
 $(T,\xi)\in S\Hen_2\subset \SLG(\CC^4)$. 
 
 \paragraph{Multiplicativity of super-Weil representations and super-Gaussian transforms.}
 Considered as a $\ZZ/2$-graded Hilbert space, the space $\Wc = \Wc^{(2)} = L_2(\RR^2)$ of the genus $2$ super-Weil representation, is the topological tensor product $\Wc^{(1)} \wh\otimes \Wc^{(1)}$,
 where $\Wc^{(1)}=L_2(\RR))$ is the space of the genus $1$ super-Weil representation. 
 
 Recall \eqref {eq:emb-G-i-2} the emdeddings $\eps_i: \Gb^{(i)} \to \Gb = \OMp(1|4)$,
 where $\Gb^{(i)}$ is the copy of $\OMp(1|2)$ corresponding to the subspace
 $\CC p_i \oplus \CC  q_i$ which is equal to $\gen^{(i)}_\1$, the odd component
 of the (complexified) Lie algebra of $\Gb^{(i)}$.  As these odd components {\em commute}
 (in the ordinary, non-super sense) in the Weil representation, we obtain
 the following multiplicativity property. To formulate it, denote 
  $\varpi_1$ and $\varpi_2$ the Weil representations of $\OMp(1|2)$
 and $\OMp(1|4)$ respectively. 
 
 \begin{prop}
 With respect to the above identification $\Wc^{(2)} = \Wc^{(1)}\wh\otimes\Wc^{(1)}$,
 we have for any  $g\in \OMp(1|2)$ (i.e., $g$ being a point with values
 in some superalgebra): 
 \[
 \varpi_2(\eps_1(g)) = \varpi_1(g)\otimes 1, \quad \varpi_2(\eps_2(g)) = 
 (-1)^\F\otimes \varpi_1(g).
 \]\qed
 \end{prop}
\noindent  Here $(-1)^\F\otimes\varpi_1(g)$ sends $u_1\otimes u_2\mapsto u_1\otimes \varpi_1(g) (u_2)$  (naive, non-super  product rule). 

\vskip .2cm

Further, the diagonal factorization of super-Gaussians 
(Proposition \ref {prop:sgauss-fac-2}) implies the following multiplicativity
property of super-Gaussian transforms. Let use the notation $\SG_1, \SG_2$
for such transforms in genus $1$ and $2$ respectively.

\begin{prop}
For $u_1, u_2\in\Wc^{(1)} = L_2(\RR)$ we have
\[
\bigl((S\delta)^* \SG_2 (u_1\otimes u_2) \bigr) \bigl( (\tau_1, \xi_1), (\tau_2, \xi_2)
\bigr) \,=\, \SG_1(u_1) (\tau_1, \xi_1)\,  \SG_1(u_2) (\tau_2, \xi_2).
\]
\qed
\end{prop}

 \paragraph{The $\RR$-holonomic system in terms of the super-Weil
 representation.} 
 
 Similarly to  \S \ref {subsec:super-weil-1}
 \ref {par:Weil-compar-gen1}, we can interpret harmonic superforms
 $\Phi\in\SF_{1/2}^\DD$ as corresponding, via the super-Gaussian transform,
 to elements $u=u(x)$, $x\in \RR^2$ of some distributional completion of $\Wc$. 
 The new phenomenon here is the harmonicity condition $\DD\Phi=0$
 not present in genus $1$. 
 The $\RR$-holonomic system of Theorem \ref {theta:susy-2} on $\Phi$
 corresponds, under $\SG$, to the system on $u$
 \[
  e^{2\pi \, \i \, x_\nu} u(x) = u(x), \quad e^{\del/\del x_\nu} 
    u(x) \bigl( = u(x+\bee_\nu)\bigr) = u(x),\quad \nu=1,2,
 \]
 which gives $u(x)= c\cdot \sum_{\bn\in \ZZ^2} \delta(x-\bn)$ so
 $\SG(u) = c\cdot\theta (T).$. 
 
   \vfill\eject
   
   \section{Higher genus: generalized supercomformal symmetry}\label{sec:higher-gen}
   
   \subsection{($\Nc=1$) Lagrangian super-Grassmannian and super-Siegel plane}
   \label{subsec:SLG-sSieg}
   
    We now  consider   an arbitrary  genus $n\geq 1$.  
    The background material
    on various supergroups and their homogeneous spaces presented
    in this section, is rather elementary and was treated
    in detail for $n=2$ in \S \ref {subsec:SLG-gen2}. Therefore
    we omit most of the proofs. 
   
    \paragraph{Generalized (super) conformal structure on the (super)  Lagrangian Grassmannian.} 
  
  Let $E=(\CC^{2n}, \omega)$ be a symplectic vector space and $\LG = \LG(E)$
  the Grassmannian of Lagrangian subspaces in $E$.  It carries 
  a differential-geometric structure  which is analogous (but not identical,  for $n>2$),  to a conformal structure.  Such structures for general Hermitian symmetric spaces
  were studied in  \cite{gonch-symm, gonch-selecta}. 
  
  \vskip .2cm
  
  To describe this structure in more detail, denote by $\Lambda $
  the tautological rank $n$ bundle on $\LG$.  It is embedded into the
  trivial bundle $\wt E = E\otimes\Oc_\LG$. 
  We have then the natural identification
   \cite[Prop.2.1] {guillemin}
  \be\label{eq:TLG-gen}
  T\LG \,\= \, S^2 \Lambda^*. 
  \ee
  In other words, let $L\subset E$ be a Lagrangian subspace and $[L]\in\LG$
  be the corresponding point. 
  By \eqref {eq:TLG-gen},  tangent vectors to $\LG$ at $[L]$  can be seen
  as quadratic forms on $L$. 
 Therefore, for $1\leq r\leq n-1$ each space $T_{[L]} \LG$ carries
  a  conic subvariety (cone, for short)
 \[
 T_{[L]}^{\leq r} \LG \subset T_{[L]}\LG = S^2 L^*
 \]
 formed by quadratic forms  on $L$ on rank $\leq r$. The system of these cones 
 is the analog of a conformal structure. For $n=2$ we have only one
 cone, that of  degenerate quadratic forms which is given by one quadratic
 equation ($\det(A)=0$ for a symmetric $2\times 2$ matrix $A$), so
 we get a conformal structure in the usual  sense. 
 
\vskip .2cm

 Consider also  the supervector space $Y =  \CC^{0|1}\oplus E$
 of super-dimension $2n|1$. It has a symplectic form $\eta$ composed on $\omega$ on $E$
 and the standard symmetric form on $\CC$:
  \be\label{eq:eta-gen-n}
   \eta\bigl( \zeta, v), ( \zeta', v')\bigr) \,=\,  \zeta\cdot\zeta' + \omega(v,v'),  
   \quad x, x'\in E, \,\,\, \zeta, \zeta'\in \CC. 
   \ee
    The ($\Nc=1$)  {\em Lagrangian super-Grassmannian}
 of $E$ is  the supermanifold  $\SLG=\SLG(E)$ parametrizing isotropic subspaces
 $L\subset Y$ of purely even dimension $n|0$. It carries the tautological bundle
 $S\Lambda$ of rank $n|0$ embedded into the trivial
 bundle $\wt Y = Y\otimes\Oc_\SLG$. 
 
 \begin{prop}
 (a) $\SLG$ is a smooth complex supermanifold of dimension 
 $({n(n+1)\over 2}|n)$, whose underlying even manifold is $\SLG_\red = \LG$. 
 
 \vskip .2cm
 
 (b) The tangent bundle of $\SLG$ is identified as
 \[
 T\SLG \,\= \, \Hom(S\Lambda, \wt Y/S\Lambda) \times_{\Hom(S\Lambda, (S\Lambda)^*)} \Hom^+(S\Lambda, S\Lambda^*),
 \]
 where 
  $\Hom^+(S\Lambda, (S\Lambda)^*)\= S^2((S\Lambda)^*)$
 is the space of self-adjoint maps $S\Lambda\to S\Lambda^*$.
 
 \vskip .2cm
 
 (c) The subbundle $\Tc = \Hom(S\Lambda, (S\Lambda^)\perp/S\Lambda)\subset
 T\SLG$ of rank $(0|n)$ is a SUSY structure and    
 $\Pi\Tc|_\LG = S^2\Lambda^*$. The restriction of the Frobenius pairing
 of $\Tc$ to $\LG$ is the isomorphism \eqref {eq:TLG-gen}. \qed
 \end{prop}
 
 \paragraph{The (ortho) symplectic group and its Lie algebra.}
 We fix a symplectic basis $\{p_i, q_i\}_{i=1}^n$ in $E$ so that 
 $\omega(p_i, q_i)=1$.  The Gram matrices of $\omega$ and $\eta$ 
 with respect to this basis (and the obvious basis of $\CC^{0|1}$)
 are
  \[
  \Omega = \bpm 0&-1_n\\ 1_n&0
  \epm,
  \quad
  H = \begin{pmatrix} 1&0&0
 \\
 0&0&-1_n
 \\ 0&1_n&0
 \end{pmatrix}.  
 \]
 The symplectic Lie group $\Sp(2n) = \Sp(E)$ and the orthosymplectic 
 Lie supergroup $\OSp(1|2n) = \Sp(Y)$ are the groups of symmetries
 of $(E,\omega)$ and $(Y, \eta)$,
 acting on $\LG$ and $\SLG$ respectively.  Their points $g$
 are represented by (super) matrices
 (Latin variables are even and Greek ones are odd)
 \be\label{eq:OSP(1|2n)-matr}
 \begin{gathered}
 \Sp(2n) \,\, \ni \,\, g = \bpm A&B \\C&D\epm,
 \quad  A,B,C,D \in\Mat_{n\times n},
  \quad g^t\Omega g= 1_{2n},
 \\
 \OSp(1|2n) \,\,\ni \,\, g =
 \begin{pmatrix}
 s&\alpha&\beta 
 \\
 \gamma& A & B
 \\
 \delta & C & D
 \end{pmatrix}, \quad 
 g^T Hg = 1_{2n+1}, \quad 
 \begin{matrix} s\in\Mat_{1\times 1}, & A,B,C,D \in\Mat_{n\times n},
 \\
 \quad \alpha,\beta\in \Mat_{1\times n}, & \gamma,\delta\in\Mat_{n\times 1}.  
 \end{matrix}
\end{gathered}
 \ee
 Here $g^T$ is the super-transpose \eqref {eq:super-trans}. Note that
 \[
 \OSp(1|2n)_\red = \O(1) \times \Sp(2n) =  \{\pm 1\}\times \Sp(2n).
 \] 
% as  the conditions imply that $s_\red^2=1$. 
 We denote by $\gen = \osp(1|2n)$ the Lie superalgebra of $\OSp(1|2n)$.
 It has $\gen_\0 = \sp(2n)$, the Lie algebra of $\Sp(2n)$, while
 $\gen_\1 = E=\CC^{2n}$, with the adjoint action of $\gen_\0$ on $\gen_\1$
 being the standard action of $\sp(2n)$ on $E$. 
 Further, as a $\gen_\0$-module, $\gen_\0\= S^2 E$ and the bracket
 $\gen_\1\otimes \gen_\1\to \gen_\0$ is the symmetrization
 $E\otimes E\to S^2E$.

  We have the infinitesimal action of $\gen_\0$  on $\LG$ 
  resp.  of $\gen$ on $\SLG$ by vector fields
 $D_y$ with $y\in \gen_\0$ resp. $y\in \gen$. 
 
 \begin{prop}
 The  subbundle $\Tc\subset T\SLG$ defining the SUSY structure
  is  spanned by the odd
 vector fields $D_y, y\in \gen_\1$. 
 \qed
 \end{prop}
 
 \paragraph{The open cells in $\LG$ and $\SLG$.}\label{par:cells-gen-n}
 The Grassmannian $\LG$ 
 has  two distinguished points
 \[
 L_q=\bigoplus_{i=1}^n\, \CC q_i , \,\,\, L_p= \bigoplus_{i=1}^n\, \CC p_i \,\,\, \subset \,\,\,   E\subset Y
 = E\oplus \CC^{0|1}.
 \]
 We define the open Schubert cells $\LG^\circ\subset\LG$ resp. 
 $\SLG^\circ\subset\SLG$ to consist of isotropic subspaces $L\subset E$
 resp. $L\subset Y$ which are transverse to $L_p$. Denoting by 
 $x_i, y_i$ and $\zeta$ the coordinates in $Y$ dual to the basis vectors
 $q_i, p_i$ and  $1\in \CC^{0|1}$, points of $\LG^\circ$ resp. $\SLG^\circ$
 have the form
 \[
 \begin{gathered}
 L_T = \bigl\{ (x,y)^t: \bigl|\,  y_i = \sum_{j=1}^n  t_{ij} x_j\bigr\}, \quad \text{resp.} 
 \\
 L_{T,\xi} = \bigl\{ (\zeta, x,y)^t \bigl| \, y_i = \sum_{j=1}^n t_{ij} x_j, \,\,\zeta = 
 \xi x = \sum_{j=1}^n \xi_i x_i\bigr\},
 \end{gathered}
 \]
where $T=\|t_{ij}\|\in \Mat_{n\times n}(\CC)$ and $\xi = (\xi_1,\cdots, \xi_n)$
is a row vector of odd variables. The condition for $L_T$ to be Lagrangian
is that $T$ is symmetric: $t_{ij}=t_{ji}$. The condition for $L_{T,\xi}$
to be isotropic is that
\be\label{eq:rels-T=xi-n}
\sigma_{ij} (T, \xi) \, := \, t_{ij}-t_{ji}-\xi_i \xi_j =0,\quad  1\leq i,j \leq n. 
\ee
Thus $\LG^\circ$ is identified with 
$\Sym_{n\times n}(\CC)\= \CC^{n(n+1)\over2}$
while $\SLG^\circ$ is given in $\CC^{n^2 |n}$ by the equations
 $\sigma_{ij}(T,\xi)=0$.  In other words,  $(T, \xi)$ are dependent coordinates on
 $\SLG^\circ$. 
 One can identify $\SLG^\circ$
 with $\CC^{{n(n+1)\over 2} | n}$ by using idependent coordinates
 \be
 v= \|v_{ji}\|, \,\,\, v_{ij} = {1\over 2}( t_{ij}+ t_{ji}) = v_{ji}, \quad \xi = 
 (\xi_1,\cdots,\xi_n). 
 \ee
 Let  a point $g\in \Sp(2n)$  resp.  $g\in\OSp(1|2n)$ be given in the form 
 \eqref{eq:OSP(1|2n)-matr}
  The action of $g$ on $\LG^\circ$ resp.  
 on $\SLG^\circ$ is given, in our coordinates,  by the standard Siegel-type birational  formulas,
 compare \eqref{eq:SSiegel-bir-2}:
 \be\label {eq:SSiegel-bir-n}
 \begin{gathered}
 g(T) = (AT+B) (CT+D)^{-1} \quad \text {resp.}
 \\
 g(T,\xi) = \bigl( (\gamma\xi + AT+B)(\delta\xi +CT+D)^{-1}, \,
(s\xi + \alpha T+\beta)(\delta\xi+CT+D)^{-1}\bigr). 
 \end{gathered} 
 \ee
  We define the ``spacetime derivatives"
 \be
 \del_{ij} \,=\,{\del\over \del t_{ij}} + {\del\over \del t_{ji}} \,=\,
 {\del\over \del v_{ij}} + {\del\over \del v_{ji} }\, =\,  \del_{ji} , 
 \ee
 which can be understood in $3$ compatible ways:
 \begin{itemize}
 \item As vector fields on $\LG^\circ = \Sym_{n\times n}(\CC)$.
 
 \item As vector fields on $\CC^{n^2|n}$ tangent to $\SLG^\circ$ 
 (annihilating the equations $\sigma_{ij}$)
 and so giving vector fields on $\SLG^\circ$.
 
 \item As vector fields on $\SLG^\circ$ given in independent
 coordinates $(v,\xi)$. 
 \end{itemize}
 We further define the ``spinor derivatives''
 \be\label{eq:spin-der-n}
 D_i \,=\,{\del\over \del \xi_i} - \sum_{j=1}^n \xi_j {\del\over \del t_{ij}},
 \quad i=1,\cdots, n
 \ee
 (vector fields on $\CC^{n^2|n}$  tangent to $\SLG^\circ$)
 which in the independent coordinates $(v,\xi)$ have the form
 \[
 D_i  \,=\,{\del\over \del \xi_i} - \sum_{j=1}^n \xi_j \del_{ij}. 
 \]
 They satisfy
 \[
 [D_i, D_j]_+ = -2\del_{ij}, \quad [D_i, \del_{jk}]_- = 0. 
 \]
 We also define the derivations
 \be
 D'_i = \sum_{j=1}^n  t_{ji} D_i, \quad i=1,\cdots, n. 
 \ee
 It is convenient to use the vector and matrix notation
 \[
 \begin{gathered}
 \xi = (\xi_1, \cdots, \xi_n), \quad {\del\over\del\xi} \,=\,
 \biggl({\del\over\del\xi_1}, \cdots, {\del\over\del\xi_n}\biggr), \quad
 {\del\over \del T} = \left\| {\del\over\del t_{ij}}\right\|, 
 \\
 D = (D_1, \cdots, D_n) = {\del\over\del\xi} - \xi\cdot {\del\over\del T},\quad
 D'= (D'_1, \cdots, D'_n)  = D\cdot T. 
 \end{gathered}
 \]
 Differentiating \eqref {eq:SSiegel-bir-n}, we get:
 \begin{prop}
 The vector fields $D_y, y\in\gen_\1$ on $\SLG$ coming from the infinitesimal
 action of $\gen$, have, in the restriction to $\SLG^\circ$,  the form
 $D_{p_i} = D_i, D_{q_i}=D'_i$. \qed
 \end{prop}
 We denote by  $\PP = \Stab(L_p)\subset \OSp(1|2n)$ the parabolic subgroup stabilizing $L_p$,
 by $\NN$ the unipotent radical of $\PP$ and by $\MM$ the Levi subgroup/quotient. 
 Thus $\MM= \GL(n) = \Aut(\L_p)$ is purely even. Let $\pen, \nen, \men$ be the Lie superalgebras
 of the supergroups $\PP, \NN, \MM$. With respect to the identification $\gen = E \oplus S^2E$, we have
 \be\label{eq:m-n-p-gen}
 \begin{gathered}
 \pen = L_p \oplus (L_p\cdot E), \quad \nen  = L_p \oplus S^2 L_p, \quad \men = L_p \cdot L_q\subset L_p\cdot E, \quad \men \= \End(L_p),
 \\
 \pen_\0 = L_p\cdot E, \quad \nen_\0 = S^2 L_p, \quad \men_\0=\men. 
 \end{gathered} 
 \ee
 
 \paragraph{The (super) Siegel plane.} \label{par;super-siegel-n}
 We denote by $\Hen_n\subset \LG^\circ$ the genus $n$ Siegel plane
 consisting of $T\in\Sym_{n\times n}(\CC)$ with $\Im(T)>0$. Let $S\Hen_n\subset \SLG^\circ$ be the
 super-thickening of $\Hen_n$ induced from the  thickening  $\LG^\circ$ given by $\SLG^\circ$. 
 The action 
 \eqref{eq:SSiegel-bir-n} restricted to the group $\Sp(2m,\RR)$, preserves $\Hen_n$ and restricted
 to the group $\OSp(1|2n, \RR)$, preserves $S\Hen_n$. This gives identifications
 \[
 \begin{gathered}
 \SLG = \OSp(1|2n)/\PP, \quad S\Hen_n = \OSp(1|2n, \RR)/\bigl( \{\pm 1\}
 \times  U(n)\bigr),
 \\
 \LG = \Sp(2n)/P, \quad P=\PP_\red^{(e)}, \quad \Hen_n = \Sp(2n,\RR)/U(n),
 \end{gathered}
 \]
 where $P=\PP_\red^{(e)}$ i is the stabilizer of $L_p$ is $\Sp(2n)$ which is the unit  connected component
 of the $2$-component group $\PP_\red$. 
 
 \paragraph{Young diagrams and partitions.  Forms and superforms. 
 }\label{par:young-diag-part}
 We will refer to dominant integer weights for $\GL_n$,
 i.e., to  non-increasing sequences $\alpha = (\alpha_1\geq\cdots\geq\alpha_n)\in\ZZ^n$ 
 as {\em $\ZZ$-Young diagrams with $n$ rows} and denote the set formed by them 
 by $\ZYD_n$. 
  An $\alpha\in\ZYD_n$ defines an  irreducible representation $\rho_\alpha$ of $\GL(n)$ or, what is
  the same, the {\em Schur functor}  $V\mapsto \Sigma^\alpha V$ on $n$-dimensional $\CC$-vector spaces. Note that
  \be\label{eq:def-alpha-vee}
 ( \Sigma^\alpha V)^* \,\= \, \Sigma^\alpha (V^*) \, \= \,
 \Sigma^{\alpha^\vee} (V), \quad \text{where} \quad 
 \alpha^\vee = (-\alpha_n, \cdots, -\alpha_1). 
  \ee
 We denote the space of representation $\rho_\alpha$ as $\Sigma^\alpha = \Sigma^\alpha \CC^n$.

 \vskip .2cm
 
 As usual, non-negative 
 $\alpha = (\alpha_1 \geq\cdots, \geq \alpha_n\geq 0)$
 are called Young diagrams or {\em partitions}.  For such $\alpha$ 
 we omit eventual zeroes at the end writing, e.g.,
 $(1,1)$ for $(1,1, 0\cdots, 0)$ so that $\Sigma^{1,1}(\CC^n)=\Lambda^2(\CC^n)$. 
 We also  write
 $|\alpha|=\sum \alpha_i$ and denote by $\alpha^*$ the
 dual Young diagrams whose rows correspond to the columns of
 $\alpha$.
 
 \vskip .2cm
 
 We will also use the {\em Frobenius notation} for partitions, writing
 $\alpha = (a_1,\cdots, a_s |b_1,\cdots, b_s)$ 
 with $a_1>\cdots > a_s\geq 0$ and $b_1>\cdots > b_s\geq 0$. 
 Here $a_i$ is the number of boxes in the $i$th {\em arm} of $\alpha$,
 i.e., in the segment strictly  to the right of the $i$th diagonal cell or,
 in other words, $a_i = \alpha_i-i$ as long as $\alpha_i-i\geq 0$.
 Similarly, $b_i$ is the number of boxes in the $i$th {\em leg}
 of $\alpha$, i.e., in the segment strictly below the $i$th
 diagonal cell or, in other words, $b_i = \alpha^*_i-i$ as long
 as $\alpha^*_i-i\geq 0$.  Thus,
 \[
 (a_1,\cdots, a_s |b_1,\cdots, b_s)^* = (b_1,\cdots, b_s | a_1,\cdots, a_s).
 \]
 For $\alpha\in \ZYD_n$ we denote
 \be
 \Lc_\alpha = \Sigma^\alpha \Lambda, \quad S\Lc_\alpha = \Sigma^\alpha S\Lambda
 \ee
 the homogeneous vector bundles on $\LG$ and $\SLG$
 respectively obtained by applying $\Sigma^\alpha$
 to the tautological bundles $\Lambda$ and $S\Lambda$ on 
 these (super) manifolds. In particularly, we
  consider the spaces
 \[
 \F_\alpha = H^0(\Hen_n, \Lc_\alpha), \quad \SF_\alpha = 
 H^0(S\Hen_n, S\Lc_\alpha)
 \]
 of  sections of the above equivariant bundles over the (super) Siegel plane
 and call their elements  {\em (super) Siegel forms of weight}  
 $\alpha\in \ZYD_n$. Thus $\Sp(2n, \RR)$ acts on $\F_\alpha$ and
 $\OSp(1|2n, \RR)$ acts on $\SF_\alpha$. 
 
 \vskip .2cm 
 
  Explicitly, as  vector bundles, $\Lc_\alpha$ and $S\Lc_\alpha$
 are trivial with fiber $\Sigma^\alpha$, so their sections
 are represented by  $\Sigma^\alpha$-valued functions
  and the
  action of the corresponding group in the sections are given by:
  \be\label{eq:action-on-sforms-n}
  \begin{gathered}
  \Lc_\alpha \,\ni \, \phi: \Hen_n\to\Sigma_\alpha, \quad 
   g^* (\phi)(T) = \rho_\alpha^{-1}(CT+D)( \phi(g(T))). 
   \\
S\Lc_\alpha\,\ni \, \Phi: S\Hen_n\to \Sigma^\alpha, \quad 
(g^*\Phi)(T,\xi) = \rho_\alpha^{-1} (\delta\xi+CT+D) \Phi(g(T,\xi)). 
  \end{gathered} 
  \ee

\paragraph{Half-integer  weights.}

 Irreducible representations of the Lie algebra $\gl_n$ are labelled by 
  {\em $\CC$-Young diagrams with $n$ rows}, i.e., by sequences $\alpha = (\alpha_1,
 \cdots\alpha_n)$ of complex numbers with  each $\alpha_i-\alpha_{i+1}$ being  a nonnegative integer. 
 We denote the set of such $\alpha$ by $\CYD_n$.  For $r\in \CC$
 we denote $r^n := (r,\cdots, r)\in \CYD_n$
  (with $r$ repeated $n$ times). 
  
  \vskip .2cm
 
 We further denote $\Z2YD_n$
 the subset of $\CYD_n$ consisting of $\alpha$ with $\alpha_i\in {1\over 2} \ZZ$ and call its
 elements {\em ${1\over 2}\ZZ$-Young diagrams  with $n$ rows}. 
 Let $\wt M = \wt \GL(n)\to \GL_n$ be  the double cover consisting of pairs $(g, \phi)$ with $g\in\GL_n$, $\phi\in \CC^*$
 such that $\phi^2=\det(g)$, cf. \eqref{eq:wt-GL2} for $n=2$. 
 We also denote by $\wt U(n)$ the preimage of
   $U(n)$ in $\wt \GL(n)$. 
 Elements $\alpha\in \Z2YD_n$ 
  parametrize 
 irreducible representations $\rho_\alpha: \wt M  \to \Aut(\Sigma^\alpha)$
   We denote
 \[
 \Z2YD{}_n^- \,\,=\,\,\Z2YD{}_n \,\setminus \ZYD_n \,\, = \,\, \bigl({1\over 2}\bigr)^n 
 + \ZYD_n
 \]
 the set of strictly half-integer weights. 
 
 \vskip .2cm
 
 The Lagrangian Grassmannian $\LG=\LG(\CC^{2n})$ carries a
$G=\Sp(2n)$-equivariant $\CC^*$-gerbe $\Gc_{1/2}$ formed by determinations
of  square roots of the line bundle $\Lc_1 = \Oc_{LG}(-1)$, see Example \ref {ex:frac-powers}(b). This square roots becomes  
 a canonically defined
  $\Gc_{1/2}$-twisted line bundle $\Lc_{1/2}$ on $\LG$. Alternatively,
representing $\LG$ as $G/P$, $P=\Stab(L_p)$, we have the central extension
\[
1\to H=\{\pm 1\} \lra \wt P \lra P\to 1
\]
coming from the above  $H$-extension $\wt M=\wt \GL(n)$ of the Levi subgroup
$M=\GL(n)$,  see \S
\ref {subsec:SLG-sSieg} \ref {par:young-diag-part}. Denoting by $\sigma: H\to\CC^*$
the nontrivial character, we get the $\CC^*$-gerbe $\Gc_\sigma$ on
$\LG$ described in Example  \ref {exas:ourTDO}(a) and which is identified
with $\Gc_{1/2}$. 

 \vskip .2cm
 
 As on earlier occasions, for any $\alpha\in \Z2YD_n^-$ we have 
 $\Gc_{1/2}$-twisted vector bundles $\Lc_\alpha$ on $\LG(E)$
 and $S\Lc_\alpha$ on $\SLG(E)$. A useful alternative
 way to define $S\Lc_\alpha$ is to 
 use the  central 
extensrion
\[
1\to H=\{\pm 1\} \lra \wt \PP \lra\PP \to 1, \quad \PP = \Stab(L_p)\subset
\GG = \OSp(1|2n),  
\]
extending $\wt P$ and the representation  $\rho_\alpha$ of $\wt \PP$
extending that of $\wt P$. 

\vskip .2cm

We denote by $\Mp(2n)\to \Sp(2n,\RR)$ resp. 
$\OMp(1|2n)\to \OSp(1|2n, \RR)$ the $2$-fold covering
whose points are pairs $\wt g= (g,\phi)$ where $g$ is a point
of $\Sp(2n,\RR)$ resp. $\OSp(1|2n,\RR)$ as in 
\eqref {eq:OSP(1|2n)-matr} and
\[
\begin{gathered}
\phi = \phi(T) \in \Oc(\Hen_n)\text{ is such that } \phi(T)^2 = 
\det(CT+D), \quad \text{resp.}
\\
\phi = \phi(T,\xi) \in \Oc(S\Hen_n) \text{ is such that } \phi(T)^2 = 
\det(\delta\xi + CT+D).
\end{gathered}
\]
Note that a choice of $\wt g$ lifting $g$ defines a lift of the
$\GL(n)$-valued function
\[
\begin{gathered}
T\,\mapsto \,\, CT+D, \quad \Hen_n \lra \GL(n), \quad \text{resp.}
\\
(T,\xi) \,\, \mapsto \,\, \delta\xi + CT+D, \quad S\Hen_n \lra \GL(n)
\end{gathered} 
\]
to a $\wt \GL(n)$-valued function which we denote $(CT+D)\, \wt{}$
resp. $(\delta\xi + CT+D)\,\wt{}$. 
We call $\Mp(2n)$ and $\OMp(1|2n)$ the {\em metaplectic}
and the {\em ortho-metaplectic} group of genus $n$ respectively. 

%%%BELOW COMMENTED OUT
\iffalse
Note that
\[
\OMp(1|2n)_\red \,=\, \O(1)\times \Mp(2n) = \{\pm 1\} \times\Mp(2n).
\]
Note further  that the preimage of $U(n)\subset \Sp(2n,\RR)$ in 
$\Mp(2n)$
is $\wt U(n)$ and the preimage of $\{\pm 1\}\times U(n)$ in
$\OMp(1|2n)$ is $\{\pm 1\}\times\wt U(n)$. 
\fi
%%%ABOVE COMMENTED OUT

\vskip .2cm

After restricting to $\Hen_n\subset \LG$ the gerbe $\Gc_{1/2}$ is canonically
trivialized and each
 $\Gc_{1/2}$-twisted bundle $\Lc_\alpha$, $\alpha\in \Z2YD_n^-$,
becomes an honest $\Mp(2n)$-equivariant vector bundle on $\Hen_\alpha$
and similarly $S\Lc_\alpha$ becomes an $\OMp(1|2n)$-equivariant
vector bundle on $S\Hen_n$. Extending the previous notation,
we denote by $\F_\alpha$ resp. $\SF_\alpha$ the space of  its global sections. 
 
 We can describe $\F_\alpha$ resp.   $\SF_\alpha$ for any $\alpha\in \Z2YD^-_n$
  as the space of $\Sigma^\alpha$-valued functions on $\Hen_n$
 resp. $S\Hen_n$ with the action of $\Mp(2n)$ resp. $\OMp(1|2n)$
  given by:
  \be\label{eq:action-on-sforms-n-half}
  \begin{gathered}
  \F_\alpha \,\ni \, \phi: \Hen_n\lra\Sigma_\alpha, \quad 
   \wt g^* (\phi)(T) = \rho_\alpha^{-1}((CT+D)\,\wt{}\, )( \phi(g(T))). 
   \\
\SF_\alpha\,\ni \, \Phi: S\Hen_n\lra  \Sigma^\alpha, \quad 
(\wt g^*\Phi)(T,\xi) = \rho_\alpha^{-1} ((\delta\xi+CT+D)\, \wt{}\, ) \Phi(g(T,\xi)). 
  \end{gathered} 
  \ee

 \vskip .2cm
 
 For $y\in \gen$ we denote by $D_y^{(\alpha)}$ the
 infinitesimal action of $y$ in $\SF_\alpha$ consider as a
 differential operator in the space of $\Sigma^\alpha$-valued
 functions on $S\Hen_n$. For $r\in {1\over 2}\ZZ$
 we abreviate $D_y^{(r^n)}$ to $D_y^{(r)}$.

 Differentiating \eqref  {eq:action-on-sforms-n}
 \eqref {eq:action-on-sforms-n-half}, we find:
 
 \begin{prop}
 For  $r\in {1\over 2}\ZZ$ we have 
 \[
 D^{(r)}_{p_i} = D_i, \quad D^{(r)}_{q_i} = D'_i - r \xi_i.  \qed
  \]
 \end{prop}

 \vfill\eject
 
 \subsection {The super-Weil representation and its BGG resolution}
 \label{subsec:sweil-BGG}
 
 \paragraph{ The (super-)Weil at the group and Lie algebra levels.}
 \label{par:sweil-n} 
 In the same way as discussed for $n=1,2$ in \S\S  {}  {}
 \ref{subsec:super-weil-1}\ref{par:weil-group-1} and \ref {subsec;super-weil-2} \ref{par:weil-group-2}, we have the Weil representation
 $\varpi$ of
 $\Mp(2n)$ in the space $\Wc =L_2(\RR^n)$ whose elements
 will be denoted $u(x)$, $x=(x_1,\cdots, x_n)^t$. As a projective
 representation of $\Sp(2n, \RR)$, $\varpi$ is uniquely defined
 by the values
 \be\label{eq:weil-expl-n}
 \varpi \left( \bpm 1_n&0 \\ C &1_n \epm\right)  u \,=\, e^{-\pi  \i \, x^t C x} \cdot u,\quad \text{while} \quad 
 \varpi  \left(\bpm 0& -1_n\\ 1_n&0 \epm\right) u  \,=\, \wh u, 
 \ee
 where  $C\in \Sym_{n\times n}(\RR)$ and $\wh u$ is the Fourier transform
 of $u$. The representation $\Wc$ splits into the sum $\Wc_\0\oplus
 \Wc_\1$ of the spaces of even and odd functions, each of
 which is an irreducible representation of $\Mp(2n)$.
 We extend the action of $\Sp(2n,\RR)$ to  that of 
 $\OSp(1|2n, \RR)_\red = \{\pm 1\}\times \Sp(2n,\RR)$ by making $-1\in\{\pm 1\}$  act by
 the parity operator, i.e., to send $f(x)\mapsto f(-x)$. 
 
 \vskip .2cm
 
 At the level of Lie algebras, $\varpi$ differentiates to a representation
 of $\gen_\0=\sp(2n)$ in $C^\oo(\RR^n)$ 
 formed by the operators
 \be\label{eq:xxdd-n}
 e_{\mu\nu } ={1\over 4\pi \, \i} \cdot  {\del^2\over\del x_\mu \del x_\nu}
 = e_{\nu \mu}, \quad 
 f_{\mu\nu } = - \pi \i \,  x_\mu x_\nu = f_{\nu \mu}, \quad h_{\mu\nu }
 = - x_\mu{\del\over\del x_\nu }
 - {1\over 2} \delta_{\mu\nu}, \quad 1\leq \mu, \nu \leq n
 \ee
which, again, preserve the decomposition into even and odd functions.
This action extends to a super-representation of $\gen = \osp(1|2n)$
by defining the action of $p_\nu, q_\nu\in \gen_\1 = E$ by
 \be\label{eq:weil-odd-n}
 p_\nu \mapsto {1\over\sqrt{-4\pi \i }}\, \,  {\del \over \del x_\nu}, \quad q_\nu \mapsto \sqrt{-\pi \i } \,\, x_\nu, \quad \nu=1,\cdots, n. 
 \ee
 Combining the action of $\OSp(1|2n,\RR)_\red$ and the compatible
 action of $\gen$ on smooth vectors as explained in
 \S \ref{subsec:super-weil-1} \ref{par:grsweil-1},
 we make $\Wc$ into 
  an irreducible representation of the Lie supergroup $\OSp(1|2n, \RR)$.
  
  \vskip .2cm
  
   Let $\Heis_n = \Heis(E, {1\over 2} \omega)$ be the Heisenberg algebra
 associated to the symplectic space $E=\CC^{2n}$ with
  the form ${1\over 2} \omega$.
 So it is generated by $y\in E$ subject to the relations
 $
 yz -zy = {1\over 2} \omega(y,z) \cdot 1
 $.  
 
 \begin{prop}\label{prop:sweil-heis-n}
 The formulas \eqref{eq:weil-odd-n} \eqref{eq:xxdd-n}
   define a  (surjective) homomorphism 
 of associative algebras $U(\gen) \to \Heis_n$.  \qed
 \end{prop}

 %\vfill\eject
  
\paragraph{ The parabolic categories $\Oc$ for Lie superalgebras.} 
\label{par:parab-cat-O}
Recall \cite{serganova-qreductive} that a complex algebraic
supergroup $G$ is called {\em quasi-reductive}, if the underlying
ordinary algebraic group $G_\red$ is reductive. 
It it known \cite{fioresi-shu}  that such $G$ are given by root data, in particular,
we can speak about parabolic subgroups in $G$, their Levi
quotients etc. 

\vskip .2cm

Let $G$ be a reductive algebraic Lie supergroup, $P\subset G$
 a parabolic subgroup
with Levi quotient $P\to M$.  Denote by $\bg, \bp, \bm$ the Lie superalgebras of
$G, P, M $. Recall that we can  lift $M$ to a subgroup of $G$ and $\bm$
to a Lie subalgebra in $\bg$. Let also $N$ be the unipotent radical of $P$ and
$\bn = \Lie(N)$ the nilpotent radical of $\bp$, so 
$\bp = \bm \oplus\bn$. 

\vskip .2cm

We denote by $\Oc^{\bg,\bp}$ the parabolic category 
$\Oc$ associated to $(\bg,\bp)$.
Explicitly, see, e.g. \cite[Ch.9]{humphreys}  it consists of $\bg$-modules $V$ which are:
\begin{itemize}
\item Finitely generated over $U(\bg)$.

\item Considered as a $U(\bm)$-module, split into a direct sum of finite-dimensional
simple modules.

\item Locally $U(\bn$)-finite, i.e., the $\bn$-module $U(\bn)v$ is finite-dimensional
for any $v\in V$. 
\end{itemize} 
It is known that $\Oc^{\gen,\pen}$ is a Noetherian and Artinian abelian category with
enough injectives and projectives. 

\vskip .2cm

Given a finite-dimensional representation $\rho: \bm\to\gl(A)$
we denote
\be\label{eq:Verma} 
\V^{\bg, \bp}_\rho = \Ind_\bp^\bg\, A =  U(\bg) \otimes_{U(\bp)}  A
\ee
the (induced, parabolic) {\em Verma module} corresponding to $\rho$. It is an object of
$\Oc^{\bg,\bp}$.

 \begin{exas}\label{exas:parab-O-n}
 (a) We will be initially interested in two cases:
 \begin{itemize}
 \item [(a1)]  The super-case when $\bg = \gen = \osp(1|2n)$ and $\bp =\pen$ is the stabilizer of
 $L_p$, so that $\bn = \nen$, $\bm=\men = \gl(n)$,
  see \eqref {eq:m-n-p-gen}.
  
  \item [(a2)]  The even case $\bg = \gen_\0=\sp(2n)$ while
  $\bp=\pen_\0$, $\bn=\nen_\0$, $\bm = \men = \men_\0$. 
 \end{itemize}
 
 In these cases any  $\alpha \in \CYD_n$  gives an irreducible
  representation
 $\Sigma^\alpha$ of $\gl(n) = \men = \men_\0$.  
 and we denote the corresponding 
  (super) Verma modules by 
 \[
 \V_\alpha = \V_\alpha^{\gen_\0, \pen_\0} = \Ind_{\pen_\0}^{\gen_\0} \, \Sigma^\alpha \,\in 
 \Oc^{\gen_\0, \pen_\0}, \quad
  \SV_\alpha = \V_\alpha^{\gen, \pen} = 
   \Ind_\pen^\gen\,  \Sigma^\alpha \,\in \,
   \Oc^{\gen, \pen}. 
 \]
(b)  The {\em polynomial super-Weil representation } 
 $\CC[x] = \CC[x]_\0 \oplus \CC[x]_\1$
 with the $\gen$-action given by \eqref{eq:xxdd-n}
 \eqref{eq:weil-odd-n} is an object of $\Oc^{\gen, \pen}$. 
 Each of the summands $\CC[x]_\0$, $\CC[x]_\1$ is an object of
 $\Oc^{\gen_\0, \pen_\0}$. 
    \end{exas}
 
 %\vfill\eject
 
 \paragraph{ The BGG resolutions: statement of results.}
Let $m\geq 0$  and $d$ be integers. Following Weyman \cite{weyman-book},
we denote by $Q_d(m)$ the set of partitions $\alpha$ of $m$
such that in the Frobenius notation
 $\alpha = (a_1,\cdots, a_s| b_1,\cdots, b_s)$ we have $b_j= a_j+d$
 for all $j$. Thus $Q_0(m)$ is the set of self-dual partitions:
 $\alpha = \alpha^*$.  Partitions from $Q_d(m)$ will be called 
 {\em $d$-unbalanced}. 
 
 \begin{thm}\label{thm:BGG}
 (a) Let $r\in \{0,1\}$ and $\ol r$ be the image of $r$  in $\ZZ/2$. There is a resolution of $\gen_\0 = \sp(2n)$-modules
\[
\K_r^\bullet =\bigl\{ \cdots \to \K_r^{-1}\to \K_r^0\bigr\}
\buildrel \on{qis}\over  \to \CC[x]_{\ol r}
\]
with
\[
\K_r^{-m} = \bigoplus_{\alpha = (b_1,\cdots, b_s|b_1,\cdots, b_s) \atop \sum b_i = m, 
\, |\alpha| \equiv r \text{ mod } 2} \V_{\alpha^\vee -\hn}.
\]
%Here $\alpha$ runs over  self-dual partitions
 %$(b_1,\cdots, b_s|b_1,\cdots, b_s) $ 
 %(in the Frobenius notation).  
 \vskip .2cm
 
 (b)  There is  a resolution   of $\gen= \osp(1|2n)$-modules
\[
\SK^\bullet = \bigl\{  \cdots \to \SK^{-2} \to \SK^{-1} \to \SK^0\bigr\} \buildrel
\on{qis}\over  \to \CC[x]
 \]
 with
 \[
 \SK^{-m} = \bigoplus_{\alpha \in Q_1(2m)} \SV_{\alpha^\vee - \hn}. 
  \]
   Here $\alpha^\vee$ is defined by \eqref{eq:def-alpha-vee}. 

 \end{thm} 
 The  rightmost parts of the resolutions of the theorem have the form
 
 \be\label{eq:3-res-end}
 \begin{gathered}
 \cdots \lra \V_{(-{1\over 2}, -{1\over 2}, \cdots, -{1\over 2}, -2{1\over 2}, -2{1\over 2})}
 \lra \V_{(-{1\over 2})^n} \lra \CC[x]_\0 \to 0, 
 \\
 \cdots \lra    V_{ (-{1\over 2}, -{1\over 2}, \cdots, -{1\over 2}, -1{1\over 2}, -2{1\over 2}) }   
  \lra V_{(-{1\over 2}, -{1\over 2}, 
 \cdots,  - {1\over 2}, -1{1\over 2})} \lra \CC[x]_\1 \to 0,
 \\
 \cdots \lra \SV_{(-{1\over 2}, -{1\over 2}, \cdots, -{1\over 2}, -1{1\over 2},
 -1{1\over 2})}
\lra \SV_{(-{1\over 2})^n} \lra \CC[x]  \to 0.
 \end{gathered}
 \ee
 The ``generator'' parts of the resolutions (the terms of degree $0$)
 correspond to $\alpha = 0$, $\alpha =1$ 
 (i.e., $\alpha = (1,0,\cdots, 0)$) and $\alpha=0$
 respectively. To explain them, note that 
   $\CC[x]_\0$ is generated, as a $\gen_\0$-module,
 by the vector $1$ which is annihilated by $\nen_\0$ and is an 
 eigenvector on $\men_\0 = \gl(n)$ with the character
 $y\mapsto {1\over 2} \tr(y)$. This follows from \eqref {eq:xxdd-n}.
 Similarly, $\CC[x]$ is generated as a $\gen$-module  by the same $1$
 which is  now annihilated by $\nen$ and is an eigenvector of $\men=\men_\0$ as before. Finally, $\CC[x]_\1$ is generated,
 as a $\gen_\0$-module, by the space of linear polynomials
 which is annihilated by $\nen_\0$ and is a representation of
 $\men_\0$ associated to $\alpha =  (-{1\over 2}, -{1\over 2}, \cdots,
 -{1\over 2}, -1{1\over 2})$.  This gives surjections from the
 corresponding 
 (super) Verma modules to (the components of) $\CC[x]$,
 i.e., the rightmost morphisms in  \eqref {eq:3-res-end}. 
 
 \vskip .2cm
 
 The ``relation'' parts of the resolution (terms with $m=1$)
 correspond to the diagrams $\alpha = (2,2)$, $\alpha = (2,1)$
 and $\alpha = (1,1)$ respectively. Further terms of the resolution
 can have contributions from several $\alpha$.
 Note that the resoltuions are finite since 
  the contributions from diagrams with $>n$ rows are zero. 
 
 \vskip .2cm
 
 \begin{rems}\label{rems:durfee}
  (a) The partitions $\alpha$ appearing  in part (a) of the theorem
 are self-dual: $\alpha=\alpha^*$. 
 The homological degree $m=\sum b_i$ to  such $\alpha$ contributes,
  is not a function of $|\alpha|$.
 In fact, the length $s$ of the Frobenius notation is the size
of the {\em Durfee square} of $\alpha$ which is the maximal
square inscribed into $\alpha$ considered as a Young diagram,
see Fig. \ref{fig:durfee}. 

 \begin{figure}[h]
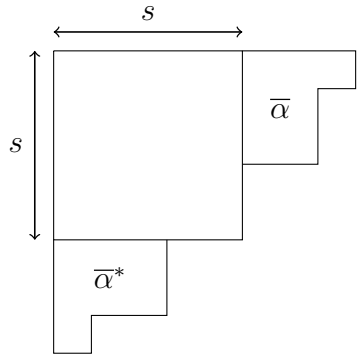
 
   \centering
   \btp[scale=0.5]
   \draw (-3,3) -- (2,3) -- (2, -2) -- (-3, -2) -- (-3,3); 
   \draw (2,3) -- (5,3) -- (5,2) -- (4,2) -- (4,0) -- (2,0); 
  \draw (-3, -2) -- (-3, -5) -- (-2, -5) -- (-2, -4) -- (0, -4) -- (0, -2); 
  
  \draw [<->, line width = 0.6] (-3, 3.5) -- (2,3.5); 
    \draw [<->, line width = 0.6] (-3.5, 3) -- (-3.5, -2); 

\node at (-0.5, 4)  {$s$}; 
  \node at (-4, 0.5)  {$s$}; 
  \node at (3,1.5){$\ol\alpha$}; 
  
 \node at (-1.5, -3){$\ol\alpha^*$}; 
   
   \etp
   \caption{The Durfee square and the remainder diagrams of a self-dual diagram.}
   \label{fig:durfee}
   \end{figure}
   
   \noindent The diagrams on the right and the bottom of the Durfee square
   are referred to as the {\em remainder diagrams} of $\alpha$.
   For $\alpha$ self-dual, the remainder diagrams are dual to
   each other and denoted $\ol\alpha$ and $\ol\alpha^*$. The 
   homological degree is then 
   %\be\label{eq:degree-durfee}
   \[
  m= \sum b_i \,= \,{s\choose 2} + |\ol\alpha| \, = \, {|\alpha| -s \over 2}. 
   \]
   (b) In contrast, for $1$-unbalanced partitions $\alpha$ appearing in part (b) of the theorem,  $|\alpha|$ is always even and 
   $m=|\alpha|/2$. 
 
 \end{rems}
   
 \noindent      Further, 
 the following  is classical, see, e.g.,
 \cite[Prop.2.3.9]{weyman-book}. 
 
 \begin{prop}
 For any finite-dimensional $\CC$-vector space $L$ we have
isomorphisms of $\GL(L)$-modules
 \[
 \begin{gathered}
 \Lambda^m(\Lambda^2(L)) \,\= \bigoplus_{\alpha\in Q_{1}(2m)} \Sigma^\alpha L, 
 \quad\quad 
 \Lambda^m(S^2(L)) \,\=  \bigoplus_{\alpha\in Q_{-1}(2m)}
  \Sigma^\alpha L.
 \qed
 \end{gathered} 
 \]
\end{prop}
\noindent This  allows us to write the $(-m)$th term of the super-resolution of
$\CC[x]$ as
\[
\SK^{-m} =  \Ind_{\pen}^{\gen}\, 
 \bigl( \Lambda^m(\Lambda^2 (\CC^n)^*) \otimes \Sigma^{-\hn}\bigr). 
\]
In a forthcoming paper with  V. Serganova  we will  give a direct construction
of the resolution in this explicit form and thus give an alternative proof and
a generalization  of
Theorem \ref {thm:BGG}(b).

% \vfill\eject
 
 \paragraph{Discussion and plan of the proof.}\label{par:plan-proof}
  Part (a) of Theorem
 \ref {thm:BGG} is in fact known. It is a particular case of the results of
\cite[\S3]{enright-hunziker} on BGG resolutions of the Wallach
representations, of which the 
$\CC[x]_i$, $i=\0, \1$,  are examples. See also \cite[Cor.6.8]{enright-hunzuker-pruett}. 

\vskip .2cm

This part is also related to the works 
\cite {weyman-resol, reiner-roberts} determining minimal
 free resolutions, over
the coordinate ring $\CC[\Sym_{n\times n} (\CC)] =  S^\bullet
(S^2(\CC^n))$ of the space of symmetric matrices,
of the  quotient ring $\CC[Z]$, where $Z$ is the variety of
matrices of rank $\leq 1$ and of the natural module ``of odd functions''
over $\CC[Z]$. 
The relation comes from the fact that like any induced Verma,  each $\V_\alpha$ is free
 over $U(\ol\nen_\0)$, where $\ol\nen_\0$ is the nilpotent
 subalgebra opposite to $\nen_\0$. As $\ol\nen_\0 = S^2(L_q) 
 \= S^2(\CC^n)$ is abelian, $U(\ol\nen_0)$ is the coordinate
 ring of symmetric matrices. 
 
 \vskip .2cm

 However,  the formulation
given in \cite{enright-hunziker, enright-hunzuker-pruett} 
uses notation somewhat different  from Theorem  \ref{thm:BGG}(a).
So before proving part (b) we
  explain the conceptual origin of 
 the resolutions of part (a). In fact, they are 
 consequences of an equivalence of categories constructed
by Enright and Shelton \cite{enright-shelton} and  of the ``standard'' parabolic BGG resolution of the trivial module, which we recall first. 

 \vskip .2cm
 
Part (b) will then be obtained from part (a) by
applying another equvialence
of categories due to Gorelik \cite {gorelik-typical} and connecting
representations of the Lie superalgebra $\gen = \osp(1|2n)$
and of its even part $\gen_\0  = \sp(2n)$.

% \vfill\eject
 
 \paragraph{Parabolic BGG for a regular block.} We start by
 recalling  the standard theory. 
 For this,  suppose
 that $\bg$ is a purely even reductive Lie algebra, and so are
 $\bp, \bn, \bm$. For definiteness let us specify that the
 category $\Oc^{\bg. \bp}$ consist of {\em ungraded} modules
 (formally, considering $\bg$ as a Lie superalgebra leads to the
 category of $\ZZ/2$-graded modules). 
 We choose a Cartan subalgebra
  $\bh\subset \bm$, which is then also a Cartan subalgebra in $\bg$. 
  Let $W$ be the Weyl group of $\gen$ and $l: W\to \ZZ_+$
  the length function. 
 Further,  let $W_\bm\subset W$ be the Weyl group of $\bm$ and 
$W^\bp\subset W$ be the set of minimal length representatives
of cosets from $W/W^\bm$. 
  
  \vskip .2cm

Let $\Zc(\bg)$ be the center of $U(\bg)$. 
For any character $\chi: \Zc(\bg)\to\CC$ we denote
by $\Oc^{\bg, \bp}_\chi\subset \Oc^{\bg, \bp}$ the $\chi$-{\em block} 
consisting of  generalized $\chi$-eigenvectors of $\Zc(\bg)$, i.e., of $v$  such that
for any $z\in \Zc(\bg)$ there is $n$ such that $(z-\chi(z)\cdot 1)^n v=0$. 
It is standard  that $\Oc^{\bg,\bp}$ splits
  into the
direct sum of the  $\Oc^{\bg, \bp}_\chi$. 

\vskip .2cm

By Chevalley's theorem $\Spec\,  \Zc(\bg)$ is identified with
 $W\backslash \bh^*$, where $W$ acts on
$\bh^*$, the space of weights, by the {\em dot action}
\[
w\cdot\lambda = w(\lambda+\rho)-\rho,
\]
with  $\rho$ being  the half-sum of the positive roots. So characters of $\Zc(\bg)$ are labelled
by $W$-orbits on $\bh^*$ with respect to this action. Given 
$\lambda\in\bh^*$, we will
denote by $\chi_\lambda: \Zc(\bg)\to\CC$ the character  corresponding to the orbit
$W\cdot\lambda$ and
write $\Oc^{\bg,\bp}_\lambda$ for $\Oc^{\bg,\bp}_{\chi_\lambda}$. 

\vskip .2cm

 Recall, see, e.g., \cite[\S1.8]{humphreys} that a weight $\lambda\in\hen^*$ is called
{\em (dot-)regular}, if $|W\cdot\lambda|=W$, i.e., $(\lambda+\rho, \alpha^\vee)\neq 0$
for any root $\alpha$,  and {\em (dot-)singular} otherwise. 
We will omit the prefix ``dot''.  
We will call the block 
$\Oc^{\bg,\bp}_\lambda$
regular or singular, if $\lambda$ is such. 

\vskip .2cm

 It is known that the  categories $\cal O^{\bg, \bp}_\lambda$
for all integral regular $\lambda$ are equivalent (Jantzen-Zuckerman translation principle \cite[\S7.8]{humphreys})
and so we will denote any such category $\Oc^{\bg,\bp}_\reg$. 
For example, the weight $0$ is regular and so we can think of 
$\Oc^{\bg,\bp}_\reg$
as $\Oc^{\bg,\bp}_0$. 

\vskip .2cm

Let $\lambda\in\bh^*$ be a dominant integral weight for $\bg$.
It is regular. It is also a dominant integral weight for $\bm$ 
and so we have the irreducible representations $\Sigma^\lambda_\bm$,
$\Sigma^\lambda_\bg$ of $\bm$ and $\bg$ respectively with
highest weight $\lambda$. We denote 
$\V_\lambda = \V_\lambda^{\bg,\bp} = \Ind_\bp^\bg \, \Sigma^\lambda_\bm$ the parabolic Verma module associated to
$\lambda$. Then we have the standard parabolic BGG
resolution in the regular block $\Oc_\lambda^{\bg,\bp}$:
\be\label{eq:BGG-stand}
\begin{gathered}
\Vc_\lambda^\bullet \,=\,\bigl\{ \cdots\lra  \Vc_\lambda^{-1} 
\lra \Vc_\lambda^0\bigr\} \buildrel \qis\over\lra  \Sigma^\lambda_\bg,
\\
\text{where} \quad 
\Vc_\lambda^{-m} \,=\bigoplus_{w\in W^\bp, \,\,  l(w)=m} 
\V_{w(\lambda+\rho_\bn)-\rho_\bn}
\end{gathered}
\ee
and $\rho_\bn$ is the half-sum of roots of $\bh$ in $\bn$. 
For different dominant integral $\lambda$, the $\Vc^\bullet_\lambda$
correspond to each other under the Jantzen-Zuckerman equivalences.
 Alternatively, we have 
  the  Beilinson-Bernstein equivalence  
%\be\label{eq:BB-eq}
\[
\Oc^{\gen,\pen}_\reg\,  \= \, \Perv(G/P, \Sc)
\]
with the category of perverse sheaves on $G/P$ smooth w.r.t. the stratification $\Sc$
by $B$-orbits (Schubert cells). Under this equivalence, Verma modules
$\V_{w(\lambda+\rho_\bn)-\rho_\bn}$ correspond to
constant sheaves on the Schubert cells and each $\Vc^\bullet_\lambda$
corresponds to the Cousin resolution of $\ul\CC_{G/P}$ by such
sheaves.

 %\vfill\eject
 
 \paragraph{ The Enright-Shelton reduction for $\sp(2n)$.}
 The two weights $-\hn$  and $(-{1\over 2}, \cdots, -{1\over 2}, -1{1\over 2})$
of $\gen_\0=\sp(2n)$ corresponding to $\CC[x]_\0$ and $\CC[x]_\1$
are  singular (and  not  integral either). They correspond, however,
to the same block  
  of 
$\Oc^{\gen_\0,\pen_\0}$ containing $\CC[x]_\0$ or $\CC[x]_\1$.
We denote this block $\Oc^{\gen_\0, \pen_\0}_{-\hn}$. It is 
 not equivalent to $\Oc^{\gen_\0,\pen_\0}_\reg$. 

\vskip .2cm

 Enright and Shelton
\cite{enright-shelton} studied singular parabolic categories 
$\Oc^{\bg,\bp}_\lambda$ for classical Hermitian
symmetric pairs (i.e., the cases when $\bp$ is a maximal parabolic subgroup in a classical $\bg$).
Their main result is that for singular $\lambda$ the category 
$\Oc^{\bg, \bp}_\lambda$
is identified with the direct sum of one or two coplies of the category $\Oc^{\bg', \bp'}_\reg$
for some other Hermitian symmetric pair $(\bg', \bp')$
with the root system of $\bg'$ being a subset of that of $\bg$. 
See \cite[Th.5.2.1]{enright-hunzuker-pruett} for a concise formulation omitting some details.

\vskip .2cm

In the case $\bg=\gen_\0 = \sp(2n), \bp=\pen_\0$ and 
$\lambda=-\hn$ the Enright-Shelton procedure
gives $\bg'=\so(2n)$ and  $\bp'$ being  the stabilizer of a maximal
($n$-dimensional) isotropic plane,  see 
\cite[Ex.6.4]{enright-hunzuker-pruett}, the case $k=1$. 
We denote this stabilizer by $\qen(n)$
and the corresponding parabolic subgroup $Q(n)\subset \SO(2n)$. 
The result of \cite{enright-shelton} in this case is as follows.

 \begin{thm}\label{thm:enright}
 We have an equivalence of categories
 \[
 \Ec: \Oc^{\gen_\0, \pen_\0}_{-\hn} \lra \Oc^{\so(2n), \qen(n)}_0 \oplus 
 \Oc^{\so(2n), \qen(n)}_0,
 \]
 with the properties:
 \begin{itemize}
 \item [(a)] 
 The image under $\Ec$ of a   Verma module  in $ \Oc^{\sp(2n), \pen(n)}_{-\hn} $
  is a  Verma module in one of the copies
 of  $ \Oc^{\so(2n), \qen(n)}_0 $ (paired with $0$ in the other copy). 
 
 \item [(b)] $\Ec(\CC[x]_\0) = (\CC, 0)$ and $\Ec(\CC[x]_\1)= (0,\CC)$. \qed
 \end{itemize} 
  \end{thm}
  
  %BELOW COMMENTED OUT
  \iffalse
  \begin{defi}\label{def:ES-invo}
  In particular, permuting the direct summands in the target  of 
  $\Ec$ gives an involutive self-equivalence of 
  $\Oc^{\gen_\0, \pen_\0}_{-\hn}$ which we call the
  {\em Enright-Shelton equivalence}. 
  \end{defi}
  \fi
  %ABOVE COMMENTED OUT

 The theorem implies that a resolution of $\CC[x]_i$, $i=\0, \1$ in terms of Verma modules in the block 
 $\Oc^{\, \gen_\0, \pen_\0}_{-\hn} $ 
 can be  obtained by applying  the inverse $\Ec^{-1}$
 to the BGG resolution \eqref {eq:BGG-stand}
  of $\CC$ in $\Oc^{\so(2n), \qen(n)}_0 $ (case $\lambda=0$). 
 The bijection (induced by $\Ec$)  of the posets labeling  Verma modules in 
 the source and target of   $\Ec$  
 was described in the proof of  \cite[Th.5.2.1]{enright-hunzuker-pruett}
  or in Proposition
 3.2(b) of \cite[Prop.3.2(b)]{enright-crelle},
  and one can verify that 
 applying the inverse of this bijection
 gives precisely the resolution claimed in  part (a) of Theorem \ref{thm:BGG},
 thus proving this part.

% \vfill\eject
 
 \paragraph{The Gorelik equivalence: resolution for $\osp(1|2n)$.}
 Part (b) of Theorem \ref{thm:BGG} will be
 deduced from  part (a) by  an equivalence of categories
  constructed by Gorelik
\cite{gorelik-typical} which we now recall. 

\vskip .2cm

Recall the notation $\gen = \osp(1|2n)$, $\gen_\0=\sp(2n)$, $\pen\subset \gen$ etc..
Let $\Zc(\gen_\0)\subset U(\gen_\0)$, $\Zc(\gen)\subset U(\gen)$ be the centers
of the universal algebras (understood in the super sense). 
Note that neither of $\Zc(\gen_\0), \Zc(\gen)$ is contained in the other. 

\vskip .2cm

We further denote by  $\wt \Zc(\gen)\subset  U(\gen)$ be the  ghost center  of $U(\gen)$, 
see \S \ref {subsec:DO-super} \ref{par:super-general}. So as a vector space
$\wt \Zc(\gen)=  \Zc(\gen) \oplus \Ac$, where $\Ac$ is the  anticenter. 
More explicitly, in our case  $\Ac= Q\cdot \Zc(\gen)$,
where $Q\in  \Ac$ is the particular element called {\em Casimir's ghost}, or {\em sCasimir} in 
\cite {arnaudon, musson-center}. 
  For $n=1$, the element $Q$ was described in \S \ref {subsec:evenheis-1} \ref{par:ghostcas-1}. 
  Thus for any $n$ the algebra
$\wt \Zc(\gen)$ is a free $\Zc(\gen)$-module of rank $2$, and $\Spec\,  \wt \Zc(\gen)\to \Spec\, \Zc(\gen)$ is a 
$2$-fold ramified covering, as $Q^2\in \Zc(\gen)$. 

\vskip .2cm

We want to compare appropriate blocks in  the parabolic $\Oc$-categories $\Oc^{\gen,\pen}$
and $\Oc^{\gen_\0, \pen_\0}$. 
  Note that that objects of $\Oc^{\gen, \pen}$ are
$\ZZ/2$-graded representations (as $\gen$ is a Lie superalgebra),
while  we assume that objects of  $\Oc^{\gen_\0, \pen_\0}$ are   ungraded. 

\vskip .2cm

Let $\phi = \phi_{-\hn}:  \Zc(\gen)\to\CC$ be the character by which  $\Zc(\gen)$ acts on the irreducible
 $\gen$-module $\CC[x]$
and $\Oc^{\gen, \pen}_{-\hn}$ be the corresponding block of $\Oc^{\gen, \pen}$.

\vskip .2cm

Let $I = I_{-\hn}=\Ker(\phi)$
be the maximal ideal corresponding to $\phi$, a point of $\Spec\,  \Zc(\gen)$. 
The character $\phi$ (or, equivalently, the maximal ideal $I$) is 
{\em weakly atypical} in the terminology of \cite[\S4]{gorelik-typical}, that is,  each Weyl group translation
of the highest weight $\hn$
of $\CC[x]$ lies on just one among the  hyperplanes $(-, \beta -\rho)=0$ corresponding to positive odd roots $\beta$.
See, e.g.,  \cite{malcom} for detailed discussion. This implies, in particular,  that $I$ is the image of a ramification point 
$\wt I$ of
the $2$-fold covering $\Spec\, \wt \Zc (\gen)\to\Spec\, \Zc(\gen)$, and  $\wt I$ is the unique preimage of $I$. 
As an ideal, $\wt I  = I  + Q \Zc(\gen)$, see  \cite[\S4]{gorelik-typical}. 
Let $\wt\chi: \wt \Zc(\gen)\to \CC$ be  the  character with kernel $\wt I$.

\vskip .2cm

At the same time, let $\chi = \chi_{-\hn}: \Zc(\gen_\0)\to\CC$ be the (same) character by 
which $\Zc(\gen_\0)$ acts on 
the irreducible $\gen_\0$-modules $\CC[x]_i$, $i=\0, \1$ and $\Oc^{\gen_\0, \pen_\0}_{-\hn}$ be the corresponding
block of $\Oc^{\gen_\0, \pen_\0}$. 
Considerations of   \cite[\S4.5]{gorelik-typical} imply that  $\chi$ is a {\em perfect mate} for $\wt\chi$ in
the terminology of {\em loc. cit} and so   \cite[Th.4.3]{gorelik-typical} gives:

\begin{thm}\label{thm:gorelik}
(a) The functors
\[
\begin{gathered}
\Phi: \Oc^{\gen, \pen}_{-\hn} \lra \Oc^{\gen_\0, \pen_\0}_{-\hn}, \quad \wt N \mapsto (\wt N_\0)_\chi,
\\
\Psi: \Oc^{\gen_\0, \pen_\0}_{-\hn}\lra  \Oc^{\gen, \pen}_{-\hn}, \quad N \mapsto (\Ind_{\gen_\0}^\gen N)_{\wt\chi}
\end{gathered} 
\]
are mutuallly quasi-inverse equivalences of categories. Here the subscript $\chi$ or $\wt\chi$
means the generalized eigenspace corresponding  to the character $\chi$ or $\wt\chi$. 
\vskip .2cm

(b) The functors $\Phi$ and $\Psi$ take induced Verma modules to induced Verma modules. 

\qed
\end{thm}

\noindent Note that $\Oc^{\gen, \pen}_{-\hn}$ has an involutive self-equivalence $\Pi$ given by the change of
parity of super-representations. It corresponds to the  involution on $ \Oc^{\gen_\0, \pen_\0}_{-\hn}$
(the category of ungraded representations) 
 coming from the  interchange of the two components in the direct sum identification given by
 Theorem \ref {thm:enright}.

 \vskip .2cm
 
 In particular $\CC[x], \Pi\CC[x]\in\Oc^{\gen, \pen}_{-\hn}$
 correspond to $\CC[x]_\0, \CC[x]_\1\in  \Oc^{\gen_\0, \pen_\0}_{-\hn}$ respectively. 
 This means that  a resolution of $\CC[x]$ in $\Oc^{\gen, \pen}_{-\hn}$ can be obtained by applying $\Psi$
 to the BGG resolution $\K_\0^\bullet$ of $\CC[x]_\0$ in $ \Oc^{\gen_\0, \pen_\0}_{-\hn}$ given
 by Theorem \ref{thm:BGG}(a). So to prove part (b) of Theorem \ref{thm:BGG} we just need
 to identify the Verma modules in $\Psi(\K_\0^\bullet)$ explicitly.  
It is convenient to proceed as follows.

\paragraph{Reformulation via Kostant homology.}\label{par:kostant}

 Recall the decomposition $E=L_p\oplus L_q$, so that $\nen = L_p\oplus S^2 L_p$. 
 Denote for short $L=L_q \= L_p^*$. Then $\ol\nen =L\oplus S^2 L$ is the 
 radical of the opposite parabolic $\ol\pen  = L\oplus L\cdot E$. 
 Then any (parabolic, induced) Verma module is free over $U(\ol\nen)$:
 \[
 \SV_\alpha = \Sigma^\alpha L \otimes_\CC U(\ol\nen), 
 \]
 equivariantly with respect to the Levi $\men = \gl_n$. So the $(-m)$th term
 of the image under $\Psi$ of the resolution from Theorem \ref{thm:BGG}(a) is 
 \[
 \Psi(\K^{-m}_\0) \, \= \, \Ind_{\pen}^\gen \,\, \bigl(H_m(\ol\nen, \CC[x])\bigr), 
 \]
where $H_m=H_m^\Lie$ is the Lie (super)algebra homology. 

\vskip .2cm

Now   the only supercommutator in $\ol\nen$ is 
 g the symmetric pairing $L\otimes L\to S^2(L)$.  Considered as a $\ol\pen$-module, i.e., as a $\gl_n$-equivariant 
 $\ol\Nen$-module, 
 $\CC[x] \= S^\bullet (L) \otimes \CC_{-\hn}$ where the $\ol\nen$-action on
 $S^\bullet(L)$ is given by multiplication:
 \[
 L\otimes S^p(L) \to S^{p+1}(L), \quad S^2(L)\otimes S^p(L) \to S^{p+2}(L). 
 \]
Since tensoring with the $1$-dimensional representation $\CC_{-\hn}$ can be
performed at the very end,  Theorem \ref{thm:BGG}(b)
reduces to the following statement.

\begin{prop}\label{prop:Hm-Nen}
For any $m\geq 0$ we have an identification of $GL(n) = GL(L)$-modules
\[
H_m(\ol\nen, S^\bullet(L)) \, \= \, \bigoplus_{\beta\in Q_1(2m)}\Sigma^\beta(L)
\, = \, \Lambda^m(\Lambda^2(L)). 
\]
\end{prop}

\paragraph{Computation of the character.} To prove Proposition \ref {prop:Hm-Nen}, we
note that all the Verma
$\gen_\0$-modules in the resolution $\Kc_\0^\bullet$ correspond to different
representations of $\gl_n$. So the same must hold for the Verma $\gen$-modules
in $\Psi(\Kc_\0^\bullet)$. By the reformulation of \S \ref{par:kostant}, 
we conclude that  each representation 
of $\GL(n)$ enters the total space  $H_\bullet(\ol\Nen, S^\bullet(L))$ 
no more than  once. So the total representation content of this $H_\bullet$
can be determined by computing the Euler characteristic
\[
[H_\bullet(\ol\nen, S^\bullet(L))] \, := \, 
\sum_m (-1)^m [H_m(\ol\nen, S^\bullet(L))] \, \in \,
\wh K_0(\Rep (\GL_n)) \,= \, \ZZ[[y_1, \cdots, y_n]]^{S_n}
\]
in the completed Grothendieck ring  of representations  of $\GL(n)$ which is the ring
of symmetric power series. 
This Euler characteristic is equal to that of the  chain (Chevalley-Eilenberg)
complex
\[
[C_\bullet(\ol\nen, S^\bullet(L))] \, = \, 
[S^\bullet(L) \otimes\Lambda^\bullet(S^2(L)) \otimes S^\bullet(L)]. 
\]
Here each $S^p(L)$ in the first appearance on the RHS has homological
degree $p$ as it comes from the exterior algebra of $\Pi L \oplus S^2(L)$,
so the Euler characteristic of the first $S^\bullet(L)$  gives
$ \sum_p (-1)^p [S^p(L)] = \prod(1+y_i)^{-1}$. The character 
$[\Lambda^\bullet(S^2(L))]$  is $\prod_{i\leq j} (1-y_i y_j)$,
as again, $\Lambda^p(\Lambda^2)$ is in homological degree $p$.
Finally, the second $S^\bullet(L)$ (the coefficient module)
is situated in homological degree $0$, so its character
is $\sum_p [S^p(L)] = \prod (1-y_i)^{-1}$. This gives
\[
[C_\bullet(\ol\nen, S^\bullet(L))]  = {\prod_{i\leq j} (1-y_i y_j)
\over \prod_i (1+y_i) \prod_i (1-y_i)} \, = \,
 {\prod_{i\leq j} (1-y_i y_j)
\over   \prod_i (1-y_i^2)} \, = \, \prod_{i<j} (1-y_iy_j) \,=\,
[\Lambda^\bullet(\Lambda^2(L))]. 
\]
This shows that as a $\GL(L)$-module  (possibly without respecting the
homological grading)
\be\label{eq:w/ograding}
H_\bullet(\ol\nen, S^\bullet(L)) \, \= \, \Lambda^\bullet(\Lambda^2(L)). 
\ee

\paragraph{Determination of the grading.}  Notice that $\ol\nen$ is a central
extension of abelian Lie superalgebras
\[
0\to  S^2(L) \lra \ol\nen \lra \Pi L \to 0 \quad  \quad  (\Pi = \text{change of parity}),
\]
so we have a
Hochschild-Serre
spectral sequence 
\[
E^1_{pq} = C_q(\Pi L,  H_p(S^2(L), S^\bullet(L))) = 
S^q(L)\otimes  H_p(S^2(L), S^\bullet(L)) \, \Rightarrow \,
 H_{p+q}(\ol\nen, S^\bullet(L)). 
\]
This sequence is $\GL(L)$-equivariant. 
In virtue of \eqref{eq:w/ograding},  to prove Proposition \ref{prop:Hm-Nen}
(and thus Theorem \ref{thm:BGG}(b)), it suffices to 
 establish the following.

\begin{prop}\label{prop:SS-once}
Let $\beta\in Q_1(2m)$ be a $1$-unbalanced diagram. Then $\Sigma^\beta(L)$
enters the term $E^1$ above exactly once (and so cannot be killed
 by the differentials and
remains present in $E^\oo$). In addition, the summand $E^1_{pq}$ where it enters, 
has  $p+q=m$., so $\Sigma^\beta(L) \subset H_m(\ol\nen, S^\bullet(L))$. 
\end{prop}

\noindent {\sl Proof:} The homology 
$H_\bullet(S^2(L), S^\bullet(L)) =
 \on{Tor}_\bullet^{S^\bullet(S^2(L))}(\CC, S^\bullet(L))$,
 being the Kostant homology of $\sp(2n)$, is directly found from the BGG resolution of
 Theorem \ref{thm:BGG}(a). More classically,  as discussed in \S \ref {par:plan-proof} above,
 this homology is known as the 
 space of generators of the minimal resolution of the coordinate ring of symmetric
 matrices of rank $\leq 1$  as well as of the natural ``dualizing'' module over it,
 see \cite [Th.3.19]{weyman-resol}, \cite[Th.6.3.1] {weyman-book} or 
 \cite[Th.1.2]{reiner-roberts}. The result (following also from  Theorem \ref{thm:BGG}(a)) is that
 \[
 H_p(S^2(L), S^\bullet(L)) \,=\, \bigoplus_{\alpha  = (b_1,\cdots, b_s|
 b_1, \cdots b_s) \atop \sum b_i = p} \Sigma^\alpha(L), 
 \]
the sum over self-dual diagrams $\alpha$ with $(|\alpha|-s)/2=p$, see
 Remark \ref {rems:durfee}(a). 
Recall the  Pieri formula,
see, e.g., \cite[Cor.2.3.5]{weyman-book}.

\begin{prop}
For any Young diagram $\alpha$, self-dual or not, we have a 
multiplicity-free decomposition
\[
S^q(L) \otimes\Sigma^\alpha(L) \, \= \,  \bigoplus_\beta \, \Sigma^\beta(L)
\]
where $\beta$ runs over Young diagrams  obtained from $\alpha$ by adding 
$q$ new cells in such a way that no two new cells lie in the same column.
\qed
\end{prop}

Given a self-dual $\alpha = (b_1,\cdots, b_s|
 b_1, \cdots b_s)$, we form the $1$-unbalanced diagram
 $\beta=\beta(\alpha) =  (b_1,\cdots, b_s|
 b_1+1, \cdots b_s+1)$ by adding one new cell to each leg of $\alpha$.
 By the Pieri formula,
 $\Sigma^{\beta(\alpha)}(L) $ is contained in
  $S^s(L) \otimes \Sigma^\alpha(L)= E^1_{p, s}$ where $p=(|\alpha|-s)/2$,
  so it contributes to the homological degree $|\beta(\alpha)|/2$,
  as claimed in Proposition \ref{prop:SS-once}. 
  
  \vskip .2cm
  
  It remains to prove that each such $\Sigma^{\beta}(L)$, $\beta = \beta(\alpha)$,
  appears in the $E^1$-term only once.  In other words, if $\Sigma^\beta(L)\subset
  S^q(L)\otimes\Sigma^\gamma(L)$ for some other self-dual diagram $\gamma$,
  then $\gamma = \alpha$ and $q=s$. Indeed, $\beta$ has the same arm lengths
  as $\alpha$. So if a self-dual $\gamma$ is contained in $\beta$, the arms
  of $\gamma$ must be contained in the arms of $\alpha$ and so 
  $\gamma\subset
  \alpha$ by self-duality.  Now, if $\gamma\neq\alpha$, then some leg of $\gamma$,
  denote this leg by $\Len$, must be shorter than the corresponding leg of $\alpha$,
  denote that leg by  $\Len'$. Since $\beta$ is obtained from $\alpha$ by adding
  one cell to each leg of $\alpha$, including $\Len'$, it is obtained from $\gamma$
  by adding, among
  other things, at least two new cells to the leg $\Len$. This contradicts the
  Pieri formula, since the  leg $\Len$ is a part of a column.
  
  Proposition  \ref{prop:SS-once} and hence Proposition \ref{prop:Hm-Nen} and
  Theorem \ref{thm:BGG} are proved. 
  
    \vfill\eject

\subsection{The super Siegel Laplacian. Even Heisenberg relations}
\label{par:super-lap-n}

\paragraph{The BGG resolutions at the level of differential operators.}
Recall the  setup and notation of \S \ref {subsec:SLG-sSieg}, that is: 
\[
\begin{gathered}
\GG=\OSp(1|2n,\CC) \,\, \supset \,\, \PP = \Stab(L_p), \quad \SLG(E)=\GG/\PP,
\\
\gen = \osp(1|2n) \supset \pen = \Lie(\PP),
\\
G=\Sp(2n, \CC) \,\,\supset \,\, P=\Stab(L_p), \quad \LG(E)=G/P,
\\
\gen_\0 = \Lie(G)=\sp(2n) \supset \pen_\0=\Lie(P),
\\
\men = \men_0 = \gl_n \text{ is the Levi in both } \pen \text{ and } \pen_\0.
\end{gathered}
\]
Further, the Levi subgroup $M\subset P$ is $\GL(n)$
while $\MM\subset \PP$ is $\O(1)\times \GL(n) = \{\pm 1\}\times \GL(n)$. 
Accordingly, the double covers $\wt \MM$ and $\wt\PP$ are disconnected
while $\wt M=\wt\GL(n)$ and $\wt P$ are connected. 

\vskip .2cm
 
Consider the block $\Oc_{-\hn}^{\gen_\0, \pen_\0}$ of the parabolic category
$\Oc$ for $\gen_\0=\sp(2n)$. The $\hen$-weights on module from this block
are all strictly half-integral, so the $\pen_\0$-action integrates to a unique
action of $\wt P$, with $H=\{\pm1 \}\subset \wt P$ acting by the nontrivial character
$\sigma$. 
This defines a full  embedding $\Oc_{-\hn}^{\gen_\0, \pen_\0} \subset \Oc^{\gen_\0,
\wt P, \sigma}$, see \S \ref {subsec:App-frac} for notation. 

\vskip .2cm

Similarly, if $V\in \Oc_{-\hn}^{\gen, \pen}$, then the action of $\pen$
integrates to that of the connected component $\PP^{(e)}$ of the $2$-component
group $\wt\PP = \{\pm 1\}\times\PP$ and we define the action of $(-1)$ in
the first factor to be given by the parity operator $(-1)^\F$. In this way be
get a full embedding $\Oc_{-\hn}^{\gen, \pen}
\subset \Oc^{\gen, \wt\PP,  \sigma}$. 

\vskip .2cm

In particular the BGG resolutions $\K^\bullet_r$ of $\CC[x]_r$,
$r = \0, \1$ of Theorem \ref {thm:BGG} can be considered as resolutions in
$\Oc^{\gen_\0, \wt P,\sigma}$ and the resolution $\SK^\bullet$ of $\CC[x]$
can be considered as a resolution on $\Oc^{\gen, \PP, \sigma}$. Therefore we can
apply the localization functors $\sharp_{G, \wt P, \sigma}$ resp.
$\sharp_{\GG, \wt \PP, \sigma}$ given by \S  \ref {subsec:App-frac}  and 
Proposition \ref {prop:harsish-twisted-n} implies the following.

\begin{cor}\label{cor:BGG-D-mod}
We have equivariant rersolutions of twisted $\Dc$-modules
\[
(\K_r^\bullet)^\sharp_\LG \buildrel \qis \over\lra (\CC[x]_r)^\sharp_\LG, 
\,\, r=\0, \1, \quad
(\SK^\bullet)^\sharp_\SLG \buildrel \qis\over\lra \CC[x]^\sharp_\SLG
\]
on $\LG$ and $\SLG$ respectively, with
\[
\begin{gathered}
(\K^{-m}_r)^\sharp_\LG \,,=\, 
\bigoplus_{\alpha = (b_1,\cdots, b_s|b_1,\cdots, b_s) \atop \sum b_i = m, 
\, |\alpha| \equiv r \text{ mod } 2}  \Dc_\LG \otimes_{\Oc_\LG} \Lc_{\alpha^\vee -\hn},
\\
(\SK^{-m})^\sharp_\SLG \,= 
 \bigoplus_{\alpha \in Q_1(2m)} \Dc_\SLG \otimes_{\Oc_\SLG} S\Lc_{\alpha^\vee - \hn}.
\end{gathered}
\]
\qed
\end{cor}

The identifications of the terms of these resolutions imply, by 
\eqref {eq:3-res-end}, that 
\be\label{eq:BGG-last-Dmod}
\begin{gathered}
(\CC[x]_\0)^\sharp_\LG \, \= \, \Mc_\Delta \, := \, 
\Coker \bigl\{ \Dc_\LG \otimes_{\Oc_\LG} \Lc_{(-{1\over 2}, -{1\over 2}, \cdots, -{1\over 2}, -2{1\over 2}, -2{1\over 2})} \buildrel \Delta^\vee \over 
 \lra \Dc_\LG \otimes_{\Oc_\LG} \Lc_{(-{1\over 2})^n} \bigr\}, 
 \\
 (\CC[x]_\1)^\sharp_\LG \, \= \, \Mc_\dirac \, := \, 
 \Coker \bigl\{ \Dc_\LG \otimes_{\Oc_\LG} \buildrel \dirac^\vee \over 
 \lra \Dc_\LG \otimes_{\Oc_\LG} \Lc_{(-{1\over 2}, -{1\over 2}, 
 \cdots,  - {1\over 2}, -1{1\over 2})} \bigr\}, 
 \\
 \CC[x]^\sharp_\SLG \, = \, \Mc_\DD \, := 
  \,
 \Coker\bigl\{ \Dc_\SLG\otimes_{\Oc_\SLG} 
  S\Lc _{(-{1\over 2}, -{1\over 2}, \cdots, -{1\over 2}, -1{1\over 2},
 -1{1\over 2})} \buildrel ^L\DD^\vee \over
\lra \Dc_\SLG \otimes_{\Oc_\SLG}  S\Lc_{(-{1\over 2})^n}. 
  \bigr\} 
 \end{gathered} 
\ee

Here we introduced the notations $\Mc_\Delta, \Mc_\dirac, \Mc_\DD$ for the
cokernels and $\Delta, \dirac, ^L\DD$ for the morphisms to suggest
the Laplace and Dirac operators for the reason that we now explain. 

\paragraph{The last arrows: Laplace, Dirac and super-Laplace.} 
To understand the morphisms in \eqref {eq:BGG-last-Dmod}, consider the
induced morphisms on the sheaves of solutions $\ul\Hom_\Dc (-, \Oc)$: 
\be\label{eq:lap-dirac-n}
\begin{gathered}
\Delta: \Lc_{\hn}  \lra \Lc_{(2,2) + \hn}
 \quad \dirac: \Lc_{(1)+\hn} \lra \Lc_{(2,1) + \hn}
 \\
^L\DD: S\Lc_\hn \lra S\Lc_{(1,1)+\hn}.
\end{gathered}
\ee
These morphisms are given by invariant
differential operators (between $\Gc_{1/2}$-twisted
bundles) which we call the {\em Siegel Laplacian}, {\em Siegel Dirac operator}
and {\em Siegel super-Laplacian} respectively. 

\vskip .2cm

Let us  restrict to the Siegel plane $\Hen_n\subset \LG$ with coordinates
$T=\|t_{ij}\|$ and the super Siegel plane $S\Hen\subset \SLG$ with coordinates
$(T, \xi)$, see \S  \ref  {subsec:SLG-sSieg}  \ref {par:cells-gen-n} and \ref {par;super-siegel-n}. Then the $\Lc_\alpha$, resp. $S\Lc_\alpha$, $\alpha\in\Z2YD_n^-$ become
honest vector bundles equivariant w.r.t. $\Mp(2n)$ resp. $\OMp(1|2n)$
so that  $\Delta, \dirac, \DD$ are invariant differential operators. Let us
determine their coordinate form explicitly. 

\vskip .2cm

As $\Lc_\hn$ and $S\Lc_\hn$ are now identified, as bundles, with $\Oc_{\Hen_n}$
and $\Oc_{S\Hen_n}$ respectively, we can think of sections of $\Lc_{\beta+\hn}$
as functions $\Hen_n\to \Sigma^\beta$ and of sections of
$S\Lc_{\beta+\hn}$ as functions $S\Hen_n \to\Sigma^\beta$. The relevant nonzero
$\beta$ in \eqref {eq:lap-dirac-n} are $\beta = (2,2), (1), (2,1), (1,1)$. 
We recall that the corresponding representations of $\GL(n)$
are described as follows:
\begin{itemize}
\item $\Sigma^{(1)}(\CC^n) = \CC^n, \quad \Sigma^{(1,1)}(\CC^n)=
\Lambda^2 \CC^n$.

\item $\Sigma^{(2,1)}\CC^n$ is the space of polynomials in entries of an indeterminate
$2\times n$ matrix $\bpm x_1& \cdots &x_n \\ y_1 & \cdots & y_n
\epm$
(of which the rows transform by the vector representations of $\GL(n)$)
spanned by the polynomials $d_k^{ij} = x_ix_j y_k - x_i x_k y_j$.

\item $\Sigma^{(2,2)}\CC^n$ is the space of polynomials in entries
of an indeterminate symmetric $n\times n$ matrix $\|x_{ij}\|$ 
spanned by the $2\times 2$ minors
\[
d_{ij}^{kl} \,=\, \left| \begin{matrix} x_{ij} & x_{il} \\ x_{kj} & x_{kl}
\end{matrix}\right|, \quad i<k, \,\, j<l. 
\]
\end{itemize}
With these identifications, we have:

\begin{prop}
(a) The operator $\Delta$ takes a function $\phi = \phi(T)$ to the
collection of the functions
\[
\Delta_{ij}^{kl} (\phi) \,=\, \left| \begin{matrix} \del_{ij} & \del_{il} \\ \del_{kj} & \del_{kl}
\end{matrix}\right| (\phi), \quad i<k, \,\, j<l,
\]
representing a $\Sigma^{(2,2)}(\CC^n)$-valued function $\Delta \phi$.

\vskip .2cm
(b) The operator $\dirac$ takes a vector function $\psi(T) = (\psi_1(T),\cdots,
\psi_n(T))$ to the collection of the functions
\[
\dirac_k^{ij} (\psi) \,=\, \del_{ij} \psi_k - \del_{ik} \psi_j,
\]
representing a $\Sigma^{(2,1)}(\CC^n)$-valued function $\dirac \psi$. 

\vskip .2cm

(c) The operator $^L\DD$ takes a (super) function $\Phi = \Phi(T, \xi)$ into the
collection of the (super) functions
\[
^L \DD_{ij} (\Phi) \,=\, (-1)^F \DD_{ij} (\Phi), \quad 
\DD_{ij}(\Phi) = 
 [D_i, D_j]_- \, \Phi \,=\, 
 (D_i D_j - D_j D_i) \Phi, \quad i<j,
\]
representing a $\Lambda^2\CC^n$-valued (super) function $^L\DD \Phi$.
Here $D_i$ is the spinor derivative \eqref {eq:spin-der-n}. 
\end{prop}

\noindent {\sl Proof:} Let $\ol\nen = L_q\oplus S^2L_q$ be the milpotent
radical of $\ol \pen\subset\gen$, the parabolic opposite to $\pen$. As
a $U(\ol\nen_\0)$-module, $\V_\alpha \= \Sigma^\alpha\otimes U(\ol\nen_\0)$
and, as $U(\ol\nen)$-module, $\SV_\alpha \= \Sigma^\alpha\otimes U(\ol\nen)$
are free on the space of generators $\Sigma^\alpha$, cf. \S \ref {subsec:sweil-BGG}
\ref {par:kostant}. So the BGG resolution $\K^\bullet_r$, $r=\0, \1$, considered
as complexes of $U(\ol\nen_\0)$-modules, si the minimal free
resolution of $\CC[x]_r$ over $U(\ol\nen_\0)$ (such a resolution is unique
up to a unique isomorphism, see, e.g., \cite{weyman-book}). Similarly
$\SK^\bullet$ gives the minimal free resolution of $\CC[x]$ over $U(\ol\nen)$.

\vskip .2cm

We are interested in the last parts of these resolutions, i.e., in the minimal
presentations of our modules by generators and relations. The last
morphism of one of these resolutions, i.e.,  the explicit form of the relations
as expressions in the generators, translates directily into the shape
of $\Delta, \dirac, ^L\DD$ as differential operators. Using the basis
$q_1, \cdots, q_n$ of $L_q$, let $\{q_{ij} = q_j\cdot q_j\}$ (the symmetric product)
be the corresponding basis of $S^2L_q = \ol\nen_\0$. After an easy rescaling,
we can assume that $q_i\in \ol\nen_\1$ acts on $\CC[x]$ by multiplication with
$x_i$ and $q_{ij}$ acts by multiplication with $x_i x_j$. 

\vskip .2cm

Now, as a $U(\ol\nen_\0) = \CC[q_{ij}]$-module, $\CC[x]_\0$ is generated by the
vector $1$ subject to the relations
\[
\left| \begin{matrix} q_{ij} & q_{il} \\ q_{kj} & q_{kl}
\end{matrix} \right| \cdot 1 = 0 , \quad 
 i<k,\,\, j<l  \quad \bigl(
\text {as } 
\left| \begin{matrix} x_i x_j & x_i x_l \\ x_k x_j & x_k x_l
\end{matrix} \right| =0\bigr). 
\]
 This shows part (a) of the proposition.
 Similarly, $\CC[x]_\1$ is generated by the linear forms $x_1, \cdots, x_n$
 subject to the relations
 \[
 q_{ij} x_k - q_{ik} x_j = 0, \quad \bigl(\text{as} \quad  x_i x_j x_k - x_i x_k x_j=0\bigr), 
 \]
 givinng (b). Finally, as a $U(\ol\nen)$-module, $\CC[x]$ is again generated
 by $1$ which is now subject to the relations
 \[
 (q_iq_j-q_j q_i)\cdot 1 = 0 \quad \bigl(\text{as} \quad x_ix_j-x_jx_i=0\bigr),
 \]
 giving (c). Here, as in  Remark \ref{rem:left-right-D},
 the twist by $(-1)^\F$ appears from the difference between 
 eft and right invariant differential operators on $\NN^-$.  
 
  \qed
 
 \vskip .2cm
 
 Applying now the solution functor $\ul\Hom_\Dc (-, \Oc)$
 to the full resolutions of Corollary \ref {cor:BGG-D-mod}, we obtain
 complexes of differential operators
 \be\label{eq:BGG-DO}
\begin{gathered}
 0\to \ul\Ker(\Delta) \to  \Lc_{\hn} \buildrel \Delta\over\to  \Lc_{(2,2) + \hn} \to\cdots\to
 \bigoplus_{\alpha = (b_1,\cdots, b_s|b_1,\cdots, b_s) \atop \sum b_i = m, 
\, |\alpha| \equiv 0 \text{ mod } 2} \Lc_{\alpha + \hn} \to\cdots, 
\\
0\to  \Ker(\dirac) \to \Lc_{(1)+\hn} 
\buildrel\dirac\over\to \Lc_{(2,1) + \hn} \to\cdots\to 
\bigoplus_{\alpha = (b_1,\cdots, b_s|b_1,\cdots, b_s) \atop \sum b_i = m, 
\, |\alpha| \equiv 1 \text{ mod } 2}
 \Lc_{\alpha + \hn} \to\cdots, 
 \\
0\to  \ul\Ker(\DD) \to S\Lc_\hn \buildrel ^L\DD\over \to S\Lc_{(1,1)+\hn}\to
\cdots\to \bigoplus_{\alpha\in Q_1(2m)} S\Lc_{\alpha + \hn} \to\cdots
 \end{gathered}
\ee
on $\Hen_n$ resp. $S\Hen_n$. 

\begin{prop}
The complexes \eqref {eq:BGG-DO} are exact everywhere on $\Hen_n$ resp.
$S\Hen_n$ thus giving right resolutions of the sheaf kernels of $\Delta, \dirac,
\DD$ by equivariant vector bundles. 
\end{prop}

\noindent{\sl Proof:} Because the resolutions of Corollary \ref {cor:BGG-D-mod}
consist of locally free $\Dc$-modules, the $m$-th cohomology sheaves of the
complexes \eqref {eq:BGG-DO}  are identified as
\[
\ul\Ext^m_{\Dc_{\Hen_n}} (\Mc_\Delta, \Oc_{\Hen_n}), 
\quad \ul\Ext^m_{\Dc_{\Hen_n}} (\Mc_\dirac, \Oc_{\Hen_n}), \quad
\ul\Ext^m_{\Dc_{S\Hen_n}} (\Mc_\DD, \Oc_{S\Hen_n})
\]
respectively. 
Now, these sheaves are $\Mp(2n)$ resp. $\OMp(1|2n)$-equivariant, since the $\Dc$-modules 
 in questions 
are. So our statement  for the first two complexes
is a consequence of  the following general fact,
see    \cite[Th. 2.4.2+Rem.3]{kash-thesis}. This fact  implies that the $\ul \Ext$-sheaves in question must
vanish. 

\begin{prop}\label{prop:kashiwara} 
Let $X$ be a complex manifold and $\Mc$ a coherent left $\Dc_X$-module.
For any $m\geq 1$ the sheaf  $\ul\Ext^m_{\Dc_X}(\Mc, \Oc_X)$
is supported on an analytic subset of codimension $\geq 1$.
\qed
\end{prop}
The statement for the third complex, with $\ul\Ext^m$ over $\Dc_{S\Hen_n}$
follows from the theorem by applying the Penkov equivalence
(Theorem \ref {thm:penkov} above)
between $\Dc$-modules over $S\Hen_n$ and $\Hen_n$.  \qed

\paragraph{Component analysis of the super-Laplacian.}
Consider.a (local) holomorphic section $\Phi$ of $S\Lc_{\hn}$ on $S\Hen_n$
which we represent in components as a super-function
\[
\Phi(x,\xi) = \phi(x) + \sum_{k=1}^n \,\, \sum_{1\leq i_1 <\cdots < i_k \leq n}
\phi_{i_1\cdots i_k}(x) \, \xi_{i_1} \cdots \xi_{i_k}
 \]
using the independent coordinates $(x=\|x_{ij}\|, \xi = (\xi_1, \cdots, \xi_n))$. 
The following fact generalizes the component interpretation of the 
vanishing of $3d$ $\Nc=1$
super-conformal
Laplacian in \S \ref {subsec:super-Lap-gen2}
 \ref {par:super-lap-2}.

     \begin{prop}\label{eq:DD-comp-n}
     The equation $\DD\Phi=0$ is  equivalent
     to the following three conditions: 
   
     \begin{enumerate}
     \item[(1)]  [Constraints] $\Phi$ is at most linear in $\xi$, i.e., 
    has the form $\phi(x) + \sum_{i=1}^n
     \phi_i (x)\xi_i$ with $\phi,   \phi_i$ being ordinary holomorphic functions.

     \item[(2)]  $\phi$ satisfies the (Siegel-) Laplace equation $\Delta \phi=0$    
    
    \item[(3)]   The vector function  $\psi(x)  = (\phi_1(x) \cdots, \phi_n(x))$
     satisfies the  (Siegel-) Dirac equation  $\dirac \psi=0$. 
     
     \end{enumerate}
     
     \end{prop}

     \vskip .2cm
     
     \noindent{\sl Proof:}  For the components of $\DD$ we have explicitly:
     \be\label{eq;square-ij}
    \DD_{ij} =  [D_i, D_j]_- \,=\, 2  {\del^2\over\del \xi_i \del \xi_j}
     -  2\sum_k\, \,  \xi_k \, \wt \dirac_k^{ij} + 2 \sum_{k, l} \,\, \xi_k \xi_l \,
     \wt  \Delta_{ij}^{kl},
     \ee
     where 
      \[
    \wt  \dirac_k^{ij} = {\del\over\del\xi_i} {\del\over\del_{jk}} - 
     {\del\over\del\xi_j} {\del\over\del_{ik}}, \quad 
   \wt  \Delta_{ij}^{kl}  = \left\| \begin{matrix} \del_{ij}&\del_{il}\\ \del_{kj} &\del_{kl}
     \end{matrix}\right\|
     \]
     are  versions of the components of $\dirac$ and $\Delta$, now
     considered as operators in the space of all superfunctions. 
     
      \vskip .2cm 
     
     Suppose $\Phi =\Phi(x,\xi)$ is annihilated by all the $\DD_{ij}$. 
Let us write $F$ as a sum
\[
\Phi = \Phi^{(0)} + \Phi ^{(1)} + \cdots + \Phi^{(n)},
\]
where $\Phi^{(p)}$ is homogeneous in $\xi$ of degree $p$. 
 Similarly, we write $\DD_{ij}(\Phi)^{(p)}$
     for the degree $p$ component of $\DD_{ij}(\Phi)$. 
     The three summands in \eqref {eq;square-ij} change the degree in $\xi$ by
     $-2$, $0$ and $+2$ respectively.  
     
     \vskip .2cm
     
     Now, $\DD_{ij}(\Phi)^{(0)} = \phi_{ij}(x)$, as only the first summand in
      \eqref {eq;square-ij}  contributes. So we must have all $\phi_{ij}=0$,
     i.e., $\Phi^{(2)}=0$.
     
     \vskip .2cm
     
     Next,  
     \[
     \DD_{ij}(\Phi)^{(1)}\, = \, \sum_{1\leq k\leq n}
     g_k(x), \quad \text{where} \quad
    g_k(x) =  - 2  \eps_{ijk} \,   \phi_{ijk}(x)  + 2 \bigl( \del_{jk} ( \phi_i) - 
    \del_{ik} (\phi_j) \bigr).
     \]
    Here  $\eps_{ijk}$ is   totally antisymmetric,
         equal to $+1$,  if $i<j<k$ (and vanishing unless $i,j,k$ are distinct). 
 Suppose  that $i,j,k$ are distinct. Then we have three
     quantities
     \[
     a= \del_{jk}(\phi_{i}),  \quad b= \del_{ik} (\phi_{j}), \quad
     c=\del_{ij}  (\phi_k)
     \]
     and the condition $g_k(x)=0$ means
     that
     \[
     a-b = \pm (b-c) = \pm (a-c)  \quad \bigl( = \pm 2 \phi_{ijk}(x)\bigr), 
     \]
     which implies that $a=b=c$ and so $\phi_{ijk}(T)=0$. This proves that 
     $\Phi^{(3)}=0$,
     and  $\Phi^{(1)}$ satifies $\dirac_k^{ij}=0$ for all $i,j,k$ by   \eqref {eq;square-ij},
     giving the Siegel-Dirac equation  as claimed in (3). 
     
     \vskip .2cm
     
     Further, since $\Phi^{(2)}=0$, only the first and third summands in 
      \eqref {eq;square-ij} wil contribute to
      \[
      \DD_{ij}(\Phi)^{(2)} = \sum_{1\leq k<l\leq n}  g_{kl}(x) \xi_k \xi_l, \quad
      g_{kl}(x) = 2\eps_{ijkl} \phi_{ijkl}(x) + 2 (\Delta_{ij}^{kl} (\phi)(x).
      \]
     Here, simillarly to the above. $\eps_{ijkl}\in\{0, \pm 1\}$ is totally antisymmetric
     and equal to $1$ for $i<j<k<l$.  
     
     \vskip .2cm
     
     To analyze the consequences of $ \DD_{ij}(\Phi)^{(2)}=0$,
     consider the space $\Lambda^2\CC^n$ with basis $e_{ij}= e_i\wedge e_j$,
     $1\leq i<j\leq n$, with $e_i$ forming a basis of $\CC^n$. Consider further
     the tensor product $\Lambda^2\CC^n\otimes\Lambda^2\CC^n$ with
     basis $e_{ij}\otimes e_{kl}$.  As a $\GL(n)$-module, we have a
     decomposition 
     \[
     \Lambda^2\CC^n\otimes\Lambda^2\CC^n \,=\, \Sigma^{2,2} \CC^n
     \oplus \Lambda^4\CC^n. 
     \]
     As mentioned earlier, the 
      $2\times 2$ minors of a symmetric $n\times n$ matrix  
     of indeterminates span the representation  $\Sigma^{2,2}$. 
   Thus for each $T$ we have
    \[
    \sum_{  k<l}
    \Delta_{ij}^{kl}(\phi) e_{ij}\otimes e_{kl} \,\in\, \Sigma^{2,2}\CC^n, 
    \quad \text{while} \quad  \sum_{  k<l} \eps_{ijkl}\,  \phi_{ijkl}\,\in \, \Lambda^4\CC^n. 
    \]
     So vanishing of  $\DD_{ij}(\Phi)^{(2)}$ implies that 
     all $\phi_{ijkl}=0$, i.e., $\Phi^{(4)}=0$, and all $\Delta_{ij}^{kl}(\phi)=0$,
     as claimed in (2). 
     
     \vskip ,2cm
     
     Now, since $\Phi^{(3)}=\Phi^{(4)}=0$, the vanishing of all the
      $\DD_{ij}(\Phi)^{(3)}$
     implies at once that $\Phi^{(5)}=0$, as only the first term in 
      \eqref {eq;square-ij}   will contribute. In a similar way we establish
       that all further $\Phi^{(p)}=0$, i.e., condition (1) of the
       proposition holds and  the equations $\Delta\phi=0$, $\dirac \psi=0$
       are already established. 
       
       \vskip .2cm
       
       In the reverse direction it is clear, using 
         \eqref {eq;square-ij}, 
        that a function $\Phi(x,\xi)$
       satisfying the conditions (1-3) of the proposition is annihilated
       by all the $\DD_{ij}$. The proposition is proved. 
       
       \qed
       
       \vskip .2cm
       
       In fact, the argument just given can be upgraded to a more conceptual statement
       about the $\Dc$-modules $\Mc_\Delta, \Mc_\dirac, \Mc_\DD$   whose
       solution sheaves  are the kernels of the corresponding operators. 
      Recall (Theorem \ref{thm:penkov}) that we have the 
        Penkov equivalence $\eps^*: \Dc_{S\Hen_n}\Mod \to \Dc_{\Hen_n}\Mod$
induced by the embedding
 $\eps: \Hen_n\hra S\Hen_n$. 
 
 \begin{prop}\label{prop:DD-lap-dir}
 We have an ($\Mp(2n)$-equivariant) isomorphism of $\Dc_{\Hen_n}$-modules
 \[
 \eps^*(\Mc_\DD) \= \Mc_\Delta \oplus \Mc_\dirac.
 \]
 \end{prop}   
 
 \noindent{\sl Proof:} The argument of Proposition \ref {eq:DD-comp-n}
 does not use  the assumption 
 that $\Phi$ is a function on $S\Hen_n$ and applies equally
  to solutions  $\Phi$ with values in  any $\Dc_{S\Hen_n}$-module.
  It identifies the functors on the category of $\Dc_{\Hen_n}$-modules
  represented by  $\Hom$ from the LHS and  the RHS of the proposed isomorphism
  and so  implies   its existence. \qed 
  
  \paragraph{The super-Gaussian transform.}  \label{par:sgauss-n}   
   If $n>2$, then  $\DD\Phi=0$ is a system of several scalar differential
   equations  with constant coefficients
   on a superfunction $\Phi(x,\xi)$.
   Here $(x,\xi)$ can be allowed to be  any point of $\SLG^\circ = 
   \Sym_{n\times n}(\CC) \times \CC^{0|n}$. 
    The analysis of 
    \ref{subsec:super-Lap-gen2}\ref {par:super-gauss-2}
    extends  from the case $n=2$ to that of arbitrary $n$ without changes giving the following generalization of
    Proposition \ref {prop:sGauss-2}.

  \begin{prop}\label{prop:sGauss-n}
  (a) The simple exponential solutions of $\DD\Phi=0$ are given by the {\em super-Gaussians}
  \[
 \theta_x(T, \xi) =   e^{\pi \i\,\,  x^t T x}  (1+\sqrt{-\pi \i}\, \,  \xi\cdot x)  = 
 e^{ \pi \i\,\,  x^t T x + \sqrt{-\pi \i}\, \,  \xi\cdot x},\quad 
 x = (x_1,
 \cdots  x_n)^t \in \CC^n. 
\]
 Here $\xi\cdot x = \sum_{i=1}^n \xi_i x_i$, as $\xi$ is a row vector.
 The exponent of $\theta_x$ is equal to $\pi \i \, x\cdot x^t = \| \pi \i x_\mu x_\nu\|$ considered
 as a point of $\Sym_{n\times n}(\CC)$.  
 
 \vskip .2cm
 
 (b) In addition to $\DD\theta_x=0$,
 the super-Gaussians satisfy the differential equations
 \[
 \begin{gathered}
 \sqrt{-4\pi \i} \,\,   (D^{(1/2)}_{p_\nu} \, \theta_x)(T, \xi) =- 2 \pi \i\,\,  x_\nu\,  \theta_x(T, \xi),\quad \nu=1,
 \cdots, n, 
 \\
\sqrt{-4\pi \i} \, \,  \,  (D^{(1/2)}_{q_\nu}\, \theta_x) (T, \xi) = -{\del\over\del x_\nu} \theta_x(T, \xi), \quad \nu=1,\cdots, n. 
 \end{gathered} 
 \]
 (c) Let $v = \|v_{ij}\|\in \Sym_{n\times n}(\CC)$. Let $\Phi$ be a holomorphic solution of $\DD\Phi=0$
 defined near $v$ and let $k\geq 0$. There exists a finite linear combination 
 $\Phi' = \sum_{i=1}^N c_i \theta_{x^{(i)}}(T,\xi)$, $x^{(i)}\in \CC^n$, such that $\Phi-\Phi'$ vanishes
 at $v$ together with all the partial derivatives of order $\leq k$. \qed

  \end{prop}

Recall that $\Wc=L_2(\RR^n)$ is the space of the Weil representation of
$\OMp(1|2n)$. Similarly to \eqref{eq:FT-2} we define 
 the {\em Fourier transform automorphism}
 \be\label{eq:FT-n}
 \FT: \OMp(1|2n) \lra \OMp(1|2n), \quad g\mapsto \wt S g \wt S^{-1},
 \ee
 where $\wt S\in \Mp(2n)$ is any preimage of 
 $S = \bpm0 & -1_n\\ 1_n & 0\epm$. 
 Further, extending  \eqref {eq:SG-2}, we define the {\em super-Gaussian transform} 
     by
      \be
 \SG=\SG_n: \Wc \lra  \SF_{\hn}, \quad u(x) \mapsto \SG(u)(T, \xi) = \int_{x\in \RR^n} 
 u(x) \theta_x(T, \xi) dx
 \ee
 and obtain, from Proposition \ref {prop:sGauss-n}(b):
  \begin{prop}\label{prop:SG=BW-n}
 $\SG$ is an injective $\FT$-twisted morphism of representations of $\OMp(1|2n)$ taking values 
   in $\Ker(\DD)\subset \SF_{\hn}$.\qed
 \end{prop}
 
\paragraph{Even Heisenberg relations on $\ul\Ker(\DD)$ and $\Mc_\DD$. } 
Since $\DD$ is an invariant differential operator, the  operators
$D_y^{(1/2)}$, $y\in \gen$ 
describing  the action
of 
$\gen = \osp(1|2n)$  in $S\Lc_{\hn}$, 
 define endomorphisms of $\Mc_\DD$ as a $\Dc_{S\Hen_n}$-module
 and of $\ul\Ker(\DD)$ as a sheaf. That is, we have homomorphisms of 
 associative algebras
 \be\label{eq:ug-act-MD}
 U(\gen) \lra \End_{\Dc_{S\Hen_n}}(\Mc_\DD), \quad U(\gen) \lra
  \End_\Sh(\ul\Ker(\DD)).
 \ee
 Recall (\S  \ref{subsec:sweil-BGG}\ref{par:sweil-n})
  the Heisenberg algebral $\Heis_n = \Heis(E, {1\over 2} \omega)$  
 associated to the symplectic space $E=\CC^{2n}$ with
  the form ${1\over 2} \omega$.

  \begin{prop}\label{prop:even-Heis-gen-n}
  The operators $D_y^{(1/2)}$, $y\in\gen_\1=E$ on $\Mc_\DD$ and $\ul\Ker(\DD)$
   satisfy even Heisenberg relations, i.e., 
 \[
 [D_y^{(1/2)}, D_z^{(1/2)}]_-   \, = \,  {1\over 2} \omega(y,z) \cdot 1, \quad y,z\in E. 
 \]
 So the morphisms  \eqref{eq:ug-act-MD} descend to
  morphisms of associative algebras
 \[
 \Heis_n \lra \End_{\Dc_{S\Hen_n}}(\Mc_\DD), \quad \Heis_n \lra \End_\Sh(\ul\Ker(\DD)).
 \]
   \end{prop}
 
 \noindent{\sl Proof:} By Proposition 
 \ref {prop:sweil-heis-n}, the action of $U(\gen)$ on $\CC[x]$
 factors through $\Heis_n$. Now, $\Mc_\DD = \CC[x]^\sharp_\SLG$
 is obtained from $\CC[x]$ by a functorial construction.  Therefore
 the action of $U(\gen)$ on $\Mc_\DD$ factors through $\Heis_n$.
 The statement for $\ul\Ker(\DD) = \ul\Hom_{\Dc_{S\Hen_n}}(\Mc_\DD, \Oc)$
 follows. \qed

 \begin {rem}\label{rem:trunc-lapl-n}

 In the terminology of \cite{KK-Takei, yoshida}, Proposition
 \ref{prop:even-Heis-gen-n} means that 
  the
 operators $D_y^{(1/2)}$ in $S\Lc_\hn$ form a Jacobi structure relatively to 
 $\Mc_\DD$. 
 In our approach they are realized (in components) as
 as 
 $2^n\times 2^n$ matrix differential operators.
 The operators used in 
   \cite{yoshida} are  $(n+1)\times (n+1)$ matrix operators
   which can be obtained from our $D_y^{(1/2)}$ by truncation, i.e.,  by
   imposing the constraints of Proposition 
 \ref{eq:DD-comp-n} (a) beforehand rather than obtaining them
 as consequences at the end. Our untruncated $D_y^{(1/2)}$
 satisfy more transparent commutation relations than the truncated
 ones.  They also fit into the  conceptual interpretation of $\Mc_\DD$
 as $\CC[x]^\sharp_\SLG$ which allows us to 
  prove  Proposition
 \ref{prop:even-Heis-gen-n}  without any computations.    
 \end{rem}

  \vfill\eject
  
  \subsection{The general Riemann theta from SUSY point of view}
  
  \paragraph{The general Riemann theta.}
  Similarly to the genus $2$ case, the general  Riemann theta function  is defined
 \cite{mumford}  by 
 \be
 \Theta(z, T) = \sum_{\bn\in \ZZ^n} \exp\bigl( \pi \i \,  \bn^t T \bn  + 2\pi \i \bn^t \cdot  z\bigr),
 \quad T\in \Hen_n \,\, z\in \CC^n. 
 \ee
 Here both $\bn$ and $z=(z_1,\cdots,  z_n)^t$ are considered as column vectors. It satisfies the
``heat equations''
\be\label{eq:heat-gen-n}
\del_{ij} \Theta(z,T) = {-\i \over 2\pi} \,\, {\del^2 \Theta \over \del z_i \del z_j},
\quad 1\leq i,j\leq n, \quad \del_{ij} := {\del\over \del t_{ij}} + {\del\over \del t_{ji}}, 
\ee
and the quasi-periodicity conditions
\[
\Theta( z+\bm, T) =\Theta(z, T), \quad \Theta (z+T\bm, T) = 
\exp\bigl( - \pi \i \, \bm^t T \bm - 2\pi \i\,  \bm^t z\bigr), \quad \bm\in \ZZ^n. 
\]
The Thetanullwert is the function
\be
\theta(T) = \Theta (0,T) = \sum_{\bn\in \ZZ^n} \exp(\pi \i \, \bn^t T \bn\bigr), \quad T\in\Hen_n. 
\ee
  It  satisfies the   periodicity
 \[
 \theta(T+U) = \theta(T), \quad U = \|u_{ij}\|\in \Sym_{n\times n}(\ZZ), \,\, 
 u_{ii}\in 2\ZZ
 \]
  and modularity 
  \be\label{eq:modul-gen2}
  \theta (-T^{-1}) =  e^{\pi \, \i \, n \over 4} \cdot \det(T)^{1/2}\, \theta(T). 
  \ee

  \paragraph{Characterization  of the general $\theta(T)$ 
   by super-differential equations.}
  Let us consider $\theta(T)$ as  Siegel superform  $\theta(T,\xi)$ independent on $\xi$,
  i.e., a section of $S\Lc_{1/2}$ on $S\Hen_n$. 
  Recall that $\gen =\osp(1|2n)$ acts on $S\Hen_n$,  so the basis vectors 
  $p_\nu, q_\nu\in E=\gen_\1$, $\nu=1, \cdots, n$ give rise to differential operators 
  $D_{p_\nu}^{(1/2)}$,
  $D_{q_\nu}^{(1/2)}$ of order $1/2$ on $S\Lc_{1/2}$.

   \begin{thm}\label{theta:susy-n}
   (a) 
     $\theta(T,\xi)$ is a solution of  the following  system of differential equations:
      \[
 \DD \Phi=0, \quad e^{\sqrt{-4\pi\, i} \, D^{(1/2)}_{p_\nu}} \Phi = \Phi, \quad e^{\sqrt{-4\pi\, i} \, D^{(1/2)}_{q_\nu}} \Phi = \Phi, \,\,\,\nu=1,\cdots n, 
 \]
 
   (b) 
 Any local holomorphic section  $\Phi=\Phi(T, \xi)$ of  $S\Lc_{\hn}$ on $S\Hen_n$ satisfying the system
 in (a),  is a constant multiple of $\theta(T)$. 
 \end{thm}
 
 \begin{rem}
 Written in components, the above  system becomes a system of
 $2^n\times 2^n$ matrix differential equations.  The system of
 $(n+1)\times (n+1)$ matrix equations considered in
 \cite{yoshida} can be obtained from it by truncation,
 i.e., by imposing the constraints of  
 Proposition 
 \ref{eq:DD-comp-n} (a) at the very outset. 
 Our approach
 elucidates the supersymmetric meaning of that $(n+1)\times (n+1)$ system. 
 See Remarks  \ref{rem:trunc-1-gen2} and \ref {rem:trunc-lapl-n} above. 
  \end{rem}
  
    Part (a) of the theorem is proved right away by expressing $\theta(T)$
  in terms of super-Gaussians:
  \be \label {eq:theta-super-sum-n}
 \theta(T)  = \sum_{\bn\in \ZZ^n} \theta_\bn(T, \xi). 
 \ee
 By Proposition \ref {prop:sGauss-n}(a), each summand $\theta_\bn$
 satisfies $\DD\theta_n = 0$, thus giving  $\DD\theta = 0$.
 Further, 
  for any $\bn\in \ZZ^n$ we have, by part (b) of Proposition \ref {prop:sGauss-n},
   \[
   \begin{gathered}
   e^{\sqrt{-4\pi \i}\,  D^{(1/2)}_{p_\nu}} \, \theta_\bn = e^{-2\pi \i n_\nu} \theta_\bn = \theta_\bn, 
   \\
   e^{\sqrt{-4\pi \i} \, D^{(1/2)}_{q_\nu}}\,  \theta_\bn = \theta_{\bn-\bee_\nu},
   \end{gathered}
   \]
   where $\bee_\nu$ is the standard basis vector of $\ZZ^n$. 
   So the  second equation of Theorem  \ref{theta:susy-n}(a) 
    holds term by term  w.r.t. 
    \eqref{eq:theta-super-sum-n} while the third  one holds
   by shifting the terms in the sum.  
   
      We also note that the equation $\DD\theta=0$, i.e., (since $\theta$ is
   independent on $\xi$) $\Delta\theta=0$ also follows directly
   from the heat equations \eqref {eq:heat-gen-n}. 
   
   \paragraph{Diagonal factorization.} To prepare the ground for the
   proof of Theorem \ref {theta:susy-n}(b),  we study the relation 
   between  theta-functions,
    super-Gaussians and super-Weil representations
     for different values of $n$, To emphasize
   the dependence on $n$, we write $\theta^{(n)}(T), T\in \Hen_n$
   for the genus $n$ Thetanullwert and similarly $\theta^{(n)}_x(T, \xi)$,
   $x\in \CC^n$, $(T, \xi)\in S\Hen_n$,  for the super-Gaussian
   and $\Wc^{(n)}$ for the super-Weil representation of
   $\OMp(1|2n)$. 
   \vskip .2cm
 
Let $m,n\geq 1$. We have the diagonal embedding
\[
\delta_{m,n}: \Hen_m\times \Hen_n \lra \Hen_{m+n}, 
\quad (T, T') \mapsto \bpm T& 0 \\ 0 & T'\epm.
\]
We see directly that the Thetanullwerts factorize under such embeddings:
\[
(\delta_{m,n} ^* \theta^{(m+n)}) (T, T') \,=\, \theta^{(m)}(T) \cdot
\theta^{(n)}(T'). 
\]
Consider also the following (non-symmetric) super-version:
\be\label{eq:sdelta-m-n}
\begin{gathered}
S\delta_{m,n}: S\Hen_m \times S\Hen_n \lra S\Hen_{m+n}, \quad 
\bigl( (T, \xi), (T', \xi)\bigr) \mapsto 
\biggl( 
\bpm
 T&0 \\
\| \xi'_i \xi_j\| & T'
\epm,  (\xi, \xi')\biggr),
\\
\text{where} \quad 
(\xi, \xi') \, := \, (\xi_1, \cdots, \xi_m, \xi'_1, \cdots, \xi'_n)
\end{gathered} 
\ee
and  the off-diagonal block consists of the $\xi'_i \xi_j$,  $i=1,\cdots, n$,
$j=1, \cdots, m$. We see that the
relations \eqref {eq:rels-T=xi-n} for $(T, \xi)$ and $(T', \xi')$ imply
the validity of such relations for the target ot $S\delta_{m,n}$. 
As for super-Gaussians, we have: 

\begin{prop}\label{prop:sgauss-mult-n}
For  $x\in \CC^m, x'\in \CC^n$ let $(x,x')\in \CC^{m+n}$ be the concatenation
of $x$ and $x'$ similarly to  the second line of \eqref {eq:sdelta-m-n}. 
Then
\[
((S\delta_{m,n})^* \theta^{(m+n)}_{(x,x')}) \bigl((T, \xi), (T', \xi')\bigr) \,=\,
\theta^{(m)}_x(T, \xi) \cdot \theta^{(n)}_{x'} (T', \xi'). 
 \]
\end{prop}
 
 \noindent{\sl Proof: } Diirect calculation.\qed
 
 \vskip .2cm
 
 The non-symmetry in the definition of $S\delta_{m,n}$ corresponds
 to the non-symmetry of the multiplication of the 
 super-Gaussians which are non-homogeneous and so their product
 depends  on the otder. 
 
 \vskip .2cm
 
 At the level of super-Weil representations, consider also the map
 \[
 \boxtimes_{m,n}: \Wc^{(m)}\otimes \Wc^{(n)} \lra \Wc^{(m+n)},
 \quad f\otimes f' \mapsto f\boxtimes f', \quad
 (f\boxtimes f') (x, x') = f(x) f'(x'). 
 \]
 These maps are also non-symmetric in the sense that 
 $\boxtimes_{m,n}$ is not obtained from $\boxtimes_{n,m}$ by
 applying the standard symmetry isomorphism
 $\Wc^{(m)}\otimes \Wc^{(n)} \to \Wc^{(n)} \otimes \Wc^{(m)}$
 of super-vector spaces. Indeed, this isomorphism involves
 the Koszul sign rule while $f(x) f'(x')$ is the usual commutative
 product of functions. 
 Proposition \ref{prop:sgauss-mult-n} implies the following
 multiplicativity property of super-Gaussian transforms.
 
 \begin{prop}
 For $f\in\Wc^{(m)}$, $f'\in\Wc^{(n)}$ we have
 \[
\bigl( (S\delta_{m,n})^*  \SG_{m+n}(f\boxtimes_{m,n} f') \bigr)
\bigl( (T, \xi), (T', \xi')\bigr) \,=\,
\SG_m(f) (T, \xi) \cdot \SG_n(f') (T', \xi'). \qed
 \]
 \end{prop}
 
 By composing the diagonal embeddings, we get the
 embeddings of the product of $n$  copies of the
 (super) Lobachevsky plane
 \be\label{eq:diag-emb-fully}
 \begin{gathered}
 \delta_{1,\cdots, 1}: \Hen^n = \Hen \times \cdots \times \Hen \lra \Hen_n,
 \quad (\tau_1, \cdots, \tau_n) \,\mapsto \, \diag (\tau_1, \cdots, \tau_n),
 \\ S\delta_{1,\cdots, 1}: S\Hen\times\cdots \times S\Hen
 \lra S\Hen_n, \quad \bigl((\tau_1, \xi_1), \cdots, (\tau_n,\xi_n)\bigr) 
 \, \mapsto \, 
 \bpm \tau_1 & 0 \cdots &0
 \\
  \|\xi_i \xi_j\| & \ddots 0 & 0
 \\
 & &\tau_n
 \epm,
 \end{gathered} 
 \ee
 where the below-diagonal entries are $\xi_i \xi_j$, $i>j$,  and 
 the above-diagonal
 entries are $0$. They satisfy
 \be
 \begin{gathered} 
( \delta_{1,\cdots, 1}^* \theta^{(n)} )(\tau_1, \cdots, \tau_n) \,=\,
\theta_J (\tau_1)\, \cdots \,  \theta_J(\tau_n), 
\\
\bigl( (S\delta_{1, \cdots, 1})^*  \theta_x\bigr)  
\bigl( (\tau_1, \xi_1), \cdots, (\tau_n, \xi_n)\bigr) 
= \theta^{(1)}_{x_1}(\tau_1, \xi_1)\,  \cdots \, \theta^{(1)}_{x_n} (\tau_n, \xi_n),
\end{gathered}
 \ee
 where $\theta_J(\tau) = \theta^{(1)}(\tau)$ is the genus $1$ (Jacobi)
 theta and $x=(x_1, \cdots, x_n)^t\in\CC^n$. 
  \vfill\eject
  
  \subsection {The Koszul complex for the general Riemann theta}
  Here we prove Theorem  \ref {theta:susy-n}(b). 
  Our analysis extends that of \S \ref {subsec-koszul-2}. 
  
  \paragraph{The Koszul complex of pluri-harmonic superforms.} 
  \label{par:koszul-pluri}
  The sheaf $\uKer(\DD)\subset S\Lc_\hn$ is preserved by the
  operators $D_y^{(1/2)}$, $y\in \gen_\1$ and so by their exponentials
  which are well-defined differential operators of infinite order.
  By Proposition \ref {prop:even-Heis-gen-n}, the operators
   \be\label{eq:AB-nu-genn}
 A_\nu = e^{\sqrt{-4\pi \i} D_{p_\nu}^{(1/2)}} -1, \quad
 B_\nu =  e^{\sqrt{-4\pi \i} D_{q_\nu}^{(1/2)}} -1, \quad  \nu=1,\cdots n, 
 \ee
acting on $\uKer(\DD)$, commute with each other. 
So we can form  the corresponding Koszul complex  of sheaves
 \[
 \Cc^\bullet = \uKer(\DD) \otimes_\CC \Lambda^\bullet (\CC^{2n}) \,\, = \,\,
 \bigl\{ \uKer(\DD) \buildrel (A_\nu,  B_\nu) \over
 \lra
 \uKer(\DD)^{\oplus 2n} \lra 
  \uKer(\DD)^{\oplus  {2n\choose 2}} \lra
 \cdots \lra \uKer(\DD)\bigr\},
 \]
situated in the degrees $[0,2n]$.  
Thus $\ul H^0(\Cc^\bullet)$
 is the sheaf of solutions of the system of Theorem \ref {theta:susy-n}. 
 As in  \S\S  {} \ref {subsec:koszul-1} and   \ref {subsec-koszul-2}, 
 we deduce that theorem
 from the next more general
 statement and the two propositions implying it.

 \begin{thm}\label{thm:koszul-LD-n}
 We have $\ul H^0(\Cc^\bullet) \= \ul\CC_{\Hen_n}$, the constant
 sheaf generated by $\theta(T)$ while $\ul H^i(\Cc^\bullet)=0$ for $i>0$. 
 \end{thm}

  \begin{prop}\label{prop:koszul-lc-n}
 The sheaves $\ul H^i(\Cc^\bullet)$ are locally (and hence globally)
  constant on $\Hen_n$. 
 \end{prop}
 
 \begin{prop}\label{prop:koszul-gl-n}
  Consider the diagonal embedding
   $\delta = \delta_{1,\cdots 1}:\Hen^n\to\Hen_n$ from \eqref {eq:diag-emb-fully}
 and the corresponding 
 sheaf-theoretic restriction $\delta^{-1}\Cc^\bullet$, a complex on
 $\Hen^n$.  Then 
 $\ul H^0(\delta^{-1}\Cc^\bullet) \= \ul\CC_{\Hen^n}$
 is the constant sheaf spanned by  the restriction 
 \[
 \theta(\delta(\tau_1, \cdots \tau_n )) = \theta_J(\tau_1)
 \cdots \theta_J(\tau_n), 
 \]
 while $\ul H^i(\delta^{-1}\Cc^\bullet) = 0$ for $i>0$. 
  \end{prop}
  
  \paragraph{The $\Dc^\oo$-module version of the Koszul complex.} 
  More fundamentally, Proposition \ref {prop:even-Heis-gen-n}
  says that the $D_y^{(1/2)}$ act on the $\Dc$-module $\Mc_\DD$
  by endomorphisms
    and satisfy
  even Heisenberg relations as such. 
  
  \vskip .2cm
  
  So we form the $\Dc_{S\Hen_n}^\oo$-module
  $\Mc_\DD^\oo = \Dc_{S\Hen_n}^\oo \otimes {\Dc_{S\Hen_n}}\Mc_\DD$. 
  As the $D_y^{(1/2)}$ have effective order $\leq 1/2$ with respect to the natural 
  filtration of $\Mc_\DD$, 
   Proposition \ref{prop:exp-endo}  gives  well defined exponentials 
 endomorphisms 
 $e^{z D_y^{(1/2)}} \in \End_{\Dc^\oo_{S\Hen_n}}(\Mc_\DD^\oo)$
 for any $z\in \CC$ and so can define  $A_\nu, B_\nu \in \End_{\Dc^\oo_{S\Hen_n}}(\Mc_\DD^\oo)$
 by   \eqref {eq:AB-nu-genn}. 
 The even Heisenberg relations imply that  the $2n$ endomorphisms $A_\nu, B_\nu$  
 commute with each other  and so we can form the associated Koszul complex
 of $\Dc_{S\Hen_2}^\oo$-modules
 \[
 \Nc^\bullet = \Mc_\DD^\oo \otimes_\CC \Lambda^{-\bullet}(\CC^{2n}),
 \]
 situated in degrees $[-2n,0]$.  
 Similarly to Proposition \ref {prop:C-as-rhom},
 the complex $\Cc^\bullet$ from
 \S \ref {par:koszul-pluri} is then recovered as
 \be
 \Cc^\bullet \= \ul{R\Hom}_{\Dc_{S\Hen_n}^\oo}(\Nc^\bullet, \Oc_{S\Hen_n}) 
 \ee
 due to:
  \begin{itemize}
 \item [(1)] Flatness of $\Dc^\oo_{S\Hen_n}$
 over $\Dc_{S\Hen_n}$ (reduced by the Penkov equivalence, to the
 purely even case of 
  \cite[\S 8.2.15]{bjork}).
  
  \item[(2)] Vanishing of 
  $\ul \Ext^p_{\Dc_{S\Hen_n}}(\Mc_\DD, \Oc_{S\Hen_n})$   for $p>1$,
  which follows by homogeneity from Proposition
  \ref {prop:kashiwara}. 
\end{itemize}

 \paragraph{Characteristic variety of the $\Dc^\oo$-module Koszul complex.} 
 Similarly to \S \ref{subsec-koszul-2}\ref{par:Koszul-D-gen2}, we use
 the resolution $(\SK^\bullet)^\sharp_\SLG$
 of $\Mc_\DD$ from Corollary \ref{cor:BGG-D-mod} and Eq. \eqref{eq:BGG-last-Dmod}
 to deduce, in virtue of  Lemma \ref{lem:perf-complex}, that $\Nc^\bullet$
 is quasi-isomorphic to a finite complex of free $\Dc_{S\Hen_n}^\oo$-modules.
 So applying Theorem \ref{thm:Char-SS},  we see that 
  Proposition \ref {prop:koszul-lc-n}  is a consequence of the next fact
 which we prove in this paragraph. 
 \begin{prop}\label{prop:Ch(N)-genn}
  The characteristic variety of $\Nc^\bullet$ is contained in  $\ul 0_{\Hen_n}$, the zero section of $T^*\Hen_n$. 
  \end{prop}
  
  To prove this,
   we use the coordinates $T=\|t_{ij}\|$ on $\Hen_n$ and 
   the independent coordinates  $(x_{ij}, \xi_i)$ on $S\Hen_n$.  Let 
   $\lambda_{ij}$ be
 the momenta coordinates on $T^*\Hen_n$ corresponding to the vector fields $\del_{ij}$. 
 We arrange them into a symmetric matrix $\lambda = \|\lambda_{ij}\|\in \Sym_{n\times n}(\CC)$. Let $Z\subset \Sym_{n\times n}(\CC)$
 be the variety of matrices of rank $\leq 1$. 
 Putting 
   \[
  \ul Z \,=\, \Hen_n \times Z \,\,\subset \,\, \Hen_n \times \Sym_2(\CC) \,=\, T^*\Hen_n,
  \]
  we obtain:
  \[
  \Ch(\Mc_\DD) = \ul Z, \quad \Ch(\Nc^\bullet) \subset \ul Z.
  \]
  Here, the first equality follows from 
  Proposition \ref{prop:DD-lap-dir} relating $\Mc_\DD$ to the
  (Siegel versions of the) Laplace and Dirac operators, while the second
  inclusion follows 
  since $\Nc^\bullet$ is constructed out of copies of $\Mc_\DD^\oo$.

  \vskip .2cm
   Let $Z_{A_i}, Z_{B_i }\subset T^*\Hen_n$ be the loci of non-invertibility of $A_i, B_i $
 in $\Ec_{S\Hen_n}^\RR$. Then 
\[
\Ch(\Nc^\bullet) \,\subset \, \ul Z\cap Z_A \cap Z_B, \quad \text{where}
\quad Z_A := \bigcap_{i=1}^n Z_{A_i}, \,\,\,
Z_B := \bigcap_{i=1}^n Z_{B_i }. 
\]
 This is because a Koszul complex is  exact whenever one of the commuting operators
giving rise to it is invertible. 

\vskip .2cm

Now, in components $\Oc_{S\Hen_n} = \Oc_{\Hen_n} \otimes \Lambda [\xi]$,
where $\Lambda [\xi] = \Lambda [\xi_1, \cdots, \xi_n] \= \CC^{2^n}$
and so the $D_y^{(1/2)}$, $y\in \gen_\1$ as well as $A_i ,B_i $
are $2^n\times 2^n$ matrix differential operators. We view
$\Ec^\RR_{S\Hen_n}$ as $\Mat_{2^n}(\Ec^\RR_{\Hen_n})$. 

\vskip .2cm

In this realization, each $D_{p_i}^{(1/2)} = D_i$ is a matrix
operator with constant coefficients.  To estimate $Z_{A_i}$ we, 
as in $n=2$ case 
(\S \ref {subsec-koszul-2} \ref {par:char-var-gen2}), we apply
Proposition \ref{thm:aoki-inv} to $C=\sqrt{4\pi \i} \, D_i$, so $A_i = e^C-1$. 
We decompose $\Lambda[\xi]$
by degree as $\Lambda^0 \oplus \Lambda^1 \oplus \cdots\oplus \Lambda^n$
and represent $D_i$ as an $(n+1)\times (n+1)$ block matrix
with respect to this decomposition. Let
$D_{i,p,q}: \Oc_{\Hen_n} \otimes \Lambda^q \to \Oc_{\Hen_n}\otimes\Lambda^p$
be the $(p,q)$th matrix element of $D_i$ with respect to this decomposition.
As $ D_i = {\del/\del \xi_i} - \sum_j \xi_j \del_{ij}$, the only nonzero matrix
elements are the
\[
\begin{gathered}
D_{i, q-1,q} = {\del/\del \xi_i}: \,\, \Lambda^q \lra \Lambda^{q-1}, \text{ of order } 0, \text{ and} 
\\
D_{i, q+1, q} = -\sum_j \xi_j \del_{ij}: \,\, \Oc_{\Hen_n}\otimes \Lambda^q 
\lra \Oc_{\Hen_n}\otimes \Lambda^{q+1}, \text{ of order } 1. 
\end{gathered} 
\]
So putting, in the context of Proposition  \ref{thm:aoki-inv},  $r_p = p/2$, $p=0,\cdots, n$
and $\rho=1/2$ we get the condition $\ord (D_{i,p,q}) \leq r_p-r_q +\rho$ and,
moreover, 
\[
\sigma_{D_{i,p,q}}^{r_p-r_q+\rho} (\lambda) = \sigma_{D_{i,p,q}} (\lambda)
\]
i.e., the highest symbol equals the total symbol. This means that the block matrix
$\|\sigma_{D_{i,p,q}}^{r_p-r_q+\rho} (\lambda)\|$ coincides with $\sigma_{D_i}(\lambda)$.
Therefore the condition for invertibility of $A_i$ in $\Mat_{2^n}(\Ec^\RR_{\Hen_n})$
 at $(T,\lambda)$ given
by Proposition  \ref{thm:aoki-inv} is that no eigenvalue of $\sqrt{4\pi \i} \, \sigma_{D_i}(\lambda)$ (as an endomorphism
of $\Lambda[\xi]$) lies in $\i\,\RR$.

\begin{prop}
The only eigenvalues of $\sigma_{D_i}(\lambda)$    are $\pm \sqrt{-\lambda_{ii}}$. 
\end{prop}

\noindent{\sl Proof:} As in the proof of Proposition
\ref {prop:ev-sigma-Di-2}, we reduce to the $1$-dimensional
case. That is, consider first  $\CC^2 = \Lambda[\xi_1]$, the exterior
algebra in one variable and the $\CC$-endomorphism $\sigma(\alpha) = 
{\del/\del \xi_1} -\alpha \xi_1$ of it  depending on $\alpha\in \CC$. 
 Its eigenvalues are  $\pm\sqrt{-\alpha}$. Now decompose
 \[
 \Lambda[\xi_1,\cdots, \xi_n] \,=\, 
 \bigoplus_{p=0}^{n-1}\,\,\,  \bigoplus_{2\leq i_1 <\cdots < i_p\leq n}
 \Lambda[\xi_1] \cdot \xi_{i_1}\cdots \xi_{i_p}
 \]
  as a direct sum of several copies of $\Lambda[\xi_1]$ and write
 \[
 \sigma_{D_1}(\lambda) \,=\, \biggl({\del\over\del \xi_1} - \lambda_{11}
 \xi_1\biggr) - \lambda_{12} \xi_2 - \cdots - \lambda_{1n} \xi_n. 
 \]
 Here the part in brackets is $\sigma(\lambda_{11})$ operating
 in each copy of $\Lambda[\xi_1]$ while the other summands increase
 the grading $p$ and so 
represent a matrix   strictly block-triangular  with respect to the ordering
(by $p$) of the summands of the decomposition. So 
the other summands do not affect the eigenvalues and
we obtain
the stamenent for $\sigma_{D_1}(\lambda)$. The case of
any other $\sigma_{D_i}$ is similar. \qed
 
 \vskip .2cm
 
 As in  \S \ref {subsec:koszul-1} \ref {prop:char-var-gen1}, we now   invoke Theorem  \ref {thm:aoki-inv},
to find
\[
\ul Z \cap Z_A \,\,\subset \,\, \bigl\{ (T, \lambda) \bigl| \,\, \lambda\in Z, \,\, \lambda_{ii} \in \i \RR_-, \,\,\, i=1,\cdots, n \bigr\}. 
\]
 Next, $Z_B$ is obtained from $Z_A$ by the change change of variables given by
 the modular transformation $S: T \mapsto -T^{-1}$ of $\Hen_n$
 and therefore  $Z_B = (dS)^t(Z_A)$. 
 So it is enough to prove:
  \begin{lem}
 The intersection
 \[
  \ul Z \cap Z_A \cap (dS)^t(Z_A) \,=\, \ul 0_{\Hen_n}
 \]
 is the zero section of $T^*\Hen_n$. 
 \end{lem}
 
 \noindent{\sl Proof:} Explicitly, the differential of $S$ at a point  $T_0\in\Hen_n$  is found as the linear operator 
 \[
 \begin{gathered}
 d_{T_0} S \,=\, - d_{T_0}(T^{-1}) \,=\, T_0^{-1} (dT)  T_0^{-1},  \text{ which takes } 
  \\
 \bigl(X\in \Mat_{n\times n}(\CC) = T_{T_0}\Hen_n\bigr)
  \mapsto T_0^{-1} X T_0^{-1}= S(T_0) X S(T_0). 
  \end{gathered} 
 \]
 Since $T_0^{-1}$ is symmetric, the  matrix of this  operator in the basis
 of matrix units, is symmetric as well, so we can computationally identify $(dS)^*$
 with $dS$. 

\vskip .2cm

 Now recall the  $2:1$ parametrization
 \[
 \rho: \CC^n \lra Z, \quad  y= (y_1, \cdots, y_n)^t \mapsto y\cdot y^t = 
 \|y_i y_j\|. 
 \]
  If $(T_0,\lambda) \in \ul Z\cap Z_A$, then 
 $y=\rho^{-1}(\lambda)$,
 defined uniquely up to sign, must have $y_i^2 = \lambda_{ii} \in \i \RR_-$,
 so $y= \pm \sqrt{-\i} \, z$, $z\in \RR^n$. 
 
  If, at the same time,  $(T_0, \lambda) \in \ul Z \cap Z_B$, then
  by the above $(d_{T_0}S)(\lambda) = -\i w\cdot  w^t$ for some $w\in \RR^n$ as well. But  
  \[
   (d_{T_0}S)(\lambda) = -\i S(T_0) z \cdot z^t S(T_0) =
    -\i (S(T_0)z) \cdot (S(T_0) z)^t
   \]
 which implies, by the $2:1$ property of $\rho$, that $w=\pm S(T_0) z$, so
 $S(T_0)z\in \RR^n$. 
This means that the real matrix $\Im (S(T_0))$  
annihilates $z\in \RR^n$
which implies $z=0$ since $\Im( S(T_0))$ is positive definite.
This proves the lemma and   Proposition \ref {prop:Ch(N)-genn}. \qed

\paragraph{The Cauchy problem for the  super-diagonal embedding.}  Here we prove
Proposition \ref {prop:koszul-gl-n}. In addition to $\delta$, 
consider  also the super-diagonal embedding
$S\delta = S\delta_{1,.\cdots, 1}:  (S\Hen)^n \to S\Hen_n$ 
from \eqref {eq:diag-emb-fully} and the $\Dc$-module inverse image
$(S\delta)^* \Mc_\DD$. 

\begin{prop}
The subvariety $\delta(\Hen^n) \subset \Hen_n$ is non-characteristic
for $\Mc_\DD$.
\end{prop}

\noindent{\sl Proof:} 
 By definition \cite{kash-thesis}, non-characteristicity means that
 the conormal bundle $T^*_{\Hen^n} \Hen_n \subset T^*\Hen_n$
 intersects $\Ch(\Mc_\DD) = \ul Z$ 
 only along the zero section. Using the matrix
 coordinates $\lambda$ in $T^*\Hen_n$ as  before,
  we find that $T^*_{\Hen^n }\Hen_n$ is given by the condition $\lambda_{ii}=0$
  for all $i$. At the same time. 
 $\ul Z$ is given by the condition
  $\on{rk}(\lambda)\leq 1$ which implies that $\lambda_{ij} = z_i z_j$
 for some $z_1, \cdots, z_n\in \CC$. So $\lambda_{ii}=0$ means that $z_i^2=0$
 and so $z_i=0$
 for all $i$, i.e., $\lambda=0$. \qed
 
 \vskip .2cm
 
 The proposition implies (Cauchy-Kowalewska theorem for non-characteristic embeddings 
 \cite[\S2.3] {kash-thesis})
 that
 \[
 \delta^{-1} \ul\Hom_{\Dc_{S\Hen_n}}(\Mc_\DD, \Oc_{S\Hen_n}) \,
 \=
 \ul\Hom_{\Dc_{S\Hen^n}} ((S\delta)^*\Mc_\DD, \Oc_{S\Hen^n}).
 \]
 Recall that the $\Dc$-module pullback $(S\delta)^*$ is defined at the level of
 $\Oc$-modules and  $\Mc_\DD$, being homogeneous, is flat over $\Oc_{S\Hen_n}$.
 Therefore the full derived functor $L(S\delta)^*\Mc_\DD$ reduces to 
 $(S\delta)^*\Mc_\DD$. Further, by the same homogeneity argument
 (and Proposition   \ref {prop:kashiwara}) both $\ul \Hom$-sheaves in the above identification
 are the same as the corresponding $\ul{R\Hom}$. 
 
 \vskip .2cm
 
 Therefore the complex $\delta^{-1}\Cc^\bullet$ of Proposition 
 \ref {prop:koszul-gl-n} is identified as follows:
 \[
 \delta^{-1}\Cc^\bullet = 
 \ul{R\Hom}_{\Dc_{S\Hen^n}^\oo}(\Nc_\delta^\bullet, \Oc_{S\Hen^n}),   
 \]
 where $\Nc^\bullet_\delta$ is the Koszul complex of $\Dc_{S\Hen^n}^\oo$-modules
 associated to the operators  $A_\nu, B_\nu$ acting in the Cauchy
  restricted $\Dc^\oo$-module
 $ ((S\delta)^* \Mc_\DD)^\oo = (S\delta)^* \Mc_\DD^\oo$.
 So to identify the cohomology sheaves of $\delta^{-1}\Cc^\bullet$ as claimed
 in Proposition 
 \ref {prop:koszul-gl-n}, it suffices to identify it with the  external tensor product
 of $n$ copies of the genus $1$ complex from \S \ref {subsec:koszul-1}. 
 This is done by the next statement. 
 
 \begin{prop}
(a) The $\Dc_{S\Hen^n}$-module $(S\delta)^*\Mc_\DD$
is identified with $\Dc_{S\Hen^n}  = \Dc_{S\Hen} \boxtimes \cdots \boxtimes  \Dc_{S\Hen}$. 

\vskip .2cm

(b) Under this identification, the operators $D_{p_i}^{(1/2)}$ and $D_{q_i}^{(1/2)}$
are identified with
\[
(-1)^\F \boxtimes\cdots \boxtimes (-1)^\F \boxtimes D_p^{(1/2)} \boxtimes 1\boxtimes
\cdot  \boxtimes 1
\text{ and  } 
 (-1)^\F \boxtimes\cdots \boxtimes (-1)^\F \boxtimes D_q^{(1/2)} \boxtimes 1\boxtimes
\cdot  \boxtimes 1
\]
respectively, where $D_p^{(1/2)}, D_q^{(1/2)}$ are the genus $1$ operators
acting in forms of weight $1/2$ on $S\Hen$, see
Theorem \ref {thm:1d-theta-charact}. 
 \end{prop}
 
 \noindent {\sl Proof:} Let $X$ be a complex algebraic supermanifold, 
 $z\in X$ be a $\CC$-point and  $I_z\subset \Oc_X$ the
 corresponding maximal ideal. For a quasi-coherent $\Oc_X$-module $\Fc$
 we call $\Fc_z := \Fc/I_z\Fc$ the {\em fiber} of $\Fc$ at $z$.
 
 \vskip .2cm
 
 By definition, the $\Dc_{S\Hen_n}$-module $\Mc_\DD$ is the restriction, from $\SLG$
 to $S\Hen_n$, of the ($\Gc_{1/2}$-twisted)  $\Dc_\SLG$-module $\CC[x]^\sharp_\SLG$,
 the twisted homogenization (\S \ref{subsec:App-frac}) of $\CC[x]$,  the polynomial Weil 
 representation.  Thus the fiber $(\Mc_\DD)_z$ of $\Mc_\DD$
 considered as an $\Oc$-module,  at $z\in \Hen_n$, is identified with $\CC[x]$ uniquely
 up to action of   $\wt \PP$, the $2$-fold covering of
 the parabolic subgroup $\Stab(L_p) \subset \OSp(1|2n)$.  Let us refer to $\CC[x]$ as
 the {\em typical fiber} of $\Mc_\DD$. 
 
 \vskip .2cm
 
 Now, the $\Dc$-module pullback appearing in the Cauchy-Kowalewska theorem,
 is performed at the level of the underlying $\Oc$-modules, so it inherits
 the fibers in the above sense.  We will show that the pullback $(S\delta)^*$
 in question inherits also the natural 
 tensor product decomposition of the typical fiber $\CC[x]$.
 That is,  multiplication of polynomials defines an isomorphism of super vector spaces
 \be\label{eq:C[x]-tensor}
 \CC[x_1] \otimes \cdots \otimes \CC [x_n] \lra \CC[x] = \CC[x_1, \cdots, x_n]. 
 \ee
 For $i=1,\cdots, n$ let 
 $E_i  = \CC\, p_i \oplus \CC\, q_i \subset E =\CC^{2n}$ and $Y_i = \CC^{0|1}\oplus E_i
 \subset Y=\CC^{0|1}\oplus E$. 
 Each $\CC[x_i]$ is the polynomial Weil representation of $\gen^{(i)} = \osp(Y_i) 
 \= \osp(1|2)$ and so gives the twisted $\Dc$-module 
 $\Mc_\DD^{(i)} = \CC[x_i]^\sharp_{\PP^{1|1}_i}$ on
 its  super Lagrangian Grassmannian $\PP^{1|1}_i := \SLG(Y_i) \= \PP^{1|1}$.
 Let us  restrict $\Mc_\DD^{(i)}$  to $S\Hen(i)\subset \PP^{1|1}_i$, the $i$-th copy of
 the super-Lobachevsky plane which we identify with the $i$th factor in the
 source of $S\delta$.  We see that
 \[
 (S\delta)^* \Mc_\DD \, \= \, \Mc^{(1)}_\DD \boxtimes \cdots\boxtimes
 \Mc^{(n)}_\DD, 
 \]
since the identification \eqref {eq:C[x]-tensor} of the typical fiber extends
to the identification of fibers  everywhere over the diagonal locus. 

\vskip .2cm

Next, the (pluri) harmonicity condition  $\DD\Phi=0$  is vacuous in genus $1$,
 so $\Mc^{(i)}_\DD$
is the  locally free $\Dc_{S\Hen}$-module $\Dc(S\Lc_{1/2}, \Oc_{S\Hen})$
which, by our trivialization of $S\Lc_{1/2}$, we identify with  the globally
free $\Dc$-module $\Dc_{S\Hen}$.
This gives part (a) of the proposition. 

\vskip .2cm

To see part (b), we  notice that the action of $p_i, q_i \in \gen_\1$ on
$\CC[x]$ in the RHS of \eqref{eq:C[x]-tensor} is written on the LHS
in the pattern similar to what is claimed in (b). That is,
the action of $p_i$ is written as
\[
(-1)^\F \otimes\cdots \otimes (-1)^\F \otimes p_i \otimes 1\cdots \otimes 1
\]
and similarly for $q_i$. 
The appearance of $(-1)^\F$ reflects the fact that these  (odd) operators for
different $i$ commute rather than anticommute.  Since the action of
$D_{p_i}^{(1/2)}$ resp.  $D_{q_i}^{(1/2)}$ on the homogenization comes
from the action of $p_i$ resp. $q_i$ on the module, 
we deduce the statement of (b). \qed

\vskip .2cm
 
 This finishes the proof of Proposition  \ref {prop:koszul-gl-n} and so of 
 Theorems \ref{thm:koszul-LD-n} and \ref {theta:susy-n}.

    \vfill\eject
  
  \appendix
  \section {Morphisms between Verma modules and invariant
  differential operators: homogenization}\label{app:A}
  
  The relation between BGG resolutions and complexes of
  equivariant differential operators goes back to the very origins
  of  the BGG theory \cite{BGG} and has been discussed in
  detail in \cite [\S VI.5]{faltings} \cite[\S 11.2]{baston}. Here
  we present a more general approach which includes arbitrary
  objects of the appropriate category $\Oc$, not just Verma modules. 
  This appendix logically depends on Ch. \ref {sec:susy-source}  and
  will be used in subsequent chapters. 
  
  \subsection{ Harish-Chandra modules and equivariant $\Dc$-modules}
  \label{subsec:A-harish}

Let $G$ be an affine  complex  algebraic supergroup with Lie
 superalgebra $\gen$
and  $K\subset G$ an  (also complex) algebraic
sub-(super)group with Lie superalgebra $\ken$. 
Recall that a {\em (Harish-Chandra) $(\gen, K)$-module} is
a (super) vector space $V$ with structure of a  $\gen$-module
 and of a $K$-module 
 (in the sense of  algebraic supergroups, i.e., a coaction of the 
 Hopf (super) algebra $\Oc(K)$)  which
are compatible.
 That is, the  infinitesimal $\ken$-action obtained from the $K$-action
equals the restriction of the $\gen$-action.
We denote by 
$K\Mod$ and 
$(\gen, K)\Mod$ the categories of $K$-modules and $(\gen, K)$-modules. 

\vskip .2cm

Consider the homogeneous (super) space $X=G/K$
and the projection $p: G\to X$.  Let $\Oc_X\Mod^G$ resp.  $\Dc_X\Mod^G$
be the category of $G$-equivariant quasi-coherent $\Oc_X$-modules
resp.  $\Dc_X$-modules. 

\vskip .2cm

Each  $K$-module
$V$ gives rise to a $G$-equivariant vector bundle
(posibly of infinite rank) $V^\sharp_X = G\times_K V$
on $X$
with fiber over the distinguished point $1\cdot K$ equal to $V$.
In other words, $V^\sharp_X$ is the $G$-equivariant
 quasi-coherent sheaf on $X$
whose sections over  any Zariski open $U\subset X$ are given by
\be\label{eq:def-sharp}
\Gamma(U, V^\sharp_X)\,=\,
\bigl\{ f: p^{-1}(U) \to V \bigl| \, f(gk) = k^{-1} f(g), \quad \forall \, g\in p^{-1}(U),
\, k\in K\bigr\}. 
\ee
This  obviously defines a fully faithful exact functor
\[
\sharp_K: K\Mod \lra \Oc_X\Mod^G, \quad V\mapsto V^\sharp_X,
\]
  which we call
the {\em homogenization} (forming a homogeneous bundle). 

\begin{prop}\label{prop:harish1}
(a) If $V$ is a $(\gen, K)$-module, then $V^\sharp_X$ has a natural structure
of a $G$-equivariant  $\Dc_{X}$-module. It is $\Dc_X$-coherent, if $V$ is
finitely generated.

\vskip .2cm

(b) The correspondence $V\mapsto V^\sharp_X$ gives  a fully faithful  
exact functor
\[
\sharp_{\gen, K}: (\gen, K)\Mod \lra \Dc_X\Mod^G. 
\]
\end{prop}

The functor in (b) will also be called homogenization. 

\vfill\eject

\noindent {\sl Proof:} Only part (a) needs a proof, as (b) then follows
from the corresponding  properties of $\sharp_K$. 
To prove (a), 
we start with the case $K=\{1\}$, so $X=G$. 
Let $\Dc(G) = H^0(G, \Dc_G)$ be the ring of global differential operators
on $G$.  Then $U(\gen) = \Dc(G)^G$ is realized as the subring of
invariants with respect to the $G$-action by right translations. Note that
product of differential operators induces an isomorphism of $\Oc_G$-modules
\be\label{eq:D_Gsemi}
\Oc_G\otimes_\CC U(\gen) \lra \Dc_G.
\ee
 This expresses $\Dc_G = \Oc_G \# U(\gen)$
as the {\em smash product} of $\Oc_G$ and the Hopf algebra $U(\gen)$
acting on it \cite[p.155]{sweedler}.

 If $V$ is a vector space (a $K=\{1\}$-module), then    
  $V^\sharp_G = \Oc_G\otimes V$. 
If, moreover,  $V$ is a $\gen$-module, then
we can form the smash product module
$\Oc_G\#  V$ over the smash product algebra 
$\Dc_G = \Oc_G \# U(\gen)$. This smash product module considered
,as an $\Oc_G$-module, is just
 $\Oc_G\otimes_\CC V =  V^\sharp_G$. This is the equivariant
  $\Dc$-module
 structure claimed in (a) for $K=\{1\}$. Alternatively, \eqref {eq:D_Gsemi}
 implies that
 \[
 V^\sharp_G = \Dc_G \otimes_{U(\gen)} V
 \]
 is just the induced module with $G$-action coming from the first factor. 
 In particular, it is coherent, if $V$ is finitely generated. This proves
 (a) for $K=\{1\}$. 
 
 \vskip .2cm
  
Let now $K$ be arbitrary.  
We claim that the $\Dc_G$-module $V^\sharp_G$ is the  pullback 
$p^* \Mc$ 
of a (unique, up to unique isomorphism) $\Dc_X$-module 
$\Mc$ isomorphic as an $\Oc_X$-module, to $V^\sharp_X$. 
Recall that left $\Dc$-module pullback is the same as $\Oc$-module
pullback, i.e., 
\[
p^* \Mc  = \Oc_G \otimes_{p^{-1} \Oc_X} p^{-1}  \Mc
\]
  By definition
\eqref {eq:def-sharp}, 
$p^{-1} V^\sharp_X$
consists of $V$-valued functions on $G$ satisfying the right
$K$-equivariance condition. 
 So we have an embedding  of $p^{-1} V^\sharp_X$ into  $V^\sharp_G$ as
a subsheaf, generating it over $\Oc_G$. Further, it is immediate that
this subsheaf is 
  annihilated by vector fields along
the fibers of $G\to X$. This implies that $V^\sharp_G$ is the pullback,
in the sense of $\Dc$-modules, of a uniquely defined $\Dc_X$-module
structure on $V^\sharp_X$. \qed

\subsection{The integrable category $\Oc$. Morphisms of Verma modules}\label{subsec:A-verma}
Consider now the situation of \S \ref{subsec:sweil-BGG}
\ref {par:parab-cat-O}, so 
  $G$ is a quasi-reductive  complex  algebraic  supergroup, $P$ a parabolic subgroup
with Levi subgroup/quotient $M$ and unipotent radical $N$, and $\bg, \bp, \bm, \bn$ the corresponding
Lie superalgebras.  We take $K=P$ and 
 denote $X=G/P$. 
 
 \vskip .2cm
 
 Let $\Oc^{\bg, P}\subset \Oc^{\bg, \bp}$ be the subcategory formed by
 $\bg$-modules $V$ such that the $\bp$-action integrates to a
 $P$-action. This is equivalent to requiring that the $\bm$-action
 integrates to an $M$-action, as $P = M\ltimes N$ and the $\bn$-action
 always integrates to an $N$ action by the local nilpotency condition. 
 [ADDRESS THE CASE WHEN  $M$ IS DISCONNECTED?]
 Thus $\Oc^{\bg, P}$ is a full subcategory in $(\bg, P)\Mod$
 and we can apply Proposition \ref {prop:harish1} to obtain
 a fully faithful exact functor
 \[
 \sharp_{gen, P}: \Oc^{\gen, P} \lra \Dc_X\Mod^G. 
 \]

\noindent Let $\rho: M\to GL(A)$ be a finite-dimensional 
irreducible representation considered
as a representation of $P$ via $P\to M$. Then, on one hand, we can form
the corresponding homogeneous vector
bundle 
  $\Lc^{G,P}_\rho = A^\sharp_X \in \Oc_X\Mod^G$    of finite rank
  on $X=G/P$.
  
  \vskip .2cm
  
  On the other hand, considering $A$ as a representation of $\bm$
  and $\bp$, we can form 
   the Verma module $\V^{\bg, \bp}_\rho$.
   Note that the $\bm$-action on $V^{\bg, \bp}_\rho$ integrates to
   an action of $M$, as it comes from $A$ which is an $M$-module. 
   Therefore $V^{\bg, \bp}_\rho$ is a
  $P$-module and so it  is a $(\gen, P)$-module. 
  Therefore we can form the 
   $G$-equivariant  $\Dc_{X}$-module $(\V^{\bg, \bp}_\rho)^\sharp_X$.

\begin{prop}\label{prop:loc-verma}
 We have an isomorphism in $\Dc_X\Mod^G$:
\[
(\V^{\bg, \bp}_{\rho})^\sharp_X  \, \= \,  
\Dc_{X} \otimes_{\Oc_{X}} \Lc^{G,P}_\rho. 
\]
 
\end{prop}

\noindent{\sl Proof:}  Each  side of the proposed isomorphism
are obtained by applying two natural induction procedures:
On the left we have  algebraic induction from $\pen$ to $\gen$
followed by  geometric induction to an equivariant quasicoherent
sheaf on $X$ (an $\Oc$-module which turns out to be a $\Dc$-module).
 On the right we have,  first, the geometric
induction to an equivariant coherent sheaf $(\Oc$-module)
on $X$ followed
the algebraic induction from $\Oc$-modules to $\Dc$-modules. 
The equivalence of both these $2$-step induction procedures follows
because their results are given by the same universal property. \qed

\vskip .2cm

This implies the following classical corollary going back to \cite{BGG}, see, e.g., 
\cite[Prop.VI.5.1]{faltings} for the particular
case of $G=\Sp(2n, \CC)$ or \cite[\S 11.2]{baston} in general. 

\begin{cor}\label{prop:Verma-DO}
For any finite-dimensional irreducible representations $\rho_1, \rho_2$ of 
$M$ we have
a canonical identification
\[
\Hom_{\Oc^{\bg, P}} (\V^{\bg, \bp}_{\rho_1}, \V^{\bg, \bp}_{\rho_2}) \,\= \,  
\Dc\bigl((\Lc^{G,P}_{\rho_2})^\vee, (\Lc^{G,P}_{\rho_1})^vee\bigr)^G, 
\]
where the RHS is the space of $G$-invariant differential operators
between the dual  equivariant vector bundles
$(\Lc^{G,P}_{\rho_i})^\vee = \Lc^{G,P}_{\rho_i^\vee}$. 
\qed
\end{cor} 

\subsection{ The fractional case}\label{subsec:App-frac}
We now generalize the homogenization procedure to twisted  bundles.

\vskip .2cm

In the situation of \S \ref {subsec:A-harish}, assume that we are
given a central extension
\be\label{eq:ext-K-A}
1\to H\lra \wt K \lra K\to 1
\ee
by a purely even finite group $H$ and a character $\sigma: H\to \CC^*$. 
Let $\wt K\Mod_\sigma\subset \wt K\Mod$ be the full subcategory
of representations on which $H$ acts via $\sigma$. 

\vskip .2cm

On the other hand, we can apply the construction of Example 
\ref {exas:ourTDO}(a) to the right $K$-torsor $p: Y=G\to X=G/K$
and obtain a $\CC^*$-gerbe $\Gc_\sigma = \sigma_*\delta(Y)$ on $X$. 
This gerbe is $G$-equivariant and we can consider the following categories:
\begin{itemize}
\item $\Oc_X\Mod^G_\sigma$, the category of $G$-equivariant
$\Gc_\sigma$-twisted quasi-coherent $\Oc_X$-modules,
 i.e., of data
which to each local object of $\Gc_\sigma$ over $U\subset X$  associate a
quasi-coherent
$\Oc_U$-module with standard compatibility properties, see \S 
\ref {subsec:twisted}. 

\item $\Dc_X\Mod^G_\sigma$, the  the category of $G$-equivariant
$\Gc_\sigma$-twisted quasi-coherent $\Dc_X$-modules, i.e., of data
which to each local object of $\Gc_\sigma$ over $U\subset X$  similarly
 associate a
quasi-coherent
$\Dc_U$-module.
\end{itemize}
\noindent
Given $V\in \wt K\Mod_\sigma$, we associate to it an object $V^\sharp_X
\in \Oc_X\Mod^G_\sigma$ as follows. Given an object of $\Gc_\sigma(U)$,
i.e., a right $\wt K$-torsor $\wt p: \wt Y_U \to U$ lifting the $K$-torsor
$p^{-1}(U)\to U$, the sheaf $V^\sharp_X(\wt Y)$ on $U$ has, for any
$U'\subset U$
 \[
\Gamma(U', V^\sharp_X)\,=\,
\bigl\{ f: \wt p^{-1}(U') \to V \bigl| \, f(z \wt k) = \wt k^{-1} f(g), \quad \forall \, 
z\in\wt  p^{-1}(U'),
\, \wt k\in \wt K\bigr\}, 
\]
 generalizing  \eqref {eq:def-sharp} and we have:

\begin{prop}
The correspondence $V\mapsto V^\sharp_X$ extends to a fully faithful
exact functor 
\[
\sharp_{\wt K,\sigma}: \wt K\Mod_\sigma \lra \Oc_X\Mod^G_\sigma. 
\]
\qed
\end{prop}

\noindent Further, we can speak about Harish-Chandra  $(\bg, \wt K)$-modules,
i.e., of $\bg$-modules $V$ equipped with a compatible structure
of a $\wt K$-module (note that $\Lie(\wt K)\= \Lie(K)\subset \bg$). 
We  denote by 
  $(\bg, \wt K)\Mod$  their category  and by by 
  $(\bg, \wt K)\Mod_\sigma$ the full subcategory formed by $V$ on
  which $\wt K$ acts by $\sigma$. 
  
   Proposition \ref {prop:harish1} immediately extends to the following.

  \begin{prop}\label{prop:harish-tw}
  For any $V\in (\bg, \wt K)\Mod_\sigma$ the $\Gc_\sigma$-twisted
  qoasicoherent sheaf $V^\sharp_X$ has a natural structure of 
  of a $\Gc_\sigma$-twisted $\Dc_X$-module and we get a
  fully faithful exact functor
  \[
  \sharp_{\bg, \wt K,\sigma}: (\bg, \wt K)\Mod_\sigma \lra \Dc_X\Mod^G_\sigma.
  \]
  \qed
  \end{prop}
  
   We will call $\sharp_{\wt K,\sigma}$ {\em twisted homogenization}. 

\vskip .2cm

 Let us specialize to the situation o \S \ref {subsec:A-verma}, so $K=P\subset G$
  is a parabolic subgroup in a quasi-reductive algebraic supergroup $G$.  
  We assume that the central extension $\wt P$ of $P=M\rtimes N$ comes from
  an extension $\wt M$ of $M$; note that $\wt M$ is reductive. 
  
  Let $\Oc^{\bg, \wt P, \sigma}$ be the category formed by $\bg$-modules
  $V\in\Oc^{\bg, \bp}$ together with a compatible action of $\wt P$ so that $A$
  acts by $\sigma$. If $\wt P$ is connected, the existence of a compatible action
  is a condtition   on $V$ as an object of $Oc^{\bg, \bp}$. Indeed, $\Lie(\wt P) = \bp$,
  $\Lie(\wt M)=\bm$ and  so the compatible action exists if and only if the action of $\bm$
  integrates to that of $\wt M$. Proposition \ref {prop:harish-tw} gives then the
  exact fully faithful  functor
  \[
  \sharp_{\bg, \wt P, \sigma}: \Oc^{\bg, \wt P, \sigma} \lra \Dc_X\Mod^G_\sigma.
  \]
  For a  finite-dimensional representation  $\rho: \wt M\to \GL(B)$   in $\wt M\Mod_\sigma$ let
  $\V^{\bg, \wt P}_\rho = \Ind_\bp^\bg( B)$ be the corresponding induced Verma module.
  Then we have a twisted version of Proposition \ref {prop:loc-verma} and Corollary
  \ref {prop:Verma-DO}.
  
  \begin{prop}\label{prop:harsish-twisted-n}
  (a) The $\Gc_\sigma$-twisted $\Dc_X$-module  $(\V^{\bg, \wt P}_\rho )_X^\sharp$
  is identified with $\Dc_X \otimes_{\cal O_X} \Lc^{G,P}_\rho$,
  the tensor product of $\Dc_X$ and the $\Gc$-twisted vector bundle $\Lc^{G,P}_\rho := B^\sharp_X$.
  
  \vskip .2cm
  
  (b) For any two finite-dimensional $\rho_1,\rho_2\in M\Mod_\sigma$ we have
  \[
\Hom_{\Oc^{\bg, \wt P,\sigma }} (\V^{\bg, \wt P}_{\rho_1}, \V^{\bg, \wt P}_{\rho_2}) \,\= \,  
\Dc\bigl((\Lc^{G,\wt P}_{\rho_2})^\vee, (\Lc^{G,\wt P}_{\rho_1})^\vee\bigr)^G, 
\]
\qed
  \end{prop}

\vfill\eject

	\vskip 2cm
	
{\small	Kavli  IPMU (WPI),  UTIAS, The University of Tokyo, Kashiwa, Chiba 277-8583, Japan, email: {\tt mikhail.kapranov@protonmail.com}

}

\end{document}